\documentclass[11pt]{article}
\usepackage[centertags,leqno]{amsmath}
\usepackage{amssymb,amsfonts}
\usepackage{mathrsfs}
\usepackage{stmaryrd}
\usepackage{marvosym}
\usepackage{color}
\usepackage{tcolorbox}
\usepackage{graphicx}
\usepackage[colorlinks=true]{hyperref}

\newtheorem{theorem}{Theorem}[section]

\newtheorem{corollary}[theorem]{Corollary}
\newtheorem{example}[theorem]{Example}

\newtheorem{lemma}[theorem]{Lemma}
\newtheorem{proposition}[theorem]{Proposition}
\newtheorem{remark}[theorem]{Remark}

\definecolor{violet}{rgb}{0.5,0,0.5}
\definecolor{orange}{cmyk}{0,0.3,0.7,0}

\definecolor{poznCol}{rgb}{0.5,0,1}
\newcommand{\poznamka}[1]%
{\noindent\bf\textcolor{poznCol}{$\maltese$...#1...$\maltese$}}

\newcommand{\vertiii}[1]%
{{\left\vert\kern-0.25ex\left\vert\kern-0.25ex\left\vert #1
  \right\vert\kern-0.25ex\right\vert\kern-0.25ex\right\vert}}

\newtheorem{definition}{Definition}[section]

\newtheorem{hypo}{Hypothesis}

\newtheorem{hypos}{Hypotheses}

\newcommand{\proof}{{\em Proof. }}
\newcommand{\qed}{\rule{2mm}{2mm}}

\newcommand{\Square}{$\sqcap$\hskip -1.5ex $\sqcup$}

\numberwithin{equation}{section}
\newcommand{\RR}{\mathbb{R}}
\newcommand{\CC}{\mathbb{C}}
\newcommand{\NN}{\mathbb{N}}
\newcommand{\ZZ}{\mathbb{Z}}
\newcommand{\ii}{\mathrm{i}}
\newcommand{\ee}{\mathrm{e}}

\newcommand{\eqdef}{\stackrel{{\mathrm {def}}}{=}}
\newcommand{\eps}{\varepsilon}
\renewcommand{\colon}{:\,}

\newcommand{\RE}{\mathop{\Re\mathfrak{e}}}  
\newcommand{\IM}{\mathop{\Im\mathfrak{m}}}  

\newcommand{\td}{\,\mathrm{d}}       

\begin{document}

\baselineskip=16pt

\begin{center}
\begingroup\large\bf
  Space\--Time Analyticity of Weak Solutions to\\
  Semi\-linear Parabolic Systems with Variable Coefficients%
\footnote{
Research supported in part by
Deutsche Forschungs\-gemeinschaft (D.F.G., Germany) \\
under Grant {\#} TA~213/16--1.
}
\\
\endgroup
\end{center}

\vspace{0.5cm}
\begin{center}
       Falko {\sc Baustian}\\
       and\\
       Peter {\sc Tak\'a\v{c}}
\\
\vspace*{2.0mm}
\begingroup\small
       Institut f\"ur Mathematik, Universit\"at Rostock,\\
       Ulmenstra{\ss}e~69, Haus~3, D-18051 Rostock, Germany,\\
       {\it e\/}-mail: {\tt peter.takac@uni-rostock.de}\\
\vspace*{2.0mm}
\color{blue}
       Web: {\tt
       https://www.mathematik.uni-rostock.de/struktur/lehrstuehle/%
             angewandte-analysis/}
\color{black}
\endgroup
\end{center}


\begin{center}
\today
\\
\color{blue}
\color{black}
\end{center}

\vspace{0.1cm}

\baselineskip=12pt
%
\begin{abstract}
\begingroup\footnotesize
Analytic smooth solutions of
a general, strongly parabolic semi\-linear Cauchy problem
of $2m$-th order in $\mathbb{R}^N\times (0,T)$
with analytic coefficients (in space and time variables)
and analytic initial data (in space variables) are investigated.
They are expressed in terms of holomorphic continuation of
global (weak) solutions to the system valued in
a suitable Besov interpolation space of $B^{s;p,p}$\--type
at every time moment $t\in [0,T]$.
Given $0 < T'< T\leq \infty$,
it is proved that any $B^{s;p,p}$\--type solution
$u\colon \mathbb{R}^N\times (0,T)\to \CC^M$
with analytic initial data possesses a bounded holomorphic continuation
$u(x + \mathrm{i}y, \sigma + \mathrm{i}\tau)$ into a complex domain in
$\mathbb{C}^N\times \mathbb{C}$ defined by
$(x,\sigma)\in \mathbb{R}^N\times (T',T)$, $|y| < A'$ and $|\tau | < B'$,
where $A', B'> 0$ are constants depending upon~$T'$.
The proof uses the extension of a weak solution to
a $B^{s;p,p}$-type solution in a domain in
$\mathbb{C}^N\times \mathbb{C}$, such that this extension
satisfies the Cauchy\--Riemann equations.
The holomorphic extension is obtained with a help from
holomorphic semigroups and maximal regularity theory for
parabolic problems in Besov interpolation spaces of
$B^{s;p,p}$\--type imbedded (densely and continuously) into
an $L^p$-type Lebesgue space.
Applications include {\em risk models\/}
for European options in Mathematical Finance.
\endgroup
\end{abstract}

\vspace{0.1cm}

\par\vspace*{0.2cm}
\noindent
\begin{tabular}{ll}
{\bf Keywords:}
& Space-time analyticity, parabolic PDE;\\
& holomorphic semigroup, Besov space;\\
& maximal regularity, Hardy space;\\
& holomorphic continuation to a complex strip;\\
& European option, bilateral counterparty risk.
\end{tabular}

\vfill
\par\vspace*{0.2cm}
\noindent
\begin{tabular}{lll}
{\bf 2020 Mathematics Subject Classification:}
& Primary   & 35B65, 35K10;\\
& Secondary & 32D05, 91G40;\\
\end{tabular}

\newpage
\tableofcontents

\vfill\eject
\baselineskip=15pt
\section{Introduction}
\label{s:Intro}

In this article we investigate analyticity (in space and time variables)
of strict $L^p$-type solutions
\begin{math}
  \mathbf{u} = (u_1,\dots,u_M)\colon \mathbb{R}^N\times (0,T)\to \CC^M
\end{math}
(or $\CC^M$)
of the classical Cauchy problem for a strongly parabolic system of
$M$ (coupled) semi\-linear partial differential equations
of order $2m$ ($m\geq 1$ -- an integer)
with {\em analytic\/} coefficients and with {\em analytic\/} initial data
$\mathbf{u}_0$ belonging to the real interpolation space
$\mathbf{B}^{s;p,p}(\RR^N)$, such that the function
$\mathbf{u}\colon [0,T]\to \mathbf{B}^{s;p,p}(\RR^N)$
is continuous.
Here,
$\mathbf{B}^{s;p,p}(\RR^N) = [ B^{s;p,p}(\RR^N) ]^M$
where $B^{s;p,p}(\RR^N)$ denotes the Besov space
\begin{equation*}
    B^{s;p,p}(\RR^N)\eqdef
    \left( L^p(\RR^N) ,\, W^{2m,p}(\RR^N) \right)_{s/(2m),p}
  = \left( L^p(\RR^N) ,\, W^{2m,p}(\RR^N) \right)_{ 1 - (1/p) ,\, p }
\end{equation*}
with $1 < p < \infty$, $p > 2 + \frac{N}{m}$, and
$s = 2m\left( 1 - \frac{1}{p}\right)\in (0,2m)$.
It is defined by {\it real\/} interpolation, e.g., in
{\sc R.~A.\ Adams} and {\sc J.~J.~F.\ Fournier}
\cite[Chapt.~7]{AdamsFournier}, {\S}7.6--{\S}7.23, pp.\ 208--221,
{\sc A.\ Lunardi}
\cite[Chapt.~1]{Lunardi}, {\S}1.2.2, pp.\ 20--25, or in
{\sc H.\ Triebel} \cite[Chapt.~1]{Triebel},
{\S}1.2--{\S}1.8, pp.\ 18--55.
Since the Besov space $B^{s;p,p}(\RR^N)$ is not imbedded into
the Hilbert space $L^2(\RR^N)$ whenever $2 < p < \infty$,
we find it convenient to consider {\em strict\/} $L^p$-type solutions
$\mathbf{u}\colon \mathbb{R}^N\times (0,T)\to \CC^M$
having the {\em maximal regularity property\/} 
(cf.\
{\sc A.\ Ashyralyev} and {\sc P.~E.\ Sobolevskii}
\cite[Chapt.~3, pp.\ 21--36]{Ashyr-Sobol-94} and
{\sc J.\ Pr\"uss} \cite{Pruess})
rather than weak $L^2$\--type solutions treated in
{\sc P.\ Tak\'a\v{c}} \cite{Takac}
for the corresponding linear partial differential equation,
but with arbitrary {\em non\-smooth\/} initial data
$\mathbf{u}_0\in \mathbf{L}^2(\RR^N)$.
Consequently, we will be able to apply the classical theory of
linear and semilinear evolutionary problems of parabolic type
in a Besov space as presented, e.g., in
{\sc H.~Amann} \cite[Chapt.~III, {\S}4, pp.\ 128--191]{Amann-Birkh},
{\sc Ph.\ Cl\'ement} and {\sc S.\ Li} \cite{Clem-Li},
{\sc A.\ Lunardi} \cite[Chapt.~7, pp.\ 257--289]{Lunardi},
{\sc M.\ K\"ohne}, {\sc J.\ Pr\"u{ss}}, and {\sc M.\ Wilke}
\cite{Koehne-Pruess-W}, and
{\sc H.~Tanabe} \cite[Chapt.\ 5--6, pp.\ 117--229]{Tana-79}.
Our Cauchy problem has the following general form for
a semilinear $2m^{\mathrm{th}}$\--order parabolic problem,
\begin{equation}
\label{e:Cauchy}
\left\{
\begin{alignedat}{2}
    \frac{\partial\mathbf{u}}{\partial t}
  + \mathbf{P}
    \Bigl( x,t, \frac{1}{\ii} \frac{\partial}{\partial x} \Bigr)
    \mathbf{u}
& = \mathbf{f}
    \left( x,t;
    \left( \frac{ \partial^{|\beta|}\mathbf{u} }{\partial x^{\beta}}
    \right)_{|\beta|\leq m}
    \right)
&&  \quad\mbox{ for } (x,t)\in \mathbb{R}^N\times (0,T) \,;
\\
  \mathbf{u}(x,0)&= \mathbf{u}_0(x)
&&  \quad\mbox{ for }\; x\in \mathbb{R}^N \,.
\end{alignedat}
\right.
\end{equation}

Here,
\begin{math}
  {\partial} / {\partial x} =
  ( {\partial}/{\partial x_1} ,\dots , {\partial}/{\partial x_N} )
\end{math}
stands for the spatial gradient and
\begin{math}
  \xi\mapsto \mathbf{P}(x,t,\xi)
\end{math}
is a~polynomial of order $2m$ in the variable
$\xi = (\xi_1,\dots,\xi_N)\in \RR^N$ (or $\CC^N$);
its coefficients are $M\times M$ matrices (real or complex)
which are assumed to be real analytic (jointly) in both variables
$x\in \RR^N$ and $t\in (0,T)$.
Also the nonlinearity
\begin{math}
  (x,t;X)\mapsto \mathbf{f}(x,t;X)
\end{math}
(a reaction function valued in $\RR^M$ or $\CC^M$)
is assumed to be analytic in all variables
$x\in \RR^N$, $t\in (0,T)$, and
\begin{math}
  X = \left( X_{\beta}\right)_{|\beta|\leq m} \in \RR^{M\tilde{N}}
\end{math}
(or $\CC^{M\tilde{N}}$),
where we have substituted
\begin{math}
  X_{\beta} = \frac{ \partial^{|\beta|}\mathbf{u} }{\partial x^{\beta}}
  \in \RR^M
\end{math}
(or $\CC^M$)
for the (mixed) partial derivative of $\mathbf{u}$
with a multi\--index
$\beta = (\beta_1,\dots,\beta_N)$ $\in (\ZZ_+)^N$ of order
$|\beta| = \beta_1 + \dots + \beta_N$, $|\beta|\leq m$.
Here,
$\ZZ_+ = \{ 0,1,2,\dots\}$ and
the Euclidean dimension of the $m$-jet $X$ equals to $M\tilde{N}$ with
\begin{equation}
\label{e:N^tilde}
    \tilde{N}
  = \sum_{k=0}^m \sum_{|\beta| = k}
    \genfrac{(}{)}{0pt}0{k}{\beta}
  = \sum_{k=0}^m N^k
    \quad\mbox{ where }\quad
    \genfrac{(}{)}{0pt}0{k}{\beta}
  \eqdef \frac{k!}{ \beta_1 !\, \beta_2 !\,\dots\, \beta_N ! } \,.
\end{equation}
As usual, $\RR^N$ and $\CC^N$, respectively, denote
the $N$-dimensional real and complex Euclidean spaces,
$\ii = \sqrt{-1}$, and $M,N\in \NN$ where
$\NN = \{ 1,2,3,\dots\}$.
We have identified
\begin{math}
  X_{\beta} = \frac{ \partial^{|\beta|}\mathbf{u} }{\partial x^{\beta}}
  \equiv \mathbf{u}(x,t)
\end{math}
for $\beta = (0,0,\dots,0)$ of order $|\beta| = 0$.
%

As already indicated,
we impose certain standard {\em strong ellipticity\/} and
{\em analyticity\/} hypotheses on the coefficients of
the partial differential operator
\begin{math}
    \mathbf{P}
    \left( x,t, \frac{1}{\ii} \frac{\partial}{\partial x} \right)
\end{math}
and on the reaction function $\mathbf{f}(x,t;X)$ as well.
Assuming that
$\mathbf{u}_0\in \mathbf{B}^{s;p,p}(\RR^N)$ ($p > 2 + \frac{N}{m}$)
possesses a complex analytic extension to a strip
$\mathfrak{X}^{(\kappa_0)}$ of constant width in $\CC^N = \RR^N + \ii\RR^N$
and the first\--order partial derivatives
\begin{equation*}
  \frac{\partial}{\partial t} \mathbf{f}(x,t;X)
  \quad\mbox{ and }\quad
  \frac{\partial}{ \partial X_{\beta} }\,
  \mathbf{f}(x,t;X) \,,
    \quad\mbox{ for $|\beta|\leq m$, }
\end{equation*}
are locally uniformly bounded for
\begin{math}
  (x,t;X)\in \mathbb{R}^N\times (0,T)\times \mathbb{R}^{M\tilde{N}} ,
\end{math}
in this work we show that the (unique) strict ($L^p$-type) solution
$\mathbf{u} = \mathbf{u}(x,t)$ of problem \eqref{e:Cauchy}
is {\em real analytic\/} in $(x,t)\in \RR^N\times (0,T)$.
Notice that the latter condition
(local boundedness of all first\--order partial derivatives
\begin{math}
  { \partial \mathbf{f} } / { \partial X_{\beta} }
)
\end{math}
is equivalent with
$X\mapsto \mathbf{f}(x,t;X)$
being locally uniformly Lipschitz\--continuous.

This analyticity claim is motivated by the standard formula for
the solution of the Cauchy problem for the heat equation in $\RR^N$
(with the Laplace operator $\Delta$, i.e.,
\begin{math}
    \mathbf{P}
    \left( x,t, \frac{1}{\ii} \frac{\partial}{\partial x} \right)
  = - \Delta \,,
\end{math}
$\mathbf{f}(x,t;X) = \mathbf{0}$, and $M=1$);
see e.g.\
{\sc F.\ John} \cite{John}, Chapt.~7, Sect.~1, eq.\ (1.11), p.~209.
The heat equation case has been significantly generalized in
{\sc P.\ Tak\'a\v{c}} et al.\ \cite[Theorem 2.1, p.~429]{TakacBoller},
where only the leading coefficients of the operator
\begin{math}
    \mathbf{P}
    \left( x,t, \frac{1}{\ii} \frac{\partial}{\partial x} \right)
\end{math}
are assumed to be constant, but it is required that
$\mathbf{u}_0\in \mathbf{L}^{\infty}(\RR^N) = [L^{\infty}(\RR^N)]^M$.
In our present work, the {\em analyticity\/} hypothesis on the initial data
$\mathbf{u}_0$ resembles more to a non\-local version of
the classical Cauchy\--Kowalewski theorem
({\sc F.\ John} \cite{John}, Chapt.~3, Sect.\ 3(d), pp.\ 73--77).
We will show that, under this analyticity hypothesis on
$\mathbf{u}(\,\cdot\,,0) = \mathbf{u}_0$, if a solution
$\mathbf{u}\colon \RR^N\times [0,T)\to \CC^M$ exists,
then it must be analytic in $\RR^N\times (0,T)$.
We are able to specify also the {\em domain of analyticity\/}
in terms of a {\em complex analytic extension\/}.
The restriction on the initial data
$\mathbf{u}_0\in \mathbf{B}^{s;p,p}(\RR^N)$, with the conditions
$p > 2 + \frac{N}{m}$ and
$s = 2m\left( 1 - \frac{1}{p}\right)\in (0,2m)$,
allows us to take advantage of (the continuity of)
the Sobolev\-(-Besov) imbedding
\begin{math}
  B^{s;p,p}(\RR^N) \hookrightarrow
  C^m(\RR^N) \cap W^{m,\infty}(\RR^N) ;
\end{math}
see, e.g.,
{\sc R.~A.\ Adams} and {\sc J.~J.~F.\ Fournier}
\cite[Chapt.~7]{AdamsFournier}, Theorem 7.34(c), p.~231.
This more restrictive condition on the initial data $\mathbf{u}_0$
enables us to work with an $m$-jet
\begin{math}
  X = \left( X_{\beta}\right)_{|\beta|\leq m} \in \CC^{M\tilde{N}}
\end{math}
whose all components
\begin{math}
  X_{\beta} = \frac{ \partial^{|\beta|}\mathbf{u} }{\partial x^{\beta}}
  \in \CC^M
\end{math}
are bounded continuous functions of $(x,t)\in \RR^N\times [0,T)$;
thus, each $X_{\beta}(\,\cdot\,,t)$ ($|\beta|\leq m$)
belongs to
$\mathbf{L}^{\infty}(\RR^N)$ at every time $t\in [0,T)$.
Consequently, we can apply the Banach fixed point theorem to
problem \eqref{e:Cauchy} in a way similar to
\cite[Theorem 2.1, p.~429]{TakacBoller}.
For instance, in a typical second\--order parabolic problem
(i.e., eq.~\eqref{e:Cauchy} with $m=1$)
we can allow for a reaction function
\begin{math}
  \mathbf{f}
    \left( x,t; \mathbf{u},\, \frac{ \partial\mathbf{u} }{\partial x}
    \right)
\end{math}
depending on $\mathbf{u}$ and its gradient
${ \partial\mathbf{u} } / {\partial x}$
(${}\equiv \ii\, D_x\mathbf{u}$), besides the independent variables
$x\in \RR^N$ and $t\in (0,T)$.

The main contribution of our present article is that
we are able to {\em remove\/} the hypothesis that
the leading coefficients must be constant, in analogy with
{\sc P.\ Tak\'a\v{c}} \cite[Theorem 3.3, p.~59]{Takac}
where the corresponding linear system is treated.
In contrast to \cite[Proposition A.4, p.~446]{TakacBoller},
this means that we {\em cannot\/} calculate the Green function for
the Cauchy problem with the leading coefficients only,
\begin{equation}
\label{e:Cauchy-lead}
\left\{
\begin{alignedat}{2}
    \frac{\partial\mathbf{u}}{\partial t}
  + (-1)^m \sum_{|\alpha| = 2m} \mathbf{P}^{(\alpha)} (x,t)\,
    \frac{ \partial^{|\alpha|} \mathbf{u} }{ \partial x^{\alpha} }
& = \mathbf{0}
&&  \quad\mbox{ for } (x,t)\in \mathbb{R}^N\times (0,T) \,;
\\
  \mathbf{u}(x,0)&= \mathbf{u}_0(x)
&&  \quad\mbox{ for } x\in \mathbb{R}^N \,,
\end{alignedat}
\right.
\end{equation}
and then simply take advantage of the variation\--of\--constants formula
\cite[eq.\ (3.22), p.~437]{TakacBoller}
in order to obtain the solution of the original problem \eqref{e:Cauchy}.
Fortunately, the methods from \cite{Takac},
based on a~priori $L^2$-type estimates combined with
the Cauchy\--Riemann equations,
are applicable also to our semilinear system
\eqref{e:Cauchy} provided that already the initial data $\mathbf{u}_0$
are analytic.
Here, each
$\mathbf{P}^{(\alpha)} (x,t)$ is an $M\times M$ matrix and recall that
\begin{math}
  { \partial^{|\alpha|} \mathbf{u} } / { \partial x^{\alpha} } =
  \frac{ \partial^{|\alpha|} \mathbf{u} }
       { \partial x_1^{\alpha_1} \,\dots\, \partial x_N^{\alpha_N} }
\end{math}
denotes the (mixed) partial derivative of
$\mathbf{u}\colon \mathbb{R}^N\times (0,T)\to \CC^M$
with a multi\--index
$\alpha = (\alpha_1,\dots,\alpha_N)$ $\in (\ZZ_+)^N$ of order
$|\alpha| = \alpha_1 + \dots + \alpha_N$.
This means that, for the semilinear parabolic Cauchy problem
\eqref{e:Cauchy}, we {\em do not improve\/} the regularity properties of
(in general) non\-smooth initial data to {\em analytic\/} regularity
as time passes by (for $t\in (0,T)$).
We show only that the analytic regularity of the initial data
$\mathbf{u}_0$ (at $t=0$)
is {\em preserved\/} for all times $t\in (0,T)$.
In contrast, analytic regularity of the initial data is {\em not assumed\/}
in \cite{Takac, TakacBoller}.

As in \cite{Takac, TakacBoller},
our method is based on the simple fact that a function
$u\colon \mathbb{R}^N\times (0,T)\to \RR$ (or $\CC$)
is real analytic if and only if it has
a {\em holomorphic\/} (i.e., complex analytic) extension
$\tilde{u}\colon \Omega\to \CC$ to some complex domain $\Omega$ such that
\begin{math}
  \mathbb{R}^N\times (0,T)\subset \Omega\subset
  \mathbb{C}^N\times \mathbb{C} ,
\end{math}
i.e.,
$u = \tilde{u}\vert_{ \mathbb{R}^N\times (0,T) }$,
the restriction of $\tilde{u}$ to $\mathbb{R}^N\times (0,T)$.
If the domain $\Omega$ is fixed then the holomorphic extension
$\tilde{u}$ of $u$ to $\Omega$ is always unique, see e.g.\
{\sc F.\ John} \cite{John}, Chapt.~3, Sect.\ 3(c), pp.\ 70--72.
Thus, in order to show that the weak solution
$\mathbf{u} = \mathbf{u}(x,t)$ of problem \eqref{e:Cauchy}
is real analytic in $\mathbb{R}^N\times (0,T)$,
it suffices to construct a holomorphic extension
$\tilde{\mathbf{u}}$ of $\mathbf{u}$ to some complex domain $\Omega$
\begin{math}
  (
  \mathbb{R}^N\times (0,T)\subset \Omega\subset
  \mathbb{C}^N\times \mathbb{C}
  ) .
\end{math}
Due to the uniqueness (of a holomorphic extension),
we often drop the tilde ``$\,\overset{\sim}{\phantom{u}}\,$''
in the notation for the (unique) holomorphic extension.
Analogous ideas
(holomorphic extension, uniqueness, and
 Bergman and Szeg\H{o} spaces of holomorphic functions)
were used earlier in {\sc N.\ Hayashi}
\cite{Hayashi_Duke90, Hayashi_Duke91, Hayashi_SIAM91, Hayashi_IHP92}.

Instead of using the Green function method (cf.~\cite{TakacBoller}),
we establish the existence of solutions to
the Cauchy problem \eqref{e:Cauchy}
in a complex parabolic domain
$\mathfrak{X}^{(r)}\times [0,T)$ in $\CC^N\times \CC$
with initial data $\mathbf{u}_0$ from
a space of holomorphic functions whose domain
$\mathfrak{X}^{(r)} = \RR^N + \ii Q^{(r)}$
is a {\em tube\/} in $\CC^N$ with {\em base\/} $Q^{(r)} = (-r,r)^N$,
for some $0 < r < \infty$, see
{\sc P.\ Tak\'a\v{c}} \cite[eq.~(21), p.~58]{Takac}.
The (complex) analyticity in space is then verified by means of
the Cauchy\--Riemann equations, whereas the (complex) analyticity in time
is obtained from the properties of holomorphic semigroups
in the Besov space
$\mathbf{B}^{s;p,p}(\RR^N) = [ B^{s;p,p}(\RR^N) ]^M$.
Our use of the Cauchy\--Riemann equations already at the initial time
$t=0$ requires that $\mathbf{u}_0$ be (complex) {\em analytic\/} in
$\mathfrak{X}^{(r)}$.

In order to provide a quick, nontechnical hint to our approach,
we now give an illustrative weaker version of our main result,
Theorem~\ref{thm-Main} in Section~\ref{s:Main},
for a single equation in one space dimension ($M=N=1$),
\begin{equation}
\label{e:Cauchy:M=N=1}
\left\{
\begin{alignedat}{2}
    \frac{\partial u}{\partial t}
  = a(x,t) \frac{\partial^2 u}{\partial x^2}
& + b(x,t) \frac{\partial u}{\partial x}
  + c(x,t) u
  + f\left( x,t; u ,\, \frac{\partial u}{\partial x} \right)
\\
&   \quad\mbox{ for } (x,t)\in \mathbb{R}^1\times (0,T) \,;
\\
  u(x,0) = u_0(x)
&   \quad\mbox{ for } x\in \mathbb{R}^1 \,.
\end{alignedat}
\right.
\end{equation}

We begin with the complexifications of the spatial and temporal variables,
$x\in \RR^1$ and $t\in (0,T)$, respectively:
Given any real numbers $0 < r < \infty$ and $0 < T'\leq T < \infty$,
we introduce the complex domains
\begin{align}
\nonumber
  \mathfrak{X}^{(r)} &\eqdef
  \{ z = x + \ii y\in \CC\colon |y| < r\}
  = \RR + \ii (-r,r) \,,
\\
\label{loc:T^(infty)}
  \Delta_{\vartheta} &\eqdef
  \{ t = \varrho\ee^{\ii\theta} \in \CC\colon
     \varrho > 0\mbox{ and } \theta\in (-\vartheta, \vartheta) \} \,,
    \quad \vartheta = \arctan (r/T') \,,
\\
\label{loc:T^(T')}
  \Delta_{\vartheta}^{(T')} &\eqdef
    \Delta_{\vartheta} \cap \{ t\in \CC\colon 0 < \RE t < T'\}
\\
\nonumber
& = \{ t = \varrho\ee^{\ii\theta} \in \CC\colon
       |\theta| < \vartheta \,\mbox{ and }\,
       0 < \varrho < T'/ \cos\theta \} \,,
\\
\intertext{and}
\label{glob:T^(r)}
  \Delta_{\vartheta}^{T',T} &\eqdef
  \Delta_{\vartheta}^{(T)}
    \cap \{ t\in \CC\colon |\IM t| < T'\cdot \tan\vartheta \}
  = \bigcup_{ 0\leq \xi\leq T-T'}
    \left( \xi + \Delta_{\vartheta}^{(T')} \right)
\\
\nonumber
& = \bigcup_{ 0\leq \xi\leq T-T'}
    \{ \xi + t'\in \CC\colon t'\in \Delta_{\vartheta}^{(T')} \}
  = [ 0, T-T'] + \Delta_{\vartheta}^{(T')}
\end{align}
with the angle $\vartheta\in (0, \pi / 2)$ given by
$\tan\vartheta = r/T'$.
Of course, if $T = T'$ then
$\Delta_{\vartheta}^{T',T} = \Delta_{\vartheta}^{(T')}$
is an open triangle.
Clearly, we have
\begin{equation*}
    \Delta_{\vartheta}^{T',T}
  = \bigcup_{ 0 < r\leq T'}
    \mathfrak{T}_{ r\cdot \cot\vartheta ,\, T}^{(r)}
  = \bigcup_{ 0 < r\leq T'}
    [ ( r\cdot \cot\vartheta ,\, T) + \ii (-r,r) ]
\end{equation*}
where
\begin{align}
\label{def:T^(r)}
  \mathfrak{T}_{T',T}^{(r)} &\eqdef
  \{ t = \sigma + \ii\tau\in \CC\colon
     T'< \sigma < T\mbox{ and } |\tau| < r\}
  = (T',T) + \ii (-r,r) \,.
\end{align}
We set
$\mathfrak{T}_{0,T}^{(r)} = (0,T) + \ii (-r,r)$
if $T'= 0$.
The closures in $\CC$ of $\mathfrak{X}^{(r)}$,
$\Delta_{\vartheta}$, $\Delta_{\vartheta}^{(T')}$,
$\Delta_{\vartheta}^{T',T}$, and
$\mathfrak{T}_{T',T}^{(r)}$ are denoted by
$\overline{\mathfrak{X}}^{(r)}$,
$\overline{\Delta}_{\vartheta}$,
$\overline{\Delta}_{\vartheta}^{(T')}$,
$\overline{\Delta}_{\vartheta}^{T',T}$, and
$\overline{\mathfrak{T}}_{T',T}^{(r)}$, respectively.

\begin{figure}
\label{fig:1}
\centering
\includegraphics[width=0.7\textwidth]{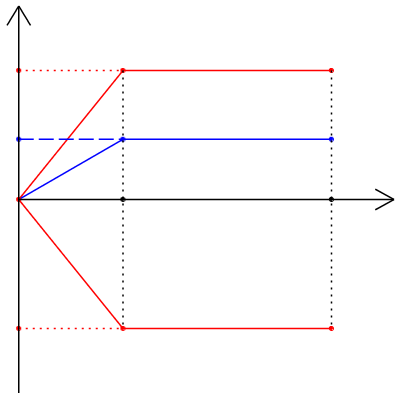}
\caption{
The triangle $\Delta_{\vartheta}$ starting at the origin has been defined
in \eqref{loc:T^(infty)}.
Its shifts to the right create the region $\Delta_{\vartheta}^{(T')}$
in \eqref{loc:T^(T')}.
}
\end{figure}

The Banach space of all continuous ($B^{s;p,p}(\RR)$\--valued) functions
$u\colon [0,T]\to B^{s;p,p}(\RR)$ is denoted by
$C\left( [0,T]\to B^{s;p,p}(\RR)\right)$;
it is endowed with the natural supremum norm
\hfil\break
\begin{math}
  \vertiii{u}_{ L^{\infty}(0,T) }\eqdef
  \sup_{t\in [0,T]} \| u(\,\cdot\,,t)\|_{ B^{s;p,p}(\RR) } < \infty .
\end{math}
%

\begin{theorem}\label{thm-nontech}
$(M=N=1)$
Let\/ $p > 2 + \frac{1}{m}$,
$s = 2m\left( 1 - \frac{1}{p}\right)$, and\/ $0 < T < \infty$.
Assume that there are some constants
$A,B > 0$ such that
all coefficients $a$, $b$, and\/ $c$ and the partial derivative
${\partial a}/{\partial x}$
are bounded, continuously differentiable functions in
the Cartesian product\/
\begin{math}
  \overline{\mathfrak{X}}^{(A)} \times
  \overline{\mathfrak{T}}_{0,T}^{(B)} ,
\end{math}
with $\RE a\geq \mathrm{const} > 0$, and all\/
$a$, $b$, and\/ $c$ are holomorphic in
$\mathfrak{X}^{(A)} \times \mathfrak{T}_{0,T}^{(B)}$.
Furthermore, let us assume that
the first\--order time derivatives of all functions
$a$, $b$, $c$, and\/
${\partial a}/{\partial x}$ are bounded in
\begin{math}
  \overline{\mathfrak{X}}^{(A)} \times
  \overline{\mathfrak{T}}_{0,T}^{(B)} .
\end{math}
Finally, assume that\/
$f$ is holomorphic in
$\mathfrak{X}^{(A)} \times \mathfrak{T}_{0,T}^{(B)} \times \CC^2$,
$f = f(x,t;u,\eta)$ where $\eta = {\partial u}/{\partial x}$,
with all functions
$f$, ${\partial f}/{\partial t}$, ${\partial f}/{\partial u}$, and\/
${\partial f}/{\partial\eta}$ being locally bounded in
$\mathfrak{X}^{(A)} \times \mathfrak{T}_{0,T}^{(B)} \times \CC^2$,
and it satisfies
\begin{equation*}
  \int_{-\infty}^{\infty} |f(x + \ii y,t;0,0)|^p \,\mathrm{d}x
  \leq K^p \hspace{05pt} (K = \mathrm{const} < \infty)
  \quad\mbox{ for all $y\in [-A,A]$ and
              $t\in \overline{\mathfrak{T}}_{0,T}^{(B)}$. }
\end{equation*}

{\rm (i)}$\;$
Then, given any\/ $u_0\in B^{s;p,p}(\RR)$,
the Cauchy problem \eqref{e:Cauchy:M=N=1}
possesses a unique weak solution
$u\in C\left( [0,T_0]\to B^{s;p,p}(\RR)\right)$
defined on a (possibly shorter) time interval
$[0,T_0]\subset [0,T]$ of some positive length $T_0\in (0,T]$.

{\rm (ii)}$\;$
Furthermore, if\/ $u_0\colon \RR^1\to \CC$
possesses a (unique) holomorphic extension to a complex strip
$\mathfrak{X}^{(r_0)}\subset \CC$, $0 < r_0\leq A$, denoted by\/
$u_0\colon \mathfrak{X}^{(r_0)}\to \CC$ again, such that\/
\begin{equation}
\nonumber
  \mathfrak{N}^{(r_0)} (u_0) \eqdef
    \sup_{ y\in [-r_0, r_0] } \| u_0(\,\cdot\, + \ii y)
    \|_{ B^{s;p,p}(\RR) } < \infty
\end{equation}
holds {\em (cf.\ ineq.~\eqref{sup_y:u_0} below)\/},
then also any (global) weak solution
$u\in C\left( [0,T]\to B^{s;p,p}(\RR)\right)$
can be (uniquely) extended to a {\em holomorphic\/} function in
$\mathfrak{X}^{(A')} \times \Delta_{\vartheta}^{T',T}$,
denoted again by $u(x + \ii y,t)$, where all numbers\/
$A'\in (0,A]$, $T'\in (0,T]$, and\/
$\vartheta\in (0, \pi / 2)$ are sufficiently small,
$T'\cdot \tan\vartheta\leq B$, and\/
\begin{math}
  u(\,\cdot\, + \ii y,\,\cdot\,)\colon
  t \,\longmapsto\, u(\,\cdot\, + \ii y,t) \colon
  \overline{\Delta}_{\vartheta}^{T',T} \,\longrightarrow\, B^{s;p,p}(\RR)
\end{math}
is continuous for every fixed $y\in [-r_0, r_0]$ together with\/
\begin{equation*}
  \sup_{ t\in \Delta_{\vartheta}^{T',T}}
  \mathfrak{N}^{(A')} ( u(\,\cdot\, + \ii y,t) )
  = \sup_{ (y,t)\in [-r_0, r_0]\times \Delta_{\vartheta}^{T',T} }
    \| u(\,\cdot\, + \ii y,t) \|_{ B^{s;p,p}(\RR^N) } < \infty \,.
\end{equation*}
In particular, the extension $u$ is {\em holomorphic\/} in
$\mathfrak{X}^{(A')} \times \mathfrak{T}_{T',T}^{(B')}$
with\/
$B'= T'\cdot \tan\vartheta\leq B$,
where $T'> 0$ and\/ $\vartheta > 0$ are small enough.
\end{theorem}
\par\vskip 10pt

We remark that
\begin{math}
  \mathfrak{T}_{T',T}^{(B')} \subset \Delta_{\vartheta}^{T',T}
  \subset \mathfrak{T}_{0,T}^{(B)} ,
\end{math}
by $B'= T'\cdot \tan\vartheta\leq B$.
If $0 < T_0 < T$ in Part~{\rm (i)} of this theorem,
then we have to replace $T$ by $T = T_0$ in Part~{\rm (ii)}.
Notice that the condition that both (continuous) partial derivatives
${\partial f}/{\partial u}$ and
${\partial f}/{\partial\eta}$ are locally bounded in
$\mathfrak{X}^{(A)} \times \mathfrak{T}_{0,T}^{(B)} \times \CC^2$
is equivalent with
$(u,\eta)\mapsto f(x,t;u,\eta)$
being locally Lipschitz\--continuous.

\par\vskip 10pt
\begin{remark}\label{rem-weak-classic}\nopagebreak
\begingroup\rm
It follows easily from Theorem~\ref{thm-nontech}
that every {\em weak\/} solution
$u\in C\left( [0,T]\to B^{s;p,p}(\RR)\right)$
to the Cauchy problem \eqref{e:Cauchy:M=N=1}
is {\em classical\/} in the sense that it is of class $C^{\infty}$
over the open set $\RR\times (0,T)$ and verifies
{\rm eq.}~\eqref{e:Cauchy:M=N=1} pointwise and the initial condition
$u(\,\cdot\,,0) = u_0\in B^{s;p,p}(\RR)$
in the $B^{s;p,p}(\RR)$-limit, i.e.,
\begin{math}
  \left\| u(\,\cdot\,,t) - u_0\right\|_{ B^{s;p,p}(\RR) } \to 0
\end{math}
as $t\to 0+$.
For the special case of the reaction function
$f\colon (u,\eta)\mapsto f(x,t;u,\eta)$
being linear, a {\em weak\/} solution
$u\in C\left( [0,T]\to L^2(\RR)\right)$
to the (linear) Cauchy problem \eqref{e:Cauchy:M=N=1}
is defined e.g.\ in
 {\sc L.~C.\ Evans} \cite{Evans-98}, Chapt.~7, {\S}1.1, p.~352, or
 {\sc J.-L.\ Lions} \cite{Lions-61}, Chapt.~IV, {\S}1, p.~44, or
 \cite{Lions-71}, Chapt.~III, eq.~(1.11), p.~102.
In that case the initial condition
$u(\,\cdot\,,0) = u_0\in L^2(\RR)$
holds in the sense of the $L^2(\RR)$-limit,
\begin{math}
  \left\| u(\,\cdot\,,t) - u_0\right\|_{ L^2(\RR) } \to 0
\end{math}
as $t\to 0+$.
The main reason why we prefer to work with the notion of
a {\em weak solution\/} as opposed to a {\em classical solution\/}
of the Cauchy problem \eqref{e:Cauchy:M=N=1}
is the fact that already a weak solution is {\em unique\/}.
The uniqueness of a weak solution is an important technical argument
in our proofs of Theorem~\ref{thm-nontech} and
Theorem~\ref{thm-Main} (Section~\ref{s:Main}).

In fact, we work sometimes also with
the so\--called {\em mild solutions\/}
to the Cauchy problem \eqref{e:Cauchy:M=N=1} that make sense in
$C\left( [0,T]\to L^2(\RR)\right)$; cf.\
{\sc P.\ Tak\'a\v{c}} \cite[Sect.\ 3 and~4]{Takac},
even though we use them in the Besov space
$B^{s;p,p}(\RR)$ in place of $L^2(\RR)$.
Mild solutions do not require any additional regularity knowledge
as they are defined by
the well\--known variation\--of\--constants formula
({\sc A.\ Pazy} \cite[{\S}5.7, p.~168]{Pazy}).
Thus, they are even ``weaker'' than the weak solutions, but
in our situation one can easily verify that every mild solution is also
a weak solution to problem \eqref{e:Cauchy:M=N=1} and vice versa;
see e.g.\
{\sc J.~M.\ Ball} \cite{J.M.Ball} (or \cite[Theorem on p.~259]{Pazy}).

The same remarks apply also to the more general Cauchy problem
\eqref{e:Cauchy}.
\hfill\Square
\endgroup
\end{remark}
\par\vskip 10pt

This article is organized as follows.
We introduce some basic notation (mostly complex domains)
in Section~\ref{s:Notation}.
Our main analyticity result, Theorem~\ref{thm-Main},
supplemented by an additional explanation in Proposition~\ref{prop-Main},
is stated in Section~\ref{s:Main}.
Their proofs are gradually built up in
Sections \ref{s:Cauchy_abstr} through~\ref{s:proof-Main}:
First, the Cauchy problem in $\RR^N\times (0,T)$
is treated as an abstract initial value problem
in Section~\ref{s:Cauchy_abstr}.
There, an important abstract a~priori $B^{s;p,p}$\--type estimate
is established in Theorem~\ref{thm-Clem-Li_glob}.
Analyticity in time for this abstract problem is proved
in Section~\ref{s:Analyt_abstr} (Theorem~\ref{thm-Analytic_glob}).
Then, in Section~\ref{s:Analyt_space},
we treat analyticity in space for
the semilinear parabolic Cauchy problem in $\RR^N\times (0,T)$, see
Proposition~\ref{prop-smooth_ext},
provided the initial data are already analytic (in space).
We show that the analyticity in space is preserved for all times in $[0,T]$.
We combine the time and space analyticity results from
Sections \ref{s:Analyt_abstr} and~\ref{s:Analyt_space}
into Theorem~\ref{thm-smooth_ext}
in Section~\ref{s:Analyt_(x,t)}.
This theorem is still only ``local'' in time.
Our main analyticity result, Theorem~\ref{thm-Main},
is proved in Section~\ref{s:proof-Main},
together with Proposition~\ref{prop-Main} and, in particular,
the ``regularity'' estimates
\eqref{eq:du/dt_Hardy^2} and~\eqref{eq:u_Hardy^2}.
Section~\ref{s:Appl} treats an application to
a {\em Risk Model in Mathematical Finance\/}.
Finally, Section~\ref{s:Hist} contains
some historical remarks and comments concerning the analyticity of
solutions to linear and semi\-linear elliptic and parabolic systems
and its applications to relevant classical problems.


\section{Notation}
\label{s:Notation}

We stick to the classical notation 
$\NN = \{ 1, 2, 3,\dots\}$,
$\ZZ_+ = \{ 0, 1, 2, 3,\dots\} = \NN\cup \{ 0\}$, and
$\ZZ = \{ 0,\pm 1,\pm 2,\pm 3,\dots\} = \ZZ_+\cup (-\ZZ_+)$,
together with
$\RR = (-\infty,\infty)$, $\RR_+ = [0,\infty)$, and
$\CC = \RR + \ii\RR\cong \RR^2$.
Typically, we denote by
$x = (x_1, x_2, \dots, x_N)$ and $y = (y_1, y_2, \dots, y_N)$
points in $\RR^N$ and by
$z = (z_1, z_2, \dots, z_N)$ points in $\CC^N$.
We often write $\zeta = \xi + \ii\eta$ for
$\zeta\in \CC$ and $\xi,\eta\in \RR$, i.e.,
$\RE\zeta = \xi$ and $\IM\zeta = \eta$
are the real and imaginary parts of $\zeta\in \CC$, respectively.
Similarly,
$z = x + \ii y$ for $z\in \CC^N$ and $x,y\in \RR^N$,
or equivalently $z_i = x_i + \ii y_i$ ($i=1,2,\dots,N$)
for $z_i\in \CC$ and $x_i, y_i\in \RR$, i.e.,
$\RE z = x$ and $\IM z = y$.
Hence, we identify
$\CC^N = \RR^N\oplus \ii\RR^N$
(or simply $\CC^N = \RR^N + \ii\RR^N$)
as vector spaces over the field $\RR$ and thus consider
$\RR^N$ to be a (vector) subspace of $\CC^N$.
We use a bar $(\bar{\phantom{z}})$ to denote
the complex conjugate $\bar{\zeta}$ of a number $\zeta\in \CC$.
The complex conjugate of a vector $z\in \CC^N$ is denoted by
$\bar{z} = (\bar{z}_1, \bar{z}_2, \dots, \bar{z}_N)$.
Similarly,
the complex conjugate function of a complex\--valued function $f(z)$
(for $f\colon \CC^N\to \CC$, for instance)
is denoted by $\bar{f}(z)\equiv \overline{f(z)}$.
Furthermore, we denote by
$(z,w) = \sum_{i=1}^N z_i\bar{w}_i$
the standard Euclidean inner product of $z,w\in \CC^N$ and by
$|z| = \left( \sum_{i=1}^N |z_i|^2\right)^{1/2}$
the induced (Euclidean) norm of $z\in \CC^N$.
We will often use the sum ($\ell^1$-) and
the maximum ($\ell^{\infty}$-) norms of $z\in \CC^N$, respectively:
\begin{eqnarray*}
& |z|_1 = \sum_{i=1}^N |z_i| \quad\mbox{ and }\quad
  |z|_{\infty} = \max_{1\leq i\leq N} |z_i| \,.
\end{eqnarray*}
Finally, we write
$z\cdot w = \sum_{i=1}^N z_i w_i$
for the {\em bilinear\/} product of $z,w\in \CC^N$,
which is not to be confused with the inner product
$(z,w) = \sum_{i=1}^N z_i\bar{w}_i$ if $w\not\in \RR^N$
(which is {\em sesquilinear\/}).
The Euclidean ($\ell^2$-) norm of $z\in \CC^N$ is abbreviated as
\begin{math}
  |z|\equiv |z|_2 = \sqrt{(z,z)}
  = \left( \sum_{i=1}^N |z_i|^2\right)^{1/2} .
\end{math}

The vector space (over the field $\RR$) of
all real\--valued (square) $M\times M$ matrices
$\mathbf{A} = ( a_{ij} )_{i,j=1}^M$ is denoted by $\RR^{M\times M}$.
Similarly, the vector space (over the field $\CC$) of
all complex\--valued $M\times M$ matrices
$\mathbf{A}$ is denoted by $\CC^{M\times M}$.

Given $r\in (0,\infty)$, we denote by
\begin{math}
  Q^{(r)} = (-r,r)^N = \{ y\in \RR^N\colon |y|_{\infty} < r\}
\end{math}
the $N$-dimensional open cube in $\RR^N$ (centered at the origin)
with side lengths $2r$, and by
$\overline{Q}^{(r)} = [-r,r]^N$ its closure.

In order to formulate our main hypotheses,
given $r,T\in (0,\infty)$ and $T'\in [0,T)$,
we introduce the following complex domains for
the complexifications of the spatial and temporal variables,
$x\in \RR^N$ and $t\in (0,T)$, respectively:
\begin{align}
\label{def:X^(r)}
  \mathfrak{X}^{(r)} &\eqdef
  \{ z = x + \ii y\in \CC^N\colon |y|_{\infty} < r\}
  = \RR^N + \ii Q^{(r)} \,,
\\
\nonumber
  \mathfrak{T}_{T',T}^{(r)} &\eqdef
  \{ t = \sigma + \ii\tau\in \CC\colon
     T'< \sigma < T\mbox{ and } |\tau| < r\} \,;
  \quad\mbox{ see \eqref{def:T^(r)}. }
\end{align}
The former, $\mathfrak{X}^{(r)}$,
is a {\em tube\/} (often called a {\em strip\/})
in $\CC^N$ with {\em base\/} $Q^{(r)}$
and the latter, $\mathfrak{T}_{T',T}^{(r)}$, is an open rectangle
in the complex plane~$\CC$.
Notice that $\mathfrak{T}_{T',T}^{(r)}$ contains the interval $(T',T)$.
The remaining temporal domains,
$\Delta_{\vartheta}$ and $\Delta_{\vartheta}^{(T')}$,
have been introduced in
{\rm eqs.}\ \eqref{loc:T^(infty)} and \eqref{loc:T^(T')}, respectively,
with the angle $\vartheta\in (0, \pi / 2)$ given by
$\tan\vartheta = r/T'$.

Our techniques will use holomorphic semigroups in an open {\em sector\/}
$\Delta_{\vartheta}\subset \mathbb{C}$
defined in \eqref{loc:T^(infty)},
with a given angle $\vartheta\in (0, \pi / 2)$,
but often locally in time in an open {\em triangle\/}
$\Delta_{\vartheta}^{(T')}\subset \mathbb{C}$
defined in \eqref{loc:T^(T')}, where $0 < T'< \infty$.
Their respective closures in $\CC$ are denoted by
$\overline{\Delta}_{\vartheta}$ and
$\overline{\Delta}_{\vartheta}^{(T')}$;
both contain the origin $0\in \CC$.
Finally, for $0 < T'\leq T < \infty$ we recall the definition of
the temporal domain
$\Delta_{\vartheta}^{T',T}$ introduced in {\rm eq.}~\eqref{glob:T^(r)}
with the angle $\vartheta\in (0, \pi / 2)$ given by
$\tan\vartheta = r/T'$.
Its closure in $\CC$ is denoted by
$\overline{\Delta}_{\vartheta}^{T',T}$.
Clearly,
\begin{math}
    \xi + \Delta_{\vartheta}^{(T')}
  = \{ \xi + t'\in \CC\colon t'\in \Delta_{\vartheta}^{(T')} \} .
\end{math}
Of course, if $T = T'$ then
$\Delta_{\vartheta}^{T',T} = \Delta_{\vartheta}^{(T')}$
is an open triangle.

Throughout this article we work with complex\--valued functions; hence,
all Banach and Hilbert spaces of functions we consider are complex
(over the field $\CC$).
We work with the standard inner product in $L^2(\RR^N)$ defined by
$(u,v)_{L^2}\eqdef \int_{\RR^N} u\bar{v} \,\mathrm{d}x$
for $u,v\in L^2(\RR^N)$.
The induced norm is abbreviated by
$\| u\|_{L^2}\equiv \| u\|_{ L^2(\RR^N) }$.
We warn the reader that we identify the dual space
$\mathscr{H}^{\prime} = \left( L^2(\RR^N)\right)^{\prime}$
of the complex Hilbert space $\mathscr{H} = L^2(\RR^N)$
with $\mathscr{H}$ itself by means of
the (complex) Riesz representation theorem which yields
an {\em anti\--linear\/} isomorphism of $\mathscr{H}$ onto
$\mathscr{H}^{\prime}$
(cf.\ 
 {\sc R.~A.\ Adams} and {\sc J.~J.~F.\ Fournier}
 \cite[Chapt.~2]{AdamsFournier}, Theorem 2.44, p.~47).

The following notation is taken from
{\sc S.~G.\ Krantz} \cite[Chapt.~0]{Krantz}.
Given a domain $\Omega$ in $\mathbb{R}^r$, we denote by
$C^k(\Omega)$ ($k\in \ZZ_+$)
the vector space of all $k$-times continuously differentiable functions
$f\colon \Omega\to \CC$ and by
$C^k(\overline{\Omega})$ the vector space of all
$f\colon \overline{\Omega}\to \CC$ such that
$f|_{\Omega}\in C^k(\Omega)$ and each partial derivative
\begin{math}
  \frac{ \partial^{|\alpha|} f }{\partial x^{\alpha}}
\end{math}
of $f$ ($\alpha\in (\ZZ_+)^r$) of order $|\alpha|\leq k$
can be extended to a continuous function on~$\overline{\Omega}$.
Of course, $f|_{\Omega}$ stands for the restriction of $f$ to $\Omega$
and all partial derivatives are taken in the real variable sense
($x\in \RR^r$).
In case
$\Omega\subset \CC^r = \RR^r\oplus \ii\RR^r\cong \RR^{2r}$ ($r\in \NN$)
is a complex domain, the (mixed) partial derivative
\begin{math}
  \frac{ \partial^{|\alpha| + |\beta|} f }%
       { \partial x^{\alpha}\, \partial y^{\beta} }
\end{math}
of $f$ ($\alpha, \beta\in (\ZZ_+)^r$) of order
$|\alpha| + |\beta|\leq k$
is taken in the real variable sense, where
$x,y\in \RR^r$ in $(x,y)\cong z = x + \ii y\in \CC^r$.
The vector spaces $C^k(\Omega)$ and $C^k(\overline{\Omega})$
are defined analogously, with $|\alpha| + |\beta|\leq k$.
If $\Omega$ is bounded then
$C^k(\overline{\Omega})$ is a Banach space endowed with
a maximum\--type norm.
If $\Omega$ is not bounded in general then we denote by
$C^0_{\mathrm{unif}}(\overline{\Omega})$
the vector space of all uniformly continuous functions
$f\colon \overline{\Omega}\to \CC$, and by
$C^0_{\mathrm{bdd}}(\Omega) = C^0(\Omega)\cap L^{\infty}(\Omega)$
the vector space of all bounded continuous functions
$f\colon \Omega\to \CC$.
Alternatively, if $\Omega$ is not bounded
(in $\RR^r$ or $\CC^r$ endowed with the Euclidean metric $d$)
then we may consider
a {\em compactification\/} $\widetilde{\Omega}$ of $\Omega$, i.e.,
a compact metric space $\widetilde{\Omega}$
(endowed with a metric $\widetilde{d}$)
such that there is a homeomorphism
$j\colon \Omega\to \widetilde{\Omega}$ of $\Omega$
onto a dense subset $j(\Omega)$ of $\widetilde{\Omega}$,
such that for any pair of sequences
$\{ x_n\}_{n=1}^{\infty}$, $\{ y_n\}_{n=1}^{\infty}$
$\subset \Omega$ we have
$d(x_n,y_n)\to 0$ if and only if
$\widetilde{d}(j(x_n), j(y_n)) \to 0$ as $n\to \infty$.
Clearly, by identifying $\Omega$ with the subset
$j(\Omega)$ of $\widetilde{\Omega}$ we identify
$\Omega$ with a dense subset of $\widetilde{\Omega}$.
Hence, we can identify
$C^0(\widetilde{\Omega})$ with a vector subspace of
\begin{math}
  C^0_{\mathrm{bdd}}(\Omega) \cap
  C^0_{\mathrm{unif}}(\overline{\Omega}) .
\end{math}
As a simple example, we may take
$\widetilde{\Omega}$ to be the one\--point compactification of
a domain $\Omega\subset \RR^r$ (or $\Omega\subset \CC^r$).
In particular, we have
$\widetilde{\mathbb{R}}^r = \mathbb{R}^r\cup \{ \infty\}$
and
$\widetilde{\mathbb{C}}^r = \mathbb{C}^r\cup \{ \infty\}$
endowed with the metric $\widetilde{d}$ defined in the next section
(Section~\ref{s:Main}, eq.~\eqref{e:d^tilde}).

Finally, if $\Omega\subset \CC^r$ is a complex domain, we denote by
$A(\Omega)$ the Fr\'echet space of all holomorphic functions
$f\colon \Omega\to \CC$
endowed with the (complete metrizable) topology of uniform convergence
on compact subsets of~$\Omega$.
As usual, we abbreviate
the {\em Cauchy\--Riemann\/} partial differential operators
\begin{equation}
\label{def:Cauchy-Riem}
  \frac{\partial}{\partial z_i} \eqdef
  \frac{1}{2}
    \left( \frac{\partial}{\partial x_i}
    - \ii\,\frac{\partial}{\partial y_i}
    \right)
  \quad\mbox{ and }\quad
  \frac{\partial}{\partial\bar{z}_i} \eqdef
  \frac{1}{2}
    \left( \frac{\partial}{\partial x_i}
    + \ii\,\frac{\partial}{\partial y_i}
    \right) \,.
\end{equation}
%

\begin{remark}\label{rem-Cauchy-Riem}\nopagebreak
\begingroup\rm
We will often use the following classical fact; see e.g.\
{\sc F.\ John} \cite[Theorem, p.~70]{John}
or
{\sc S.~G.\ Krantz} \cite[Definition~II, p.~3]{Krantz}:
Let $\Omega\subset \CC^N$ be a complex domain ($N\geq 1$).
A continuously differentiable function $h\colon \Omega\to \CC$
(in the real variable sense, partially with respect to
 $x_i, y_i\in \RR$; $i=1,2,\dots,N$)
is holomorphic (i.e., complex analytic) {\em if and only if\/}
it verifies the {\em Cauchy\--Riemann\/} equations in $\Omega$, i.e.,
\begin{math}
  \;
  {\partial h}/{\partial\bar{z}_i} = 0
  \,\mbox{ in $\Omega$; for all }\, i=1,2,\dots,N .
\end{math}
\hfill\Square
\endgroup
\end{remark}
\par\vskip 10pt

\section{Statement of the main result}
\label{s:Main}

Let us abbreviate the differential operators (and derivatives)
\begin{align*}
& D_x\eqdef \frac{1}{\ii} \frac{\partial}{\partial x}
  = \left( \frac{1}{\ii} \frac{\partial}{\partial x_i}
    \right)_{i=1}^N \,,\quad
  \partial_t\eqdef
    \frac{\partial}{\partial t} \,,
    \quad\mbox{ and }
\\
& D_x^{\alpha}\eqdef \ii^{- |\alpha|}\,
    \frac{ \partial^{|\alpha|} }{ \partial x^{\alpha} }
  = \ii^{- |\alpha|}\,
    \frac{ \partial^{|\alpha|} }%
         { \partial x_1^{\alpha_1} \,\dots\, \partial x_N^{\alpha_N} }
  \quad\mbox{ for } \alpha = (\alpha_i)_{i=1}^N\in (\ZZ_+)^N \,.
\end{align*}
We assume that the operator
\begin{equation}
\label{def:P}
\begin{aligned}
&   \mathbf{P}(x,t,D_x)
  = \textstyle\sum_{|\alpha|, |\beta|\leq m} D_x^{\alpha}
    \left( \mathbf{P}^{\alpha\beta} (x,t)\, D_x^{\beta} \right)
\\
& \equiv
    \sum_{|\alpha|, |\beta|\leq m} \ii^{- |\alpha| - |\beta|}\,
    \frac{ \partial^{|\alpha|} }{ \partial x^{\alpha} }
    \left( \mathbf{P}^{\alpha\beta} (x,t)\,
    \frac{ \partial^{|\beta|} }{ \partial x^{\beta} }
    \right) \,,
  \quad\mbox{ for }\; (x,t)\in \mathbb{R}^N\times (0,T) \,,
\end{aligned}
\end{equation}
is a linear partial differential operator of order $2m$
in {\em divergence form\/} with the coefficients
\hfil\break
$\ii^{- |\alpha| - |\beta|} \mathbf{P}^{\alpha\beta} (x,t)$
indexed by $\alpha, \beta\in (\ZZ_+)^N$ with
$|\alpha|\leq m$ and $|\beta|\leq m$, where each
\begin{math}
  \mathbf{P}^{\alpha\beta} (x,t) = ( P^{\alpha\beta}_{jk} )_{j,k=1}^M
\end{math}
is an $M\times M$ matrix with real (or complex) entries
$P^{\alpha\beta}_{jk} = P^{\alpha\beta}_{jk}(x,t)$.
The reader is referred to
{\sc A.\ Friedman} \cite[Part.~1, Sect.~12, pp.\ 32--37]{Friedman-69}
or
{\sc F.\ John} \cite[Chapt.~6, Sect.~2, pp.\ 190--195]{John}
for general facts about such operators.

Let us abbreviate the product domain
\begin{math}
  \Omega = \mathfrak{X}^{(r_0)}\times \Delta_{\vartheta_0}^{T_0,T}
  \subset \CC^N\times \CC ,
\end{math}
with some
$r_0\in (0,\infty)$, $0 < T_0\leq T < \infty$, and
$\vartheta_0\in (0, \pi / 2)$;
the closure of $\Omega$ in $\CC^N\times \CC$ is denoted by
$\overline{\Omega}$.
We introduce also a compactification $\widetilde{\Omega}$ of~$\Omega$,
\begin{equation*}
    \widetilde{\Omega}
  = \widetilde{ \mathfrak{X} }^{(r_0)}
  \times \overline{\Delta}_{\vartheta_0}^{T_0,T}
  = \left( \widetilde{\RR}^N + \ii \overline{Q}^{(r)} \right)
  \times \overline{\Delta}_{\vartheta_0}^{T_0,T}
  \cong \widetilde{\RR}^N\times \overline{Q}^{(r)}
  \times \overline{\Delta}_{\vartheta_0}^{T_0,T} \,,
\end{equation*}
where
$\widetilde{\RR}^N = \RR^N\cup \{ \infty\}$
denotes the one\--point compactification of $\RR^N$
endowed with the standard metric
\begin{equation}
\label{e:d^tilde}
  \widetilde{d}(x,y) = \widetilde{d}(y,x) =
\left\{
\begin{alignedat}{2}
  \frac{|x-y|}{1 + |x-y|}
& \quad\mbox{ if }\; x,y\in \RR^N \,;
\\
  1\qquad
& \quad\mbox{ if }\; x\in \RR^N \,,\ y = \infty \,;
\\
  0\qquad
& \quad\mbox{ if }\; x = y = \infty \,.
\end{alignedat}
\right.
\end{equation}
Hence,
%
\begin{math}
        \widetilde{ \mathfrak{X} }^{(r_0)}
      = \widetilde{\RR}^N + \ii \overline{Q}^{(r)}
  \cong \widetilde{\RR}^N\times \overline{Q}^{(r)}
\end{math}
%
is a compactification of
$\mathfrak{X}^{(r_0)} \subset \CC^N$.

We assume that the operator $\mathbf{P}$ and the function $\mathbf{f}$
satisfy the following hypotheses in the product domain
\begin{math}
  \Omega = \mathfrak{X}^{(r_0)}\times \Delta_{\vartheta_0}^{T_0,T}
\end{math}
defined in eqs.\ \eqref{glob:T^(r)} and \eqref{def:X^(r)}:

\begin{hypos}\nopagebreak
\begingroup\rm
\begin{enumerate}
\setcounter{enumi}{0}
\renewcommand{\labelenumi}{{\bf (H\arabic{enumi})}}
\item
\makeatletter
\def\@currentlabel{{\bf H\arabic{enumi}}}\label{hy:UnifEll}
\makeatother
For each pair $\alpha, \beta\in (\ZZ_+)^N$ with
$|\alpha|\leq m$ and $|\beta|\leq m$,
the entries
$P^{\alpha\beta}_{jk} \colon \overline{\Omega}\to \CC$
($j,k=1,2,\dots,M$)
of the coefficient
$\mathbf{P}^{\alpha\beta} = ( P^{\alpha\beta}_{jk} )_{j,k=1}^M$
belong to
$C^1(\overline{\Omega})\cap L^{\infty}(\Omega)\cap A(\Omega)$.
Moreover, we assume that also all partial derivatives
\begin{math}
  \frac{ \partial^{|\alpha'|} }{ \partial x^{\alpha'} }
  P^{\alpha\beta}_{jk} (x,t)
\end{math}
of order $|\alpha'|\leq |\alpha|$ ($\alpha'\in (\ZZ_+)^N$)
are in $C^1(\overline{\Omega})$.
The entries $P^{\alpha\beta}_{jk}$ of the {\em leading coefficients\/}
($|\alpha| = |\beta| = m$)
are assumed to belong also to
\begin{math}
  C^0_{\mathrm{unif}}(\overline{\Omega})
\end{math}
besides being in
$C^1(\overline{\Omega})\cap L^{\infty}(\Omega)\cap A(\Omega)$.

\begin{quote}
This is the case if
the entries $P^{\alpha\beta}_{jk}$ of the leading coefficients
($|\alpha| = |\beta| = m$)
belong also to
\begin{math}
  C^0(\widetilde{\Omega})
\end{math}
besides being in
$C^1(\overline{\Omega})\cap L^{\infty}(\Omega)\cap A(\Omega)$.
This claim follows directly from
\begin{math}
  C^0(\widetilde{\Omega}) \subset
  C^0_{\mathrm{unif}}(\overline{\Omega})
\end{math}
which means that any continuous function
$f\colon \widetilde{\Omega}\to \CC$
is uniformly continuous on $\widetilde{\Omega}$
with the restriction $f|_{\Omega}$ to $\Omega$
being uniformly continuous and, thus,
$f|_{\Omega}\colon \Omega\to \CC$
possesses a (unique) continuous extension
$\widetilde{f}\colon \overline{\Omega}\to \CC$
to $\overline{\Omega}$ which turns out to be
uniformly continuous, as well.
\end{quote}
\item
\makeatletter
\def\@currentlabel{{\bf H\arabic{enumi}}}\label{hy:Ellipt}
\makeatother
The operator $\mathbf{P}$ is
{\em strongly elliptic\/} in $\overline{\Omega}$, i.e.,
there exists a constant $c\in (0,\infty)$ such that the inequality
\begin{equation}
\label{e:Ellipt}
  \RE \biggl( \sum_{j,k=1}^M \sum_{|\alpha| = |\beta| = m}
    P^{\alpha\beta}_{jk}(z,t)\,
    \xi^{\alpha + \beta}\, \eta_k\bar{\eta}_j
      \biggr)
  \geq c\, |\xi |^{2m}\, |\boldsymbol{\eta}|^2
\end{equation}
holds for all $(z,t)\in \overline{\Omega}$ and for all
$\xi = (\xi_1,\dots,\xi_N)\in \RR^N$
and
$\boldsymbol{\eta} = (\eta_1,\dots,\eta_M) \in \CC^M$,
where
\begin{math}
  \xi^{\alpha + \beta} =
  \xi_1^{\alpha_1 + \beta_1} \,\dots\, \xi_N^{\alpha_N + \beta_N}
\end{math}
and
$\alpha = (\alpha_1,\dots,\alpha_N) \in (\ZZ_+)^N$,
$\beta = (\beta_1,\dots,\beta_N) \in (\ZZ_+)^N$.
\item
\makeatletter
\def\@currentlabel{{\bf H\arabic{enumi}}}\label{hy:f_Hardy^2}
\makeatother
The components
\begin{math}
  f_j\colon \overline{\Omega}\times \CC^{M\tilde{N}}\to \CC
\end{math}
($j=1,2,\dots,M$)
of the reaction function
$\mathbf{f} = (f_1,\dots,f_M)$ belong to
\begin{equation*}
  C^1( \overline{\Omega}\times \CC^{M\tilde{N}} ) \cap
  A( \Omega\times \CC^{M\tilde{N}} ) \,.
\end{equation*}
(Recall that the integer $\tilde{N}$ is defined in
 {\rm eq.}~\eqref{e:N^tilde}.)
Moreover, we assume that, for every bounded subset
$\Sigma\subset \CC^{M\tilde{N}}$,
their first\--order time derivatives
\begin{math}
  \frac{\partial}{\partial t} f_j(x,t;X)
\end{math}
together with their first\--order partial derivatives
\begin{equation*}
  \frac{\partial}{ \partial X_{\beta,k} }\, f_j(x,t;X) \,,
    \quad\mbox{ for $|\beta|\leq m$ and $j,k=1,2,\dots,M$, }
\end{equation*}
with respect to the components $X_{\beta,k}$ of the vector
\begin{math}
    X_{\beta}
  = ( X_{\beta,1}, \dots, X_{\beta,M} )
  = \frac{ \partial^{|\beta|}\mathbf{u} }{\partial x^{\beta}}
  \in \CC^M
\end{math}
(or $\CC^M$)
are {\em uniformly bounded\/} on the set $\Omega\times \Sigma$.
Finally, we assume that the function\/
\begin{math}
  \mathbf{f}\colon
  \overline{\Omega}\times \CC^{M\tilde{N}}\to \CC^M
\end{math}
satisfies
\begin{equation}
\label{e:f_Hardy^2}
  \int_{\RR^N} |\mathbf{f} (x + \ii y,t; \vec{\mathbf{0}})|^p
  \,\mathrm{d}x
  \leq K^p
  \quad\mbox{ for all $y\in \overline{Q}^{(r_0)}$ and
              $t\in \overline{\Delta}_{\vartheta_0}^{T_0,T}$, }
\end{equation}
where $K\in (0,\infty)$ is a constant and
\begin{math}
  \vec{\mathbf{0}}\eqdef \left( 0\right)_{|\beta|\leq m}
  \equiv (0,\dots,0)\in \CC^{M\tilde{N}} .
\end{math}
\end{enumerate}
%
\endgroup
\end{hypos}
\par\vskip 10pt

\begin{remark}\label{rem-f_Hardy^2}\nopagebreak
\begingroup\rm
The local boundedness condition in Hypothesis~\eqref{hy:f_Hardy^2}
on the first\--order partial derivatives
$\;$
\begin{math}
  \frac{\partial}{\partial t} f_j(x,t;X)
\end{math}
and
\begin{math}
  \frac{\partial f_j}{ \partial X_{\beta,k} } (x,t;X)
\end{math}
$\;$
on the set $\Omega\times \Sigma$
will be used later
(cf.\ ineq.~\eqref{e:|df/dX|<C} in Section~\ref{s:Analyt_space})
in the following equivalent form:

Given a bounded subset
$\Sigma\subset \CC^{M\tilde{N}}$ and indices $j=1,2,\dots,M$,
there is a constant
$C_j\equiv C_j(\Sigma) \in (0,\infty)$
such that the following inequalities,
\begin{align}
\label{e:f_j'}
  \left| \frac{\partial f_j}{\partial t} (x,t;X)\right|
  \leq C_j(\Sigma)
  \qquad\mbox{ and }\qquad
  \left|
    \frac{\partial f_j}{ \partial X_{\beta,k} } (x,t;X)
  \right| \leq C_j(\Sigma)
\\
\nonumber
  \quad\mbox{ for all $(x,t)\in \Omega$ and for all }\,
  X = \left( X_{\beta}\right)_{|\beta|\leq m} \in \Sigma \,,
\end{align}
hold for all
$\beta = (\beta_1,\dots,\beta_N) \in (\ZZ_+)^N$ with
$|\beta|\leq m$ and for all $k=1,2,\dots,M$.
\hfill\Square
\endgroup
\end{remark}
\par\vskip 10pt

A simple, but more restrictive alternative to formulate
{\rm Hypotheses} \eqref{hy:UnifEll}, \eqref{hy:Ellipt}, and\/
\eqref{hy:f_Hardy^2}
is to replace
\begin{math}
  \Omega = \mathfrak{X}^{(r_0)}\times \Delta_{\vartheta_0}^{T_0,T}
\end{math}
by a larger, but simpler product domain
\begin{math}
  \Omega_0 = \mathfrak{X}^{(r_0)}\times \mathfrak{T}_{0,T}^{(\tau_0)}
\end{math}
(defined in eqs.\ \eqref{def:T^(r)} and \eqref{def:X^(r)})
with $\tau_0 = T_0\cdot \tan\vartheta_0$; hence,
$\Omega\subset \Omega_0$, thanks to
\begin{math}
  \Delta_{\vartheta_0}^{T_0,T} \subset \mathfrak{T}_{0,T}^{(\tau_0)} .
\end{math}

Let us recall our abbreviation
\begin{math}
  X = \left( X_{\beta}\right)_{|\beta|\leq m} \in \RR^{M\tilde{N}}
\end{math}
(or $\CC^{M\tilde{N}}$)
with
$\tilde{N} = \sum_{k=0}^m N^k$ from eq.~\eqref{e:N^tilde}
and make it more precise as follows:
When dealing with {\it complex\/} (partial) derivatives of the function
$f_j$ with respect to the variable $X_{\beta,k}$,
we prefer to replace $X_{\beta,k}$ by the complex variable
$Z_{\beta,k} = X_{\beta,k} + \ii Y_{\beta,k}\in \CC$
($X_{\beta,k}, Y_{\beta,k}\in \RR$);
$j,k=1,2,\dots,M$, and write
\begin{math}
  \frac{\partial f_j}{ \partial Z_{\beta,k} } (x,t;Z)
\end{math}
in place of
\begin{math}
  \frac{\partial f_j}{ \partial X_{\beta,k} } (x,t;X) .
\end{math}

The {\em strong ellipticity\/} inequality \eqref{e:Ellipt}
can be improved as follows; cf.\
{\sc P.\ Tak\'a\v{c}} \cite[Remark 3.1, p.~57]{Takac}:

\begin{remark}\label{rem-hypos}\nopagebreak
\begingroup\rm
In a smaller domain
\begin{math}
  \Omega'= \mathfrak{X}^{(r_0)}\times \Delta_{\vartheta_0'}^{T_0,T}
  \subset \Omega ,
\end{math}
with some number
$\vartheta_0'\in (0,\vartheta_0]$ (small enough),
inequality \eqref{e:Ellipt}
holds in the following qualitatively stronger form, cf.\
{\sc S.\ Agmon} \cite{Agmon}, Theorem 7.12, ineq.\ (7.21) on p.~87:
\begin{equation}
\nonumber
\tag{\ref{e:Ellipt}$'$}
  \RE \biggl( \ee^{\ii\theta}\cdot
    \sum_{j,k=1}^M \sum_{|\alpha| = |\beta| = m}
    P^{\alpha\beta}_{jk}(z,t)\,
    \xi^{\alpha + \beta}\, \eta_k\bar{\eta}_j
      \biggr)
  \geq c'\, |\xi |^{2m}\, |\boldsymbol{\eta}|^2
\end{equation}
for all
$\theta\in [-\vartheta_0', \vartheta_0']$,
for all $(z,t)\in \overline{\Omega}'$, and for all
$\xi = (\xi_1,\dots,\xi_N)\in \RR^N$
and
$\boldsymbol{\eta} = (\eta_1,\dots,\eta_M) \in \CC^M$,
where $c'\in (0,c]$ is a constant.
Recall that
$0 < \vartheta_0'\leq \vartheta_0$ and $\Omega'\subset \Omega$.
Consequently, without loss of generality, we may remove the prime $(')$
from both $\vartheta_0'$ and $c'$ in (\ref{e:Ellipt}$'$)
and assume that
\begin{equation}
\label{e:Ellipt'}
  \RE \biggl( \ee^{\ii\theta}\cdot
    \sum_{j,k=1}^M \sum_{|\alpha| = |\beta| = m}
    P^{\alpha\beta}_{jk}(z,t)\,
    \xi^{\alpha + \beta}\, \eta_k\bar{\eta}_j
      \biggr)
  \geq c\, |\xi |^{2m}\, |\boldsymbol{\eta}|^2
\end{equation}
for all
$\theta\in [-\vartheta_0, \vartheta_0]$
and for all $(z,t)\in \overline{\Omega}$,
$\xi = (\xi_1,\dots,\xi_N)\in \RR^N$,
and
$\boldsymbol{\eta} = (\eta_1,\dots,\eta_M) \in \CC^M$,
where $c>0$ is a constant.
We prefer to use
inequality \eqref{e:Ellipt'} in place of \eqref{e:Ellipt}.
\hfill\Square
\endgroup
\end{remark}
\par\vskip 10pt

The {\em G\r{a}rding inequality\/} (in the whole space $\RR^N$) below
is an important consequence of inequality \eqref{e:Ellipt'};
see e.g.\
{\sc S.\ Agmon} \cite[Theorem 7.6, p.~78]{Agmon}:

\begin{corollary}\label{cor-hypos}
{\rm (G\r{a}rding's inequality)}
Under both 
{\rm Hypotheses} \eqref{hy:UnifEll} and\/ \eqref{hy:Ellipt},
there exist some constants $c_1$ and\/ $c_2$,
$c_1 > 0$ and\/ $0\leq c_2 < \infty$, such that\/
\begin{equation}
\label{e:Garding}
\begin{aligned}
&   \RE
    \biggl[ \ee^{\ii\theta}\cdot
           \sum_{|\alpha| = |\beta| = m} \int_{\RR^N}
    \overline{ D_x^{\alpha} \mathbf{w} }
    \cdot  \mathbf{P}^{\alpha\beta} (x + \ii y,t)\, D_x^{\beta} \mathbf{w}
           \,\mathrm{d}x
    \biggr]
\\
& \geq c_1 \textstyle\sum_{|\alpha| = m}
    \| D_x^{\alpha} \mathbf{w}\|_{ L^2(\RR^N) }^2
  - c_2\, \| \mathbf{w}\|_{ L^2(\RR^N) }^2
\end{aligned}
\end{equation}
holds for all\/
$\mathbf{w}\in W^{m,2}(\RR^N)$ and for all\/
$\theta\in [-\vartheta_0, \vartheta_0]$,
$y\in \overline{Q}^{(r_0)}$, and\/
$t\in \overline{\Delta}_{\vartheta_0}^{T_0,T}$.
\end{corollary}
\par\vskip 10pt

\proof
The reader is referred to
{\sc S.\ Agmon} \cite[Theorem 7.6, pp.\ 78--86]{Agmon} for a proof.
We remark that the proof of G\r{a}rding's inequality
(\cite[Lemma 7.9, p.~81]{Agmon})
requires the {\em uniform equi\-continuity\/} of the leading coefficients
$\mathbf{P}^{\alpha\beta} (x + \ii y,t)$
as functions of $x\in \RR^N$ parametrized by
$y\in \overline{Q}^{(r_0)}$ and
$t\in \overline{\Delta}_{\vartheta_0}^{T_0,T}$, where
\begin{math}
  \mathbf{P}^{\alpha\beta} (x + \ii y,t)
  = ( P^{\alpha\beta}_{jk} )_{j,k=1}^M
\end{math}
for $|\alpha| = |\beta| = m$.
This is guaranteed by our Hypothesis~\eqref{hy:UnifEll}
that all
$P^{\alpha\beta}_{jk}$ ($|\alpha| = |\beta| = m$)
belong to
$C^0_{\mathrm{unif}}(\overline{\Omega})$
as functions of $(z,t) = (x + \ii y,t)$ $\in \overline{\Omega}$.
\qed
\par\vskip 10pt

In order to give a natural lower estimate on the domain of holomorphy
(i.e., the domain of complex analyticity)
of a weak solution $\mathbf{u}$ to the Cauchy problem \eqref{e:Cauchy},
we introduce a few more subsets of $\CC^N\times \CC$
(cf.\ \cite[p.~58]{Takac} or \cite[p.~428]{TakacBoller}):

We use the subdomain 
\begin{math}
  \Gamma^{(T')}_T (r',\vartheta') =
  \mathfrak{X}^{(r')}\times \Delta_{\vartheta'}^{T',T}
\end{math}
of the (larger) domain
\begin{equation}
\label{e:Gamma}
  \Omega = \Gamma^{(T_0)}_T (r_0,\vartheta_0) \eqdef
  \mathfrak{X}^{(r_0)}\times \Delta_{\vartheta_0}^{T_0,T}
  \subset \CC^N\times \CC
\end{equation}
defined above before
{\rm Hypotheses} \eqref{hy:UnifEll}, \eqref{hy:Ellipt}, and\/
\eqref{hy:f_Hardy^2}.
The three constants $T'\in (0,T_0]$,
$r'\in (0,r_0]$, and $\vartheta'\in (0,\vartheta_0]$
used below will be specified later (in Theorem~\ref{thm-Main}).

We recall the {\em tube\/}
\begin{math}
  \mathfrak{X}^{(r_0)} = \RR^N + \ii Q^{(r_0)}
\end{math}
(called often a {\em strip\/})
in $\CC^N$ with {\em base\/}
$Q^{(r_0)}\subset \RR^N$ defined in {\rm eq.}~\eqref{def:X^(r)}
and recall from {\rm eq.}~\eqref{loc:T^(T')} also the definition of the set
$\Delta_{\vartheta}^{(T')}$ in $\CC$.
In formula \eqref{glob:T^(r)}
we employ the time translation of
$\Delta_{\vartheta}^{(T')}$ by $r$ units, i.e., the set
$r + \Delta_{\vartheta}^{(T')}$, to define the union
$\Delta_{\vartheta}^{T',T}$ of such translations for
$0\leq r\leq T-T'< \infty$.
It is easy to see that, for
$0 < s\leq T < \infty$ and $0 < \vartheta_0 < \pi / 2$, we have
\begin{equation}
\label{eq:Lambda}
\begin{split}
& \Delta_{\vartheta_0}^{s,T}
  = \Delta_{\vartheta_0}^{(s)} \cup
    \left( [s,T)
         + \ii ( -s\cdot \tan\vartheta_0 ,\, s\cdot \tan\vartheta_0 )
    \right)
  =
\\
&   \left\{ t = \sigma + \ii\tau\in \CC\colon
    0 < \sigma < T \,\mbox{ and }\,
    |\tau| < r\cdot \tan\vartheta_0 \,\mbox{ where }\,
    r = \min\{ \sigma, s\} \right\} \,.
\end{split}
\end{equation}
Recall that the closure of $\Delta_{\vartheta_0}^{s,T}$
(and $\Delta_{\vartheta_0}^{(s)}$, respectively)
in $\CC$ is denoted by
$\overline{\Delta}_{\vartheta_0}^{s,T}$
(and $\overline{\Delta}_{\vartheta_0}^{(s)}$).
Given any $r\in [0,T)$, we observe that the (real) time $r$ section of
$\Delta_{\vartheta_0}^{s,T}$ is given by
\begin{equation*}
\begin{split}
    \{ t\in \Delta_{\vartheta_0}^{s,T} \colon \RE t = r \}
  = r + \ii ( -r'\cdot \tan\vartheta_0 ,\, r'\cdot \tan\vartheta_0 )
    \subset \CC
    \quad\text{ where }\, r' = \min \{ r, s\} \,.
\end{split}
\end{equation*}
These sets in the complex plane $\mathbb{C}$ have already been sketched
in Figure~\ref{fig:1} above and will be sketched more precisely
in Figure~$2$ below.

The Cartesian product
\begin{math}
  \Gamma^{(s)}_T (r_0,\vartheta_0)
  = \mathfrak{X}^{(r_0)} \times \Delta_{\vartheta_0}^{s,T}
\end{math}
defined in eq.~\eqref{e:Gamma} above
is our most important complex analyticity domain in $\CC^N\times \CC$.
Recall that
\begin{math}
  B^{s;p,p}(\RR^N) =\hfil\break
  \left( L^p(\RR^N) ,\, W^{2m,p}(\RR^N) \right)_{ 1 - (1/p) ,\, p }
    \hookrightarrow L^p(\RR^N)
\end{math}
with $p > 2 + \frac{N}{m}$.

Our main result is as follows:

\begin{theorem}\label{thm-Main}
Let\/ $m,M,N\geq 1$, $p > 2 + \frac{N}{m}$, $0 < T < \infty$,
and assume that all\/
{\rm Hypotheses} \eqref{hy:UnifEll}, \eqref{hy:Ellipt}, and\/
\eqref{hy:f_Hardy^2}
are satisfied in the product domain
\begin{math}
  \Omega = \Gamma^{(T_0)}_T (r_0,\vartheta_0)
         = \mathfrak{X}^{(r_0)}\times \Delta_{\vartheta_0}^{T_0,T}
  \subset \CC^N\times \CC
\end{math}
with some constants\/
$0 < r_0 < \infty$, $0 < T_0\leq T$, and\/ $0 < \vartheta_0 < \pi / 2$;
cf.\ {\rm eqs.}\ \eqref{glob:T^(r)} and\/ \eqref{def:X^(r)}.

{\bf (i)}$\;$
Then, given any\/ $t_0\in [0,T)$ and any initial value
$\mathbf{u}_0\in \mathbf{B}^{s;p,p}(\RR^N)$ at time $t = t_0$,
there is a number\/ $T_1\in (t_0,T]$,
depending on $t_0$ and $\mathbf{u}_0$, such that
the Cauchy problem \eqref{e:Cauchy} on the (local) time interval\/
$[t_0,T_1]\subset [0,T]$ with the initial condition
$\mathbf{u}(\,\cdot\,,t_0) = \mathbf{u}_0$ in
$\mathbf{B}^{s;p,p}(\RR^N)$ possesses a unique weak solution
\begin{math}
  \mathbf{u}\in
  C\left( [t_0,T_1]\to \mathbf{B}^{s;p,p}(\RR^N)\right) .
\end{math}

{\bf (ii)}$\;$
Furthermore, given any initial data
$\mathbf{u}_0\in \mathbf{B}^{s;p,p}(\RR^N)$ at time $t = 0$,
then any (global) weak solution
\begin{math}
  \mathbf{u}\in C\left( [0,T]\to \mathbf{B}^{s;p,p}(\RR^N)\right)
\end{math}
of the Cauchy problem \eqref{e:Cauchy}, {\em if it exists\/},
is always unique and it possesses
a unique (temporal) extension to the space\--time domain
$\RR^N\times \overline{\Delta}_{\vartheta'}^{T',T}$,
denoted again by $\mathbf{u}$, such that
the $\mathbf{B}^{s;p,p}(\RR^N)$\--valued function
\begin{math}
  \mathbf{u}\colon \overline{\Delta}_{\vartheta'}^{T',T}
            \mapsto \mathbf{B}^{s;p,p}(\RR^N)
\end{math}
is continuous in
$\overline{\Delta}_{\vartheta'}^{T',T}$ and its restriction to
$\Delta_{\vartheta'}^{T',T}$ is {\em holomorphic\/}, provided the numbers
$T'\in (0,T_0]$ and $\vartheta'\in (0,\vartheta_0]$ are small enough.
This temporal extension is a unique weak solution to
the Cauchy problem \eqref{e:Cauchy} in
$\RR^N\times \overline{\Delta}_{\vartheta'}^{T',T}$.

{\bf (iii)}$\;$
If, in addition to\/
$\mathbf{u}_0\in \mathbf{B}^{s;p,p}(\RR^N)$,
$\mathbf{u}_0$ possesses a (unique) {\em holomorphic\/} extension
\begin{math}
  \tilde{\mathbf{u}}_0\colon \mathfrak{X}^{(\kappa_0)}\to \CC^M
\end{math}
from $\RR^N$ to the complex domain
$\mathfrak{X}^{(\kappa_0)} \subset \CC^N$,
for some $\kappa_0\in (0,r_0]$, such that the function
\begin{equation*}
  \tilde{\mathbf{u}}_0 (\,\cdot\, + \ii y)
  \colon x\longmapsto \tilde{\mathbf{u}}_0 (x + \ii y)
  \colon \RR^N\longrightarrow \CC^M
\end{equation*}
belongs to $\mathbf{B}^{s;p,p}(\RR^N)$ 
for each $y\in Q^{(\kappa_0)}$ and
\begin{equation}
\label{sup_y:u_0}
  \mathfrak{N}^{(\kappa_0)} (\tilde{\mathbf{u}}_0) \eqdef
    \sup_{ y\in Q^{(\kappa_0)} }
    \| \tilde{\mathbf{u}}_0 (\,\cdot\, + \ii y)
    \|_{ B^{s;p,p}(\RR^N) } < \infty \,,
\end{equation}
then also any (global) weak solution
\begin{math}
  \mathbf{u}\in C\left( [0,T]\to \mathbf{B}^{s;p,p}(\RR^N)\right)
\end{math}
of the Cauchy problem \eqref{e:Cauchy}, {\em if it exists\/},
possesses a unique extension to the space\--time domain
\begin{math}
  \mathfrak{X}^{(r')}\times \overline{\Delta}_{\vartheta'}^{T',T} ,
\end{math}
denoted by $\tilde{\mathbf{u}}$, with the following properties,
provided the numbers $T'\in (0,T_0]$, $r'\in (0,\kappa_0]$, and
$\vartheta'\in (0,\vartheta_0]$ are small enough:
\begin{quote}
{\bf ($\mathbf{iii}_1$)}$\;$
\begin{math}
  \tilde{\mathbf{u}} (\,\cdot\, + \ii y,\, t)
  \in \mathbf{B}^{s;p,p}(\RR^N)
\end{math}
holds for all\/
$(y,t)\in Q^{(r')}\times \overline{\Delta}_{\vartheta'}^{T',T}$,
together with
\begin{equation}
\label{sup_(y,t):tilde_u}
  \sup_{ t\in \overline{\Delta}_{\vartheta'}^{T',T} }
  \mathfrak{N}^{(r')} \left( \tilde{\mathbf{u}} (\,\cdot\, ,\, t)\right)
  = \sup_{ t\in \overline{\Delta}_{\vartheta'}^{T',T} }
    \sup_{ y\in Q^{(\kappa_0)} }
    \| \tilde{\mathbf{u}} (\,\cdot\, + \ii y,\, t)
    \|_{ B^{s;p,p}(\RR^N) } < \infty \,,
\end{equation}
%
{\bf ($\mathbf{iii}_2$)}$\;$
the function
\begin{equation*}
  \tilde{\mathbf{u}}\colon (y,t) \,\longmapsto\,
  \tilde{\mathbf{u}} (\,\cdot\, + \ii y,\, t)
  \colon Q^{(r')}\times \overline{\Delta}_{\vartheta'}^{T',T}
  \,\longrightarrow\, \mathbf{B}^{s;p,p}(\RR^N)
\end{equation*}
is continuous in the space\--time variable
\begin{math}
  (y,t)\in Q^{(r')}\times \overline{\Delta}_{\vartheta'}^{T',T} ,
\end{math}
and\/
\hfill\break
{\bf ($\mathbf{iii}_3$)}$\;$
$\tilde{\mathbf{u}}$ is {\em holomorphic\/} in the complex domain
\begin{equation*}
  \Gamma^{(T')}_T (r',\vartheta') =
  \mathfrak{X}^{(r')}\times \Delta_{\vartheta'}^{T',T}
  \subset
  \mathfrak{X}^{(\kappa_0)}\times \Delta_{\vartheta_0}^{T_0,T}
  \subset
  \Omega = \Gamma^{(T_0)}_T (r_0,\vartheta_0)
         = \mathfrak{X}^{(r_0)}\times \Delta_{\vartheta_0}^{T_0,T} \,.
\end{equation*}
\end{quote}
Finally, the extension $\tilde{\mathbf{u}}$ verifies
the partial differential equation in the Cauchy problem \eqref{e:Cauchy}
pointwise in $\Gamma^{(T')}_T (r',\vartheta')$, i.e.,
in the classical sense, and obeys the initial data as follows,
\begin{equation}
\label{init:Cauchy}
  \| \tilde{\mathbf{u}} (\,\cdot\, + \ii y,\, t)
   - \tilde{\mathbf{u}}_0 (\,\cdot\, + \ii y)
  \|_{ B^{s;p,p}(\RR^N) } \,\longrightarrow\, 0
    \quad\mbox{ as }\, t\to 0 \,,\quad t\in \Delta_{\vartheta'}^{(T')} \,,
\end{equation}
for every\/ $y\in Q^{(r')}$.
\end{theorem}
\par\vskip 10pt

In Part~{\bf (iii)}, properties
{\bf ($\mathbf{iii}_1$)} and {\bf ($\mathbf{iii}_2$)}
combined with the Sobolev\-(-Besov) imbedding\hfil\break
\begin{math}
  B^{s;p,p}(\RR^N)\hookrightarrow
  C^m(\RR^N) \cap W^{m,\infty}(\RR^N)
\end{math}
guarantee the continuity and boundedness of the function
\begin{math}
  \tilde{\mathbf{u}}\colon
  \mathfrak{X}^{(r')}\times \overline{\Delta}_{\vartheta'}^{T',T}
  \to \CC^M .
\end{math}

Our condition $p > 2 + \frac{N}{m}$ is natural (and sharp)
to guarantee the continuity of the Sobolev imbedding
\begin{equation*}
  \mathbf{B}^{s;p,p}(\RR^N) = [ B^{s;p,p}(\RR^N) ]^M
  \hookrightarrow
  \mathbf{C}^m(\RR^N) \cap \mathbf{W}^{m,\infty}(\RR^N)
  = [ C^m(\RR^N) \cap W^{m,\infty}(\RR^N) ]^M
\end{equation*}
which follows from the Sobolev inequalities
and the Sobolev imbedding
\begin{equation*}
  B^{s-m;p,p}(\RR^N) \hookrightarrow C^0_{\mathrm{bdd}}(\RR^N)\eqdef
  C(\RR^N)\cap L^{\infty}(\RR^N)
  \quad\mbox{ for }\; N/p < s-m < \infty \,;
\end{equation*}
see, e.g.,
{\sc R.~A.\ Adams} and {\sc J.~J.~F.\ Fournier}
\cite[Chapt.~7]{AdamsFournier}, Theorem 7.34(c), p.~231.
In our proposition below we abbreviate
\begin{math}
  \varsigma_1(s) = \min\{ s,1\}
\end{math}
for $s\in \RR_+$.


\begin{figure}
\label{fig:2}
\centering
\includegraphics[width=0.9\textwidth]{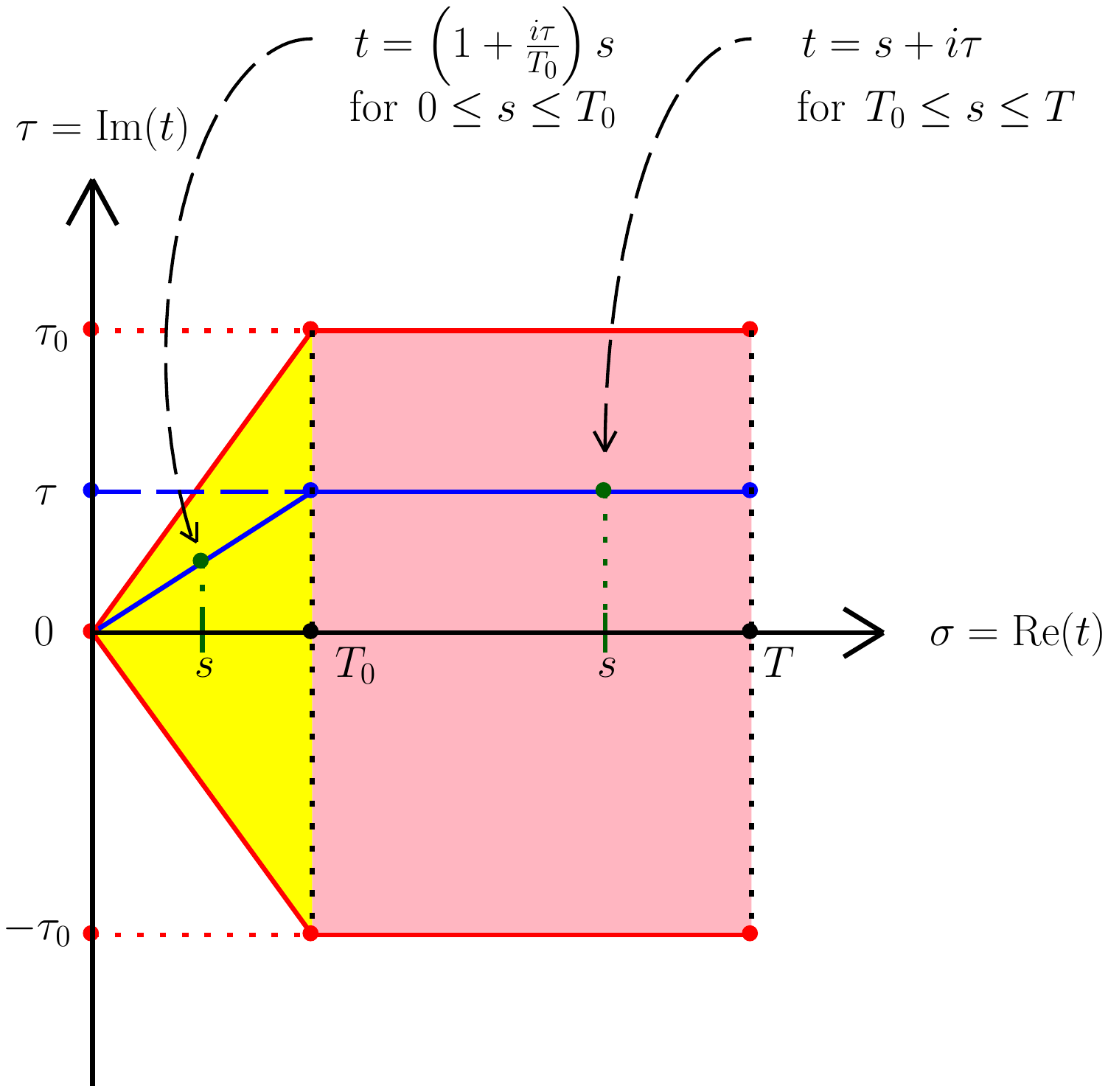}
\caption{
Here is a more detailed version of Figure~\ref{fig:1}
used in Proposition~\ref{prop-Main} above.
}
\end{figure}

This proposition follows directly from our proof of
Part~{\bf (iii)} of Theorem~\ref{thm-Main} given in
Section~\ref{s:proof-Main}.
The temporal integration path in the double space\--time integrals
\eqref{eq:du/dt_Hardy^2} and~\eqref{eq:u_Hardy^2}
is sketched in Figure~$2$ above.

\begin{proposition}\label{prop-Main}
In the situation of\/ {\rm Theorem~\ref{thm-Main}} above,
there exist constants\/ $c_0 > 0$ and\/ $C_0 > 0$
{\em depending solely\/} on
$\kappa_0$, $\vartheta_0$, $K$, $T'$, $r'$, $\vartheta'$ and
the supremum norm
\begin{equation*}
  \vertiii{\mathbf{u}}_{ L^{\infty}(0,T) }\eqdef
  \| \mathbf{u}\|_{ C\left( [0,T]\to \mathbf{B}^{s;p,p}(\RR^N)\right) }
  = \sup_{0\leq t\leq T}
         \| \mathbf{u}(\,\cdot\,,\, t)\|_{ B^{s;p,p}(\RR^N) }
  \hspace*{10pt} (< \infty)
\end{equation*}
of the (global) weak solution
\begin{math}
  \mathbf{u}\in C\left( [0,T]\to \mathbf{B}^{s;p,p}(\RR^N)\right) ,
\end{math}
such that the following estimate holds for all pairs\/
\begin{math}
  (y,t)\in Q^{(r')}\times \overline{\Delta}_{\vartheta'}^{T',T}
\end{math}
with $t = \sigma + \ii\tau$ ($\sigma, \tau\in \mathbb{R}$):
\begin{equation}
\label{eq:du/dt_Hardy^2}
\begin{aligned}
& \int_0^{\sigma}
    \int_{\mathbb{R}^N}
    \left| \partial_t \tilde{\mathbf{u}}
  \left( x + \ii y ,\, s + \ii\varsigma_1(s/T') \tau
  \right)
    \right|^p \,\mathrm{d}x \,\mathrm{d}s
\\
& {}
  + c_0\sum_{|\alpha|\leq 2m} \int_0^{\sigma}
    \int_{\mathbb{R}^N}
    \left| D_x^{\alpha} \tilde{\mathbf{u}}
  \left( x + \ii y ,\, s + \ii\varsigma_1(s/T') \tau
  \right)
    \right|^p \,\mathrm{d}x \,\mathrm{d}s
  \leq C_0 \,.
\end{aligned}
\end{equation}
Similarly, there are constants\/ $c_0'> 0$ and\/ $C_0'> 0$
{\em depending solely\/} on the constants\/ $c_0$ and\/ $C_0$,
such that the following estimate holds for all pairs\/
\begin{math}
  (y,t)\in Q^{(r')}\times \overline{\Delta}_{\vartheta'}^{T',T}
\end{math}
with $t = \sigma + \ii\tau$ ($\sigma, \tau\in \mathbb{R}$):
\begin{equation}
\label{eq:u_Hardy^2}
\begin{aligned}
& \left\Vert
    \tilde{\mathbf{u}}
  \left( \,\cdot\, + \ii y ,\, \sigma + \ii\varsigma_1(\sigma/T') \tau
  \right)
  \right\Vert_{ B^{s;p,p}(\RR^N) }^p
\\
& {}
  + c_0'\sum_{|\alpha|\leq 2m} \int_0^{\sigma}
    \int_{\mathbb{R}^N}
    \left| D_x^{\alpha} \tilde{\mathbf{u}}
  \left( x + \ii y ,\, s + \ii\varsigma_1(s/T') \tau
  \right)
    \right|^p \,\mathrm{d}x \,\mathrm{d}s
  \leq C_0' \,.
\end{aligned}
\end{equation}
\end{proposition}
\par\vskip 10pt

%
%

\begin{remark}\label{rem-Main-substi}\nopagebreak
\begingroup\rm
(See Figure~$2$ above).
Notice that, in {\rm eqs.}
\eqref{eq:du/dt_Hardy^2} and~\eqref{eq:u_Hardy^2} above,
the temporal argument in the function
%
\begin{math}
    D_x^{\alpha} \mathbf{u}
  \left( x + \ii y ,\, s + \ii\varsigma_1(s/T')\tau \right)
\end{math}
reads
\begin{math}
  s + \ii\varsigma_1(s/T')\tau = (1 + \ii (\tau / T')) s
\end{math}
whenever $0\leq s < T'< \infty$, whereas
\begin{math}
  s + \ii\varsigma_1(s/T')\tau = s + \ii\tau
\end{math}
holds whenever $0 < T'\leq s$ ($\leq \sigma\leq T < \infty$).
%
\hfill\Square
\endgroup
\end{remark}
\par\vskip 10pt

\section{Abstract Cauchy problem in an interpolation space}
\label{s:Cauchy_abstr}

We assume that $E = (E_0,E_1)$ is a {\em Banach couple\/}, that is,
$E_0$, $E_1$ are Banach spaces such that
$E_1$ is densely and continuously imbedded into $E_0$, i.e.,
$E_1\hookrightarrow E_0$.
We consider only complex Banach spaces over the field $\CC$.
Given a number $1 < p < \infty$, we denote by
\begin{equation*}
  E_{ 1 - \frac{1}{p} ,\, p } \equiv
  (E_0, E_1)_{ 1 - \frac{1}{p} ,\, p }
\end{equation*}
the real interpolation space between $E_1$ and $E_0$ obtained by
the {\em trace method\/} as follows, with the paremeter
$\theta = 1 - \frac{1}{p}\in (0,1)$.
We define such an interpolation space for any $\theta\in (0,1)$
below, cf.\
{\sc A.\ Lunardi} \cite[Chapt.~1]{Lunardi}, {\S}1.2.2, pp.\ 20--25.
The reader is referred to
{\sc R.~A.\ Adams} and {\sc J.~J.~F.\ Fournier}
\cite[Chapt.~7]{AdamsFournier}, {\S}7.6--{\S}7.23, pp.\ 208--221, or
{\sc H.\ Triebel}
\cite[Chapt.~1]{Triebel}, {\S}1.8, pp.\ 41--55,
for further details.
The trace spaces were originally introduced in
{\sc J.-L.\ Lions}
\cite{Lions-60, LionsPetr-61, LionsMag-61}.

Let $X_{\theta}^p$ denote the Banach space of
all Bochner\--measurable functions
$u\colon \RR_+\to E_0$ endowed with the weighted Lebesgue norm
\begin{equation}
\label{norm:X_th^p}
  \| u\|_{ X_{\theta}^p } \eqdef
  \left(
    \int_0^{\infty} \| t^{1-\theta}\, u(t)\|_{E_0}^p
    \,\frac{\mathrm{d}t}{t}
  \right)^{1/p} \equiv
  \left(
    \int_0^{\infty} \| u(t)\|_{E_0}^p
    \,\frac{\mathrm{d}t}{ t^{1 - (1-\theta) p } }
  \right)^{1/p} < \infty \,.
\end{equation}
Notice that
$X_{ 1 - \frac{1}{p} }^p = L^p(\RR_+\to E_0)$.
Analogously, we define the Banach space $Y_{\theta}^p$ of all functions
$u\in X_{\theta}^p$ with the following properties:
$u$ can be identified (by equality a.e.\ in $\RR_+$)
with a Bochner\--measurable function
$u\colon \RR_+\to E_1$ satisfying
\begin{equation}
\label{snorm:Y_th^p}
  \llbracket u\rrbracket_{ Y_{\theta}^p } \eqdef
  \left(
    \int_0^{\infty} \| t^{1-\theta}\, u(t)\|_{E_1}^p
    \,\frac{\mathrm{d}t}{t}
  \right)^{1/p} \equiv
  \left(
    \int_0^{\infty} \| u(t)\|_{E_1}^p
    \,\frac{\mathrm{d}t}{ t^{1 - (1-\theta) p } }
  \right)^{1/p} < \infty \,,
\end{equation}
and there is a function $v\in X_{\theta}^p$,
denoted by $v = u'$ in the sequel, such that the equality
\begin{equation}
\label{int:u'=v}
  u(t_2) - u(t_1) = \int_{t_1}^{t_2} v(s) \,\mathrm{d}s
  \quad\mbox{ holds in $E_0$ for all }\, 0 < t_1\leq t_2 < \infty \,.
\end{equation}
Applying H\"older's inequality to this equation,
it is easy to show that every
$u\in Y_{\theta}^p$ is $\theta$\--H\"older continuous on
any compact interval $[0,T]$ as a function valued in $E_0$, see, e.g.,
{\sc A.\ Lunardi} \cite[Chapt.~1]{Lunardi}, {\S}1.2.2, p.~20.
The Banach space $Y_{\theta}^p$ is endowed with the norm
\begin{equation}
\label{norm:Y_th^p}
  \| u\|_{ Y_{\theta}^p } \eqdef
    \llbracket u\rrbracket_{ Y_{\theta}^p }
  + \| u'\|_{ X_{\theta}^p } < \infty \,.
\end{equation}
For the special case $\theta = 1 - \frac{1}{p}$,
a useful equivalent norm on $Y_{ 1 - \frac{1}{p} }^p$ is defined by
\begin{equation}
\label{norm:Y_p^p^}
  \| u\|_{ Y_{ 1 - \frac{1}{p} }^p }^{\sharp} \eqdef
  \left( \int_0^{\infty} \| u(t)\|_{E_1}^p \,\mathrm{d}t
       + \int_0^{\infty} \| u'(t)\|_{E_0}^p \,\mathrm{d}t
  \right)^{1/p} \,.
\end{equation}
Thus,
\begin{equation*}
    Y_{ 1 - \frac{1}{p} }^p
  = L^p(\RR_+\to E_1)\cap W^{1,p}(\RR_+\to E_0)
\end{equation*}
is an abstract Sobolev space
(see {\sc H.\ Amann} \cite[Chapt.~III, {\S}1.1]{Amann-Birkh}, pp.\ 88--89).

Finally, the interpolation {\em trace space\/} is the vector space
\begin{equation*}
  E_{\theta,p}\equiv (E_0, E_1)_{\theta,p}\eqdef
  \left\{ u(0)\in E_0\colon u\in Y_{\theta}^p \right\}
\end{equation*}
of the initial values $x = u(0)\in E_0$ of all functions
$u\in Y_{\theta}^p$ endowed with the trace norm
\begin{equation}
\label{norm:E_th^p}
  \| x\|_{ E_{\theta,p} } \eqdef
  \inf\left\{ \| u\|_{ Y_{\theta}^p }\colon x = u(0)
              \;\mbox{ for some }\; u\in Y_{\theta}^p
      \right\}
\end{equation}
which makes the (linear) {\em trace mapping\/}
$\tau\colon u\mapsto u(0)\colon Y_{\theta}^p\to E_{\theta,p}$
bounded (i.e., continuous), with the operator norm $\leq 1$.
Equivalently to eq.~\eqref{norm:E_th^p}, we have
\begin{math}
  \| x\|_{ E_{\theta,p} }\leq \| u\|_{ Y_{\theta}^p }
\end{math}
for every $u\in Y_{\theta}^p$ with $u(0) = x$
and there exists a sequence
$\{ u_n\}_{n=1}^{\infty}\subset Y_{\theta}^p$
such that $u_n(0) = x$ and
\begin{math}
  \| u_n\|_{ Y_{\theta}^p } \,\longrightarrow\, \| x\|_{ E_{\theta,p} }
\end{math}
as $n\to \infty$.
It can be shown that there is a constant $c = c(\theta,p) > 0$,
depending only on $E = (E_0,E_1)$,
$\theta\in (0,1)$, and $p\in (1,\infty)$, such that
\begin{equation}
\label{inter:Hoelder}
  \| x\|_{ E_{\theta,p} } \leq c\,
  \| x\|_{E_0}^{1-\theta}\, \| x\|_{E_1}^{\theta}
  \quad\mbox{ holds for all }\, x\in E_1 \,;
\end{equation}
see, e.g.,
{\sc H.\ Triebel} \cite[Chapt.~1]{Triebel}, Theorem 1.3.3(g), p.~25,
combined with Theorem 1.8.2, pp.\ 44--45.
As an easy consequence of the definition of $E_{\theta,p}$ for
$\theta = 1 - \frac{1}{p}$, i.e., $(1-\theta) p = 1$,
one can show that the abstract Sobolev space
$Y_{ 1 - \frac{1}{p} }^p$
is continuously imbedded into the Fr\'echet space
\begin{math}
  C\left( \RR_+\to E_{ 1 - \frac{1}{p} ,\, p } \right)
\end{math}
of all continuous functions
$u\colon \RR_+\to E_{ 1 - \frac{1}{p} ,\, p }$
endowed with the (locally convex) topology of uniform convergence
on every compact time interval $[t_1,t_2]\subset \RR_+$,
\begin{equation*}
    Y_{ 1 - \frac{1}{p} }^p
  = L^p(\RR_+\to E_1)\cap W^{1,p}(\RR_+\to E_0)
  \hookrightarrow
  C\left( \RR_+\to E_{ 1 - \frac{1}{p} ,\, p } \right) .
\end{equation*}
We complete our definition by setting
$E_{\theta,p}\eqdef E_{\theta}$ if $\theta\in \{ 0,\,1\}$.

In what follows we deal with applications of 
the interpolation trace space $E_{\theta,p}$
(with $\theta = 1 - \frac{1}{p}$)
to abstract linear and nonlinear evolutionary problems of type
\begin{equation}
\label{e:abstr_Cauchy}
\left\{
\begin{alignedat}{2}
    \frac{\mathrm{d}u}{\mathrm{d}t}
  - A(t,u(t)) u(t) &= f(t,u(t)) + g(t)
  \quad\mbox{ for a.e. } t\in (0,T) \,;
\\
  u(0) &= u_0\in E_{ 1 - \frac{1}{p} ,\, p } \,.
\end{alignedat}
\right.
\end{equation}
Here, $u\colon (0,T)\to E_0$
is the unknown function valued in the Banach space $E_0$ and
$0 < T\leq \infty$.
A rigorous definition of a {\it\bfseries strict solution\/} $u$
of the initial value problem \eqref{e:abstr_Cauchy}
will be given below, in Definition~\ref{def-strict_sol}.
Essentially, we follow
{\sc Ph.\ Cl\'ement} and {\sc S.\ Li}
\cite{Clem-Li}, Section~1, pp.\ 17--18.
A closely related approach is carried out also in
{\sc M.\ K\"ohne}, {\sc J.\ Pr\"u{ss}}, and {\sc M.\ Wilke}
\cite{Koehne-Pruess-W}.

We denote by $\mathcal{L}(E_1\to E_0)$
the Banach space of all bounded (i.e., continuous) linear operators
$B\colon E_1\to E_0$ endowed with the standard operator norm
$\| B\|_{ \mathcal{L}(E_1\to E_0) }$.
Let us denote by
$I\colon E_1\hookrightarrow E_0$
the continuous imbedding of $E_1$ into $E_0$; hence,
$I\in \mathcal{L}(E_1\to E_0)$.
We identify $I$ with the identity mapping in the whole of $E_0$
and abbreviate
$\mathcal{L}(E_0)\equiv \mathcal{L}(E_0\to E_0)$.

If, for some complex number $\lambda\in \CC$, the operator
$\lambda I - B\in \mathcal{L}(E_1\to E_0)$
is invertible with an inverse denoted by
\begin{math}
  (\lambda I - B)^{-1}\colon E_0\to E_1\hookrightarrow E_0
\end{math}
such that this inverse is bounded from $E_0$ into itself, i.e.,
$(\lambda I - B)^{-1}\in \mathcal{L}(E_0\to E_0)$,
then we alternatively (equivalently) view $B$ as
a densely defined, closed linear operator
$B\colon E_0\to E_0$ with the domain $\mathcal{D}(B) = E_1$,
by the closed graph theorem, cf.\
{\sc H.~Amann} \cite[Chapt.~I, Lemma 1.1.2]{Amann-Birkh}, p.~10.
Indeed, if the graph $\mathcal{G}(B)$ of $B$ is closed in
$E_0\times E_0$, it is closed also in $E_1\times E_0$.
In this case, the norm $\|\cdot\|_{E_1}$ on $E_1$
is equivalent with the {\em graph norm\/}
\begin{equation*}
  \| x\|_{ \mathcal{D}(B) }\eqdef
  \| Bx\|_{E_0} + \| x\|_{E_0} \,,\quad x\in \mathcal{D}(B) \,,
\end{equation*}
on $\mathcal{D}(B) = E_1$.
An important class of such operators, denoted by
$\mathrm{Gen}(E)\equiv \mathrm{Gen}(E_1\to E_0)$, is formed by
all closed linear operators $B\colon E_0\to E_0$
with the domain $\mathcal{D}(B) = E_1$
that generate a strongly continuous semigroup
$\{ \ee^{tB}\colon t\geq 0\}$ on $E_0$.
We will consider only {\it\bfseries generators\/} $B$ with domain
$\mathcal{D}(B) = E_1$.
Finally, we denote by
$\mathrm{Hol}(E)\equiv \mathrm{Hol}(E_1\to E_0)$
the subset of all (infinitesimal) generators $B\in \mathrm{Gen}(E)$
that generate
a {\it\bfseries holomorphic\/} (i.e., analytic) semigroup on $E_0$.
We refer to
{\sc H.~Amann} \cite[Chapt.~I, {\S}1]{Amann-Birkh}, pp.\ 9--24,
{\sc A.~Pazy} \cite[Chapt.\ 1--2]{Pazy}, pp.\ 1--75, or
{\sc H.~Tanabe} \cite[Chapt.~3, {\S}3.1--{\S}3.4]{Tana-79}, pp.\ 51--72,
for details about strongly continuous (and holomorphic) semigroups.

Next, given an operator $B\in \mathrm{Hol}(E)$,
let us consider the following special (linear) case of
problem~\eqref{e:abstr_Cauchy}, namely,
\begin{equation}
\label{lin:abstr_Cauchy}
\left\{
\begin{alignedat}{2}
    \frac{\mathrm{d}u}{\mathrm{d}t} - B u(t) &= g(t)
  \quad\mbox{ for a.e. } t\in (0,T) \,;
\\
  u(0) &= u_0\in E_{ 1 - \frac{1}{p} ,\, p } \,.
\end{alignedat}
\right.
\end{equation}
Here,
$u_0\in E_{ 1 - \frac{1}{p} ,\, p }$ is a given initial value,
$g\in L^p((0,T)\to E_0)$ is a given function,
$1 < p < \infty$, and $0 < T < \infty$.
In analogy with our definition of the Banach spaces
$X_{\theta}^p$ and $Y_{\theta}^p$ of functions
$u\colon \RR_+\to E_0$ on the entire half\--line $\RR_+$,
endowed with the norms given by {\rm eqs.}\
\eqref{norm:X_th^p} and \eqref{norm:Y_th^p}, respectively,
we introduce the corresponding Banach spaces
$X_{\theta}^p(0,T)$ and $Y_{\theta}^p(0,T)$ of functions
$u\colon [0,T)\to E_0$ on a bounded interval $[0,T)$, $0 < T < \infty$.
Of course, in eq.~\eqref{norm:X_th^p}, the integral
$\int_0^{\infty}\dots \,\frac{\mathrm{d}t}{t}$
has to be replaced by
$\int_0^T\dots \,\frac{\mathrm{d}t}{t}$ .
It is not difficult to show that if one replaces the pair of spaces
$X_{\theta}^p$ and $Y_{\theta}^p$ by
$X_{\theta}^p(0,T)$ and $Y_{\theta}^p(0,T)$, respectively,
in the definition of the trace space $E_{\theta,p}$ and its norm in
eq.~\eqref{norm:E_th^p},
the {\em same\/} interpolation trace space is obtained.
These facts can be inferred easily from the treatment of trace spaces
in the monographs
\cite{AdamsFournier, Amann-Birkh, Lunardi, Triebel} or from
the original works by
{\sc J.-L.\ Lions} \cite{Lions-60, LionsPetr-61, LionsMag-61}.
In particular, we have the continuous imbedding
\begin{equation}
\label{e:Besov_imbed}
    Y_{ 1 - \frac{1}{p} }^p(0,T)
  = L^p((0,T)\to E_1)\cap W^{1,p}((0,T)\to E_0)
  \hookrightarrow
    C\left( [0,T]\to E_{ 1 - \frac{1}{p} ,\, p } \right) ,
\end{equation}
see, e.g.,
\cite[Chapt.~7]{AdamsFournier}, {\S}7.67, p.~255.
Thus, the (linear) trace mapping
\begin{equation*}
  \tau\colon u\mapsto u(0)\colon
  Y_{ 1 - \frac{1}{p} }^p(0,T) \to E_{ 1 - \frac{1}{p} ,\, p }
\end{equation*}
is continuous.

We say that a function
$u\colon [0,T)\to E_0$ is a {\it\bfseries strict solution\/}
of the initial value problem \eqref{lin:abstr_Cauchy} if
\begin{equation*}
  u\in Y_{ 1 - \frac{1}{p} }^p(0,T) \,,\quad \tau u\equiv u(0) = u_0 \,,
\end{equation*}
and the differential equation in \eqref{lin:abstr_Cauchy}
is satisfied with all terms in
$X_{ 1 - \frac{1}{p} }^p(0,T) = L^p((0,T)\to E_0)$.

\begin{definition}\label{def-max_regul}\nopagebreak
\begingroup\rm
An infinitesimal generator $B\in \mathrm{Hol}(E)$
of a holomorphic semigroup on $E_0$ with domain
$\mathcal{D}(B) = E_1$ is said to possess
the {\it\bfseries maximal $L^p$-regularity property\/}, symbolically
\begin{math}
  B\in \mathrm{MR}_p(E)\equiv \mathrm{MR}_p(E_1\to E_0) ,
\end{math}
if for any given initial condition
$u_0\in E_{ 1 - \frac{1}{p} ,\, p }$
and any given function $g\in L^p((0,T)\to E_0)$,
problem \eqref{lin:abstr_Cauchy} posesses a unique strict solution
$u\in Y_{ 1 - \frac{1}{p} }^p(0,T)$
that satisfies the following estimate:
\begin{quote}
\begingroup\em
\underline{There exists} a constant $M\equiv M(p,E,B,T) > 0$,
independent from $u_0$ and $g$, such that
\endgroup
\begin{equation}
\label{est:strict_sol}
\begin{aligned}
    \int_0^T \genfrac{\|}{\|}{}0{\mathrm{d}u}{\mathrm{d}t}_{E_0}^p
    \,\mathrm{d}t
  + \int_0^T \| B u(t)\|_{E_0}^p \,\mathrm{d}t
  \leq M
    \left(
      \| u_0\|_{ E_{ 1 - \frac{1}{p} ,\, p } }^p
    + \int_0^T \| g(t)\|_{E_0}^p \,\mathrm{d}t
    \right) \,.
  \quad\mbox{\Square}
\end{aligned}
\end{equation}
\end{quote}
%
\endgroup
\end{definition}
\par\vskip 10pt

We have adopted this definition of class $\mathrm{MR}_p(E)$ from
{\sc Ph.\ Cl\'ement} and {\sc S.\ Li} \cite[p.~18]{Clem-Li}
and from the monograph by
{\sc A.\ Ashyralyev} and {\sc P.~E.\ Sobolevskii}
\cite[Chapt.~1]{Ashyr-Sobol-94}, {\S}3.5, p.~28.
It may be viewed as some kind of {\em ellipticity hypothesis\/}
for the linear operator
$B\in \mathcal{L}(E_1\to E_0)$ or a {\em stability hypothesis\/}
for the linear parabolic problem~\eqref{lin:abstr_Cauchy}.
Equivalently, the abstract (linear evolutionary)
partial differential operator
\begin{equation*}
  \left( \partial_t - B ,\, \tau\right) \colon
    Y_{ 1 - \frac{1}{p} }^p(0,T) \to
    L^p((0,T)\to E_0) \times E_{ 1 - \frac{1}{p} ,\, p } \,,
\end{equation*}
defined for every
\begin{math}
  u\in Y_{ 1 - \frac{1}{p} }^p(0,T)
     = L^p((0,T)\to E_1)\cap W^{1,p}((0,T)\to E_0)
\end{math}
by
\begin{equation}
\label{abstr:strict_sol}
  \left( \partial_t - B ,\, \tau\right) \colon
  u\;\longmapsto\;
  \left( \frac{\mathrm{d}u}{\mathrm{d}t} - B u(t) \,,\; u(0) \right) \,,
\end{equation}
possesses a bounded inverse furnished by the strict solution
\begin{equation*}
  u = \left( \partial_t - B ,\, \tau\right)^{-1} (g,u_0)
      \in Y_{ 1 - \frac{1}{p} }^p(0,T)
\end{equation*}
of problem \eqref{lin:abstr_Cauchy}; cf.\
{\sc S.\ Angenent} \cite[Lemma 2.2, p.~95]{Angenent}
for the parallel interpolation case $p = \infty$ introduced in
{\sc G.\ Da Prato} and {\sc P.\ Grisvard} \cite{DaPrato-Grisv}.

\begin{remark}\label{rem-perturb}\nopagebreak
\begingroup\rm

{\rm (a)}$\;$
It is not difficult to show that the maximal $L^p$-regularity class
$\mathrm{MR}_p(E)$ is independent from a particular choice of
$T\in (0,\infty)$; see
{\sc G.\ Dore} \cite[Sect.~5, p.~310]{Dore-2000}, Corollary 5.4,
or
{\sc J.\ Pr\"uss} \cite[p.~4]{Pruess}, remarks after Corollary 1.3.
More importantly, this class is {\em independent\/} from $p\in (1,\infty)$
as well, i.e.,
$\mathrm{MR}_p(E) = \mathrm{MR}_{p_0}(E)$ holds for all
$p, p_0\in (1,\infty)$, by a classical result due to
{\sc P.~E.\ Sobolevskii} \cite{Sobol-1961, Sobol-1964}; see, e.g.,
\cite[{\S}3.1, pp.\ 343--345]{Sobol-1961}.
Further details on the independence of
$\mathrm{MR}_p(E)$ from $p\in (1,\infty)$ can be found in
{\sc A.\ Ashyralyev} and {\sc P.~E.\ Sobolevskii}
\cite[Chapt.~1]{Ashyr-Sobol-94}, {\S}3.5, Theorem 3.6 on p.~35,
{\sc G.\ Dore} \cite[Sect.~7, p.~313]{Dore-2000}, Theorem 7.1,
or
{\sc M.\ Hieber} \cite[Corollary 4.4, p.~371]{Hieber-1999},
where in ineq.~\eqref{est:strict_sol} one may take
\begin{math}
  M\equiv M(p) = p^2 (p-1)^{-1}\, M(p_0) < \infty
\end{math}
if the constant $M(p_0)\in (0,\infty)$ is known,
by \cite{Ashyr-Sobol-94}.

{\rm (b)}$\;$
We are allowed to specify the constant
$M\equiv M(p,E,B,T) > 0$ in ineq.~\eqref{est:strict_sol}
to be the {\it smallest\/} nonnegative number $M\in \RR_+$
for which ineq.~\eqref{est:strict_sol} is valid; cf.\
{\sc Ph.\ Cl\'ement} and {\sc S.\ Li}
\cite[Proposition 2.2, p.~19]{Clem-Li}.
Then, clearly, $T\mapsto M(p,E,B,T)$
is a {\it nondecreasing\/} (nonnegative) function of time
$T\in (0,\infty)$.
Indeed, if $T'\in (0,T)$ and
$g\in L^p((0,T')\to E_0)$ is arbitrary, it suffices to apply
ineq.~\eqref{est:strict_sol} with the function
\begin{equation*}
  \tilde{g}(t) = \left\{
  \begin{alignedat}{2}
    g(t) \quad\mbox{ if }\, 0\leq t\leq T' \,;\\
    0    \quad\mbox{ if }\, T' < t\leq T \,,\\
  \end{alignedat}
  \right.
\end{equation*}
in place of $g$ in order to derive ineq.~\eqref{est:strict_sol}
for $T'$ in place of~$T$ with the same constant~$M$.
Hence, $M(p,E,B,T')\leq M(p,E,B,T)$ holds for $0 < T'< T$.
It is easy to see that $M(p,E,B,T) > 0$.
(The case $M = 0$ would easily lead to a contradiction.)
In what follows we always use this optimal value of~$M$, i.e.,
$M\equiv M(p,E,B,T) > 0$.

{\rm (c)}$\;$
Simple perturbation theory for linear operators shows that the set
$\mathrm{Hol}(E)\equiv \mathrm{Hol}(E_1\to E_0)$
is open in the Banach space $\mathcal{L}(E_1\to E_0)$.
Even a more precise, relative perturbation result is valid; see
{\sc T.\ Kato} \cite[Chapt.~IX]{Kato}, {\S}2.2, Theorem 2.4 on p.~499.
A similar result can be derived for the class
$\mathrm{MR}_p(E)\equiv \mathrm{MR}_p(E_1\to E_0)$
applying the perturbation technique from either
{\sc H.\ Amann} \cite[Chapt.~III, {\S}1.6]{Amann-Birkh},
Proposition 1.6.3 on p.~97, or from
{\sc Ph.\ Cl\'ement} and {\sc S.\ Li}
\cite[Proof of Theorem 2.1]{Clem-Li}, pp.\ 19--23:
The set $\mathrm{MR}_p(E)$ is open in $\mathcal{L}(E_1\to E_0)$; see
Lemma~\ref{lem-perturb_MR} below.
Indeed, this follows from the fact that
the set of all bounded linear operators from
\begin{equation*}
  \mathcal{L}\left(
  Y_{ 1 - \frac{1}{p} }^p(0,T) \to
  L^p((0,T)\to E_0) \times E_{ 1 - \frac{1}{p} ,\, p }
    \right)
\end{equation*}
that possess a bounded inverse is open in this Banach space,
and the inverse
$\left( \partial_t - B ,\, \tau\right)^{-1}$
is a locally Lipschitz\--continuous function of
$B\in \mathcal{L}(E_1\to E_0)$, by
Lemma~\ref{lem-perturb_MR} below and
formula \eqref{inv_Neumann:B+A} thereafter, with
$B\in \mathrm{MR}_p(E)$ being fixed and
$A\in \mathcal{L}(E_1\to E_0)$ having a sufficiently small operator norm
\begin{math}
  \| A\|_{ \mathcal{L}(E_1\to E_0) }
\end{math}
depending on $B$.
\hfill\Square
\endgroup
\end{remark}
\par\vskip 10pt

Now we are ready to define a strict solution $u$
to our abstract nonlinear evolutionary problem \eqref{e:abstr_Cauchy}.
We assume that $1 < p < \infty$, $0 < T < \infty$,
$g\in L^p((0,T)\to E_0)$, $u_0\in U$ where $U$ is an open set in
$E_{ 1 - \frac{1}{p} ,\, p }$, and the mappings
\begin{align*}
& A\colon (t,v)\mapsto A(t,v)\colon [0,T]\times U\subset
    [0,T]\times E_{ 1 - \frac{1}{p} ,\, p } \;\longrightarrow\;
    \mathcal{L}(E_1\to E_0) \,,
\\
& f\colon (t,v)\mapsto f(t,v)\colon [0,T]\times U\subset
    [0,T]\times E_{ 1 - \frac{1}{p} ,\, p } \;\longrightarrow\; E_0
\end{align*}
satisfy the following ``natural'' hypotheses
(cf.\ {\sc Ph.\ Cl\'ement} and {\sc S.\ Li}
      \cite[p.~19]{Clem-Li}, (H1)--(H3)):

\begin{hypos}\nopagebreak
\begingroup\rm
\begin{enumerate}
\setcounter{enumi}{0}
\renewcommand{\labelenumi}{{\bf (C\arabic{enumi})}}
\item
\makeatletter
\def\@currentlabel{{\bf C\arabic{enumi}}}\label{hyp:A}
\makeatother
$A\colon [0,T]\times U\to \mathcal{L}(E_1\to E_0)$
is a Lipschitz continuous mapping such that
$A(t,v)\in \mathrm{MR}_p(E)$ for all $(t,v)\in [0,T]\times U$.
\item
\makeatletter
\def\@currentlabel{{\bf C\arabic{enumi}}}\label{hyp:f}
\makeatother
$f\colon [0,T]\times U\to E_0$ is a Lipschitz continuous mapping.
\hfill\Square
\end{enumerate}
\endgroup
\end{hypos}
\par\vskip 10pt

Of course, the metric on $[0,T]\times U$ is induced by
the canonical norm on $\RR\times E_{ 1 - \frac{1}{p} ,\, p }$.
It is a matter of a straight\--forward calculation to verify that
both substitution mappings,
\begin{align*}
  (v,u) &\mapsto \left[ t\mapsto A(t,v(t))u(t)\right] \colon
  C\left( [0,T]\to U\right) \times L^p((0,T)\to E_1)
  \;\longrightarrow\; L^p((0,T)\to E_0) \,,
\\
  v &\mapsto \left[ t\mapsto f(t,v(t)) \right] \colon
  C\left( [0,T]\to U\right)
  \;\longrightarrow\; L^p((0,T)\to E_0) \,,
\end{align*}
are locally Lipschitz\--continuous with values in
$L^p((0,T)\to E_0)$; see, e.g.,
{\sc Ph.\ Cl\'ement} and {\sc S.\ Li}
\cite[Proof of Theorem 2.1]{Clem-Li}, pp.\ 19--23.

\begin{remark}\label{rem-strict_sol}\nopagebreak
\begingroup\rm
In {\rm Hypothesis} \eqref{hyp:A} we did not have to assume that
$A(t,v)\in \mathrm{MR}_p(E)$ holds for {\em all\/}
$(t,v)\in [0,T]\times U$.
We could assume {\em only\/}
$A(0,u_0)\in \mathrm{MR}_p(E)$; cf.\ results to follow below
(e.g., Theorems \ref{thm-Clem-Li} and~\ref{thm-Clem-Li_glob} and
 Remark~\ref{rem-Clem-Li}).
However, the set $\mathrm{MR}_p(E)$
being open in $\mathcal{L}(E_1\to E_0)$,
$A(0,u_0)\in \mathrm{MR}_p(E)$
would imply that there are a number $T_0\in (0,T]$ and
an open neighborhood $U_0$ of $u_0$ in $E_{ 1 - \frac{1}{p} ,\, p }$,
$u_0\in U_0\subset U$, such that
$A(t,v)\in \mathrm{MR}_p(E)$ holds for all $(t,v)\in [0,T_0]\times U_0$,
by the Lipschitz continuity of~$A$.
But this statement is qualitatively the same as
$A(t,v)\in \mathrm{MR}_p(E)$ for all $(t,v)\in [0,T]\times U$
in our {\rm Hypothesis} \eqref{hyp:A}.
\endgroup
\end{remark}
\par\vskip 10pt

\begin{definition}\label{def-strict_sol}\nopagebreak
\begingroup\rm
({\sc Ph.\ Cl\'ement} and {\sc S.\ Li} \cite[p.~18]{Clem-Li}.)
Recall that $U$ is an open set in
$E_{ 1 - \frac{1}{p} ,\, p }$ and $u_0\in U$.
We say that a function
$u\colon [0,T)\to E_0$ is a {\it\bfseries strict solution\/}
of the initial value problem \eqref{e:abstr_Cauchy} if
$u\in Y_{ 1 - \frac{1}{p} }^p(0,T)$,
$u(t)\in U$ for every $t\in [0,T]$, $u(0) = u_0$, and
the differential equation in \eqref{e:abstr_Cauchy}
is satisfied with all terms (summands) in $L^p((0,T)\to E_0)$.
\hfill\Square
\endgroup
\end{definition}
\par\vskip 10pt

We recall that the Banach space
$Y_{ 1 - \frac{1}{p} }^p(0,T)$
has been introduced in eq.~\eqref{e:Besov_imbed} above.

The main result in \cite[Theorem 2.1, p.~19]{Clem-Li}
is {\it local in time\/} and reads as follows, with
{\rm Hypothesis} \eqref{hyp:A} being somewhat weakened
in the sense of our Remark~\ref{rem-strict_sol} above.

\begin{theorem}\label{thm-Clem-Li}
Let\/ $1 < p < \infty$ and\/ $0 < T < \infty$.
Let\/ $U$ be a nonempty open set in\/
$E_{ 1 - \frac{1}{p} ,\, p }$ and\/ $u_0\in U$.
Assume that both mappings\/
$A\colon [0,T]\times U\to \mathcal{L}(E_1\to E_0)$ and\/
$f\colon [0,T]\times U\to E_0$ are Lipschitz\--continuous.
If\/
$A(0,u_0)\in \mathrm{MR}_p(E)$ then there exists some time
$T_1\equiv T_1(u_0)\in (0,T]$, depending on $u_0$, such that
the abstract initial value problem \eqref{e:abstr_Cauchy}
possesses a unique strict solution
\begin{align}
\label{e_1:strict_sol}
& u\in Y_{ 1 - \frac{1}{p} }^p(0,T_1)
\\
\nonumber
& \left( {}=
  L^p((0,T_1)\to E_1)\cap W^{1,p}((0,T_1)\to E_0)
    \hookrightarrow
  C\left( [0,T_1]\to E_{ 1 - \frac{1}{p} ,\, p } \right)
  \right)
\end{align}
on the time interval\/ $[0,T_1]$.
Consequently, one has\/
$u(t)\in U$ for every\/ $t\in [0,T_1]$.
\end{theorem}
\par\vskip 10pt

This theorem is proved in \cite{Clem-Li}, Section~2, pp.\ 20--23,
using the Banach contraction principle in the closed ball
\begin{equation*}
  \Sigma_{\rho_1,T_1}^{(u_0)} =
  \left\{ v\in Y^{T_1}\colon v(0) = u_0 \quad\mbox{ and }\quad
          \| v-w\|_{ Y^{T_1} }\leq \rho_1
  \right\}
\end{equation*}
of radius $\rho_1\in (0,\infty)$ centered at the point
$w\in Y^{T_1}$ in the Banach space
\begin{equation*}
  Y^{T_1} = Y_{ 1 - \frac{1}{p} }^p(0,T_1)
  = L^p((0,T_1)\to E_1)\cap W^{1,p}((0,T_1)\to E_0) \,.
\end{equation*}
Here, the ``center'' function $w\in Y^{T_1}$
is defined to be the restriction to $[0,T_1]$ of
the unique strict solution
$\tilde{w}\in Y^{T} = Y_{ 1 - \frac{1}{p} }^p(0,T)$
to the abstract initial value problem \eqref{lin:abstr_Cauchy}
in the time interval $[0,T]$ with the linear operator
$B = A(0,u_0)\in \mathrm{MR}_p(E)$ and the right\--hand side
$g(t)$ replaced by the sum $f(t,u_0) + g(t)$,
\begin{equation}
\label{nonlin:abstr_Cauchy}
\left\{
\begin{alignedat}{2}
    \frac{ \mathrm{d}\tilde{w} }{\mathrm{d}t} - A(0,u_0) \tilde{w}(t)
& = f(t,u_0) + g(t)
  \quad\mbox{ for a.e. } t\in (0,T) \,;
\\
  \tilde{w}(0) &= u_0\in E_{ 1 - \frac{1}{p} ,\, p } \,.
\end{alignedat}
\right.
\end{equation}
Although the proof in \cite{Clem-Li} has been carried out only for
$A(t,u) = A(u)$ independent from time $t\in [0,T]$,
it is a matter of straight\--forward calculations to adapt this proof to
the case of $A(t,u)$ depending on time $t$, cf.\
\cite[p.~23]{Clem-Li}, Remark at the end of Section~2.
A detailed treatment of the latter case is presented in
{\sc J.\ Pr\"uss} \cite[pp.\ 9--13]{Pruess}, Chapt.~3,
under slightly different assumptions
(see also {\sc M.\ K\"ohne}, {\sc J.\ Pr\"u{ss}}, and {\sc M.\ Wilke}
 \cite{Koehne-Pruess-W}).

\begin{remark}\label{rem-Clem-Li}\nopagebreak
\begingroup\rm
Furthermore, one can easily conclude from the proof of Theorem 2.1 in
\cite[pp.\ 20--23]{Clem-Li}
that if
$\overline{B}_{R_0}(w_0)$ is any closed ball in the Banach space
$E_{ 1 - \frac{1}{p} ,\, p }$ of radius $R_0\in (0,\infty)$
centered at a point $w_0\in E_{ 1 - \frac{1}{p} ,\, p }$, such that
$\overline{B}_{R_0}(w_0) \subset U$ and $R_0 > 0$ is small enough,
then the constants
$\rho_1\in (0,\infty)$ and $T_1\in (0,T]$
can be chosen small enough to depend solely on $R_0$,
but not on $w_0$, provided
$u_0\in \overline{B}_{R_0}(w_0) \subset U$.
The estimates in \cite[pp.\ 20--23]{Clem-Li},
based on the Lipschitz constants for $A$ and $f$ in $[0,T]\times U$
and the estimate in \eqref{est:strict_sol},
remain valid for any
$u_0\in \overline{B}_{R_0}(w_0)$.
Thus, we have
$T_1\equiv T_1(R_0)\in (0,T]$ and
$\rho_1\equiv \rho_1(R_0)\in (0,\infty)$.
Finally, using similar estimates, cf.\
\cite[p.~22]{Clem-Li}, (2.14)--(2.17),
one can show that the (strict) solution mapping
\begin{equation*}
  u_0\longmapsto u\colon \overline{B}_{R_0}(w_0)\subset U
    \subset E_{ 1 - \frac{1}{p} ,\, p } \longrightarrow
  Y^{T_1} = Y_{ 1 - \frac{1}{p} }^p(0,T_1)
\end{equation*}
is Lipschitz continuous with a Lipschitz constant
$L\equiv L(R_0)\in (0,\infty)$ independent from
$w_0\in E_{ 1 - \frac{1}{p} ,\, p }$, such that
$\overline{B}_{R_0}(w_0) \subset U$ and $R_0 > 0$ is small enough.
This means that if
$u_1, u_2\colon [0,T_1]\to E_{ 1 - \frac{1}{p} ,\, p }$
are two strict solutions to problem \eqref{e:abstr_Cauchy}
on the time interval $[0,T_1]$,
with (possibly different) initial values
$u_1(0) = u_{0,1}$ and $u_2(0) = u_{0,2}$ in
$\overline{B}_{R_0}(w_0) \subset U$, then one has
$u_1(t), u_2(t)\in U$ for all $t\in [0,T_1]$ and
\begin{equation}
\label{e:abstr_Lip}
    \| u_1 - u_2\|_{ Y^{T_1} }
  \leq L\,
    \| u_{0,1} - u_{0,2}\|_{ E_{ 1 - \frac{1}{p} ,\, p } } \,.
\end{equation}
Combining this inequality with the continuous imbedding
\begin{math}
  Y^{T_1}\hookrightarrow
  C\left( [0,T_1]\to E_{ 1 - \frac{1}{p} ,\, p } \right)
\end{math}
in \eqref{e:Besov_imbed}, we obtain also
\begin{equation}
\label{e:abstr_Lip:0}
    \| u_1(t) - u_2(t)\|_{ E_{ 1 - \frac{1}{p} ,\, p } }
  \leq L_1\,
    \| u_{0,1} - u_{0,2}\|_{ E_{ 1 - \frac{1}{p} ,\, p } }
  \quad\mbox{ for all }\, t\in [0,T_1] \,,
\end{equation}
with another Lipschitz constant
$L_1\equiv L_1(R_0)\in (0,\infty)$.
\hfill\Square
\endgroup
\end{remark}
\par\vskip 10pt

A number of sufficient conditions that guarantee the existence of
a {\em global\/} weak solution
$\mathbf{u}\colon \RR^N\times (0,T)\to \RR^M$ $(\CC^M)$
for all times $t\in (0,T)$
to the parabolic Cauchy problem \eqref{e:Cauchy}
can be found in
{\sc H.\ Amann} \cite{Amann-1, Amann-2} for systems similar to ours.
As we do not wish to impose those kinds of restrictive growth conditions
on the reaction function $\mathbf{f}$
on the right\--hand side of eq.~\eqref{e:Cauchy},
we prefer to \underline{assume\/} the {\it existence\/} of
a fixed global strict solution
(cf.\ \eqref{e_1:strict_sol})
\begin{align}
\label{e_w:strict_sol}
& w\in Y_{ 1 - \frac{1}{p} }^p(0,T)
\\
\nonumber
& \left( {}=
  L^p( (0,T)\to E_1 )\cap W^{1,p}( (0,T)\to E_0 )
  \hookrightarrow
  C\left( [0,T]\to E_{ 1 - \frac{1}{p} ,\, p } \right)
  \right)
\end{align}
to problem \eqref{e:abstr_Cauchy}
on the whole time interval $[0,T]$, for some $T\in (0,\infty)$,
with a prescribed initial value
$w(0) = w_0\in U\subset E_{ 1 - \frac{1}{p} ,\, p }$ and such that
$w(t)\in U$ and $A(t,w(t))\in \mathrm{MR}_p(E)$ for all $t\in [0,T]$.
Then the local Theorem~\ref{thm-Clem-Li} and Remark~\ref{rem-Clem-Li}
from above may be applied on any time interval
$[t_0 ,\, t_0 + T_1]\subset [0,T]$
of sufficiently short length $T_1 > 0$ in order to obtain
unique strict solutions $u$ ``along'' the known solution $w$
to the following abstract initial value problem:
\begin{equation}
\label{e:abstr_Cauchy:t_0}
\left\{
\begin{alignedat}{2}
    \frac{\mathrm{d}u}{\mathrm{d}t}
  - A(t,u(t)) u(t) &= f(t,u(t)) + g(t)
  \quad\mbox{ for a.e. } t\in (t_0 ,\, t_0 + T_1) \,;
\\
  u(t_0) &= u_0\in E_{ 1 - \frac{1}{p} ,\, p } \,.
\end{alignedat}
\right.
\end{equation}
Here, $u_0\in \overline{B}_{R_0}(w(t_0))$ is arbitrary,
where the radius $R_0 > 0$ is small enough, as described in
Remark~\ref{rem-Clem-Li}, such that
$\overline{B}_{R_0}(w(t_0)) \subset U$.
By Theorem~\ref{thm-Clem-Li}, the strict solution
$u\colon [t_0 ,\, t_0 + T_1]\to E_{ 1 - \frac{1}{p} ,\, p }$
satisfies $u(t)\in U$ for every $t\in [t_0 ,\, t_0 + T_1]$.
The image
$w([0,T]) = \{ w(t)\colon t\in [0,T] \}$
of the solution $w$ being compact in the open set
$U\subset E_{ 1 - \frac{1}{p} ,\, p }$,
we may choose $R_0 > 0$ even smaller, such that
$\overline{B}_{R_0}(w(t)) \subset U$ holds for all $t\in [0,T]$.

In addition to these claims that follow immediately from the proof of
\cite[Theorem 2.1, pp.\ 20--23]{Clem-Li},
one can deduce from inequalities analogous to those in
\cite[p.~22]{Clem-Li}, (2.14)--(2.17),
cf.\ Remark~\ref{rem-Clem-Li} above, ineq.~\eqref{e:abstr_Lip:0},
that there exists a Lipschitz constant $L_1\in [1,\infty)$, such that if
\begin{math}
  u_1, u_2\colon\hfil\break
  [t_0 ,\, t_0 + T_1]\to E_{ 1 - \frac{1}{p} ,\, p }
\end{math}
are two strict solutions to problem \eqref{e:abstr_Cauchy:t_0}
with initial values $u_1(t_0) = u_{0,1}$ and $u_2(t_0) = u_{0,2}$ in
$\overline{B}_{R_0}(w(t_0)) \subset U$, then one has
$u_1(t), u_2(t)\in U$ and
\begin{equation}
\label{e:abstr_Lip:t_0}
    \| u_1(t) - u_2(t)\|_{ E_{ 1 - \frac{1}{p} ,\, p } }
  \leq L_1\,
    \| u_1(t_0) - u_2(t_0)\|_{ E_{ 1 - \frac{1}{p} ,\, p } }
  \quad\mbox{ for all }\, t\in [t_0 ,\, t_0 + T_1] \,.
\end{equation}
Consequently, fixing the smallest integer $m\in \NN$ such that
$m\geq T/T_1$ (${}\geq 1$),
we obtain, by ``induction'' on $k = 1,2,3,\dots,m$, first
\begin{equation}
\label{e:abstr:k.T_1}
\begin{aligned}
    \| u(t) - w(t)\|_{ E_{ 1 - \frac{1}{p} ,\, p } }
  \leq L_1^k\,
    \| u_0 - w_0\|_{ E_{ 1 - \frac{1}{p} ,\, p } }
  \leq L_1^k\cdot R_0 / L_1^k = R_0
\\
  \quad\mbox{ for all }\,
    t\in \left[ 0 ,\, \min\{ k\, T_1 ,\, T\} \right] \,,
  \quad\mbox{ whenever }\;
    \| u_0 - w_0\|_{ E_{ 1 - \frac{1}{p} ,\, p } }
  \leq R_0 / L_1^k \,,
\end{aligned}
\end{equation}
then also
$u_1(t), u_2(t)\in \overline{B}_{R_0}(w(t)) \subset U$ and
\begin{equation}
\label{e:abstr_Lip:k.T_1}
\begin{aligned}
    \| u_1(t) - u_2(t)\|_{ E_{ 1 - \frac{1}{p} ,\, p } }
  \leq L_1^k\,
    \| u_{0,1} - u_{0,2}\|_{ E_{ 1 - \frac{1}{p} ,\, p } }
  \quad\mbox{ for all }\,
    t\in \left[ 0 ,\, \min\{ k\, T_1 ,\, T\} \right] \,,
\\
  \quad\mbox{ whenever }\;
    \| u_{0,j} - w_0\|_{ E_{ 1 - \frac{1}{p} ,\, p } }
  \leq R_k\eqdef R_0 / L_1^k\hspace{5pt} ({}> 0) \,;\quad j=1,2 \,,
\end{aligned}
\end{equation}
for every $k = 1,2,3,\dots,m$.

We have thus obtained the following result, {\it global in time\/}
on an arbitrary time interval\/ $(t_0,T)$, $0\leq t_0 < T$,
with the constants $R_0\in (0,\infty)$,
$T_1\equiv T_1(R_0)\in (0,T]$, and $L_1\in [1,\infty)$
specified above in eqs.\
\eqref{e:abstr_Lip:t_0} -- \eqref{e:abstr_Lip:k.T_1}:

\begin{theorem}\label{thm-Clem-Li_glob}
Let\/ $1 < p < \infty$, $0 < T < \infty$, and\/
$g\in L^p((0,T)\to E_0)$.
Assume that\/ $U$ is a nonempty open subset of\/
$E_{ 1 - \frac{1}{p} ,\, p }$ and\/
$A$ and\/ $f$ satisfy\/
{\rm Hypotheses} \eqref{hyp:A} and\/ \eqref{hyp:f}, respectively.
Finally, assume that\/
$w\colon [0,T]\to U\subset E_{ 1 - \frac{1}{p} ,\, p }$
is a fixed {\em global strict solution\/} to problem
\eqref{e:abstr_Cauchy} satisfying \eqref{e_w:strict_sol},
with a prescribed initial value $w(0) = w_0\in U$ and such that\/
$w(t)\in U$ and $A(t,w(t))\in \mathrm{MR}_p(E)$ for all\/ $t\in [0,T]$.
Then there exist some constant\/ $R_0\in (0,\infty)$,
sufficiently small,
with the following two properties, where
$R_m = R_0 / L_1^m\in (0,R_0]$ is the constant defined in
{\rm eq.}~\eqref{e:abstr_Lip:k.T_1} above:
\begin{itemize}
\item[{\rm (i)}]
If\/ $t_0\in [0,T)$ and\/
$u_0\in \overline{B}_{R_m}(w(t_0)) \subset U$,
then the abstract initial value problem \eqref{e:abstr_Cauchy}
on the time\--interval $(t_0,T)$ with $u(t_0) = u_0$
possesses a unique strict solution
$u\in Y_{ 1 - \frac{1}{p} }^p(t_0,T)$
(cf.\ \eqref{e_1:strict_sol})
\begin{align*}
& u\in Y_{ 1 - \frac{1}{p} }^p(t_0,T)
\\
\nonumber
& \left( {}=
  L^p( (t_0,T)\to E_1 )\cap W^{1,p}( (t_0,T)\to E_0 )
    \hookrightarrow
  C\left( [t_0,T]\to E_{ 1 - \frac{1}{p} ,\, p } \right)
  \right)
\end{align*}
such that\/
$u(t)\in \overline{B}_{R_0}(w(t)) \subset U$
for every\/ $t\in [t_0,T]$.
\item[{\rm (ii)}]
If\/ $t_0\in [0,T)$ and\/
$u_1, u_2\colon [t_0,T]\to E_{ 1 - \frac{1}{p} ,\, p }$
are two strict solutions to problem \eqref{e:abstr_Cauchy}
on the time\--interval $(t_0,T)$
with initial values $u_1(t_0) = u_{0,1}$ and $u_2(t_0) = u_{0,2}$ in
$\overline{B}_{R_m}(w(t_0)) \subset U$, then one has
$u_1(t), u_2(t)\in \overline{B}_{R_0}(w(t)) \subset U$
and\/
\begin{equation}
\label{e:abstr_Lip:T}
    \| u_1(t) - u_2(t)\|_{ E_{ 1 - \frac{1}{p} ,\, p } }
  \leq L_1^m\,
    \| u_{0,1} - u_{0,2}\|_{ E_{ 1 - \frac{1}{p} ,\, p } }
  \quad\mbox{ for all }\, t\in [t_0,T] \,.
\end{equation}
\end{itemize}
\end{theorem}
\par\vskip 10pt

\section{Analyticity in time for the abstract Cauchy problem}
\label{s:Analyt_abstr}

In this section we establish a few temporal analyticity results,
Theorem~\ref{thm-Analytic_glob}
being the most important among them,
that will be used later (in Section~\ref{s:proof-Main})
in order to prove Part~{\bf (ii)} of Theorem~\ref{thm-Main}.

\subsection{Auxiliary linear perturbation results}
\label{ss:Lin_Perturb}

We begin by quoting a well\--known result: If
$B\in \mathrm{Gen}(E)$
is the generator of a holomorphic semigroup on $E_0$
with the domain $\mathcal{D}(B) = E_1$, i.e.,
$B\in \mathrm{Hol}(E)$, then so is every operator
$B_{\nu} = (1 + \ii\nu) B\colon E_1\subset E_0\to E_0$, $\nu\in \RR$,
provided $|\nu|$ is small enough, $|\nu|\leq \delta_1 < 1$; see, e.g.,
{\sc H.~Amann} \cite[Chapt.~I, {\S}1]{Amann-Birkh}, pp.\ 9--24,
{\sc A.\ Pazy} \cite[{\S}3.2, pp.\ 80--81]{Pazy}, or
{\sc H.~Tanabe} \cite[Chapt.~3, {\S}3.1--{\S}3.4]{Tana-79}, pp.\ 51--72.
A more general perturbation theorem for
generators of holomorphic semigroups is proved in 
{\sc A.\ Pazy} \cite[{\S}3.2]{Pazy}, Theorem 2.1 on p.~80.
An analogous perturbation result for the smaller class
$\mathrm{MR}_p(E)$
\begin{math}
  \,\bigl(\,
  \mathrm{MR}_p(E)\subset \mathrm{Hol}(E)\subset \mathrm{Gen}(E)
  \,\bigr)\,
\end{math}
is proved in
{\sc H.~Amann} \cite[Chapt.~III, {\S}1.6]{Amann-Birkh},
Proposition 1.6.3 on p.~97.
Since we take advantage of the latter in an essential manner,
we now give its precise formulation.

Let $1 < p < \infty$ and $0 < T < \infty$.
Given any generator $B\in \mathrm{Gen}(E)$, let us consider
the bounded linear operator
$\tilde{K}_B\colon L^1((0,T)\to E_0)\to L^{\infty}((0,T)\to E_0)$
defined by
\begin{equation}
\label{def:K_B}
\begin{aligned}
  (\tilde{K}_B g)(t)\eqdef \int_0^t \ee^{(t-s)B} g(s) \,\mathrm{d}s
    \in E_0
    \quad\mbox{ for all }\, t\in [0,T]
\\
    \,\mbox{ and for every }\, g\in L^1((0,T)\to E_0) \,.
\end{aligned}
\end{equation}
It is proved in \cite[Chapt.~III, {\S}1.5]{Amann-Birkh},
Theorem 1.5.2 on p.~95,
that if $B\in \mathrm{Hol}(E)$ and $B$ possesses
the maximal $L^p$-regularity property, i.e.,
$B\in \mathrm{MR}_p(E)\equiv \mathrm{MR}_p(E_1\to E_0)$,
then the restriction
\begin{equation*}
  K_B = \tilde{K}_B\vert_{X^T} \;\mbox{ of $\tilde{K}_B$ to }\;
  X^T\eqdef X_{ 1 - \frac{1}{p} }^p(0,T) = L^p((0,T)\to E_0)
\end{equation*}
is a bounded linear operator from the Banach space $X^T$
into another Banach space
\begin{equation*}
    Y^T\eqdef Y_{ 1 - \frac{1}{p} }^p(0,T)
  = L^p((0,T)\to E_1)\cap W^{1,p}((0,T)\to E_0)
\end{equation*}
with the operator norm
\begin{math}
  \| K_B\|_{ \mathcal{L}(X^T\to Y^T) } < \infty .
\end{math}
For the perturbed initial value problem
\begin{equation}
\label{B_A:abstr_Cauchy}
\left\{
\begin{alignedat}{2}
    \frac{\mathrm{d}u}{\mathrm{d}t} - (B+A) u(t) &= g(t)
  \quad\mbox{ for a.e. } t\in (0,T) \,;
\\
  u(0) &= u_0\in E_{ 1 - \frac{1}{p} ,\, p } \,,
\end{alignedat}
\right.
\end{equation}
the following result is established in
\cite[Chapt.~III, {\S}1.6]{Amann-Birkh}, Proposition 1.6.3 on p.~97:

\begin{lemma}\label{lem-perturb_MR}
Assume that\/
$B\in \mathrm{MR}_p(E)$ and let\/
$A\in \mathcal{L}(E_1\to E_0)$ be arbitrary with the norm
\begin{equation*}
  \| A\|_{ \mathcal{L}(E_1\to E_0) } \leq
  \gamma \Big\slash \| K_B\|_{ \mathcal{L}(X^T\to Y^T) }
    \quad\mbox{ for some }\, \gamma\in (0,1) \,.
\end{equation*}
Then also the operator\/
$B_A = B+A\in \mathcal{L}(E_1\to E_0)$ belongs to the class
$\mathrm{MR}_p(E)$ and the operator norms of the inverses of
the abstract (linear) partial differential operators
\begin{equation*}
  \left( \partial_t - B ,\, \tau\right) \,,\;
  \left( \partial_t - B-A ,\, \tau\right) \colon
    Y_{ 1 - \frac{1}{p} }^p(0,T) \to
    L^p((0,T)\to E_0) \times E_{ 1 - \frac{1}{p} ,\, p }
\end{equation*}
defined in {\rm eq.}~\eqref{abstr:strict_sol} satisfy
\begin{equation}
\label{inv:B+A}
  \Vert \left( \partial_t - B-A ,\, \tau\right)^{-1} \Vert
    \leq C\cdot (1-\gamma)^{-1}
  \Vert \left( \partial_t - B ,\, \tau\right)^{-1} \Vert \,,
\end{equation}
where $C\equiv C(p,E,T) > 0$ is a constant independent from
$A$, $B$, and\/ $\gamma$.
\end{lemma}
\par\vskip 10pt

More precisely, we have
\begin{equation}
\label{inv_Neumann:B+A}
    \left( \partial_t - B-A ,\, \tau\right)^{-1}
  = (I - K_B A)^{-1}
    \left( \partial_t - B   ,\, \tau\right)^{-1}
\end{equation}
with the operator norm of the product
\begin{equation*}
  K_B A\colon Y^T\to Y^T = Y_{ 1 - \frac{1}{p} }^p(0,T)
  = L^p((0,T)\to E_1)
\end{equation*}
bounded above by
\begin{equation*}
  \| K_B A\|_{ \mathcal{L}(Y^T\to Y^T) }
  \leq \| K_B\|_{ \mathcal{L}(X^T\to Y^T) }
  \cdot \| A\|_{ \mathcal{L}(E_1\to E_0) }
  \leq \gamma < 1 \,.
\end{equation*}
Here, $I$ stands for the identity mapping in
$\mathcal{L}(Y^T\to Y^T)$.
Hence, the Neumann series
\begin{math}
  (I - K_B A)^{-1} = \sum_{k=0}^{\infty} (K_B A)^k
\end{math}
converges absolutely in $\mathcal{L}(Y^T\to Y^T)$ and
\begin{math}
  \| (I - K_B A)^{-1}\|_{ \mathcal{L}(Y^T\to Y^T) }
  \leq (1-\gamma)^{-1} < \infty .
\end{math}

The following claims are trivial applications of this lemma:
$\mathrm{MR}_p(E)\equiv \mathrm{MR}_p(E_1\to E_0)$
is an open subset of the Banach space $\mathcal{L}(E_1\to E_0)$.
Furthermore, if $B\in \mathrm{MR}_p(E)$ and
$A\in \mathcal{L}(E_1\to E_0)$ then also
$B_{\nu A} = B + \nu A\in \mathrm{MR}_p(E)$
holds for every $\nu\in \CC$ provided $|\nu|$ is small enough,
\begin{math}
  |\nu|\leq \delta_1 =\hfil\break
  \gamma\, \| A\|_{ \mathcal{L}(E_1\to E_0) }^{-1}
           \| K_B\|_{ \mathcal{L}(X^T\to Y^T) }^{-1} < \infty \,.
\end{math}
Naturally, the special case $A = \ii B$ is of interest.

The following perturbation lemma for problem \eqref {B_A:abstr_Cauchy}
is related to
{\sc S.\ Angenent} \cite[Lemma 2.5, p.~97]{Angenent}; see also
{\sc R.\ Denk}, {\sc M.\ Hieber}, and {\sc J.\ Pr\"uss}
\cite{Denk-Hieber},
Proposition 4.3 on p.~44 and Theorem 4.4 on p.~45.

\begin{lemma}\label{lem-perturb_MR:A<B}
Assume that\/ $B\in \mathrm{MR}_p(E)$.
Then there exists a number\/
$\delta\in (0,1)$ and a constant\/
$C_{\delta}\in \RR_+$ with the following property:\/
If\/
$A\in \mathcal{L}(E_1\to E_0)$ is arbitrary with the norm
\begin{equation}
\label{e:perturb_MR:A<B}
  \| Au\|_{E_0}\leq \delta\, \| Bu\|_{E_0} + C_{\delta}\, \| u\|_{E_0}
    \quad\mbox{ for all }\, u\in E_1 \,,
\end{equation}
then also the operator\/
$B_A = B+A\in \mathcal{L}(E_1\to E_0)$
belongs to the class $\mathrm{MR}_p(E)$.
Furthermore, there exists a constant\/
\begin{math}
  \tilde{M}\equiv \tilde{M}(p,E,B,\delta,C_{\delta},T) > 0 ,
\end{math}
independent from
\begin{math}
  (g,u_0)\in L^p((0,T)\to E_0)\times \hfil\break
  E_{ 1 - \frac{1}{p} ,\, p } ,
\end{math}
such that the unique strict solution
$v = \left( \partial_t - B-A ,\, \tau\right)^{-1} (g,u_0)$
to the perturbed initial value problem \eqref{B_A:abstr_Cauchy}
satisfies the inequality
\begin{equation}
\label{est:B+A_sol}
    \int_0^T \genfrac{\|}{\|}{}0{\mathrm{d}v}{\mathrm{d}t}_{E_0}^p
    \,\mathrm{d}t
  + \int_0^T \| (B+A) v(t)\|_{E_0}^p \,\mathrm{d}t
  \leq \tilde{M}
    \left(
      \| u_0\|_{ E_{ 1 - \frac{1}{p} ,\, p } }^p
    + \int_0^T \| g(t)\|_{E_0}^p \,\mathrm{d}t
    \right) \,,
\end{equation}
whenever\/
$u_0\in E_{ 1 - \frac{1}{p} ,\, p }$ and\/
$g\in L^p((0,T)\to E_0)$.
\end{lemma}
\par\vskip 10pt

\proof
{\it Step~$1$.}$\;$
First, we prove the lemma for $[0,T]\subset \RR_+$
replaced by a sufficiently short time interval
$[t_0, t_0 + T_1]\subset [0,T]$, i.e.,
$0\leq t_0 < t_0 + T_1\leq T$ with $T_1\in (0,\infty)$ small enough.
Without loss of generality, we may assume $t_0 = 0$ and $0 < T_1\leq T$.

Let us recall our notation and the continuous imbedding
(cf.\ \eqref{e_1:strict_sol})
\begin{equation}
\label{e:Y->C}
\begin{aligned}
  Y^{T_1} = Y_{ 1 - \frac{1}{p} }^p(0,T_1)
& = L^p((0,T_1)\to E_1)\cap W^{1,p}((0,T_1)\to E_0)
\\
& \hookrightarrow
  C\left( [0,T_1]\to E_{ 1 - \frac{1}{p} ,\, p } \right) \,.
\end{aligned}
\end{equation}
It is easy to see that a function $v\in Y^{T_1}$
is a strict solution of the perturbed initial value problem
\eqref{B_A:abstr_Cauchy} on $(0,T_1)$
{\it if and only if\/} it satisfies
\begin{equation}
\label{fix:abstr_Cauchy}
\left\{
\begin{alignedat}{2}
    \frac{\mathrm{d}v}{\mathrm{d}t} - B v(t) &= A v(t) + g(t)
  \quad\mbox{ for a.e. } t\in (0,T_1) \,;
\\
  v(0) &= u_0\in E_{ 1 - \frac{1}{p} ,\, p } \,,
\end{alignedat}
\right.
\end{equation}
in the strict sense, again.
Notice that
$\tilde{g} = Av + g\in L^p((0,T_1)\to E_0)$.
We observe that problem \eqref{fix:abstr_Cauchy}
has a unique strict solution $v\in Y^{T_1}$
{\it as soon as\/} we have shown that
the affine self\--mapping
$F\colon v\mapsto \hat{v}\colon Y^{T_1}\to Y^{T_1}$, defined by
\begin{equation}
\label{v^:abstr_Cauchy}
\left\{
\begin{alignedat}{2}
    \frac{ \mathrm{d}\hat{v} }{\mathrm{d}t} - B\hat{v}(t)
& = A v(t) + g(t)
  \quad\mbox{ for a.e. } t\in (0,T_1) \,;
\\
  \hat{v}(0) &= u_0\in E_{ 1 - \frac{1}{p} ,\, p } \,,
\end{alignedat}
\right.
\end{equation}
possesses a unique fixed point $v\in Y^{T_1}$.
Obviously, such a fixed point must belong to the (closed) affine subspace
\begin{equation*}
  Y^{T_1}_{(u_0)} = \left\{ v\in Y^{T_1}\colon v(0) = u_0 \right\}
  \quad\mbox{ of the Banach space $Y^{T_1}$; }
\end{equation*}
hence,
\begin{math}
  Y^{T_1}_{(u_0)} = u_0 + Y^{T_1}_{(0)} .
\end{math}
Clearly,
\begin{math}
  Y^{T_1}_{(0)} = \left\{ v\in Y^{T_1}\colon v(0) = 0 \right\}
\end{math}
is a closed vector subspace of $Y^{T_1}$.
The former one inherits the norm from the latter.

Next, we prove that
$F\colon v\mapsto \hat{v}$ is a contraction on $Y^{T_1}_{(u_0)}$.
To this end, let $v_i\in Y^{T_1}_{(u_0)}$ be arbitrary and set
$\hat{v}_i = F(v_i)$; $i=1,2$.
The differences $z = v_1 - v_2$ and
$\hat{z} =  \hat{v}_1 - \hat{v}_2$ are in
$Y^{T_1}_{(0)}$ and, by \eqref{v^:abstr_Cauchy}, they satisfy
\begin{equation}
\label{diff:abstr_Cauchy}
\left\{
\begin{alignedat}{2}
    \frac{ \mathrm{d}\hat{z} }{\mathrm{d}t} - B\hat{z}(t)
& = A z(t)
  \quad\mbox{ for a.e. } t\in (0,T_1) \,;
\\
  \hat{z}(0) &= 0\in E_{ 1 - \frac{1}{p} ,\, p } \,.
\end{alignedat}
\right.
\end{equation}
By Remark~\ref{rem-perturb}, {\rm Part~(a)}, the operator
$B\in \mathrm{MR}_p(E)$ satisfies ineq.~\eqref{est:strict_sol}
with a constant\hfil\break
$M(p,E,B,T_1)\leq M(p,E,B,T)\equiv M_T < \infty$.
Hence, we have
\begin{equation*}
    \int_0^{T_1}
    \genfrac{\|}{\|}{}0{ \mathrm{d}\hat{z} }{\mathrm{d}t}_{E_0}^p
    \,\mathrm{d}t
  + \int_0^{T_1} \| B\hat{z}(t)\|_{E_0}^p \,\mathrm{d}t
  \leq M_T\int_0^{T_1} \| A z(t)\|_{E_0}^p \,\mathrm{d}t \,.
\end{equation*}
Now we estimate the integrand on the right\--hand side by
ineq.~\eqref{e:perturb_MR:A<B},
\begin{align}
\label{diff:strict_sol}
    \int_0^{T_1}
    \genfrac{\|}{\|}{}0{ \mathrm{d}\hat{z} }{\mathrm{d}t}_{E_0}^p
    \,\mathrm{d}t
  + \int_0^{T_1} \| B\hat{z}(t)\|_{E_0}^p \,\mathrm{d}t
  \leq M_T\int_0^{T_1}
    \left(
  \delta\, \| B z(t)\|_{E_0} + C_{\delta}\, \| z(t)\|_{E_0}
    \right)^p \,\mathrm{d}t
\\
\nonumber
  \leq 2^{p-1}\, M_T
    \left(
    \delta^p\int_0^{T_1} \| B z(t)\|_{E_0}^p \,\mathrm{d}t
  + C_{\delta}^p\int_0^{T_1} \| z(t)\|_{E_0}^p \,\mathrm{d}t
    \right) \,.
\end{align}
The integrand in the second integral on the right\--hand side
is estimated by H\"older's inequality:
\begin{align*}
  \| z(t)\|_{E_0}
& = \left\| \int_0^t
    \genfrac{}{}{}0{\mathrm{d}z}{\mathrm{d}t}(s) \,\mathrm{d}s
    \right\|_{E_0}
  \leq \int_0^t \| z'(s)\|_{E_0} \,\mathrm{d}s
\\
& \leq
    \left( \int_0^t
    \left\| \genfrac{}{}{}0{\mathrm{d}z}{\mathrm{d}t}(s)
    \right\|_{E_0}^p \,\mathrm{d}s
    \right)^{1/p}
    \left( \int_0^t \,\mathrm{d}s \right)^{1/p'}
  \quad\mbox{ for all }\, t\in [0,T_1] \,,
\end{align*}
where $p'= p/(p-1)\in (1,\infty)$.
Here, we have used $z(0) = 0\in E_0$.
Hence,
\begin{equation*}
  \| z(t)\|_{E_0}^p
  \leq t^{p/p'}\cdot
    \left( \int_0^t
    \left\| \genfrac{}{}{}0{\mathrm{d}z}{\mathrm{d}t}(s)
    \right\|_{E_0}^p \,\mathrm{d}s
    \right) \,.
\end{equation*}
After integration we thus obtain, thanks to $p/p'= p-1$,
\begin{equation}
\label{est:|z(t)|^p}
  \int_0^{T_1} \| z(t)\|_{E_0}^p \,\mathrm{d}t
  \leq \frac{1}{p}\, T_1^p
    \int_0^{T_1}
    \left\| \genfrac{}{}{}0{\mathrm{d}z}{\mathrm{d}t}(s)
    \right\|_{E_0}^p \,\mathrm{d}s \,.
\end{equation}
Of course, the same inequality is valid also for
$\hat{z}\in Y^{T_1}_{(0)}$ in place of the function $z$.
We apply the last inequality to
the right\--hand side of \eqref{diff:strict_sol},
\begin{align}
\label{rh:strict_sol}
&   \int_0^{T_1}
    \genfrac{\|}{\|}{}0{ \mathrm{d}\hat{z} }{\mathrm{d}t}_{E_0}^p
    \,\mathrm{d}t
  + \int_0^{T_1} \| B\hat{z}(t)\|_{E_0}^p \,\mathrm{d}t
\\
\nonumber
& \leq 2^{p-1}\, M_T
    \left(
    \delta^p\int_0^{T_1} \| B z(t)\|_{E_0}^p \,\mathrm{d}t
  + \frac{1}{p}\, C_{\delta}^p\, T_1^p
    \int_0^{T_1}
    \genfrac{\|}{\|}{}0{\mathrm{d}z}{\mathrm{d}t}_{E_0}^p
    \,\mathrm{d}t
    \right) \,.
\end{align}

The integrals on both sides containing the generator
$B\in \mathrm{MR}_p(E)$ are estimated as follows.
First, there are constants $c_1, C_1\in (0,\infty)$ and
$c_2, C_2\in \RR_+$ such that the inequalities
\begin{equation*}
    c_1\, \| u\|_{E_1} - c_2\, \| u\|_{E_0}
  \leq \| Bu\|_{E_0}\leq
    C_1\, \| u\|_{E_1} + C_2\, \| u\|_{E_0}
  \quad\mbox{ hold for all }\, u\in E_1 \,.
\end{equation*}
Consequently, we have also
\begin{align*}
    c_1^p\, \| u\|_{E_1}^p
& \leq 2^{p-1}
  \left( \| Bu\|_{E_0}^p + c_2^p\, \| u\|_{E_0}^p \right)
  \quad\mbox{ and }
\\
    \| Bu\|_{E_0}^p
& \leq 2^{p-1}
  \left( C_1^p\, \| u\|_{E_1}^p + C_2^p\, \| u\|_{E_0}^p \right)
  \quad\mbox{ for all }\, u\in E_1 \,.
\end{align*}
Applying these inequalities to \eqref{rh:strict_sol}, we arrive at
\begin{align*}
&   \int_0^{T_1}
    \genfrac{\|}{\|}{}0{ \mathrm{d}\hat{z} }{\mathrm{d}t}_{E_0}^p
    \,\mathrm{d}t
  + 2^{-(p-1)}\, c_1^p\int_0^{T_1} \| \hat{z}(t)\|_{E_1}^p \,\mathrm{d}t
  - c_2^p\int_0^{T_1} \| \hat{z}(t)\|_{E_0}^p \,\mathrm{d}t
\\
& \leq 2^{p-1}\, M_T\cdot
     2^{p-1}\, C_1^p\, \delta^p
     \int_0^{T_1} \| z(t)\|_{E_1}^p \,\mathrm{d}t
\\
& {}
  + 2^{p-1}\, M_T
    \left(
    2^{p-1}\, C_2^p\, \delta^p
    \int_0^{T_1} \| z(t)\|_{E_0}^p \,\mathrm{d}t
  + \frac{1}{p}\, C_{\delta}^p\, T_1^p
    \int_0^{T_1}
    \genfrac{\|}{\|}{}0{\mathrm{d}z}{\mathrm{d}t}_{E_0}^p
    \,\mathrm{d}t
    \right) \,.
\end{align*}
Finally, we estimate the integrals
$\int_0^{T_1} \| \hat{z}(t)\|_{E_0}^p \,\mathrm{d}t$
and
$\int_0^{T_1} \| z(t)\|_{E_0}^p \,\mathrm{d}t$
above by ineq.~\eqref{est:|z(t)|^p}, thus obtaining
\begin{align}
\nonumber
& \left( 1 - \frac{1}{p}\, T_1^p\, c_2^p\right)
    \int_0^{T_1}
    \genfrac{\|}{\|}{}0{ \mathrm{d}\hat{z} }{\mathrm{d}t}_{E_0}^p
    \,\mathrm{d}t
  + 2^{-(p-1)}\, c_1^p\int_0^{T_1} \| \hat{z}(t)\|_{E_1}^p \,\mathrm{d}t
\\
\label{lh<rh:strict_sol}
& \leq 2^{2(p-1)}\, C_1^p\, \delta^p\, M_T
    \int_0^{T_1} \| z(t)\|_{E_1}^p \,\mathrm{d}t
\\
\nonumber
& {}
  + \frac{2^{p-1}}{p}\, T_1^p\, M_T
    \left( 2^{p-1}\, C_2^p\, \delta^p + C_{\delta}^p \right)
    \int_0^{T_1}
    \genfrac{\|}{\|}{}0{\mathrm{d}z}{\mathrm{d}t}_{E_0}^p
    \,\mathrm{d}t \,.
\end{align}

We finish this step by choosing first
$\delta\in (0,1)$ then $T_1\in (0,T]$ small enough, such that
\begin{align*}
& 2^{2(p-1)}\, C_1^p\, \delta^p\, M_T \leq
  \frac{1}{2}\cdot 2^{-(p-1)}\, c_1^p \quad\mbox{ and }\quad
\\
&   \frac{2^{p-1}}{p}\, T_1^p\, M_T
    \left( 2^{p-1}\, C_2^p\, \delta^p + C_{\delta}^p \right)
  \leq \frac{1}{2}
  \left( 1 - \frac{1}{p}\, T_1^p\, c_2^p\right) \,,
\end{align*}
respectively, or, equivalently,
\begin{align}
\label{ineq:delta}
& 0 < \delta\leq 2^{-(3p-2)/p}\, (c_1/C_1)\, M_T^{-1/p}
    \quad\mbox{ and }\quad
\\
\label{ineq:T_1}
& T_1^p
  \left[
    2^p\, M_T\left( 2^{p-1}\, C_2^p\, \delta^p + C_{\delta}^p \right)
  + c_2^p
  \right] \leq p \,.
\end{align}
With these choices of $\delta$ and $T_1$, we obtain
\begin{equation}
\label{contract:Y_p^p^}
  \| \hat{z}\|_{Y^{T_1}}^{\flat}
  \leq \genfrac{}{}{}1{1}{2}\, \| z\|_{Y^{T_1}}^{\flat}
\end{equation}
in the new, equivalent norm
\begin{equation}
\label{norm_flat:Y_p^p^}
  \| u\|_{Y^{T_1}}^{\flat} \eqdef
  \left[
  \left( 1 - \frac{1}{p}\, T_1^p\, c_2^p\right)
    \int_0^{T_1}
    \genfrac{\|}{\|}{}0{\mathrm{d}u}{\mathrm{d}t}_{E_0}^p
    \,\mathrm{d}t
  + 2^{-(p-1)}\, c_1^p\int_0^{T_1} \| u(t)\|_{E_1}^p \,\mathrm{d}t
  \right]^{1/p}
\end{equation}
on the abstract Sobolev space
\begin{math}
  Y^{T_1} =  Y_{ 1 - \frac{1}{p} }^p(0,T_1)
  = L^p((0,T_1)\to E_1)\cap W^{1,p}((0,T_1)\to E_0) ;
\end{math}
see \eqref{norm:Y_p^p^} and below in eq.~\eqref{est:sol}.
Inequality \eqref{contract:Y_p^p^} shows that
$F\colon v\mapsto \hat{v}$ is a contraction on
$Y^{T_1}_{(u_0)}$ ($\subset Y^{T_1}$)
with the Lipschitz constant $\frac{1}{2}$ with respect to the new norm
$\|\cdot\|_{Y^{T_1}}^{\flat}$.
Consequently, problem \eqref{fix:abstr_Cauchy}
has a unique strict solution $v\in Y^{T_1}$; in fact, we have
$v\in Y^{T_1}_{(u_0)}$.

The following estimate for $v$ can be proved by the same arguments
as those used in our proof of contraction above:
There is a constant
$\Gamma\equiv \Gamma(T_1)\in (0,\infty)$,
independent from
$u_0\in E_{ 1 - \frac{1}{p} ,\, p }$ and
$g\in L^p((0,T)\to E_0)$, such that
\begin{equation}
\label{est:sol}
\begin{aligned}
  \left( \| v\|_{Y^{T_1}}^{\sharp} \right)^p
& = \int_0^{T_1} \genfrac{\|}{\|}{}0{\mathrm{d}v}{\mathrm{d}t}_{E_0}^p
    \,\mathrm{d}t
  + \int_0^{T_1} \| v(t)\|_{E_1}^p \,\mathrm{d}t
\\
& {}
  \leq \Gamma
    \left(
      \| u_0\|_{ E_{ 1 - \frac{1}{p} ,\, p } }^p
    + \int_0^{T_1} \| g(t)\|_{E_0}^p \,\mathrm{d}t
    \right) \,.
\end{aligned}
\end{equation}
Recall that
$\|\cdot\|_{Y^{T_1}}^{\sharp}$
is an equivalent norm on the Banach space
$Y^{T_1}$; cf.\ {\rm eq.}~\eqref{norm:Y_p^p^}.
In analogy with Remark~\ref{rem-perturb}, {\rm Part~(b)},
we may take the constant $\Gamma\equiv \Gamma(T_1) > 0$
in ineq.~\eqref{est:sol} above
to be the {\it smallest\/} nonnegative number $\Gamma\in \RR_+$
for which ineq.~\eqref{est:sol} is valid.
It is easy to see that $\Gamma\equiv \Gamma(T_1)\in \RR_+$
is a {\it nondecreasing\/} function of time $T_1\in (0,T]$ and $\Gamma > 0$.
The last estimate, ineq.~\eqref{est:sol}, easily implies
ineq.~\eqref{est:B+A_sol} with $T_1$ in place of~$T$.
The imbedding \eqref{e:Y->C} being continuous,
by ineq.~\eqref{est:sol}, there is another constant
$\hat{\Gamma}\equiv \hat{\Gamma}(T_1)\in [1,\infty)$,
independent from
$u_0\in E_{ 1 - \frac{1}{p} ,\, p }$ and
$g\in L^p((0,T)\to E_0)$, such that
\begin{equation}
\label{est:sol:T_1}
      \| v(T_1)\|_{ E_{ 1 - \frac{1}{p} ,\, p } }^p
  \leq \hat{\Gamma}
    \left(
      \| u_0\|_{ E_{ 1 - \frac{1}{p} ,\, p } }^p
    + \int_0^{T_1} \| g(t)\|_{E_0}^p \,\mathrm{d}t
    \right) \,.
\end{equation}
Again, similarly to $\Gamma\equiv \Gamma(T_1) > 0$
in ineq.~\eqref{est:sol} above, we may take the constant
$\hat{\Gamma}\equiv \hat{\Gamma}(T_1) > 0$
in ineq.~\eqref{est:sol:T_1} above
to be the {\it smallest\/} number $\hat{\Gamma}\in [1,\infty)$
for which ineq.~\eqref{est:sol:T_1} is valid.
It is now easy to see that also the constant
$\hat{\Gamma}\equiv \hat{\Gamma}(T_1) \geq 1$
is a {\it nondecreasing\/} function of time $T_1\in (0,T]$.

{\it Step~$2$.}$\;$
We may take $T_1 = T/m$ sufficiently small in Step~1 above,
where $m\in \NN$ is a sufficiently large positive integer.
Next, we replace the interval $[0,T_1]$ from Step~1 by
any subinterval
$[t_0, t_0 + T_1]$ $= \mathcal{J}_k = [ (k-1) T_1 ,\, k T_1 ]$
of $[0,T]$ of length $T_1$ for $k = 1,2,\dots,m$; hence,
$\bigcup_{k=1}^m \mathcal{J}_k = [0,T]$.
We make use of the existence and uniqueness of a strict solution
\begin{equation*}
  v\in Y_{ 1 - \frac{1}{p} }^p(t_0, t_0 + T_1)
  = L^p((t_0, t_0 + T_1)\to E_1)\cap W^{1,p}((t_0, t_0 + T_1)\to E_0)
\end{equation*}
of the perturbed initial value problem
\eqref{B_A:abstr_Cauchy} in every subinterval
$\mathcal{J}_k$; $k = 1,2,\dots,m$,
together with the estimates \eqref{est:sol} and \eqref{est:sol:T_1}
on $\mathcal{J}_k$, by Step~1.
Thus, from \eqref{est:sol} and \eqref{est:sol:T_1}
we get, respectively,
\begin{align}
\label{est_k:sol}
\begin{aligned}
&     \int_{ (k-1) T_1 }^{k T_1}
      \genfrac{\|}{\|}{}0{\mathrm{d}v}{\mathrm{d}t}_{E_0}^p \,\mathrm{d}t
    + \int_{ (k-1) T_1 }^{k T_1} \| v(t)\|_{E_1}^p \,\mathrm{d}t
\\
& \leq \Gamma
    \left(
      \| v( (k-1) T_1 )\|_{ E_{ 1 - \frac{1}{p} ,\, p } }^p
    + \int_{ (k-1) T_1 }^{k T_1} \| g(t)\|_{E_0}^p \,\mathrm{d}t
    \right) \,,
\end{aligned}
\\
\label{est_k:sol:T_1}
      \| v(k T_1)\|_{ E_{ 1 - \frac{1}{p} ,\, p } }^p
  \leq \hat{\Gamma}
    \left(
      \| v( (k-1) T_1 )\|_{ E_{ 1 - \frac{1}{p} ,\, p } }^p
    + \int_{ (k-1) T_1 }^{k T_1} \| g(t)\|_{E_0}^p \,\mathrm{d}t
    \right) \,.
\end{align}
We recall that $\hat{\Gamma}\geq 1$.
Consequently, iterating inequalities \eqref{est_k:sol:T_1}
for all $k = 1,2,\dots,\ell$, $1\leq \ell\leq m$, we arrive at
\begin{equation}
\label{est_ell:sol:T_1}
\begin{aligned}
      \| v(\ell T_1)\|_{ E_{ 1 - \frac{1}{p} ,\, p } }^p
& \leq \hat{\Gamma}^{\ell}
      \| u_0\|_{ E_{ 1 - \frac{1}{p} ,\, p } }^p
  + \sum_{k=1}^{\ell} \hat{\Gamma}^{\ell - k + 1}
    \int_{(k-1) T_1}^{k T_1} \| g(t)\|_{E_0}^p \,\mathrm{d}t
\\
& {}
  \leq \hat{\Gamma}^{\ell}
    \left(
      \| u_0\|_{ E_{ 1 - \frac{1}{p} ,\, p } }^p
    + \int_0^{\ell T_1} \| g(t)\|_{E_0}^p \,\mathrm{d}t
    \right) \,.
\end{aligned}
\end{equation}

Next, we sum up inequalities \eqref{est_k:sol}
for all $k = 1,2,\dots,m$, thus obtaining
\begin{equation}
\label{est_ell:sol}
\begin{aligned}
  \left( \| v\|_{Y^T}^{\sharp} \right)^p
& = \int_0^T \genfrac{\|}{\|}{}0{\mathrm{d}v}{\mathrm{d}t}_{E_0}^p
    \,\mathrm{d}t
  + \int_0^T \| v(t)\|_{E_1}^p \,\mathrm{d}t
\\
& \leq \Gamma
    \left( \sum_{k=1}^m
      \| v( (k-1) T_1 )\|_{ E_{ 1 - \frac{1}{p} ,\, p } }^p
    \right)
    + \Gamma\int_0^T \| g(t)\|_{E_0}^p \,\mathrm{d}t \,.
\end{aligned}
\end{equation}
In order to estimate the first summand on the right\--hand side
from above, we apply ineq.~\eqref{est_ell:sol:T_1}
with $\ell = k-1$ for all $k = 1,2,\dots,m$, thus arriving at
\begin{align*}
& \sum_{k=1}^m
    \| v( (k-1) T_1 )\|_{ E_{ 1 - \frac{1}{p} ,\, p } }^p
  \leq
    \left( \sum_{\ell = 0}^{m-1} \hat{\Gamma}^{\ell} \right)
      \| u_0\|_{ E_{ 1 - \frac{1}{p} ,\, p } }^p
  + \sum_{\ell = 0}^{m-1}
    \hat{\Gamma}^{\ell} \int_0^{\ell T_1} \| g(t)\|_{E_0}^p \,\mathrm{d}t
\\
& {}
  \leq m\, \hat{\Gamma}^{m-1}\,
       \| u_0\|_{ E_{ 1 - \frac{1}{p} ,\, p } }^p
  + (m-1)\, \hat{\Gamma}^{m-1}
    \int_0^{(m-1) T_1} \| g(t)\|_{E_0}^p \,\mathrm{d}t
\\
& {}
  \leq \hat{M}\, \Gamma^{-1}
    \left(
      \| u_0\|_{ E_{ 1 - \frac{1}{p} ,\, p } }^p
    + \int_0^T \| g(t)\|_{E_0}^p \,\mathrm{d}t
    \right) \,,
\end{align*}
where
$\hat{M} = m\, \hat{\Gamma}^{m-1}\, \Gamma\in [1,\infty)$
is a constant independent from~$v$.
We apply this estimate to the right\--hand side of
{\rm eq.}~\eqref{est_ell:sol} to get
\begin{equation}
\label{est_ell:sol^sharp}
\begin{aligned}
  \left( \| v\|_{Y^T}^{\sharp} \right)^p
& \leq \hat{M}\, \| u_0\|_{ E_{ 1 - \frac{1}{p} ,\, p } }^p
     + (\hat{M} + \Gamma)\, \int_0^T \| g(t)\|_{E_0}^p \,\mathrm{d}t \,.
\end{aligned}
\end{equation}

We conclude the proof by applying
ineq.~\eqref{e:perturb_MR:A<B} with
$B,A\in \mathcal{L}(E_1\to E_0)$
to the left\--hand side of ineq.~\eqref{est:B+A_sol},
\begin{align*}
&   \int_0^T \genfrac{\|}{\|}{}0{\mathrm{d}v}{\mathrm{d}t}_{E_0}^p
    \,\mathrm{d}t
  + \int_0^T \| (B+A) v(t)\|_{E_0}^p \,\mathrm{d}t
  \leq
\\
&   \int_0^T \genfrac{\|}{\|}{}0{\mathrm{d}v}{\mathrm{d}t}_{E_0}^p
  + (1 + \delta)\int_0^T \| B v(t)\|_{E_0}^p \,\mathrm{d}t
  + C_{\delta}\int_0^T \| v(t)\|_{E_0}^p \,\mathrm{d}t
\\
&   \int_0^T \genfrac{\|}{\|}{}0{\mathrm{d}v}{\mathrm{d}t}_{E_0}^p
  + (1 + \delta)\, \| B\|_{ \mathcal{L}(E_1\to E_0) }
    \int_0^T \| v(t)\|_{E_1}^p \,\mathrm{d}t
  + C_{\delta}\int_0^T \| v(t)\|_{E_0}^p \,\mathrm{d}t
\\
& {}
  \leq M_{\delta}\, \left( \| v\|_{Y^T}^{\sharp} \right)^p \,,
\end{align*}
by {\rm eq.}~\eqref{est_ell:sol},
with a constant $M_{\delta}\in (0,\infty)$ independent from~$v$.
Now we apply ineq.~\eqref{est_ell:sol^sharp}
to the last estimate to arrive at the desired inequality
\eqref{est:B+A_sol} with the constant
$\tilde{M} = M_{\delta} (\hat{M} + \Gamma) > 0$.
We have proved that the operator
$B_A = B+A\in \mathcal{L}(E_1\to E_0)$
belongs to the class $\mathrm{MR}_p(E)$.
\qed
\par\vskip 10pt

\subsection{Proof of Analyticity in Time}
\label{ss:Analytic_t}

Now we are ready to prove that any global strict solution
$w\colon [0,T]\to U\subset E_{ 1 - \frac{1}{p} ,\, p }$
to problem \eqref{e:abstr_Cauchy}
that satisfies the hypotheses of Theorem~\ref{thm-Clem-Li_glob} above
must be analytic in time $t\in (0,T)$.
Let us recall that a {\em strict solution\/}
to problem \eqref{e:abstr_Cauchy} has been introduced
in Definition~\ref{def-strict_sol}.
Indeed, below we will prove a more detailed result
on a complex analytic (i.e., holomorphic) extension of $u(t)$
from the real time interval
$(0,T)\subset \RR\subset \CC$ to the open complex domain
$\Delta_{\vartheta}^{T',T}$
which is the intersection of the (open) triangle
$\Delta_{\vartheta}^{(T')}$ with the (open) complex strip
$\mathfrak{T}^{(r)}$ defined in eqs.\
\eqref{loc:T^(T')}, \eqref{glob:T^(r)}, and \eqref{def:T^(r)},
respectively,
where $\vartheta\in (0, \pi / 2)$ is a given angle and
$0 < T'\leq T < \infty$.
Here, the constants $\vartheta\in (0, \pi / 2)$ and $T'\in (0,T]$
will be chosen sufficiently small, but positive; hence, we have
$(0,T)\subset \Delta_{\vartheta}^{T',T}$.
Finally, we denote by
$\overline{\Delta}_{\vartheta}^{T',T}$ the closure of
$\Delta_{\vartheta}^{T',T}$ in $\CC$.

In addition to {\rm Hypotheses} \eqref{hyp:A} and\/ \eqref{hyp:f},
we assume that $A$ and\/ $f$ satisfy also
the following analyticity hypotheses, respectively
(cf.\ {\sc A.\ Lunardi} \cite[Chapt.~8]{Lunardi}, {\S}8.3.3, p.~308):

\begin{hypos}\nopagebreak
\begingroup\rm
Recall that both spaces, $E_0$ and $E_1$, in the Banach couple
$E = (E_0,E_1)$ are assumed to be complex Banach spaces
(over the field $\CC$)
with $E_1\hookrightarrow E_0$ densely and continuously.
Furthermore, we assume that there are positive constants
$\vartheta_0\in (0, \pi / 2)$ and $T_0\in (0,T]$, and open sets\,
$\mathcal{U}\subset \CC$ and
$\tilde{U}\subset E_{ 1 - \frac{1}{p} ,\, p }$
containing the compact set
$\overline{\Delta}_{\vartheta_0}^{T_0,T}$ and the open set $U$,
respectively, i.e.,
\begin{math}
  \overline{\Delta}_{\vartheta_0}^{T_0,T}
  \subset \mathcal{U}\subset \CC
\end{math}
and
\begin{math}
  U\subset \tilde{U}\subset E_{ 1 - \frac{1}{p} ,\, p } ,
\end{math}
such that
\begin{enumerate}
\setcounter{enumi}{0}
\renewcommand{\labelenumi}{{\bf (C\arabic{enumi}')}}
\item
\makeatletter
\def\@currentlabel{{\bf C\arabic{enumi}'}}\label{hyp:A_anal}
\makeatother
$A\colon [0,T]\times U\to \mathcal{L}(E_1\to E_0)$
possesses a holomorphic extension
\begin{math}
  \tilde{A}\colon \mathcal{U}\times \tilde{U}
  \to \mathcal{L}(E_1\to E_0)
\end{math}
to the complex domain $\mathcal{U}\times \tilde{U}$ which satisfies
$\tilde{A}(t,v)\in \mathrm{MR}_p(E)$ for all
$(t,v)\in \mathcal{U}\times \tilde{U}$.
\item
\makeatletter
\def\@currentlabel{{\bf C\arabic{enumi}'}}\label{hyp:f_anal}
\makeatother
$f\colon [0,T]\times U\to E_0$
possesses a holomorphic extension
\begin{math}
  \tilde{f}\colon \mathcal{U}\times \tilde{U}\to E_0
\end{math}
to the complex domain $\mathcal{U}\times \tilde{U}$.
\hfill\Square
\end{enumerate}
\endgroup
\end{hypos}
\par\vskip 10pt

Again, the metric on $\mathcal{U}\times \tilde{U}$ is induced by
the canonical norm on $\CC\times E_{ 1 - \frac{1}{p} ,\, p }$.
A precise definition of a {\it holomorphic\/}
(i.e., {\it complex analytic\/}) {\it mapping\/}
\begin{math}
  \mathscr{F}\colon \mathscr{O}\subset \mathcal{X}\to \mathcal{Y}
\end{math}
from an open subset $\mathscr{O}$ of a complex Banach space
$\mathcal{X}$ into another complex Banach space $\mathcal{Y}$
is given in
{\sc K.\ Deimling} \cite[Definition 15.1, p.~150]{Deimling}
(see also \cite[Proposition 15.2, p.~150]{Deimling}).

Without assuming
{\rm Hypotheses} \eqref{hyp:A} and\/ \eqref{hyp:f},
we observe that
{\rm Hypotheses} \eqref{hyp:A_anal} and\/ \eqref{hyp:f_anal}
still guarantee the following claims, respectively:
Given any compact set $K\subset \mathcal{U}$
and any continuous function $z\colon [0,T]\to K$,
one can easily verify that both substitution mappings,
\begin{align*}
  v &\mapsto
  \left[ t\mapsto A\left( z(t), v(z(t))\right)\right]
  \colon
\\
& C( K\to \tilde{U} ) \;\longrightarrow\;
  \mathcal{L}\left( L^p((0,T)\to E_1)\to L^p((0,T)\to E_0)\right)
  \quad\mbox{ and }
\\
  v &\mapsto \left[ t\mapsto f\left( z(t), v(z(t))\right) \right]
  \colon
  C( K\to \tilde{U} )
  \;\longrightarrow\; L^p((0,T)\to E_0) \,,
\intertext{the former one meaning}
  (v,u) &\mapsto
  \left[ t\mapsto A\left( z(t), v(z(t))\right)\, u(z(t))\right]
  \colon
\\
& C( K\to \tilde{U} )\times L^p((0,T)\to E_1)
  \;\longrightarrow\; L^p((0,T)\to E_0) \,,
\end{align*}
are locally Lipschitz continuous, the former one with values in
\begin{math}
  \mathcal{L}( L^p((0,T)\to E_1)\to\hfil\break
               L^p((0,T)\to E_0) )
\end{math}
and the latter one with values in $L^p((0,T)\to E_0)$.
We will take advantage of this local Lipschitz continuity
in our proof of Theorem~\ref{thm-Analytic_glob} below.
We remark that the operator norm in
$\mathcal{L}\left( L^p((0,T)\to E_1)\to L^p((0,T)\to E_0)\right)$
of the linear substitution operator
\begin{align*}
  u \mapsto \left[ t\mapsto A\left( z(t), v(z(t))\right) u(t)\right]
    \colon L^p((0,T)\to E_1)\to L^p((0,T)\to E_0) \,,
\\
  \mbox{ with $z\in C([0,T]\to K)$ and $v\in C(K\to \tilde{U})$
         being fixed, }
\end{align*}
is bounded above by the supremum norm
\begin{align*}
  \vertiii{ A\left( z(\cdot), v(z(\cdot))\right) }_{ L^{\infty}(0,T) }
& {}
  \eqdef
  \lVert \left[ t\mapsto A\left( z(t), v(z(t))\right) \right]
  \rVert_{ C\left( [0,T]\to \mathcal{L}(E_1\to E_0)\right) }
\\
& {}
  = \sup_{0\leq t\leq T}
    \| A\left( z(t), v(z(t))\right) \|_{ \mathcal{L}(E_1\to E_0) }
  \hspace*{10pt} (< \infty) \,.
\end{align*}
%

\begin{theorem}\label{thm-Analytic_glob}
Let\/ $1 < p < \infty$, $\vartheta_0\in (0, \pi / 2)$,
$0 < T_0\leq T < \infty$, and assume that\/
$g\in \hfil\break L^p((0,T)\to E_0)$
possesses a holomorphic extension
\begin{math}
  \tilde{g}\colon \mathcal{U}\to E_0
\end{math}
to an open set $\mathcal{U}\subset \CC$ containing
$\overline{\Delta}_{\vartheta_0}^{T_0,T}$, i.e.,
\begin{math}
  \overline{\Delta}_{\vartheta_0}^{T_0,T}
  \subset \mathcal{U}\subset \CC .
\end{math}
Assume that\/ $\tilde{U}$ is a nonempty open subset of\/
$E_{ 1 - \frac{1}{p} ,\, p }$ and\/
$\tilde{A}$ and\/ $\tilde{f}$ satisfy\/
{\rm Hypotheses} \eqref{hyp:A_anal} and\/ \eqref{hyp:f_anal},
respectively, and their respective restrictions
\begin{math}
  A = \tilde{A}\vert_{[0,T]\times U}
\end{math}
and\/
\begin{math}
  f = \tilde{f}\vert_{[0,T]\times U}
\end{math}
to
\begin{math}
  [0,T]\times U\subset \RR\times E_{ 1 - \frac{1}{p} ,\, p }
\end{math}
satisfy
{\rm Hypotheses} \eqref{hyp:A} and\/ \eqref{hyp:f}
with an open set\/
$U\subset \tilde{U}\subset E_{ 1 - \frac{1}{p} ,\, p }$.
Finally, assume that\/
$w\colon [0,T]\to U\subset E_{ 1 - \frac{1}{p} ,\, p }$
is a fixed global strict solution to problem \eqref{e:abstr_Cauchy}
(hence, satisfying \eqref{e_w:strict_sol})
with a prescribed initial value $w(0) = w_0\in U$
and such that\/
$w(t)\in U$ and $\tilde{A}(t,w(t))\in \mathrm{MR}_p(E)$
for all\/ $t\in [0,T]$.

Then there exist constants\/
$\vartheta'\in (0, \vartheta_0]$ and\/
$T'\in (0,T_0]$, small enough, and a holomorphic function
\begin{math}
  \tilde{w}\colon \Delta_{\vartheta'}^{T',T}
  \to E_{ 1 - \frac{1}{p} ,\, p }
\end{math}
with the following two properties:
\begin{itemize}
\item[{\rm (a)}]
$\tilde{w}(t)\in \tilde{U}$ for every\/
$t\in \Delta_{\vartheta'}^{T',T}$ and\/
$\tilde{w}$ verifies the abstract nonlinear evolutionary problem
\eqref{e:abstr_Cauchy} in the complex domain
$\Delta_{\vartheta'}^{T',T}$, i.e.,
\begin{equation}
\label{e:abstr_Cauchy_anal}
\left\{
\begin{alignedat}{2}
    \frac{\mathrm{d}u}{\mathrm{d}t}
  - \tilde{A}(t,u(t)) u(t) &= \tilde{f}(t,u(t)) + \tilde{g}(t)
  \quad\mbox{ for every }\, t\in \Delta_{\vartheta'}^{T',T} \,;
\\
  \lim_{ t\to 0 ,\ t\in \Delta_{\vartheta'}^{T',T} }\hspace{05pt}
  u(t) &= w_0\in E_{ 1 - \frac{1}{p} ,\, p } \,.
\end{alignedat}
\right.
\end{equation}
\item[{\rm (b)}]
$\tilde{w}(t) = w(t)$ holds for a.e.\ $t\in (0,T)$.
\end{itemize}

Such a holomorphic extension
\begin{math}
  \tilde{w}\colon \Delta_{\vartheta'}^{T',T}
  \to \tilde{U}\subset E_{ 1 - \frac{1}{p} ,\, p }
\end{math}
of
$w\colon (0,T)\to U\subset E_{ 1 - \frac{1}{p} ,\, p }$
from $(0,T)$ to $\Delta_{\vartheta'}^{T',T}$ is unique.
\end{theorem}
\par\vskip 10pt

Before proceeding to the proof of this theorem,
we clarify our notation with the open sets $U$ and $\tilde{U}$ in
$E_{ 1 - \frac{1}{p} ,\, p }$ as follows:

\begin{remark}\label{rem-Analytic_glob}\nopagebreak
\begingroup\rm
We need to take advantage of our
{\rm Hypotheses} \eqref{hyp:A} and \eqref{hyp:f}
(with an open set
 $U\subset E_{ 1 - \frac{1}{p} ,\, p }$)
and
{\rm Hypotheses} \eqref{hyp:A_anal} and \eqref{hyp:f_anal}
(with another open set
 $\tilde{U}\subset E_{ 1 - \frac{1}{p} ,\, p }$)
only for the values of $v = w(t)\in U\subset \tilde{U}$
($t\in \overline{\Delta}_{\vartheta'}^{T',T}$)
near the (compact) image
\begin{math}
  K =
  \{ w(t)\in E_{ 1 - \frac{1}{p} ,\, p } \colon t\in [0,T] \}
\end{math}
of the (continuous) curve
\begin{math}
  w\colon [0,T]\to E_{ 1 - \frac{1}{p} ,\, p } .
\end{math}
Indeed,
{\rm Hypotheses} \eqref{hyp:A_anal} and \eqref{hyp:f_anal}
imply that both holomorphic extensions
\begin{math}
  \tilde{A}\colon \mathcal{U}\times \tilde{U}
  \to \mathcal{L}(E_1\to E_0)
\end{math}
and
\begin{math}
  \tilde{f}\colon \mathcal{U}\times \tilde{U}\to E_0
\end{math}
of
$A\colon [0,T]\times U\to \mathcal{L}(E_1\to E_0)$
and
$f\colon [0,T]\times U\to E_0$, respectively,
are locally Lipschitz continuous.
Consequently, the Cartesian product
$[0,T]\times K$ being compact in the complex Banach space
$\CC\times E_{ 1 - \frac{1}{p} ,\, p }$,
we use a finite open subcover by open balls to find two bounded open sets
$\mathcal{U}\subset \CC$ and
$U = \tilde{U}\subset E_{ 1 - \frac{1}{p} ,\, p }$,
such that both mappings
$\tilde{A}$ and $\tilde{f}$ are Lipschitz continuous in
$\mathcal{U}\times \tilde{U}$.
We conclude that, in our proof of Theorem~\ref{thm-Analytic_glob}
below, we may assume that
$[0,T]\subset \mathcal{U}\subset \CC$ and
$K\subset U = \tilde{U}\subset E_{ 1 - \frac{1}{p} ,\, p }$
with both $\mathcal{U}$ and $U$ being open and bounded.
In particular, if the numbers
$\vartheta_0\in (0, \pi / 2)$ and $T_0\in (0,T]$
are taken sufficiently small, then we have also
$\overline{\Delta}_{\vartheta_0}^{T_0,T} \subset \mathcal{U}$
together with $\tilde{w}(t)\in U$ for all
$t\in \overline{\Delta}_{\vartheta'}^{T',T}$,
provided
$\vartheta'\in (0, \vartheta_0]$ and $T'\in (0,T_0]$
are small enough.
Consequently,
\begin{math}
  \overline{\Delta}_{\vartheta'}^{T',T} \subset
  \overline{\Delta}_{\vartheta_0}^{T_0,T} \subset \mathcal{U} .
\end{math}
%

%
%
In order to simplify our notation, we work {\em only\/} with
the holomorphic extensions
\begin{math}
  \tilde{g}\colon \mathcal{U}\to E_0 ,
\end{math}
\begin{math}
  \tilde{A}\colon \mathcal{U}\times U
  \to \mathcal{L}(E_1\to E_0) ,
\end{math}
and
\begin{math}
  \tilde{f}\colon \mathcal{U}\times U\to E_0
\end{math}
of the mappings $g$, $A$, and $f$, respectively.
Hence, we may remove the ``tilde'' from these symbols and write simply
$g = \tilde{g}$, $A = \tilde{A}$, and $f = \tilde{f}$.
We also may and will assume that both mappings
$A$ and $f$ are Lipschitz continuous {\it in all of\/}
$\mathcal{U}\times \tilde{U}$.
\hfill\Square
%
%
\endgroup
\end{remark}
\par\vskip 10pt

In our construction of the continuous extension
\begin{math}
  \tilde{w}\colon \overline{\Delta}_{\vartheta_0}^{T_0,T} \to \CC
\end{math}
of the strict solution $w\colon [0,T]\to \CC$,
holomorphic in $\Delta_{\vartheta_0}^{T_0,T}$,
we take advantage of a factorization approach for
the complex time variable $t = \varrho\mu$ where
$\varrho\in (0,\tau_0)$ and $\mu\in \CC$ with
$|\mu - 1| < \sin\vartheta$.
The numbers
$\tau_0\in (0,T)$ and $\vartheta\in (0,\vartheta_0)$
are suitable constants.
Fixing such a constant $\mu$, we obtain a {\em mild solution\/},
\begin{math}
  \omega\equiv \omega_{\mu}\colon
  [0,\tau_0]\to U\subset E_{ 1 - \frac{1}{p} ,\, p } ,
\end{math}
of the corresponding initial value problem with the real time variable
$t = \varrho\in [0,\tau_0]$.
Of course, this solution depends on the complex parameter $\mu$
from the open disc
\begin{equation*}
  D_r(1)\eqdef \{ \mu\in \CC\colon |\mu - 1| < r\}
\end{equation*}
centered at the point $1\in \CC$ with radius $r = \sin\vartheta$.
We will complete the proof by showing that the mild solution,
$\omega$, is holomorphic with respect to $\mu$.
This factorization approach has been used earlier in
{\sc D.\ Henry} \cite[Chapt.~3, {\S}3.4]{Henry} and
{\sc A.\ Lunardi} \cite[Chapt.~8, {\S}8.3.3]{Lunardi}.

\par\vskip 10pt
{\it Proof of\/} Theorem~\ref{thm-Analytic_glob}.
Given any two numbers
$\vartheta\in (0, \pi / 2)$ and $\tau_0\in (0,\infty)$,
we define a bounded open {\em sector\/}
in the complex plane $\CC$ by
\begin{equation}
\label{loc:tau_0}
  \mathfrak{A}_{\vartheta}^{(\tau_0)} \eqdef
  \{ t = \varrho\mu\in \CC\colon
     0 < \varrho < \tau_0 \,\mbox{ and }\, \mu\in\CC \,\mbox{ with }\,
     |\mu - 1| < \sin\vartheta \}
\end{equation}
with vertex at the origin $0\in \CC$ and angle $2\vartheta$.
Its closure in $\CC$ is denoted by
$\overline{\mathfrak{A}}_{\vartheta}^{(\tau_0)}$.
Recalling our definition of the triangle
$\Delta_{\vartheta}^{(T)}$ by {\rm eq.}~\eqref{loc:T^(T')},
and setting $r = \sin\vartheta$ (hence, $0<r<1$), we deduce that
\begin{align*}
& \Delta_{\vartheta'}^{(T_1)} \subset
  \mathfrak{A}_{\vartheta}^{(\tau_0)} \subset
  \Delta_{\vartheta}^{(T_2)}
    \quad\mbox{ holds whenever }\quad
\\
& 0 < T_1\leq \tau_0 \,,\quad 0 < \vartheta'< \arctan r \,,
  \quad\mbox{ and }\quad
  (1+r) \tau_0\leq T_2 < \infty \,.
\end{align*}
Following this factorization of the complex time
$t\in \CC$ in $t = \varrho\mu$ with $\varrho\in (0,\tau_0)$ and
\begin{math}
  \mu\in D_r(1) = \{ \mu\in \CC\colon |\mu - 1| < r\} ,
\end{math}
so that $0 < 1-r < \RE\mu < 1+r$ with $r = \sin\vartheta$ $(< 1)$,
we replace the complex time
$t\in \Delta_{\vartheta}^{T',T}$
in the initial value problem \eqref{e:abstr_Cauchy_anal}
by the product $\mu t\in \CC$ with $t\in (0,\tau_0)$ and $\mu\in D_r(1)$,
where we will choose both
$\vartheta\in (0, \pi / 2)$ and $\tau_0\in (0,\infty)$
sufficiently small, so that
\begin{math}
  \mathfrak{A}_{\vartheta}^{(\tau_0)} \subset
  \Delta_{\vartheta_0}^{T_0,T}
\end{math}
holds, i.e.,
$\mu t\in \Delta_{\vartheta_0}^{T_0,T}$ for every pair
$(t,\mu)\in (0,\tau_0)\times D_r(1)$.
Hence, we must have
\begin{equation*}
  0 < \vartheta\leq \vartheta_0 \,,\quad
  r\tau_0\leq T_0\cdot \tan\vartheta_0 \,,\quad\mbox{ and }\quad
  (1+r)\tau_0\leq T \,.
\end{equation*}

Given a fixed number $\mu\in D_r(1)\subset \CC$,
we look for an unknown continuous mapping
\begin{math}
  \omega\equiv \omega_{\mu}\colon
  [0,\tau_0]\to U\subset E_{ 1 - \frac{1}{p} ,\, p } ,
\end{math}
$\omega(t)\equiv \omega_{\mu}(t) = \tilde{w}(\mu t)$,
that according to {\rm eq.}~\eqref{e:abstr_Cauchy_anal}
must be a strict solution to the following evolutionary problem
(with the tilde ``$\,\overset{\sim}{\phantom{u}}\,$'',
 marking holomorphic extensions, having been removed),
\begin{equation}
\label{e:omega_Cauchy}
\left\{
\begin{alignedat}{2}
    \frac{\mathrm{d}\omega}{\mathrm{d}t}
  - \mu\, A(\mu t,\omega(t)) \omega(t)
& = \mu\left[ f(\mu t,\omega(t)) + g(\mu t)\right]
  \quad\mbox{ for every }\, t\in (0,\tau_0) \,;
\\
  \omega(0)
& = w_0\in E_{ 1 - \frac{1}{p} ,\, p } \,.
\end{alignedat}
\right.
\end{equation}
Of course, for $\mu = 1$ we will have
$\omega(t)\equiv \omega_1(t) = w(t)$ for a.e.\ $t\in (0,T)$,
by uniqueness.
We remark that, thanks to our hypothesis
$w(t)\in U$ and $A(t,w(t))\in \mathrm{MR}_p(E)$ for all $t\in [0,T]$,
we have also
$\mu A(t,w(t))\in \mathrm{MR}_p(E)$ for all $t\in [0,T]$ and all
$\mu\in \CC$ satisfying $|\mu - 1| < \sin\vartheta$
with $\vartheta\in (0, \pi / 2)$ small enough, say,
$0 < \vartheta < \pi / 6$ in which case $|\mu - 1| < 1/2$.
This claim follows easily from the perturbation lemma,
Lemma~\ref{lem-perturb_MR:A<B}, thanks to $\mu = 1 + \nu$ with
$\nu\in \CC$ satisfying $|\nu| < \sin\vartheta$.
Even Lemma~\ref{lem-perturb_MR} would do if
$\vartheta\in (0, \pi / 2)$ were chosen sufficiently small.
Clearly, \underline{it suffices} to prove that,
for each fixed $t\in (0,\tau_0)$, the function
\begin{math}
  \mu\mapsto \omega_{\mu}(t)\colon
  D_r(1)\to E_{ 1 - \frac{1}{p} ,\, p }
\end{math}
is holomorphic.
This approach to the analyticity in time of solutions to
semilinear parabolic problems can be found, e.g., in the monographs by
{\sc D.\ Henry} \cite[Chapt.~3]{Henry}, {\S}3.4,
Theorem 3.4.4 and Corollary 3.4.6 on pp.\ 63--66, and by
{\sc A.\ Lunardi} \cite[Chapt.~8]{Lunardi}, {\S}8.3.3, p.~308.

Choosing $\vartheta\in (0, \pi / 2)$ and $\tau_0\in (0,T]$
small enough, such that
$\overline{\mathfrak{A}}_{\vartheta}^{(\tau_0)} \subset \mathcal{U}$,
and recalling $r = \sin\vartheta\in (0,1)$, we abbreviate
\begin{equation}
\label{def:F_mu}
  F(t,v,\mu)\eqdef \mu\left[ f(\mu t,v) + g(\mu t)\right]
  \quad\mbox{ for all }\,
  (t,v,\mu)\in [0,\tau_0]\times U\times D_r(1) \,.
\end{equation}
By {\rm Hypothesis} \eqref{hyp:f_anal}, the mapping
\begin{math}
  (v,\mu)\mapsto F(t,v,\mu)\colon U\times D_r(1)\to E_0
\end{math}
is holomorphic for each $t\in [0,\tau_0]$,
with all partial derivatives of $F$ with respect to $v$ and $\mu$
being continous in $[0,\tau_0]\times U\times D_r(1)$.
According to 
{\sc H.\ Amann} \cite[Chapt.~III, {\S}4.10]{Amann-Birkh},
pp.\ 180--191, and
{\sc Ph.\ Cl\'ement} and {\sc S.\ Li} \cite[Sect.~2]{Clem-Li},
p.~18, given a fixed parameter value $\mu\in D_r(1)$,
every strict solution
\begin{math}
  \omega\equiv \omega_{\mu}\in Y_{ 1 - \frac{1}{p} }^p(0,\tau_0)
\end{math}
of the initial value problem \eqref{e:omega_Cauchy}
satisfies the following integral equation for the unknown function
\begin{math}
  \omega\equiv \omega_{\mu}\in Y_{ 1 - \frac{1}{p} }^p(0,\tau_0) ,
\end{math}
\begin{align}
\label{fix:omega_Cauchy}
  \omega(t) = \mathcal{F}(t,\omega,\mu)
  \quad\mbox{ for every }\, t\in [0,\tau_0] \,,
  \quad\mbox{ with the right\--hand side equal to }
\\
\label{def:omega_Cauchy}
\begin{alignedat}{2}
    \mathcal{F}(t,v,\mu)
& \eqdef
    \ee^{\mu t A(0,w_0)}\, w_0
  + \mu\int_0^t \ee^{\mu (t-s) A(0,w_0)}\,
    \left[ A(\mu s, v(s)) - A(0,w_0)\right] v(s) \,\mathrm{d}s
\\
& {}
  + \int_0^t \ee^{\mu (t-s) A(0,w_0)}\,
    F(s,v(s),\mu)  \,\mathrm{d}s
  \quad\mbox{ for every }\, t\in [0,\tau_0]
\end{alignedat}
\end{align}
and for all $v\in Y_{ 1 - \frac{1}{p} }^p(0,\tau_0)$ satisfying
$v(t)\in U$ for every $t\in [0,\tau_0]$.
In contrast to defining a contraction mapping using
the (unique) strict solution to prove local existence in
Theorem~\ref{thm-Clem-Li}, in the case of problem~\eqref{e:omega_Cauchy}
we prefer to use the (unique) mild solution defined by
an integral representation
(variation\--of\--constants formula); cf.\
eq.~\eqref{def:omega_Cauchy} above.
The equivalence \underline{between} strict and mild solutions is treated in
{\sc J.~M.\ Ball} \cite{J.M.Ball},
{\sc D.\ Henry} \cite[Chapt.~3]{Henry}, and
{\sc A.\ Pazy} \cite[Theorem on p.~259]{Pazy}.

Clearly, \eqref{fix:omega_Cauchy}
is a fixed point equation for the unknown function
$\omega\in Y_{ 1 - \frac{1}{p} }^p(0,\tau_0)$.
Here, one can choose $\tau_0\in (0,T_1]$, where
$T_1\in (0,T]$ and $\vartheta\in (0, \pi / 2)$
are sufficiently small, such that
\begin{math}
  \overline{\mathfrak{A}}_{\vartheta}^{(T_1)} \subset \mathcal{U}
\end{math}
and the mapping
\begin{equation*}
  \Phi\equiv \Phi_{\mu}\colon v \mapsto
    \left[ t\mapsto \mathcal{F}(t,v,\mu) \right]
\end{equation*}
is a contraction in a closed ball
\begin{equation*}
  \Sigma_{\rho_1,T_1}^{(w_0)} =
  \left\{ v\in Y^{T_1}\colon v(0) = w_0 \quad\mbox{ and }\quad
          \left\| v - w\vert_{[0,T_1]}
          \right\|_{ Y^{T_1} }\leq \rho_1
  \right\}
\end{equation*}
in the Banach space
\begin{equation*}
  Y^{T_1} = Y_{ 1 - \frac{1}{p} }^p(0,T_1)
  = L^p((0,T_1)\to E_1)\cap W^{1,p}((0,T_1)\to E_0)
\end{equation*}
of radius $\rho_1\in (0,\infty)$ centered at the point $w\in Y^{T_1}$.
As usual, the function
$w\vert_{[0,T_1]} \in Y^{T_1}$
denotes the restriction to $[0,T_1]$ of the strict solution
$w\in Y^{T} = Y_{ 1 - \frac{1}{p} }^p(0,T)$
from the hypotheses of our theorem.
The proof of this contraction property follows
the same ideas and steps as the proof of
Theorem~\ref{thm-Clem-Li} taken from
{\sc Ph.\ Cl\'ement} and {\sc S.\ Li}
\cite[Theorem 2.1, p.~19]{Clem-Li}.
The reader is referred to
{\sc J.\ Pr\"uss} \cite[pp.\ 9--13]{Pruess}, Chapt.~3,
for further details.
Notice that the numbers $\rho_1\in (0,\infty)$, $T_1\in (0,T]$, and
$\vartheta\in (0, \pi / 2)$,
if chosen small enough, such that the contraction holds
with the Lipschitz constant $\frac{1}{2}$,
are independent from the particular choice of the parameter
$\mu\in D_r(1)$ since
$\overline{\mathfrak{A}}_{\vartheta}^{(T_1)}$
is a compact subset of $\mathcal{U}$; cf.\
our remarks before this proof
(Remark~\ref{rem-Analytic_glob})
that remain valid also for the compact set
\begin{math}
  \overline{\mathfrak{A}}_{\vartheta}^{(T_1)} \times K
\end{math}
in the complex Banach space
$\CC\times E_{ 1 - \frac{1}{p} ,\, p }$.
Of course,
$r = \sin\vartheta\in (0,1)$ is sufficiently small, and both
$\rho_1\in (0,\infty)$ and $T_1\in (0,T]$ must be also so small that
$v(t)\in U$ holds for all $t\in [0,T_1]$, whenever
$v\in \Sigma_{\rho_1,T_1}^{(w_0)}$.
Finally, the constants $\rho_1$, $T_1$, and $\vartheta$
can be chosen independent from $w_0\in K$,
so that one may use them in any time interval
$[t_0 ,\, t_0 + T_1]\subset [0,T]$
of sufficiently short length $T_1 > 0$;
the initial condition $w(0) = w_0\in K$ at $t=0$ is replaced by
$w(t_0)\in K$ at arbitrary time $t_0\in [0, T-T_1]$.

Next, we analyze the holomorphy properties of the fix point mapping
\begin{equation*}
  \mathcal{F}\colon [0,T_1]\times \Sigma_{\rho_1,T_1}^{(w_0)}
             \times D_r(1)\to \Sigma_{\rho_1,T_1}^{(w_0)}
\end{equation*}
defined in {\rm eq.}~\eqref{def:omega_Cauchy}
where we may take $\tau_0 = T_1$; more precisely, those of the mapping
\begin{equation*}
  (v,\mu)\mapsto \mathcal{F}(t,v,\mu)\colon
         \Sigma_{\rho_1,T_1}^{(w_0)} \times D_r(1)\to
         \Sigma_{\rho_1,T_1}^{(w_0)} \,,
\end{equation*}
for each fixed $t\in [0,T_1]$.
To begin with, for $0\leq s < t\leq T_1$, $\mu\in D_r(1)$, and
$v\in \Sigma_{\rho_1,T_1}^{(w_0)}$, we rewrite
\begin{equation*}
    A(\mu s, v(s)) - A(0,w_0)
  = \left\{ I - \left[ \lambda I - A(\mu s, v(s)) \right]
                \left[ \lambda I - A(0,w_0) \right]^{-1}
    \right\}    \left[ \lambda I - A(0,w_0) \right] \,,
\end{equation*}
where $\lambda\in (0,\infty)$ is large enough
in order to guarantee that the (bounded) linear operator\hfil\break
$\lambda I - A(0,w_0)\colon E_1\to E_0$
has a bounded inverse
$\left[ \lambda I - A(0,w_0) \right]^{-1} \colon E_0\to E_1$,
and observe that the function (integrand)
\begin{equation*}
  \mu\mapsto \ee^{\mu (t-s) A(0,w_0)}\,
             \left[ A(\mu s, v(s)) - A(0,w_0)\right] v(s)
  \colon D_r(1)\to E_0
\end{equation*}
is holomorphic and so is the integral
\begin{equation*}
  \mu\mapsto \int_0^t \ee^{\mu (t-s) A(0,w_0)}\,
    \left[ A(\mu s, v(s)) - A(0,w_0)\right] v(s) \,\mathrm{d}s
  \colon D_r(1)\to E_0 \,.
\end{equation*}
We have used here also the fact that the operator\--valued function
\begin{equation*}
  \mu\mapsto \ee^{\mu (t-s) A(0,w_0)}
  \colon D_r(1)\to \mathcal{L}(E_0\to E_0)
\end{equation*}
is holomorphic for any fixed numbers $s,t\in \mathbb{R}$ satisfying
$0\leq s < t\leq T_1$.
Similarly, the function
\begin{equation*}
  \mu\mapsto \ee^{\mu (t-s) A(0,w_0)}\, F(s,v(s),\mu)
  \colon D_r(1)\to E_0
\end{equation*}
being holomorphic, so is the integral
\begin{equation*}
  \mu\mapsto \int_0^t \ee^{\mu (t-s) A(0,w_0)}\,
    F(s,v(s),\mu) \,\mathrm{d}s
  \colon D_r(1)\to E_0 \,.
\end{equation*}
We conclude that also the sum
\begin{equation*}
  \mu\mapsto \mathcal{F}(t,v,\mu) \colon D_r(1)\to E_0
\end{equation*}
defined by eq.~\eqref{def:omega_Cauchy} with $\tau_0 = T_1$
is holomorphic for every $t\in [0,T_1]$.

Finally, from the fixed point equation \eqref{fix:omega_Cauchy}
we deduce that also the function
\begin{math}
  \omega\equiv \omega_{\mu}\colon
  [0,\tau_0]\to U\subset E_{ 1 - \frac{1}{p} ,\, p } ,
\end{math}
which is continuous thanks to
\begin{math}
  \omega\in \Sigma_{\rho_1,T_1}^{(w_0)} \subset Y^{T_1}
            = Y_{ 1 - \frac{1}{p} }^p(0,T_1) ,
\end{math}
is holomorphic in the variable $\mu\in D_r(1)$.
Although this holomorphy claim follows directly from
a well\--known result in {\sc K.\ Deimling}
\cite[Theorem 15.3, Chapt.~4, {\S}15, p.~151]{Deimling},
cf.\ also {\sc S.~G.\ Krantz} and {\sc H.~R.\ Parks}
\cite[Theorem 6.1.2, {\S}6.1, p.~118]{Krantz-implicit},
we sketch a constructive proof below for the sake of completeness.

Indeed, any standard proof of the Banach fixed point theorem for
the ({\em contractive\/}) self\--mapping
\begin{equation*}
  \Phi\equiv \Phi_{\mu}\colon v \mapsto
    \left[ t\mapsto \mathcal{F}(t,v,\mu) \right]
  \colon \Sigma_{\rho_1,T_1}^{(w_0)} \to \Sigma_{\rho_1,T_1}^{(w_0)}
\end{equation*}
shows that, given an arbitrary ``initial'' function
$\varphi_0\in \Sigma_{\rho_1,T_1}^{(w_0)}$, the iterates
\begin{equation*}
  \varphi_n = \Phi (\varphi_{n-1}) = \Phi^2 (\varphi_{n-2})
            = \,\dots\, = \Phi^{n-1} (\varphi_1)
            = \Phi^n (\varphi_0) \,;
    \quad\mbox{ for }\, n=1,2,3,\dots \,,
\end{equation*}
form a Cauchy sequence in
$\Sigma_{\rho_1,T_1}^{(w_0)}$
which converges to the unique fixed point
$\omega\equiv \omega_{\mu}$ of $\Phi$, namely,
$\varphi_n\to \omega$ in
$\Sigma_{\rho_1,T_1}^{(w_0)} \subset Y^{T_1}$ as $n\to \infty$.
The convergence is uniform for $\mu\in D_r(1)$.
Recalling the continuous imbedding \eqref{e:Y->C},
we have also
$\varphi_n(t)\to \omega(t)$ in $E_{ 1 - \frac{1}{p} ,\, p }$
as $n\to \infty$, uniformly for $t\in [0,T_1]$ and $\mu\in D_r(1)$.
Choosing
$\varphi_0 = w\vert_{[0,T_1]}$,
a function of time $t\in [0,T_1]$ which does not depend on the parameter
$\mu\in D_r(1)$, we observe that each iterate
\begin{equation*}
  \varphi_n(t) = \mathcal{F}(t, \varphi_{n-1}, \mu) \,;\quad
  n=1,2,3,\dots \,,\quad t\in [0,T_1] \,,
\end{equation*}
is a holomorphic function in the variable (parameter) $\mu\in D_r(1)$.
Applying Osgood's theorem and
the Cauchy integral formula for discs to each iterate $\varphi_n(t)$
(see e.g.\
 {\sc S.~G.\ Krantz} \cite{Krantz}, Theorem 1.2.2 (p.~24), or
 {\sc F.\ John} \cite{John}, Chapt.~3, Sect.\ 3(c), eq.\ (3.22c), p.~71),
we conclude that also the limit function
$\omega\equiv \omega_{\mu}$
is holomorphic in the variable $\mu\in D_r(1)$ and satisfies
$\Phi(\omega) = \omega$.

We have thus verified that the strict solution
$w\colon [0,T]\to U\subset E_{ 1 - \frac{1}{p} ,\, p }$
of problem \eqref{e:abstr_Cauchy}
possesses a holomorphic extension to the bounded open sector
$\mathfrak{A}_{\vartheta}^{(T_1)}$.
In fact, we have proved that this claim is valid
in any time shift of this sector by a number $t_0\in [0, T-T_1]$,
that is, in any sector
\begin{equation*}
  t_0 +
  \mathfrak{A}_{\vartheta}^{(T_1)} \eqdef
  \{ t = t_0 + \varrho\mu\in \CC\colon
     0 < \varrho < T_1 \,\mbox{ and }\, |\mu - 1| < \sin\vartheta \}
\end{equation*}
with vertex at the point $t_0\in \CC$ and angle $2\vartheta$.
We apply the last result with $t_0$ ranging from $0$ to $T-T_1$
over the interval $[0, T-T_1]$ to conclude that the function
$w\colon [0,T]\to U\subset E_{ 1 - \frac{1}{p} ,\, p }$
possesses a holomorphic extension to the bounded open set
\begin{equation*}
  \bigcup_{ t_0\in [0, T-T_1] }
  \left( t_0 + \mathfrak{A}_{\vartheta}^{(T_1)} \right)
  \subset \CC
\end{equation*}
which contains the open complex domain
$\Delta_{\vartheta'}^{T',T}$ defined in \eqref{glob:T^(r)},
whenever $T'= T_1$ and
$0 < \vartheta'\leq \arctan (\sin\vartheta)$, owing to
\begin{math}
  \Delta_{\vartheta'}^{(T_1)} \subset
  \mathfrak{A}_{\vartheta}^{(T_1)} .
\end{math}

Hence, we have proved that there are constants
$\vartheta'\in (0, \vartheta_0]$ and
$T'\in (0,T_0]$, small enough, and a holomorphic function
\begin{math}
  \tilde{w}\colon \Delta_{\vartheta'}^{T',T}
  \to E_{ 1 - \frac{1}{p} ,\, p }
\end{math}
with the desired properties {\rm (a)} and {\rm (b)}
in the conclusion of our theorem.
Since
$(0,T)\subset \Delta_{\vartheta'}^{T',T}$,
such a holomorphic function $\tilde{w}$ must be unique.
The proof is finished.
\qed
\par\vskip 10pt

\section{Analyticity in space for the Cauchy problem in
         $\RR^N\times (0,T)$}
\label{s:Analyt_space}

In the previous two sections,
Sections \ref{s:Cauchy_abstr} and~\ref{s:Analyt_abstr},
we have treated the initial value problem \eqref{e:abstr_Cauchy}
for a strict solution $u\colon (0,T)\to E_0$
with the initial condition $u_0$ in the real interpolation space
\begin{equation*}
  E_{ 1 - \frac{1}{p} ,\, p } \equiv
  (E_0, E_1)_{ 1 - \frac{1}{p} ,\, p }
  \qquad\mbox{ between $E_1$ and $E_0$, }\;
  E_1\hookrightarrow E_{ 1 - \frac{1}{p} ,\, p } \hookrightarrow E_0 \,.
\end{equation*}
By Theorem~\ref{thm-Clem-Li_glob},
such a strict solution belongs to the abstract Sobolev space
\begin{math}
  Y^T = Y_{ 1 - \frac{1}{p} }^p(0,T)
\end{math}
introduced in \eqref{e:Besov_imbed}.
Hence, we have
$u(t)\in E_{ 1 - \frac{1}{p} ,\, p }$ for every $t\in [0,T]$.

In this section we replace
the triplet of abstract (complex) Banach spaces
\begin{math}
  E_1\hookrightarrow E_{ 1 - \frac{1}{p} ,\, p }\hookrightarrow E_0
\end{math}
by the following {\em complex\/} Sobolev, Besov, and Lebesgue spaces,
respectively,
\begin{align*}
& W^{2m,p}(\RR^N)\hookrightarrow B^{s;p,p}(\RR^N)\
                 \hookrightarrow L^p(\RR^N) \,,\quad
  s = 2m\left( 1 - \genfrac{}{}{}1{1}{p}\right)\in (0,2m) \,,
  \quad\mbox{ where }
\\
&   B^{s;p,p}(\RR^N)\eqdef
    \left( L^p(\RR^N) ,\, W^{2m,p}(\RR^N) \right)_{s/(2m),p}
  = \left( L^p(\RR^N) ,\, W^{2m,p}(\RR^N)
    \right)_{ 1 - (1/p) ,\, p }
\end{align*}
is the Besov space obtained by real interpolation
(see, e.g., \cite{AdamsFournier, Lunardi, Triebel}).
We recall that
$2 + \genfrac{}{}{}1{N}{m} < p$ $< \infty$
which guarantees $(s-m)p > N$ and, thus, the Sobolev\-(-Besov) imbeddings
\begin{math}
  B^{s-m;p,p}(\RR^N) \hookrightarrow
  C^0(\RR^N)\cap L^{\infty}(\RR^N)
\end{math}
and
\begin{math}
  B^{s;p,p}(\RR^N) \hookrightarrow
  C^m(\RR^N)\cap W^{m,\infty}(\RR^N)
\end{math}
are continuous.

Throughout this section we restrict ourselves to the easiest case
of {\color{red}\em\bfseries analytic initial conditions\/}
that {\color{red}\em we are able to treat\/} in our present work:

\begin{hypo}\nopagebreak
\begingroup\rm
\begin{enumerate}
\setcounter{enumi}{0}
\renewcommand{\labelenumi}{\bf ($\mathbf{A_0}$)}
\item
\makeatletter
\def\@currentlabel{$\mathbf{A_0}$}\label{hyp:analyt}
\makeatother
We assume that the initial data
$\mathbf{u}_0\colon \RR^N\to \CC^M$,
$\mathbf{u}_0 = (u_{0,1}, u_{0,2}, \dots, u_{0,M})$,
can be extended to a holomorphic function
\begin{math}
    \tilde{\mathbf{u}}_0
  = ( \tilde{u}_{0,1}, \tilde{u}_{0,2}, \dots, \tilde{u}_{0,M} )
  \colon \mathfrak{X}^{(r)}\to \CC^M
\end{math}
in a complex strip
$\mathfrak{X}^{(r)}\subset \CC^N$ defined in
{\rm eq.}~\eqref{def:X^(r)},
for some $r\in (0,\infty)$, such that every component
$\tilde{u}_{0,j}\colon \mathfrak{X}^{(r)}\to \CC$; $j=1,2,\dots,M$,
has the following properties:
\begin{itemize}
\item[{\rm ($A_0^{\prime}$)}]
$\;$
the function
$x\mapsto \tilde{u}_{0,j}(x + \ii y)\colon \RR^N\to \CC$
is in the (complex) Besov space $B^{s;p,p}(\RR^N)$,
\item[{\rm ($A_0^{\prime\prime}$)}]
$\;$
the Besov norm
\begin{math}
  \| \tilde{u}_{0,j}(\,\cdot\, + \ii y) \|_{ B^{s;p,p}(\RR^N) }
\end{math}
is uniformly bounded for all $y\in Q^{(r)}$, and
\item[{\rm ($A_0^{\prime\prime\prime}$)}]
$\;$
\begin{math}
  y\mapsto \tilde{u}_{0,j}(\,\cdot\, + \ii y)\colon
  Q^{(r)}\to B^{s;p,p}(\RR^N)
\end{math}
is continuously (partially) differentiable with respect to the parameter
\begin{math}
  y = (y_1,\dots,y_N)\in
  Q^{(r)} = \{ y\in \RR^N\colon |y|_{\infty} < r\} \,;
\end{math}
$j=1,2,\dots,M$.
\end{itemize}

Equivalently, the function
\begin{math}
  x\mapsto \tilde{\mathbf{u}}_0(x + \ii y)\colon \RR^N\to \CC^M
\end{math}
belongs to the Cartesian product
$\mathbf{B}^{s;p,p}(\RR^N) = [ B^{s;p,p}(\RR^N) ]^M$,
its norm
\begin{math}
  \| \tilde{\mathbf{u}}_0(\,\cdot\, + \ii y) \|_{ B^{s;p,p}(\RR^N) }
\end{math}
satisfies (cf.\ eq.~\eqref{sup_y:u_0})
\begin{equation*}
  \mathfrak{N}^{(r)} (\tilde{\mathbf{u}}_0)
  = \sup_{ y\in Q^{(r)} }
    \| \tilde{\mathbf{u}}_0 (\,\cdot\, + \ii y)
    \|_{ B^{s;p,p}(\RR^N) } < \infty \,,
\end{equation*}
and it is continuously differentiable with respect to the parameter
$y\in Q^{(r)}$.
\hfill\Square
\end{enumerate}
\endgroup
\end{hypo}
\par\vskip 10pt

The ``shift'' isometry
\begin{math}
  \| \tilde{\mathbf{u}}_0(\,\cdot\, + x_0 + \ii y_0)
  \|_{ B^{s;p,p}(\RR^N) } =
  \| \tilde{\mathbf{u}}_0(\,\cdot\, + \ii y_0)
  \|_{ B^{s;p,p}(\RR^N) }
\end{math}
is obvious for all pairs
$(x_0,y_0)\in \RR^N\times Q^{(r)}$, i.e., for all complex numbers
$z_0 = x_0 + \ii y_0\in \mathfrak{X}^{(r)}$.

The restriction in {\rm Hypothesis} \eqref{hyp:analyt}
is motivated by the following approximation property
of the Sobolev and Besov spaces, see e.g.\
{\sc H.\ Triebel} \cite[Chapt.~2]{Triebel}:

\begin{remark}\label{rem-hyp:analyt}\nopagebreak
\begingroup\rm
The Fr\'echet space $\mathcal{S}(\RR^N)$ of all complex\--valued,
rapidly decreasing infinitely differentiable functions
$\varphi\colon \RR^N\to \CC$
being dense in all of the spaces
$L^p(\RR^N)$, $W^{2m,p}(\RR^N)$, $B^{s;p,p}(\RR^N)$, and $L^2(\RR^N)$,
by \cite[Chapt.~2]{Triebel}, {\S}2.3, Theorem 2.3.2 on p.~172,
we take
$\mathbf{u}_0\colon \RR^N\to \CC^M$
so smooth and rapidly decreasing near infinity that
its holomorphic extension
$\tilde{\mathbf{u}}_0\colon \mathfrak{X}^{(r)}\to \CC^M$
satisfies even the following stronger regularity condition:
The family of functions
\begin{math}
  x\mapsto \tilde{\mathbf{u}}_0(x + \ii y)\colon \RR^N\to \CC^M ,
\end{math}
parametrized by $y\in Q^{(r)}$, belongs to a bounded subset of
\begin{equation*}
  \mathbf{L}^2(\RR^N)\cap \mathbf{W}^{2m,p}(\RR^N)
  \quad\mbox{ for some $p\in \RR$, }\;
  2 + \genfrac{}{}{}1{N}{m} < p < \infty \,.
\end{equation*}
For instance, all complex linear combinations of {\it Hermite functions\/}
form a dense vector subspace $\mathcal{V}$
of the Fr\'echet space $\mathcal{S}(\RR^N)$,
by {\sc M.\ Reed} and {B.\ Simon} \cite[Chapt.~V, {\S}3]{RSimon-I},
Theorem V.13 on p.~143.
{\it Hermite functions\/} are entire complex functions
$h\colon \CC^N\to \CC$ of the form
\begin{equation*}
  h(z) = P(z_1,z_2,\dots,z_N)\cdot \exp
  \left( {}- \frac{1}{2}\sum_{i=1}^N z_i^2\right)
  \quad\mbox{ for }\; z = (z_i)_{i=1}^N = x + \ii y\in \CC^N \,,
\end{equation*}
where $P(z)$ is a complex polynomial in $N$ complex variables
$z_i\in \CC$; $i=1,2,\dots,N$, see \cite[p.~142]{RSimon-I}.
One may take functions from $\mathcal{V}$ as components
$\tilde{u}_{0,j}$ of $\tilde{\mathbf{u}}_0$; $j=1,2,\dots,M$.
Indeed, notice that
\begin{align*}
    \left|
    \exp\left( {}- \frac{1}{2}\sum_{i=1}^N z_i^2 \right)
    \right|
  = \exp\left( {}- \frac{1}{2}\sum_{i=1}^N x_i^2
                 + \frac{1}{2}\sum_{i=1}^N y_i^2 \right)
  \leq \exp\left( \frac{1}{2}\, N r^2 \right)
    \,\cdot\, \exp\left( {}- \frac{1}{2}\, |x|_2^2 \right)
\\
  \quad\mbox{ holds for all }\;
    z = x + \ii y\in \mathfrak{X}^{(r)}\subset \CC^N \,,
  \quad\mbox{ where }\quad
\\
{ \textstyle
  |x|_2 = \left( \sum_{i=1}^N |x_i|^2 \right)^{1/2}
}
    \;\mbox{ and }\;
  |y|_{\infty} = \max_{i=1,2,\dots,N} |y_i| < r \,.
\end{align*}

It is well\--known that all three vector spaces
$\mathcal{V}\subset \mathcal{S}(\RR^N)\subset L^2(\RR^N)$
are invariant under the (unitary) Fourier transformation
$\mathcal{F}\colon L^2(\RR^N)\to L^2(\RR^N)$.
(We always consider the {\em unitary\/} Fourier transformation
 $\mathcal{F}$ as described in
 {\sc E.~M.\ Stein} and {\sc G.\ Weiss} \cite[Chapt.~I]{SteinWeiss}.)
Consequently, if the Fourier transform
$\mathcal{F} u_{0,j}\colon \RR^N\to \CC$ of each component of
$\mathbf{u}_0\colon \RR^N\to \CC^M$
decays at least exponentially fast at infinity, then
the holomorphic extension of the function
$u_{0,j}\colon \RR^N\to \CC$ to a complex strip
$\mathfrak{X}^{(r)}\subset \CC^N$, for some $r\in (0,\infty)$,
is easily obtained in the form of
the {\em inverse Fourier\--Laplace transform\/}
\begin{math}
  \mathcal{F}^{-1} (\mathcal{F} u_{0,j})
  \colon \mathfrak{X}^{(r)}\to \CC
\end{math}
of $\mathcal{F} u_{0,j}$, by
the classical Paley\--Wiener\--Schwartz theory, see e.g.\
{\sc L.\ H\"ormander} \cite[Theorem 7.4.2, p.~192]{Hoerm-PDE-1} or
{\sc E.~M.\ Stein} and {\sc G.\ Weiss} \cite[Chapt.~III]{SteinWeiss},
{\S}2, pp.\ 91--101, and {\S}6.12, pp.\ 127--128.
An interested reader is referred also to
{\sc P.\ Tak\'a\v{c}} \cite[Chapt.~5]{Takac}
for a brief review of the (inverse) Fourier\--Laplace transform
that applies to our current setting.
%
\hfill\Square
\endgroup
\end{remark}
\par\vskip 10pt

With regard to later applications
(cf.\ Proposition~\ref{prop-smooth_ext} and
      Theorem~\ref{thm-smooth_ext}),
in our {\rm Hypothesis} \eqref{hyp:analyt} above
we have not specified the number $r\in (0,\infty)$
corresponding to the half\--width of the complex strip
$\mathfrak{X}^{(r)} = \RR^N + \ii Q^{(r)}$,
a tube in $\CC^N$ with the base $Q^{(r)} = (-r,r)^N$.
{\rm Hypotheses}
\eqref{hy:UnifEll}, \eqref{hy:Ellipt}, and\/ \eqref{hy:f_Hardy^2}
in Section~\ref{s:Main}
show that only the case $0 < r\leq r_0$ is useful.
We will comment on the choice of $r\in (0,r_0]$
in Remark~\ref{rem-smooth_ext}
right after Theorem~\ref{thm-smooth_ext} below.
Concerning this question of choosing (finding) a suitable half\--width
$r\in (0,r_0]$, we begin with the following observation.

\begin{remark}\label{rem-half-width^r}\nopagebreak
\begingroup\rm
The {\it Hermite functions\/}
$h\colon \CC^N\to \CC$ desribed in {\rm Remark~\ref{rem-hyp:analyt}}
above are {\it not\/} the only way for approximating the initial values
$\mathbf{u}(\,\cdot\,,t_0) = \mathbf{u}_0\in \mathbf{B}^{s;p,p}(\RR^N)$
at time $t_0 = 0$ in the Besov space within the Besov norm
$\|\cdot\|_{ B^{s;p,p}(\RR^N) }$.
In our approximation procedure we need to guarantee
the following ``uniformity'' of the half\--width of the complex strip
$\mathfrak{X}^{(r)} = \RR^N + \ii Q^{(r)}$, i.e.,
the same half\--width $r\in (0,r_0]$ for each approximating function
\begin{math}
  \tilde{\mathbf{u}}_{0,n} \colon \mathfrak{X}^{(r)}\to \CC^M \,;
\end{math}
$n=1,2,3,\dots$.
In precise analytic terms, this means that, for any given radius
$R_1\in (0,\infty)$ of the ball $B_{R_1}(0)$ in
$\mathbf{B}^{s;p,p}(\RR^N)$,
there exists a number $r_1\in (0,r_0]$
small enough, such that, whenever $r\in (0,r_1]$,
the approximating sequence of functions
$\{ \tilde{\mathbf{u}}_{0,n} \}_{n=1}^{\infty}$
has the following properties
(cf.\ {\rm Hypotheses} \eqref{hyp:analyt}):
\begin{itemize}
\item[{\rm ($A_{0,n}^{\prime}$)}]
$\;$
each function 
$x\mapsto \tilde{\mathbf{u}}_{0,n}(x + \ii y)\colon \RR^N\to \CC^M$
is in the Besov space $\mathbf{B}^{s;p,p}(\RR^N)$,
for every $y\in Q^{(r)}$,
\item[{\rm ($A_{0,n}^{\prime\prime}$)}]
$\;$
\begingroup\large\bf
the ``proximity to $\mathbf{u}_0$'' estimate
\endgroup
\begin{math}
  \| \tilde{\mathbf{u}}_{0,n}(\,\cdot\, + \ii y) - \mathbf{u}_0
  \|_{ B^{s;p,p}(\RR^N) } < R_1
\end{math}
holds for all $y\in Q^{(r)}$ and $n=1,2,3,\dots$,
\item[{\rm ($A_{0,n}^{\prime\prime\prime}$)}]
$\;$
\begin{math}
  \tilde{\mathbf{u}}_{0,n} \colon \mathfrak{X}^{(r)}\to \CC^M
\end{math}
is holomorphic for every $n=1,2,3,\dots$, and finaly
\item[{\rm ($A_{0,n}^{\mathrm{iv}}$)}]
$\;$
the restrictions
\begin{math}
  \mathbf{u}_{0,n} = \tilde{\mathbf{u}}_{0,n} \vert_{ \mathbb{R} }
\end{math}
of $\tilde{\mathbf{u}}_{0,n}$ to the real line $\RR$ satisfy
\begin{math}
  \| \mathbf{u}_{0,n} - \mathbf{u}_0 \|_{ B^{s;p,p}(\RR^N) }
  \,\to\, 0
\end{math}
as $n\to \infty$.
\hfill\Square
\end{itemize}
%
%
\endgroup
\end{remark}
\par\vskip 10pt

We keep the natural notation
$\mathbf{L}^2(\RR^N) = [ L^2(\RR^N) ]^M$ etc.\
introduced for spaces of vector\--valued functions in
the Introduction (Section~\ref{s:Intro}).
We recall the continuous Sobolev\-(-Besov) imbeddings
\begin{equation*}
  \mathcal{S}(\RR^N) \hookrightarrow
  L^2(\RR^N)\cap W^{2m,p}(\RR^N)
  \hookrightarrow B^{s;p,p}(\RR^N) \hookrightarrow
  C^m(\RR^N)\cap W^{m,\infty}(\RR^N)
\end{equation*}
and
\begin{equation*}
  W^{2m,p}(\RR^N) \hookrightarrow
  B^{s;p,p}(\RR^N) \hookrightarrow
  L^p(\RR^N)\cap L^{\infty}(\RR^N) \,,\quad
  2 + \genfrac{}{}{}1{N}{m} < p < \infty \,.
\end{equation*}
We remark that
\begin{math}
  W^{2m,p}(\RR^N) \not\subset L^2(\RR^N)
\end{math}
if $2 < p < \infty$.
From now on we identify
$\mathbf{u}_0\equiv \tilde{\mathbf{u}}_0$
and drop the tilde ``$\,\overset{\sim}{\phantom{u}}\,$''
in the (unique) holomorphic extension.

By our {\rm Hypotheses}
\eqref{hy:UnifEll}, \eqref{hy:Ellipt}, and\/ \eqref{hy:f_Hardy^2}
(cf.\ Theorem~\ref{thm-Main}),
let us set $r = r_0\in (0,\infty)$ above.
In the Cauchy problem \eqref{e:Cauchy}
we may replace the real space variable $x\in \RR^N$ by its complex shift
$z = x + x_0 + \ii y_0$ by a fixed complex vector
$z_0 = x_0 + \ii y_0\in \mathfrak{X}^{(r)}\subset \CC^N$
with any $x_0\in \RR^N$ and any $y_0\in Q^{(r)}$.
In the sequel we consider
$z_0\in \mathfrak{X}^{(r)}$ to be a parameter and
$x\in \RR^N$ an independent variable
in the Cauchy problem \eqref{e:Cauchy} spatially ``shifted'' by $z_0$,
\begin{equation}
\label{e:Cauchy+z_0}
\left\{
\begin{alignedat}{2}
    \frac{\partial\mathbf{u}}{\partial t}
  + \mathbf{P}
    \Bigl( x + z_0, t, \frac{1}{\ii} \frac{\partial}{\partial x} \Bigr)
    \mathbf{u}
& = \mathbf{f}
    \left( x + z_0, t;
    \left( \frac{ \partial^{|\beta|}\mathbf{u} }{\partial x^{\beta}}
    \right)_{|\beta|\leq m}
    \right)
\\
&   \quad\mbox{ for } (x,t)\in \mathbb{R}^N\times (0,T) \,;
\\
  \mathbf{u}(x,0)&= \mathbf{u}_0(x + z_0)
    \quad\mbox{ for } x\in \mathbb{R}^N \,.
\end{alignedat}
\right.
\end{equation}
By our hypothesis on the initial data
$\mathbf{u}_0\colon \RR^N\to \CC^M$ and its holomorphic extension
\begin{math}
  \mathbf{u}_0\equiv
  \tilde{\mathbf{u}}_0\colon \mathfrak{X}^{(r)}\to \CC^M
\end{math}
stated above, for each $z_0\in \mathfrak{X}^{(r)}$,
the ``shifted'' function
\begin{math}
  x\mapsto \mathbf{u}_0^{(z_0)}(x)\eqdef \mathbf{u}_0(x + z_0)
  \colon \RR^N\to \CC^M
\end{math}
belongs to
\begin{math}
  \mathbf{L}^2(\RR^N)\cap \mathbf{W}^{2m,p}(\RR^N)
\end{math}
where
\begin{math}
  2 + \genfrac{}{}{}1{N}{m} < p < \infty \,.
\end{math}
Consequently, we have also
$\mathbf{u}_0^{(z_0)}\in \mathbf{B}^{s;p,p}(\RR^N)$
and, thus, we may apply
the local (in time) existence and uniqueness result
(Theorem~\ref{thm-Clem-Li}) on
a short time interval
$[t_0,T_1]\subset [0,T]$ with the initial condition
$\mathbf{u}(\,\cdot\,,t_0) = \mathbf{u}_0^{(z_0)}$
in $\mathbf{B}^{s;p,p}(\RR^N)$ at time $t = t_0\in [0,T)$
to conclude that
the spatially ``shifted'' Cauchy problem \eqref{e:Cauchy+z_0}
possesses a unique weak solution
\begin{math}
  \mathbf{u}^{(z_0)}\in
  C\left( [t_0,T_1]\to \mathbf{B}^{s;p,p}(\RR^N)\right) ,
\end{math}
local in time.
Of course, the length of the time interval $[t_0,T_1]$ depends on
the shift $z_0\in \mathfrak{X}^{(r)}$; more precisely,
on its imaginary part $y_0 = \IM z_0\in Q^{(r)}$.
However, when making use of the abstract reformulation
\eqref{e:abstr_Cauchy}
of the Cauchy problem \eqref{e:Cauchy+z_0}, we must guarantee that
the \underline{values} of the (unique) strict solution
\begin{math}
  \mathbf{u}^{(z_0)}\colon [t_0,T_1]\to \mathbf{B}^{s;p,p}(\RR^N)
\end{math}
to the Cauchy problem \eqref{e:Cauchy+z_0},
the (continuous) ``shifted'' function
\begin{math}
  t\mapsto\hfil\break \mathbf{u}(\,\cdot\, + z_0,t)
  = \mathbf{u}^{(z_0)}(\,\cdot\,,t)
  \colon [t_0,T_1]\to \mathbf{B}^{s;p,p}(\RR^N) ,
\end{math}
\underline{stay} in the bounded open set
$U = \tilde{U}\subset E_{ 1 - \frac{1}{p} ,\, p }$
for all times $t\in [t_0,T_1]$
(cf.\ {\rm Remark~\ref{rem-Analytic_glob}}).
In order to avoid this technical problem,
we make the following {\it global existence\/} hypothesis, cf.\
Theorem~\ref{thm-Main}:

\begin{hypo}\nopagebreak
\begingroup\rm
\begin{enumerate}
\setcounter{enumi}{0}
\renewcommand{\labelenumi}{{\bf (T)}}
\item
\makeatletter
\def\@currentlabel{{\bf T}}\label{hyp:global}
\makeatother
The original Cauchy problem \eqref{e:Cauchy}, i.e.,
problem \eqref{e:Cauchy+z_0} with $z_0 = 0$ and the initial data
\begin{math}
  \mathbf{u}_0 = \widehat{\mathbf{u}}_0\in \mathbf{B}^{s;p,p}(\RR^N)
\end{math}
at $t = 0$, possesses a global weak solution
\hfil\break
\begin{math}
  \widehat{\mathbf{u}}\in
  C\left( [0,T]\to \mathbf{B}^{s;p,p}(\RR^N)\right) .
\end{math}
\hfill\Square
\end{enumerate}
\endgroup
\end{hypo}
\par\vskip 10pt

Now define the set
\begin{math}
  U\subset E_{ 1 - \frac{1}{p} ,\, p } = \mathbf{B}^{s;p,p}(\RR^N)
\end{math}
that appears in
{\rm Hypotheses} \eqref{hyp:A}, \eqref{hyp:f} and
                 \eqref{hyp:A_anal}, \eqref{hyp:f_anal}
as follows: First, let
\begin{equation}
\label{def:U_0}
  U_0 = \mathrm{conv}
    \biggl(
      B_{R_0}(\mathbf{0}) \cup \bigcup_{t\in [0,T]}
      B_{R_0}( \widehat{\mathbf{u}}(\,\cdot\,,t) )
    \biggr) \subset \mathbf{B}^{s;p,p}(\RR^N)
\end{equation}
be the convex hull of the union of open balls
\begin{align*}
& B_{R_0}(\mathbf{0})\eqdef
  \left\{ \mathbf{w}\in \mathbf{B}^{s;p,p}(\RR^N)\colon
       \| \mathbf{w}\|_{ B^{s;p,p}(\RR^N) } < R_0
  \right\} \subset \mathbf{B}^{s;p,p}(\RR^N)
\\
&   \quad\mbox{ and }\quad
  B_{R_0}(\mathbf{v})\eqdef \mathbf{v} + B_{R_0}(\mathbf{0})
  \;\mbox{ with }\;
  \mathbf{v} = \widehat{\mathbf{u}}(\,\cdot\,,t)
  \;\mbox{ for }\; t\in [0,T] \,,
\end{align*}
where their radius $R_0\in (0,\infty)$ is an arbitrary positive number.
Of course, the symbol ``$\mathbf{0}$'' stands for the zero function in
$\mathbf{B}^{s;p,p}(\RR^N)$.
Alternatively, we may take
$U_0 = B_{R_0}(\widehat{\mathbf{u}}_0)$ to be any open ball in
$\mathbf{B}^{s;p,p}(\RR^N)$ centered at
$\widehat{\mathbf{u}}_0 = \widehat{\mathbf{u}}(\,\cdot\,,0)$
with (sufficiently large) radius $R_0\in (0,\infty)$, such that
$\mathbf{0}\in B_{R_0}(\widehat{\mathbf{u}}_0)$
and
\begin{math}
  \widehat{\mathbf{u}}(\,\cdot\,,t)
  \in B_{R_0}(\widehat{\mathbf{u}}_0)
\end{math}
holds for every $t\in [0,T]$.
However, this choice of $R_0$ would not fit in
Example~\ref{exam-nontech} in Section~\ref{s:Appl} below.
Clearly, $U_0$ is a bounded open set in $\mathbf{B}^{s;p,p}(\RR^N)$.
From now on we take the initial values
\begin{math}
    \mathbf{u}^{(z_0)}(\,\cdot\,,t_0) = \mathbf{u}_0^{(z_0)}
  = \mathbf{u}^{(z_0)}(\,\cdot\, + z_0,t_0)
  \in \mathbf{B}^{s;p,p}(\RR^N)
\end{math}
at time $t = t_0\in [0,T)$
in the Cauchy problem \eqref{e:Cauchy+z_0}
(and similar related initial value problems)
from the set $U_0$ only, i.e.,
$\mathbf{u}_0^{(z_0)}\in U_0$.
This choice will guarantee that
the \underline{values} of the (unique) strict solution
\begin{math}
  \mathbf{u}^{(z_0)}\colon [t_0,T_1]\to \mathbf{B}^{s;p,p}(\RR^N)
\end{math}
to the ``shifted'' Cauchy problem \eqref{e:Cauchy+z_0}
\underline{stay} for all times $t\in [t_0,T_1]$
in the bounded open set
\begin{math}
  U = \tilde{U}\subset E_{ 1 - \frac{1}{p} ,\, p }
    = \mathbf{B}^{s;p,p}(\RR^N)
\end{math}
defined next, $U\supset U_0$.

We set
\begin{align}
\label{def:U}
& U = \bigcup
    \left\{ B_{R_0}(\mathbf{v})\colon \mathbf{v}\in U_0\right\}
    \subset \mathbf{B}^{s;p,p}(\RR^N) \,; \qquad\mbox{ hence, }\quad
\\
\nonumber
& U =
  \left\{ \mathbf{w}\in \mathbf{B}^{s;p,p}(\RR^N)\colon
       \| \mathbf{w} - \mathbf{v}\|_{ B^{s;p,p}(\RR^N) } < R_0
       \,\mbox{ for some }\, \mathbf{v}\in U_0
  \right\} \,.
\end{align}
One may call $U$ the {\em open $R_0$\--neighborhood\/}
of $U_0$ in $\mathbf{B}^{s;p,p}(\RR^N)$.
Also $U = \tilde{U}$ is a bounded, open, and convex set in
the complex Besov space $\mathbf{B}^{s;p,p}(\RR^N)$
and, consequently, in $\mathbf{W}^{m,p}(\RR^N)$ and in
\begin{math}
    \mathbf{C}^m_{\mathrm{bdd}}(\RR^N)
  = \mathbf{C}^m(\RR^N)\cap \mathbf{W}^{m,\infty}(\RR^N) ,
\end{math}
as well, owing to the continuous Sobolev\-(-Besov) imbeddings
\begin{equation}
\label{e:Besov-imbed}
\begin{aligned}
& B^{s;p,p}(\RR^N) \hookrightarrow W^{m,p}(\RR^N) \;\mbox{ and }\;
  B^{s;p,p}(\RR^N) \hookrightarrow
  C^m(\RR^N)\cap W^{m,\infty}(\RR^N) \,,
  \quad\mbox{ respectively, }\;
\\
& \quad\mbox{ where }\;
  \left( 1 + \genfrac{}{}{}1{N}{2m+N} \right)
  m < s = 2m\left( 1 - \genfrac{}{}{}1{1}{p} \right) < 2m
    \;\mbox{ thanks to }\;
  2 + \genfrac{}{}{}1{N}{m} < p < \infty \,;
\end{aligned}
\end{equation}
see, e.g.,
{\sc R.~A.\ Adams} and {\sc J.~J.~F.\ Fournier}
\cite[Chapt.~7]{AdamsFournier}, Theorem 7.34(a,c), p.~231.
This shows that, for any function $\mathbf{w}\in U$,
all partial derivatives 
\begin{math}
  \frac{ \partial^{|\beta|}\mathbf{w} }{\partial x^{\beta}} \,,
\end{math}
$\beta = (\beta_1,\dots,\beta_N)$ $\in (\ZZ_+)^N$, of order
$|\beta| = \beta_1 + \dots + \beta_N\leq m$,
are uniformly bounded on $\RR^N$,
\begin{equation}
\label{e:|w|<C}
  \genfrac{|}{|}{}0{ \partial^{|\beta|}\mathbf{w}(x) }{\partial x^{\beta}}
  \leq C\equiv C(U) = \mathrm{const} < \infty
  \quad\mbox{ for all }\; x\in \RR^N \,,
\end{equation}
by a constant $C\in \RR_+$ depending solely on $U$.
These partial derivatives are arguments in the reaction function
\begin{math}
  \mathbf{f}
    \left( x, t;
    \left( \frac{ \partial^{|\beta|}\mathbf{u} }{\partial x^{\beta}}
    \right)_{|\beta|\leq m}
    \right)
\end{math}
on the right\--hand side of
in {\rm eqs.}\ \eqref{e:Cauchy} and~\eqref{e:Cauchy+z_0}.
In our {\rm Hypothesis} \eqref{hy:f_Hardy^2} on $\mathbf{f}$ we take
$\Sigma\subset \CC^{M\tilde{N}}$
to be the closed polydisc
\begin{equation*}
    \Sigma = [ \overline{D}_C(0) ]^{M\tilde{N}}
    \quad\mbox{ where }\;
  \overline{D}_C(0)\eqdef \{ z\in \CC\colon |z|\leq C\}
    \;\mbox{ is a closed disc. }
\end{equation*}
This {\em restriction\/} on the values of the (unique) strict solution
to the bounded open set
$U = \tilde{U}\subset E_{ 1 - \frac{1}{p} ,\, p }$
will be used in applications to semilinear {\sc Heston\/}\--type models
in {\sl ``Mathematical Finance''\/}
treated in Section~\ref{s:Appl}.

From {\rm Hypothesis} \eqref{hy:f_Hardy^2} we deduce immediately that
each component
\begin{math}
  f_j\colon
  \overline{\Omega}\times [ \overline{D}_C(0) ]^{M\tilde{N}} \to \CC
\end{math}
of the reaction function
$\mathbf{f} = (f_1,\dots,f_M)$
is continuously differentiable (i.e., of class $C^1$)
with the time derivative
\begin{math}
  \frac{\partial}{\partial t} f_j(x,t;X)
\end{math}
and all argument first\--order partial derivatives
\begin{equation*}
  \frac{\partial f_j}{ \partial X_{\beta,k} } (x,t;X) \,,
    \quad\mbox{ for $|\beta|\leq m$ and $j,k=1,2,\dots,M$, }
\end{equation*}
being uniformly bounded on $\Omega\times \Sigma$.
Consequently, each $f_j$ is Lipschitz\--continuous with respect to
the variables $t$ and $X_{\beta,k}$, uniformly on $\Omega\times \Sigma$.

Recalling the continuous Sobolev\-(-Besov) imbeddings
\eqref{e:Besov-imbed}, i.e.,
\begin{equation*}
  B^{s;p,p}(\RR^N) \hookrightarrow W^{m,p}(\RR^N)\cap
  C^m(\RR^N)\cap W^{m,\infty}(\RR^N) \,,
\end{equation*}
and the $L^p$\--integrability condition in ineq.~\eqref{e:f_Hardy^2},
we have just proved the following lemma
(cf.\ {\rm Hypotheses} \eqref{hyp:A}, \eqref{hyp:f} and
                       \eqref{hyp:A_anal}, \eqref{hyp:f_anal}):

\begin{lemma}\label{lem-f_Lipschitz}
Assume that\/
\begin{math}
  \mathbf{f}\colon
  \overline{\Omega}\times \CC^{M\tilde{N}}\to \CC^M
\end{math}
satisfies {\rm Hypothesis} \eqref{hy:f_Hardy^2}.
Let also {\rm Hypothesis} \eqref{hyp:global} be satisfied and define
\begin{math}
  U = \tilde{U}\subset E_{ 1 - \frac{1}{p} ,\, p }
    = \mathbf{B}^{s;p,p}(\RR^N)
\end{math}
by {\rm eq.}~\eqref{def:U}.
Then the Nemytskii operator
$\mathbf{F}\colon [0,T]\times U\to E_0 = \mathbf{L}^p(\RR^N)$
defined by
\begin{equation}
\label{def:F}
  \mathbf{F}(t,\mathbf{v})(x) \eqdef
  \mathbf{f}
    \left( x, t;
    \left( \frac{ \partial^{|\beta|}\mathbf{v} }{\partial x^{\beta}}
    \right)_{|\beta|\leq m}
    \right) \,,\quad x\in \RR^N \,,
\end{equation}
for all\/ $t\in [0,T]$ and for all\/ $\mathbf{v}\in U$,
enjoys the following properties:
\begin{itemize}
\item[{\rm (a)}]
$\mathbf{F}\colon [0,T]\times U\to E_0$
is a Lipschitz continuous mapping, i.e.,
$\mathbf{F}$ satisfies {\rm Hypothesis} \eqref{hyp:f}.
\item[{\rm (b)}]
Also the substitution mapping
\begin{math}
  \mathfrak{F}\colon C([0,T]\to U)\to L^p((0,T)\to E_0)
\end{math}
defined by
\begin{equation*}
  \mathfrak{F}(\mathbf{v})(t) \eqdef
  \mathbf{F}(t, \mathbf{v}(t)) \,,\quad t\in [0,T] \,,
  \quad\mbox{ for all }\; \mathbf{v}\in C([0,T]\to U) \,,
\end{equation*}
is Lipschitz\--continuous with values in $L^p((0,T)\to E_0)$.
\end{itemize}
\end{lemma}
\par\vskip 10pt

\proof
The only claims in {\rm Parts} {\rm (a)} and {\rm (b)}, respectively,
that remain to be verified are that
$\mathbf{F}$ maps $[0,T]\times U$ into $E_0$ and
$\mathfrak{F}$ maps $C([0,T]\to U)$ into $L^p((0,T)\to E_0)$.

{\it Part}~{\rm (a)}.
For each component $F_j$ of $\mathbf{F} = (F_1,\dots,F_M)$
we derive from {\rm eq.}~\eqref{def:F}
\begin{align}
\nonumber
&   F_j(t,\mathbf{v}_1)(x) - F_j(t,\mathbf{v}_2)(x) =
\\
\label{def:F:v_1-v_2}
&
\begin{aligned}
    \sum_{|\beta|\leq m} \sum_{k=1}^M
&   \left[
    \int_0^1
  \frac{\partial f_j}{ \partial X_{\beta,k} }
    \left( x, t;
    \left(
    (1-\theta)\,
    \frac{ \partial^{|\beta|}\mathbf{v}_1 }{\partial x^{\beta}}
  + \theta\,
    \frac{ \partial^{|\beta|}\mathbf{v}_2 }{\partial x^{\beta}}
    \right)_{|\beta|\leq m}
    \right) \,\mathrm{d}\theta
    \right]
\\
& {}
  \times\,
    \frac{ \partial^{|\beta|} }{ \partial x^{\beta} }\,
    \left( v_{1,k}(x) - v_{2,k}(x)\right) \,,\quad x\in \RR^N \,;\quad
    j=1,2,\dots,M \,,
\end{aligned}
\end{align}
for all $t\in [0,T]$ and for all
$\mathbf{v}_1, \mathbf{v}_2\in U$.
Notice that the partial derivatives with respect to $x^{\beta}$
emerge from the chain rule applied to the right\--hand side of
eq.~\eqref{def:F} using the partial derivatives 
\begin{math}
  \frac{\partial}{\partial X_{\beta,k}} f_j(x,t;X)
\end{math}
with respect to the argument $X_{\beta,k}$ and the convex combination
\begin{math}
  \mathbf{w} =\hfil\break
  (1-\theta)\mathbf{v}_1 + \theta\mathbf{v}_2\in U
\end{math}
for $0\leq \theta\leq 1$, thanks to $U$ being convex.
Consequently, with this abbreviation for $\mathbf{w}$
and our choice of the constant $C\equiv C(U)$ in ineq.~\eqref{e:|w|<C},
all partial derivatives
\begin{math}
  \frac{ \partial^{|\beta|}\mathbf{w} }{\partial x^{\beta}} \,,
\end{math}
$\beta = (\beta_1,\dots,\beta_N)$ $\in (\ZZ_+)^N$, of order
$|\beta| = \beta_1 + \dots + \beta_N\leq m$,
are uniformly bounded on $\RR^N$, by ineq.~\eqref{e:|w|<C}.
By {\rm Hypothesis} \eqref{hy:f_Hardy^2},
cf.\ Remark~\ref{rem-f_Hardy^2}, all partial derivatives
\begin{math}
  \frac{\partial}{ \partial X_{\beta,k} } f_j(x,t;\,\cdot\,)
  \colon \Sigma\to \CC
\end{math}
are uniformly bounded,
\begin{align}
\label{e:|df/dX|<C}
& \genfrac{|}{|}{}0{ \partial f_j(x,t;X) }{ \partial X_{\beta,k} }
  \leq C_1\equiv C_1(C(U)) = \mathrm{const} < \infty
  \hskip 30pt \mbox{ (cf.\ ineq.~\eqref{e:f_j'}) }
\\
\nonumber
& \quad\mbox{ for all }\; (x,t)\in \Omega
  \;\mbox{ and for all }\;
  X = \left( ( X_{\beta,k} )_{|\beta|\leq m} \right)_{k=1}^M
      \in \Sigma \,,
\end{align}
by a constant $C_1\in \RR_+$ depending solely on $C(U)$.
We apply these estimates to the integrands in eq.~\eqref{def:F:v_1-v_2}
to conclude that there is a Lipschitz constant
$L\equiv L(U)\in \RR_+$ depending solely on $U$
(through the constant $C_1(C(U))\geq 0$ in ineq.~\eqref{e:|df/dX|<C} above),
such that
\begin{equation}
\label{est:F:v_1-v_2}
    \left|
      \mathbf{F}(t,\mathbf{v}_1)(x)
    - \mathbf{F}(t,\mathbf{v}_2)(x)
    \right|
  \leq L\sum_{|\beta|\leq m} \sum_{k=1}^M
    \left|
    \frac{ \partial^{|\beta|} }{ \partial x^{\beta} }\,
    \left( v_{1,k}(x) - v_{2,k}(x)\right)
    \right| \,,\quad x\in \RR^N \,,
\end{equation}
for all $t\in [0,T]$ and for all
$\mathbf{v}_1, \mathbf{v}_2\in U$.

We recall $U\subset B^{s;p,p}(\RR^N)$ and the imbeddings
in {\rm eq.}~\eqref{e:Besov-imbed}
in order to deduce from ineq.~\eqref{est:F:v_1-v_2}
that all mappings
\begin{math}
  \mathbf{F}(t,\,\cdot\,)\colon U\to E_0 = \mathbf{L}^p(\RR^N)
\end{math}
are uniformly Lipschitz\--continous (with the same Lipschitz constant)
for all $t\in [0,T]$.
Here, we single out the special case of
$\mathbf{v}_1 = \mathbf{v}\in U$ being arbitrary and
$\mathbf{v}_2 = \mathbf{0}\in U$, i.e.,
$\mathbf{v}_2(x)\equiv \mathbf{0}\in \CC^M$ for all $x\in \RR^N$.
Then ineq.~\eqref{est:F:v_1-v_2} yields
\begin{equation}
\label{est:F:v-0}
    \left|
      \mathbf{F}(t,\mathbf{v})(x)
    \right|
  \leq
    \left|
      \mathbf{F}(t,\mathbf{0})(x)
    \right|
  + L\sum_{|\beta|\leq m} \sum_{k=1}^M
    \left|
    \frac{ \partial^{|\beta|} }{ \partial x^{\beta} }\, v_k(x)
    \right| \,,\quad x\in \RR^N \,,
\end{equation}
for all $t\in [0,T]$ and for all $\mathbf{v}\in U$, where
\begin{equation*}
    \mathbf{F}(t,\mathbf{0})(x)
  = \mathbf{f}(x,t; \vec{\mathbf{0}}) \,,\quad x\in \RR^N \,,\quad
    \vec{\mathbf{0}} = \left( 0\right)_{|\beta|\leq m}
  \equiv (0,\dots,0)\in \CC^{M\tilde{N}} \,.
\end{equation*}
satisfies
$\mathbf{F}(t,\mathbf{0})\in E_0 = \mathbf{L}^p(\RR^N)$, by
the $L^p$\--integrability condition in ineq.~\eqref{e:f_Hardy^2}, i.e.,
\begin{equation*}
    \| \mathbf{F}(t,\mathbf{0}) \|_{E_0}
  = \left(
  \int_{\RR^N} |\mathbf{f} (x,t; \vec{\mathbf{0}})|^p \,\mathrm{d}x
    \right)^{1/p}
  \leq K \quad\mbox{ for all $t\in [0,T]$, }
\end{equation*}
where $K\in (0,\infty)$ is a constant.
Now it follows from ineq.~\eqref{est:F:v-0} above that also
$\mathbf{F}(t,\mathbf{v})\in E_0$ holds
for all $t\in [0,T]$ and for all $\mathbf{v}\in U$, as claimed.

{\it Part}~{\rm (b)}.
Analogous results for the mapping
\begin{math}
  \mathfrak{F}\colon C([0,T]\to U)\to L^p((0,T)\to E_0)
\end{math}
follow from those we have just proved for
$\mathbf{F}(t,\,\cdot\,)\colon U\to E_0$, $t\in [0,T]$.
Namely, the ``supremum'' (or ``maximum'') norm on the Banach space
\begin{math}
  C\left( [0,T]\to E_{ 1 - \frac{1}{p} ,\, p }\right)
\end{math}
of all continuous functions
$u\colon [0,T]\to E_{ 1 - \frac{1}{p} ,\, p }$
is defined by
\begin{equation*}
  \vertiii{u}_{ L^{\infty}(0,T) }
  \eqdef
  \| u\|_{ C\bigl( [0,T]\to E_{ 1 - \frac{1}{p} ,\, p }\bigr) }
  = \sup_{ t\in [0,T] } \| u(t)\|_{ E_{ 1 - \frac{1}{p} ,\, p } } \,.
\vspace*{-8.0mm}
\end{equation*}
\null\hfill\qed
\par\vskip 10pt

\begin{remark}\label{rem-f_Lipschitz}\nopagebreak
\begingroup\rm
Analogous result (to {\rm Lemma~\ref{lem-f_Lipschitz}} above)
holds also for the ``shifted'' Nemytskii operator
\begin{math}
  \mathbf{F}^{(z_0)}\colon [0,T]\times U\to E_0 = \mathbf{L}^p(\RR^N) ,
\end{math}
by a complex vector\/
$z_0\in \mathfrak{X}^{(r)}$, defined by
\begin{equation}
\label{def:F+z_0}
  \mathbf{F}^{(z_0)}(t,\mathbf{v})(x) \eqdef
  \mathbf{f}
    \left( x + z_0, t;
    \left( \frac{ \partial^{|\beta|}\mathbf{v} }{\partial x^{\beta}}
    \right)_{|\beta|\leq m}
    \right) \,,\quad x\in \RR^N \,,
\end{equation}
for all $t\in [0,T]$ and for all $\mathbf{v}\in U$.
Both constants, $C_1\equiv C_1(C(U))$ and $L\equiv L(U)$
in inequalities \eqref{e:|df/dX|<C} and \eqref{est:F:v_1-v_2},
respectively, are independent from the shift by
$z_0\in \mathfrak{X}^{(r)}$ in case $x\in \RR^N$ is replaced by $x + z_0$,
thanks to $(x,t)\in \Omega$ where the domain
\begin{math}
  \Omega = \Gamma^{(T_0)}_T (r_0,\vartheta_0)
         = \mathfrak{X}^{(r_0)}\times \Delta_{\vartheta_0}^{T_0,T}
  \subset \CC^N\times \CC
\end{math}
has been introduced in Section~\ref{s:Main} before
{\rm Hypotheses}
\eqref{hy:UnifEll}, \eqref{hy:Ellipt}, and\/ \eqref{hy:f_Hardy^2}
(cf.\ eq.~\eqref{e:Gamma} and Theorem~\ref{thm-Main}).
\hfill\Square
\endgroup
\end{remark}
\par\vskip 10pt

In what follows, the initial data $\mathbf{u}_0$
in the Cauchy problem \eqref{e:Cauchy+z_0} at time $t_0\in [0,T)$
have nothing to do with the initial value
\begin{math}
    \widehat{\mathbf{u}}(\,\cdot\,,0)
  = \widehat{\mathbf{u}}_0\in \mathbf{B}^{s;p,p}(\RR^N)
\end{math}
of the global weak solution
\begin{math}
  \widehat{\mathbf{u}}\in
  C\left( [0,T]\to \mathbf{B}^{s;p,p}(\RR^N)\right)
\end{math}
in our {\rm Hypothesis}~\eqref{hyp:global}, except for the restriction
$\mathbf{u}_0\in U_0$, where $U_0$ is determined by
the values of the solution
$\widehat{\mathbf{u}}(\,\cdot\,,t)$ in $\mathbf{B}^{s;p,p}(\RR^N)$
for $0\leq t\leq T$, see eq.~\eqref{def:U_0}.

We take advantage of Remark~\ref{rem-f_Lipschitz}
to recall that, given a holomorphic function
\begin{math}
  \mathbf{u}_0\equiv \tilde{\mathbf{u}}_0\colon
  \mathfrak{X}^{(r)}\to \CC^M
\end{math}
as described before {\rm Hypothesis}~\eqref{hyp:global},
the spatially ``shifted'' Cauchy problem \eqref{e:Cauchy+z_0}
possesses a unique weak solution
\begin{math}
  \mathbf{u}^{(z_0)}\equiv \mathbf{u}^{(z_0,t_0)}
  \in C\left( [t_0,T_1]\to \mathbf{B}^{s;p,p}(\RR^N)\right) ,
\end{math}
local in time, for every fixed shift $z_0\in \mathfrak{X}^{(r)}$,
locally uniformly in the complex domain,
provided its real and imaginary parts,
\begin{math}
  x_0 = \RE z_0, y_0 = \IM z_0
  \in \overline{Q}^{(r_1)} \subset Q^{(r)} ,
\end{math}
are small enough, i.e.,
$\,\max\{ |x_0|_{\infty} ,\, |y_0|_{\infty} \} \leq r_1\,$
($< r = r_0$).
Indeed, it suffices to choose $r_1\in (0,r)$ so small that
each ``shifted'' function
\begin{math}
  x\mapsto \mathbf{u}_0^{(z_0)}(x) = \mathbf{u}_0(x + z_0)
  \colon \RR^N\to \CC^M
\end{math}
(serving as the initial data at time $t = t_0\in [0,T)$),
with $z_0 = x_0 + \ii y_0$ satisfying
$\,\max\{ |x_0|_{\infty} ,\, |y_0|_{\infty} \} \leq r_1\,$,
lies in the open set $U_0$ specified in eq.~\eqref{def:U_0}
after {\rm Hypothesis} \eqref{hyp:global}, i.e.,
\begin{math}
  \mathbf{u}_0^{(z_0)} = \mathbf{u}_0(\,\cdot\, + z_0)
  \in U_0\subset E_{ 1 - \frac{1}{p} ,\, p } \,.
\end{math}
Recall that $\overline{Q}^{(r_1)}$ and
\begin{math}
  \overline{\mathfrak{X}}^{(r_1)} = \RR^N + \ii \overline{Q}^{(r_1)}
\end{math}
stand for the respective closures of the cube
$Q^{(r_1)}\subset \RR^N$ and the strip
$\mathfrak{X}^{(r_1)}\in \CC^N$.
Indeed, our choice of $r_1\in (0,r)$ small enough to guarantee
\begin{math}
  \mathbf{u}_0^{(z_0)} =
  \mathbf{u}_0(\,\cdot\, + x_0 + \ii y_0) \in U_0
\end{math}
for every $y_0\in Q^{(r_1)}$ (and for all $x_0\in \RR^N$),
is possible thanks to the closed cube
$\overline{Q}^{(r_1)}\subset \RR^N$ being compact.
For ``small'' shifts $z_0 = x_0 + \ii y_0\in \CC^N$
we introduce the complex cube
\begin{equation*}
    Q^{(r_1)}_{\CC} \eqdef Q^{(r_1)} + \ii Q^{(r_1)}
  = \{ z = x + \ii y\in \CC\colon
       \max\{ |x|_{\infty} ,\, |y|_{\infty} \} < r_1
    \}
  \subset \mathfrak{X}^{(r_1)} = \RR^N + \ii Q^{(r_1)}
\end{equation*}
and denote by $\overline{Q}^{(r_1)}_{\CC}$ its closure in $\CC^N$;
hence,
\begin{equation*}
\begin{aligned}
    \overline{Q}^{(r_1)}_{\CC}
  = \overline{Q}^{(r_1)} + \ii \overline{Q}^{(r_1)}
& {}
  = \{ z = x + \ii y\in \CC\colon
       \max\{ |x|_{\infty} ,\, |y|_{\infty} \} \leq r_1
    \}
\\
& {}
  \subset
  \overline{\mathfrak{X}}^{(r_1)} = \RR^N + \ii \overline{Q}^{(r_1)}
  \subset \mathfrak{X}^{(r)} \,,
\end{aligned}
\end{equation*}
thanks to $0 < r_1 < r$.
We will call $\overline{Q}^{(r_1)}_{\CC}$
the {\em\bfseries small shift cube\/}.

To obtain a {\it spatially holomorphic\/} extension
\begin{math}
  \mathbf{u}\equiv \tilde{\mathbf{u}}
  \colon \mathfrak{X}^{(r_1)}\times [0,T]\to \CC^M
\end{math}
of the function
\begin{math}
  \widehat{\mathbf{u}}\in
  C\left( [0,T]\to \mathbf{B}^{s;p,p}(\RR^N)\right)
\end{math}
from {\rm Hypothesis} \eqref{hyp:global}
to the spatio\--temporal strip
$\mathfrak{X}^{(r_1)} \times [0,T]$,
having a (unique) continuous extension
\begin{math}
  \mathbf{u}\colon \overline{\mathfrak{X}}^{(r_1)}\times [0,T]\to \CC^M
\end{math}
to the closed strip
$\overline{\mathfrak{X}}^{(r_1)} \times [0,T]$, such that
\begin{math}
  (y_0,t)\mapsto \mathbf{u}(\,\cdot\, + \ii y_0,t)
  \colon \overline{Q}^{(r_1)}\times [0,T]\to
  \mathbf{B}^{s;p,p}(\RR^N)
\end{math}
is continuous and satisfies
$\mathbf{u}(\,\cdot\, + \ii y_0,t)\in U$
for every $y_0\in \overline{Q}^{(r_1)}$
and for every $t\in [0,T]$, we construct $\mathbf{u}$
first {\it locally in time\/} on a short time interval
$[t_0, t_0+T_1]\subset [0,T]$ as follows, where
$T_1 > 0$ is small enough:

Let
\begin{math}
  \mathbf{u}^{(z_0)}\colon
  t\mapsto \mathbf{u}^{(z_0)}(\,\cdot\,,t)
    \equiv \mathbf{u}^{(z_0,t_0)}(\,\cdot\,,t)
\end{math}
from
\begin{math}
  C\left( [t_0, t_0+T_1]\to \mathbf{B}^{s;p,p}(\RR^N)\right)
\end{math}
be as above, such that $z_0 = \ii y_0$ and
$\mathbf{u}^{(\ii y_0)}(\,\cdot\,,t)\in U$ for every pair
$(y_0,t)\in \overline{Q}^{(r_1)}\times [t_0, t_0+T_1]$.
Given any $y_0\in \overline{Q}^{(r_1)}$, let us define
\begin{equation}
\label{def:u-ext}
  \mathbf{u}(x + \ii y_0,t)\eqdef \mathbf{u}^{(\ii y_0)}(x,t)
  \quad\mbox{ for all }\; (x,t)\in \RR^N\times [t_0, t_0+T_1] \,.
\end{equation}
Clearly,
\begin{math}
  \mathbf{u}\colon
  \overline{\mathfrak{X}}^{(r_1)} \times [t_0, t_0+T_1]\to \CC^M
\end{math}
is a well\--defined mapping; it has the following properties:

\begin{proposition}\label{prop-smooth_ext}
Let\/ $M,N\geq 1$, $0 < T < \infty$, and assume that all\/
{\rm Hypotheses} \eqref{hy:UnifEll}, \eqref{hy:Ellipt}, and\/
\eqref{hy:f_Hardy^2}
are satisfied with some constants\/
$0 < r_0 < \infty$, $0 < T_0\leq T$, and\/ $0 < \vartheta_0 < \pi / 2$.
Furthermore, assume that\/
\begin{math}
  \widehat{\mathbf{u}}\in
  C\left( [0,T]\to \mathbf{B}^{s;p,p}(\RR^N)\right)
\end{math}
is a globally defined weak solution
to the original Cauchy problem \eqref{e:Cauchy}, i.e.,
{\rm Hypothesis} \eqref{hyp:global} is valid.
Let the sets\/
\begin{math}
  U_0\subset U\subset
  E_{ 1 - \frac{1}{p} ,\, p } = \mathbf{B}^{s;p,p}(\RR^N)
\end{math}
be specified as in {\rm eqs.}\ \eqref{def:U_0} and~\eqref{def:U}.
Given any\/ $t_0\in [0,T)$, let\/
$\mathbf{u}_0\in U_0\subset \mathbf{B}^{s;p,p}(\RR^N)$
be any initial data at time $t = t_0$, such that\/
$\mathbf{u}_0$ has a holomorphic extension
\begin{math}
  \tilde{\mathbf{u}}_0\colon \mathfrak{X}^{(r)}\to \CC^M
\end{math}
as described in {\rm Hypothesis} \eqref{hyp:analyt} and\/
{\rm Remark~\ref{rem-half-width^r}}
\begingroup\rm
(we identify $\mathbf{u}_0\equiv \tilde{\mathbf{u}}_0$),
\endgroup
with
\begin{math}
  \mathbf{u}_0^{(z_0)} = \mathbf{u}_0(\,\cdot\, + z_0)
  \in U_0\subset E_{ 1 - \frac{1}{p} ,\, p }
\end{math}
whenever\/
\begin{math}
  z_0 = x_0 + \ii y_0\in \overline{Q}^{(r_1)}_{\CC}
  \subset \overline{\mathfrak{X}}^{(r_1)}
\end{math}
for some $r_1\in (0,r)$.
(We have set\/ $r = r_0$; $r_1$ may depend on $\mathbf{u}_0$.)

Then there exists a number\/ $T_1\in (0,T-t_0]$,
depending on $r_1$ and $U$, but not on $t_0$, such that
the Cauchy problem \eqref{e:Cauchy+z_0} on the (local) time interval\/
$[t_0, t_0+T_1]\subset [0,T]$ with the initial condition
$\mathbf{u}(\,\cdot\,,t_0) = \mathbf{u}_0^{(z_0)}$
possesses a unique weak solution
\begin{math}
  \mathbf{u}^{(z_0)}\equiv \mathbf{u}^{(z_0,t_0)}\in
  C\left( [t_0, t_0+T_1]\to \mathbf{B}^{s;p,p}(\RR^N)\right) ,
\end{math}
such that
$\mathbf{u}^{(z_0)}(\,\cdot\,,t)\in U$
holds for every $t\in [t_0, t_0+T_1]$.
The family $\mathbf{u}^{(z_0)}$, parametrized by
\begin{math}
  z_0 = x_0 + \ii y_0\in \overline{Q}^{(r_1)}_{\CC} ,
\end{math}
has the following properties:
\begin{itemize}
\item[{\rm (a)}]
$\mathbf{u}^{(z_0)}(x,t) = \mathbf{u}^{(\ii y_0)}(x + x_0,t)$
holds for all\/ $(x,t)\in \RR^N\times [t_0, t_0+T_1]$ and for all\/
$z_0\in \overline{Q}^{(r_1)}_{\CC}$;
consequently, even for all\/
$z_0\in \overline{\mathfrak{X}}^{(r_1)}$.
\item[{\rm (b)}]
The mapping
\begin{math}
  \mathbf{u}\colon
  \overline{\mathfrak{X}}^{(r_1)} \times [t_0, t_0+T_1]\to \CC^M
\end{math}
defined in {\rm eq.}~\eqref{def:u-ext} satisfies
\begin{math}
    \mathbf{u}(x + x_0 + \ii y_0,t)
  = \mathbf{u}^{(x_0 + \ii y_0)}(x,t)
\end{math}
for all\/ $(x,t)\in \RR^N\times [t_0, t_0+T_1]$ and for all\/
$z_0 = x_0 + \ii y_0\in \CC^N$ with $|y_0|_{\infty}\leq r_1$.
\item[{\rm (c)}]
The mapping
\begin{math}
  \mathbf{u}\colon
  \overline{\mathfrak{X}}^{(r_1)} \times [t_0, t_0+T_1]\to \CC^M
  \colon (z,t) = (x + \ii y,t)\mapsto \mathbf{u}(x + \ii y,t)
\end{math}
is continuously (partially) differentiable with respect to
all the real variables $x_i$ and $y_i$
{\rm ($i=1,2,\dots,N$)} in
$x = (x_1,\dots,x_N)$ and $y = (y_1,\dots,y_N)$ in $\RR^N$ with
$|y|_{\infty}\leq r_1$.
\item[{\rm (d)}]
For each fixed $t\in [t_0, t_0+T_1]$, the mapping
\begin{math}
  \mathbf{u}(\,\cdot\,,t)\colon \mathfrak{X}^{(r_1)}\to \CC^M
  \colon z = x + \ii y\mapsto \mathbf{u}(x + \ii y,t)
\end{math}
is holomorphic, i.e., (partially) complex differentiable with respect to
all the complex variables $z_i = x_i + \ii y_i$
{\rm ($i=1,2,\dots,N$)} in
$z = (z_1,\dots,z_N)\in \mathfrak{X}^{(r_1)}\subset \CC^N$.
\end{itemize}
\end{proposition}
\par\vskip 10pt

As for Part~{\rm (d)} in this proposition, there are
several {\em equivalent\/} definitions of
a {\em holomorphic function\/} of several complex variables
used in the literature, cf.\
{\sc S.~G.\ Krantz} \cite[Definitions I--IV, pp.\ 3--4]{Krantz}.
We adopt the most widely used definition in
\cite{Krantz}, Definition~II (p.~3) and Definition 1.2.1 (p.~24).
From this definition it is easy to derive the existence of
an absolutely convergent power series
(\cite[Definition~III, p.~3]{Krantz})
and the Cauchy formula in a polydisc
(\cite[Definition~IV, pp.\ 3--4]{Krantz}).
Nevertheless, the equivalence of \cite[Definition~I, p.~3]{Krantz}
verified in Part~{\rm (d)} in our proposition above
with \cite[Definitions II--IV, pp.\ 3--4]{Krantz}
is a deep classical result due to {\sc F.\ Hartogs}; see
\cite[Theorem 1.2.5, p.~25]{Krantz}.
However, taking also Part~{\rm (c)} into account,
we observe that also \cite[Definition~II, p.~3]{Krantz}
is verified in our proposition.

\par\vskip 10pt
{\it Proof of\/} Proposition~\ref{prop-smooth_ext}.
Recalling our remarks on the local (in time) existence and uniqueness
before this proposition, we observe that it suffices to verify
only our claims in {\rm Parts} {\rm (a)} through~{\rm (d)}.

Proof of~{\rm (a)}.$\;$
Clearly, given any fixed
\begin{math}
  z_0 = x_0 + \ii y_0\in \overline{Q}^{(r_1)}_{\CC} ,
\end{math}
both functions
\begin{equation*}
  t\mapsto \mathbf{u}^{(z_0)}(\,\cdot\,,t)
  \quad\mbox{ and }\quad
  t\mapsto \mathbf{u}^{(\ii y_0)}(\,\cdot\, + x_0,t)
  \colon [t_0, t_0+T_1]\to \mathbf{B}^{s;p,p}(\RR^N)
\end{equation*}
are weak solutions to our Cauchy problem \eqref{e:Cauchy+z_0}
on the (sufficiently short) time interval
$[t_0, t_0+T_1]$ $\subset [0,T]$ with the same initial data
\begin{math}
    \mathbf{u}_0^{(z_0)}(\,\cdot\,)
  = \mathbf{u}_0^{(\ii y_0)}(\,\cdot\, + x_0)
\end{math}
at time $t = t_0\in [0,T)$, for some $T_1\in (0,T-t_0]$.
The uniqueness for problem \eqref{e:Cauchy+z_0} now forces
\begin{math}
         \mathbf{u}^{(z_0)}(\,\cdot\,,t)
  \equiv \mathbf{u}^{(\ii y_0)}(\,\cdot\, + x_0,t)
\end{math}
for every $t\in [t_0, t_0+T_1]$ as claimed.

Part~{\rm (b)} is an immediate consequence of Part~{\rm (a)}
applied to eq.~\eqref{def:u-ext}.

Proof of~{\rm (c)}.$\;$
At the initial time $t = t_0$,
the continuous (partial) differentiability is valid by our hypotheses
on the initial data
$\mathbf{u}_0\colon \mathfrak{X}^{(r)}\to \CC^N$
viewed as a function
\begin{equation*}
  z_0 = x_0 + \ii y_0 = (x_0,y_0)\mapsto
    \mathbf{u}_0(\,\cdot\, + z_0) = \mathbf{u}_0^{(z_0)}
  \colon \mathfrak{X}^{(r)} = \RR^N\times Q^{(r)}
  \to \mathbf{B}^{s;p,p}(\RR^N)
\end{equation*}
valued in the Besov space
$E_{ 1 - \frac{1}{p} ,\, p } = \mathbf{B}^{s;p,p}(\RR^N)$;
in particular,
\begin{math}
    \mathbf{u}_0(\,\cdot\, + z_0) = \mathbf{u}_0^{(z_0)}
  \in U_0\subset E_{ 1 - \frac{1}{p} ,\, p }
\end{math}
provided
\begin{math}
  z_0 = x_0 + \ii y_0\in \overline{Q}^{(r_1)}_{\CC}
  \subset \overline{\mathfrak{X}}^{(r_1)}
\end{math}
($\subset \mathfrak{X}^{(r)}$).
We recall that $U_0$ is an open subset of
$E_{ 1 - \frac{1}{p} ,\, p }$ defined in eq.~\eqref{def:U_0}.
We may view this $C^1$\--differentiability as
(partial) differentiability with respect to
the real parameters $x_{0,i}$ and $y_{0,i}$ in the complex shift
\begin{math}
  z_0 = x_0 + \ii y_0\in \overline{Q}^{(r_1)}_{\CC} ,
\end{math}
where
$x_0 = (x_{0,1}, \dots, x_{0,N})$ and
$y_0 = (y_{0,1}, \dots, y_{0,N})$ are in $\RR^N$ with
\begin{math}
  \max\{ |x_0|_{\infty} ,\, |y_0|_{\infty} \}\leq r_1 .
\end{math}

We now briefly {\bf interrupt\/} our {\it proof of\/}
Proposition~\ref{prop-smooth_ext} to make the following remarks:
\par\vskip 10pt

\par\noindent
{\bf Remarks\/}$\;$
The kind of theory on continuous and differentiable dependence of
the solution
\begin{equation*}
  z_0\mapsto \mathbf{u}(\,\cdot\, + z_0,t)
           = \mathbf{u}^{(z_0)}(\,\cdot\,,t)
  \colon \overline{\mathfrak{X}}^{(r_1)}
  = \RR^N\times \overline{Q}^{(r_1)}
  \to \mathbf{B}^{s;p,p}(\RR^N) \,,\quad t\in [t_0, t_0+T_1] \,,
\end{equation*}
on the real parameters $x_{0,i}$ and $y_{0,i}$
in $z_0\in \overline{\mathfrak{X}}^{(r_1)}$,
that has been developed in
{\sc D.\ Henry} \cite[Chapt.~3]{Henry}, {\S}3.4, pp.\ 62--70,
or, alternatively, in
{\sc A.\ Lunardi} \cite[Chapt.~8]{Lunardi}, {\S}8.3.1, pp.\ 302--306,
can be adapted also to our setting for
the spatially ``shifted'' Cauchy problem \eqref{e:Cauchy+z_0},
with only minor changes.
We should remark that, in this approach,
the following hypotheses on $A$ and $f$ will do;
they follow from {\rm Hypotheses}
\eqref{hy:UnifEll}, \eqref{hy:Ellipt}, and\/ \eqref{hy:f_Hardy^2}
(cf.\ Lemma~\ref{lem-f_Lipschitz} and its proof):

\begin{hypos}\nopagebreak
\begingroup\rm
In analogy with
{\rm Hypotheses} \eqref{hyp:A_anal} and \eqref{hyp:f_anal}
in Section~\ref{s:Analyt_abstr}, {\S}\ref{ss:Analytic_t}, above,
let us assume that there are positive constants
$\vartheta_0\in (0, \pi / 2)$ and $T_0\in (0,T]$, and open sets\,
$\mathcal{U}\subset \CC$ and
$\tilde{U}\subset E_{ 1 - \frac{1}{p} ,\, p }$
containing the compact set
$\overline{\Delta}_{\vartheta_0}^{T_0,T}$ and the open set $U$,
respectively, i.e.,
\begin{math}
  \overline{\Delta}_{\vartheta_0}^{T_0,T}
  \subset \mathcal{U}\subset \CC
\end{math}
and
\begin{math}
  U\subset \tilde{U}\subset E_{ 1 - \frac{1}{p} ,\, p } ,
\end{math}
such that
\begin{enumerate}
\setcounter{enumi}{0}
\renewcommand{\labelenumi}{{\bf (C\arabic{enumi}'')}}
\item
\makeatletter
\def\@currentlabel{{\bf C\arabic{enumi}''}}\label{hyp:A_diff}
\makeatother
$A\colon [0,T]\times U\to \mathcal{L}(E_1\to E_0)$
possesses a continuously (Fr\'echet-) differentiable extension
(i.e., of class $C^1$)
\begin{math}
  \tilde{A}\colon \mathcal{U}\times \tilde{U}
  \to \mathcal{L}(E_1\to E_0)
\end{math}
to the complex domain $\mathcal{U}\times \tilde{U}$ which satisfies
$\tilde{A}(t,v)\in \mathrm{MR}_p(E)$ for all
$(t,v)\in \mathcal{U}\times \tilde{U}$.
\item
\makeatletter
\def\@currentlabel{{\bf C\arabic{enumi}''}}\label{hyp:f_diff}
\makeatother
$f\colon [0,T]\times U\to E_0$
possesses a continuously (Fr\'echet-) differentiable extension
\begin{math}
  \tilde{f}\colon \mathcal{U}\times \tilde{U}\to E_0
\end{math}
to the complex domain $\mathcal{U}\times \tilde{U}$.
\hfill\Square
\end{enumerate}
\endgroup
\end{hypos}
\par\vskip 10pt

Clearly, in both these hypotheses, the mappings $A$ and~$f$,
respectively, are extended from the domain
\begin{math}
  [0,T]\times U\subset
  \overline{\Delta}_{\vartheta_0}^{T_0,T}\times U
\end{math}
to the complex domain
\begin{math}
  \mathcal{U}\times \tilde{U}\subset
  \CC\times E_{ 1 - \frac{1}{p} ,\, p } \,.
\end{math}

\par\vskip 10pt
{\it Proof of\/} Proposition~\ref{prop-smooth_ext}
({\bf continued}), {\rm Part~(c)}.\hfil\break
Recall that the metric on $\mathcal{U}\times \tilde{U}$
(${}\subset \CC\times E_{ 1 - \frac{1}{p} ,\, p }$)
is induced by
the canonical norm on $\CC\times E_{ 1 - \frac{1}{p} ,\, p }$.
It is evident that
{\rm Hypotheses} \eqref{hyp:A_anal} and \eqref{hyp:f_anal} imply
\eqref{hyp:A_diff} and \eqref{hyp:f_diff}, respectively.

Applying the results from
\cite[Chapt.~3, {\S}3.4]{Henry} or
\cite[Chapt.~8, {\S}8.3.1]{Lunardi},
we now conclude that the mapping
\begin{equation*}
  z_0\mapsto \mathbf{u}(\,\cdot\, + z_0,t)
  \colon \overline{\mathfrak{X}}^{(r_1)}
  = \RR^N\times \overline{Q}^{(r_1)}
  \to \mathbf{B}^{s;p,p}(\RR^N) \,,\quad t\in [t_0, t_0+T_1] \,,
\end{equation*}
is continuously differentiable with respect to
the real parameters $x_{0,i}$ and $y_{0,i}$
in $z_0\in \overline{\mathfrak{X}}^{(r_1)}$.
The partial derivatives,
\begin{equation*}
  \frac{\partial\mathbf{u}}{\partial x_{0,i}} \equiv
  \frac{\partial\mathbf{u}}{\partial x_i} =
  \frac{\partial\mathbf{u}}{\partial z_i}
  \;\mbox{ and }\;
  \frac{\partial\mathbf{u}}{\partial y_{0,i}} \equiv
  \frac{\partial\mathbf{u}}{\partial y_i} = \ii\cdot
  \frac{\partial\mathbf{u}}{\partial z_i}
  \colon \overline{\mathfrak{X}}^{(r_1)} \times [t_0, t_0+T_1]\to \CC^M
\end{equation*}
valued in
\begin{math}
  C\left( [t_0, t_0+T_1]\to \mathbf{B}^{s;p,p}(\RR^N)\right) ,
\end{math}
are the unique weak solutions of the following Cauchy problems
derived from \eqref{e:Cauchy+z_0}
by the corresponding partial differentiation, respectively:
\begin{equation}
\label{e:Cauchy_dx}
\left\{
\begin{alignedat}{2}
&   \frac{\partial}{\partial t}
    \left( \frac{\partial\mathbf{u}}{\partial x_i} \right)
\color{blue}
  + \frac{\partial\mathbf{P}}{\partial x_i}
    \Bigl( x + z_0, t, \frac{1}{\ii} \frac{\partial}{\partial x} \Bigr)
    \mathbf{u}(x + z_0,t)
\color{black}
\\
& {}
  + \mathbf{P}
    \Bigl( x + z_0, t, \frac{1}{\ii} \frac{\partial}{\partial x} \Bigr)
    \left( \frac{\partial\mathbf{u}}{\partial x_i} \right) (x + z_0,t)
\\
& {}
  =
\color{blue}
    \frac{\partial\mathbf{f}}{\partial x_i}
    \left( x + z_0, t;
    \left( \frac{ \partial^{|\beta|}\mathbf{u} }{\partial x^{\beta}}
           (x + z_0, t)
    \right)_{|\beta|\leq m}
    \right)
\color{black}
\\
& {}
  + \sum_{|\beta|\leq m} \sum_{k=1}^M
    \frac{\partial\mathbf{f}}{\partial Z_{\beta,k}}
    \left( x + z_0, t;
    \left( \frac{ \partial^{|\beta|}\mathbf{u} }{\partial x^{\beta}}
    \right)_{|\beta|\leq m}
    \right)
    \frac{ \partial^{|\beta|} }{\partial x^{\beta}}
    \left( \frac{\partial u_k}{\partial x_i} \right) (x + z_0,t)
\\
&   \quad\mbox{ for } (x,t)\in \mathbb{R}^N\times (t_0, t_0+T_1) \,;
\\
&   \frac{\partial\mathbf{u}  }{\partial x_i} (x + z_0,0)
  = \frac{\partial\mathbf{u}_0}{\partial x_i} (x + z_0)
    \quad\mbox{ for } x\in \mathbb{R}^N \,,
\end{alignedat}
\right.
\end{equation}
and
\begin{equation}
\label{e:Cauchy_dy}
\left\{
\begin{alignedat}{2}
&   \frac{\partial}{\partial t}
    \left( \frac{\partial\mathbf{u}}{\partial y_i} \right)
\color{blue}
  + \frac{\partial\mathbf{P}}{\partial y_i}
    \Bigl( x + z_0, t, \frac{1}{\ii} \frac{\partial}{\partial x} \Bigr)
    \mathbf{u}(x + z_0,t)
\color{black}
\\
& {}
  + \mathbf{P}
    \Bigl( x + z_0, t, \frac{1}{\ii} \frac{\partial}{\partial x} \Bigr)
    \left( \frac{\partial\mathbf{u}}{\partial y_i} \right) (x + z_0,t)
\\
& {}
  =
\color{blue}
    \frac{\partial\mathbf{f}}{\partial y_i}
    \left( x + z_0, t;
    \left( \frac{ \partial^{|\beta|}\mathbf{u} }{\partial x^{\beta}}
           (x + z_0, t)
    \right)_{|\beta|\leq m}
    \right)
\color{black}
\\
& {}
  + \sum_{|\beta|\leq m} \sum_{k=1}^M
    \frac{\partial\mathbf{f}}{\partial Z_{\beta,k}}
    \left( x + z_0, t;
    \left( \frac{ \partial^{|\beta|}\mathbf{u} }{\partial x^{\beta}}
    \right)_{|\beta|\leq m}
    \right)
    \frac{ \partial^{|\beta|} }{\partial x^{\beta}}
    \left( \frac{\partial u_k}{\partial y_i} \right) (x + z_0,t)
\\
&   \quad\mbox{ for } (x,t)\in \mathbb{R}^N\times (t_0, t_0+T_1) \,;
\\
&   \frac{\partial\mathbf{u}  }{\partial y_i} (x + z_0,0)
  = \frac{\partial\mathbf{u}_0}{\partial y_i} (x + z_0)
    \quad\mbox{ for } x\in \mathbb{R}^N \,,
\end{alignedat}
\right.
\end{equation}
where the complex variable $Z_{\beta,k}$ stands for
\begin{math}
  Z_{\beta,k} =
  \frac{ \partial^{|\beta|} }{\partial x^{\beta}} u_k =
  \ii^{|\beta|}\, D_x^{\beta} u_k\in \CC .
\end{math}
This proves Part~{\rm (c)}.

Proof of~{\rm (d)}.$\;$
We take advantage of the last two equations,
\eqref{e:Cauchy_dx} and \eqref{e:Cauchy_dy},
to apply the Cauchy\--Riemann operator
${\partial} / {\partial\bar{z}_i}$
from \eqref{def:Cauchy-Riem}
to problem \eqref{e:Cauchy+z_0} in order to conclude that
the {\em Cauchy\--Riemann derivative\/}
\begin{equation*}
  \bar{\partial}_{z_{0,i}} \mathbf{u}\equiv
  \frac{\partial\mathbf{u}}{\partial\bar{z}_{0,i}} \equiv
  \frac{\partial\mathbf{u}}{\partial\bar{z}_i} =
  \frac{1}{2}
    \left( \frac{\partial}{\partial x_i}
    + \ii\,\frac{\partial}{\partial y_i}
    \right) \mathbf{u}
  \colon \overline{\mathfrak{X}}^{(r_1)} \times [t_0, t_0+T_1]\to \CC^M
\end{equation*}
is in
\begin{math}
  C\left( [t_0, t_0+T_1]\to \mathbf{B}^{s;p,p}(\RR^N)\right)
\end{math}
and obeys the following homogeneous linear Cauchy problem,
which is a simple linear combination
%
``
\begin{math}
    \frac{1}{2}\cdot   \mbox{$\eqref{e:Cauchy_dx}$}
  + \frac{\ii}{2}\cdot \mbox{$\eqref{e:Cauchy_dy}$}
  = \mbox{$\eqref{e:Cauchy_dz}$}
\end{math}
'':
\begin{equation}
\label{e:Cauchy_dz}
\left\{
\begin{alignedat}{2}
&   \frac{\partial}{\partial t}
    \left( \bar{\partial}_{z_{0,i}} \mathbf{u}\right)
  + \mathbf{P}
    \Bigl( x + z_0, t, \frac{1}{\ii} \frac{\partial}{\partial x} \Bigr)
    \left( \bar{\partial}_{z_{0,i}} \mathbf{u}\right) (x + z_0,t)
\\
& {}
  = \sum_{|\beta|\leq m} \sum_{k=1}^M
    \frac{\partial\mathbf{f}}{\partial Z_{\beta,k}}
    \left( x + z_0, t;
    \left( \frac{ \partial^{|\beta|}\mathbf{u} }{\partial x^{\beta}}
    \right)_{|\beta|\leq m}
    \right)
    \frac{ \partial^{|\beta|} }{\partial x^{\beta}}
    \left( \bar{\partial}_{z_{0,i}} u_k\right) (x + z_0,t)
\\
&   \quad\mbox{ for } (x,t)\in \mathbb{R}^N\times (t_0, t_0+T_1) \,;
\\
&   \left( \bar{\partial}_{z_{0,i}} \mathbf{u}\right) (x + z_0,0)
    = \mathbf{0} \quad\mbox{ for } x\in \mathbb{R}^N \,.
\end{alignedat}
\right.
\end{equation}
Here, we have used the fact that both operators
\begin{equation*}
  z\mapsto \mathbf{P}
    \Bigl( z, t, \frac{1}{\ii} \frac{\partial}{\partial x} \Bigr)
  \;\mbox{ and }\;
  z\mapsto \mathbf{f}
    \Bigl( z, t; (Z_{\beta})_{|\beta|\leq m} \Bigr)
  \colon \mathfrak{X}^{(r)} \to \CC^M
  \qquad\mbox{ ($r = r_0$) }
\end{equation*}
are holomorphic, i.e.,
\begin{math}
  \bar{\partial}_{z_i} \mathbf{P}
    \Bigl( z, t, \frac{1}{\ii} \frac{\partial}{\partial x} \Bigr)
  = \mathbf{0}
\end{math}
and
\begin{math}
  \bar{\partial}_{z_i} \mathbf{f}
    \left( z, t; (Z_{\beta})_{|\beta|\leq m} \right)
  = \mathbf{0} ,
\end{math}
by our {\rm Hypotheses}
\eqref{hy:UnifEll} and \eqref{hy:f_Hardy^2}, respectively.
By our choice of
\begin{math}
  \mathbf{u}_0\equiv \tilde{\mathbf{u}}_0
  \colon \mathfrak{X}^{(r)} \to \CC^M
\end{math}
being holomorphic, we have also
$\bar{\partial}_{z_i} \mathbf{u}_0(z) = \mathbf{0}$; $i=1,2,\dots,N$.
Notice that eq.~\eqref{e:Cauchy_dz} is valid only for every
$z_0\in \overline{\mathfrak{X}}^{(r_1)}$
($\subset \mathfrak{X}^{(r)}$).

By {\rm Hypotheses} \eqref{hy:UnifEll} and \eqref{hy:Ellipt},
the linear differential operator on
the left\--hand side of eq.\ \eqref{e:Cauchy_dz}, i.e.,
\begin{equation*}
    \frac{\partial}{\partial t}
  + \mathbf{P}
    \Bigl( x + z_0, t, \frac{1}{\ii} \frac{\partial}{\partial x}
    \Bigr) \,,
\end{equation*}
is uniformly parabolic of order $2m$ with smooth coefficients.
It is proved in 
{\sc R.\ Denk}, {\sc M.\ Hieber}, and {\sc J.\ Pr\"uss}
\cite[p.~67]{Denk-Hieber}, Theorem 5.7
(cf.\ also \cite{Pruess}, Theorem 2.1 (p.~8)
 and remarks thereafter (p.~9))
that, for every $z_0\in \mathfrak{X}^{(r)}$
and for every $t\in [0,T]$,
\begin{equation}
\label{def:Cauchy_dz}
  A^{(z_0)}(t)\eqdef {}-
    \mathbf{P}
    \Bigl( x + z_0, t, \frac{1}{\ii} \frac{\partial}{\partial x}
    \Bigr) \colon
  \mathbf{W}^{2m,p}(\RR^N)\to \mathbf{L}^p(\RR^N)
\end{equation}
is a bounded linear operator, i.e.,
$A^{(z_0)}(t)\in \mathcal{L}(E_1\to E_0)$,
and it possesses
the maximal $L^p$-regularity property, i.e.,
$A^{(z_0)}(t)\in \mathrm{MR}_p(E)\equiv \mathrm{MR}_p(E_1\to E_0)$.
Let us recall that
\begin{math}
  E_1 = \mathbf{W}^{2m,p}(\RR^N) \hookrightarrow
  E_0 = \mathbf{L}^p(\RR^N) .
\end{math}

Furthermore, in view of {\rm Hypothesis} \eqref{hy:f_Hardy^2},
the pointwise multiplication and differentiation operators on
the right\--hand side of eq.~\eqref{e:Cauchy_dz}
are of order $|\beta|$ ($|\beta|\leq m < 2m$)
and all have bounded continuous coefficients, by
\begin{math}
  \mathbf{u}(\,\cdot\,,t)\in U\subset \mathbf{B}^{s;p,p}(\RR^N)
\end{math}
for every $t\in [t_0, t_0+T_1]$ combined with the Sobolev imbedding
\begin{math}
  B^{s;p,p}(\RR^N) \hookrightarrow
  C^m(\RR^N)\cap W^{m,\infty}(\RR^N) \,,
\end{math}
where
$2 + \genfrac{}{}{}1{N}{m} < p < \infty$ and
$m < s = 2m\left( 1 - \genfrac{}{}{}1{1}{p} \right) < 2m$.
We denote their sum, which appears in eq.~\eqref{e:Cauchy_dz}, by
\begin{math}
  f^{(z_0)}(t)\colon
  \mathbf{B}^{s;p,p}(\RR^N)\to \mathbf{L}^p(\RR^N) \,,
\end{math}
i.e.,
\begin{math}
  f^{(z_0)}(t)\in
  \mathcal{L}\left( E_{ 1 - \frac{1}{p} ,\, p } \to E_0 \right) .
\end{math}
Here, we allow any
$z_0\in \overline{\mathfrak{X}}^{(r_1)}$.
Consequently, the mappings
\begin{equation*}
  (t,v)\mapsto A^{(z_0)}(t)
       \colon [0,T]\times U\to \mathcal{L}(E_1\to E_0)
  \quad\mbox{ and }\quad
  (t,v)\mapsto f^{(z_0)}(t) v
       \colon [0,T]\times U\to E_0
\end{equation*}
satisfy
{\rm Hypotheses} \eqref{hyp:A} and\/ \eqref{hyp:f}
for $A$ and~$f$, respectively, with
$U = E_{ 1 - \frac{1}{p} ,\, p }$,
$A^{(z_0)}(t)\in\hfil\break \mathcal{L}(E_1\to E_0)$
being independent from $v\in U$, and
$v\mapsto f^{(z_0)}(t) v$ linear in $v\in U$.

We observe that the homogeneous linear Cauchy problem
\eqref{e:Cauchy_dz}
for the Cauchy\--Riemann derivative
$\bar{\partial}_{z_{0,i}} \mathbf{u}$ of $\mathbf{u}$
takes the following abstract linear form, whenever
$z_0\in \overline{\mathfrak{X}}^{(r_1)}$:
\begin{equation}
\label{eq:Cauchy_dz}
\left\{
\begin{alignedat}{2}
  \frac{\mathrm{d}}{\mathrm{d}t}\,
    (\bar{\partial}_{z_{0,i}} \mathbf{u})
  - A^{(z_0)}(t)
    (\bar{\partial}_{z_{0,i}} \mathbf{u})
& {}
  = f^{(z_0)}(t)
    (\bar{\partial}_{z_{0,i}} \mathbf{u})
  \quad\mbox{ for a.e. } t\in (t_0, t_0+T_1) \,;
\\
  (\bar{\partial}_{z_{0,i}} \mathbf{u})(0)
& {}
  = \mathbf{0}\in E_{ 1 - \frac{1}{p} ,\, p } \,.
\end{alignedat}
\right.
\end{equation}
This abstract linear problem corresponds to
the nonlinear initial value problem \eqref{e:abstr_Cauchy}
treated in Section~\ref{s:Cauchy_abstr}.
We apply the uniqueness part of Theorem~\ref{thm-Clem-Li}
to deduce
\begin{math}
  (\bar{\partial}_{z_{0,i}} \mathbf{u}) (x,t)\equiv \mathbf{0}
\end{math}
for all $(x,t)\in \RR^N\times [t_0, t_0+T_1]$.
This implies that the mapping
$z\mapsto u_k(z,t)\colon \mathfrak{X}^{(r_1)}\to \CC$
is holomorphic in each complex variable $z_i\in \CC$,
for every fixed time $t\in [t_0, t_0+T_1]$; $k=1,2,\dots,N$.
Moreover, by Part~{\rm (c)}, all complex partial derivatives
$\partial_{z_i} u_k(\,\cdot\,,t)$
are continuous in $\mathfrak{X}^{(r_1)}$.
Finally, we take advantage of the classical fact that such a function
$u_k(\,\cdot\,,t)\colon \mathfrak{X}^{(r_1)}\to \CC$
is holomorphic
(Remark~\ref{rem-Cauchy-Riem}); see e.g.\
{\sc F.\ John} \cite[Theorem, p.~70]{John}
or
{\sc S.~G.\ Krantz} \cite[Definition~II, p.~3]{Krantz}.
Also Part~{\rm (d)} and, thus, the entire proposition is proved.
\qed
\par\vskip 10pt

\section{Space\--time analyticity for the Cauchy problem in
         $\RR^N\times (0,T)$}
\label{s:Analyt_(x,t)}

We summarize the time and space analyticity results from
the last two sections
(Sections \ref{s:Analyt_abstr} and~\ref{s:Analyt_space}),
for the mapping
\begin{math}
  \mathbf{u}\colon
  \overline{\mathfrak{X}}^{(r_1)} \times [t_0, t_0+T_1]\to \CC^M
\end{math}
defined in {\rm eq.}~\eqref{def:u-ext}, in the following theorem:

\begin{theorem}\label{thm-smooth_ext}
Let\/ $M,N\geq 1$, $0 < T < \infty$, and assume that all\/
{\rm Hypotheses} \eqref{hy:UnifEll}, \eqref{hy:Ellipt}, and\/
\eqref{hy:f_Hardy^2}
are satisfied with some constants\/
$0 < r_0 < \infty$, $0 < T_0\leq T$, and\/ $0 < \vartheta_0 < \pi / 2$.
Furthermore, assume that\/
\begin{math}
  \widehat{\mathbf{u}}\in
  C\left( [0,T]\to \mathbf{B}^{s;p,p}(\RR^N)\right)
\end{math}
is a globally defined weak solution
to the original Cauchy problem \eqref{e:Cauchy}, i.e.,
{\rm Hypothesis} \eqref{hyp:global} is valid.
Let the sets\/
\begin{math}
  U_0\subset U\subset
  E_{ 1 - \frac{1}{p} ,\, p } = \mathbf{B}^{s;p,p}(\RR^N)
\end{math}
be specified as in {\rm eqs.}\ \eqref{def:U_0} and~\eqref{def:U}.
Given any\/ $t_0\in [0,T)$, let\/
$\mathbf{u}_0\in \mathbf{B}^{s;p,p}(\RR^N)$
be any initial data at time $t = t_0$, such that\/
$\mathbf{u}_0\in U_0$ and\/
$\mathbf{u}_0$ has a holomorphic extension
\begin{math}
  \tilde{\mathbf{u}}_0\colon \mathfrak{X}^{(r)}\to \CC^M
\end{math}
as described before {\rm Lemma~\ref{lem-f_Lipschitz}}
\begingroup\rm
(we identify $\mathbf{u}_0\equiv \tilde{\mathbf{u}}_0$),
\endgroup
with
\begin{math}
  \mathbf{u}_0^{(z_0)} = \mathbf{u}_0(\,\cdot\, + z_0)
  \in U_0\subset E_{ 1 - \frac{1}{p} ,\, p }
\end{math}
whenever
\begin{math}
  z_0\in \overline{Q}^{(r_1)}_{\CC}
  \subset\overline{\mathfrak{X}}^{(r_1)}
\end{math}
for some $r_1\in (0,r)$, cf.\ {\rm Hypothesis} \eqref{hyp:analyt}.
(We have set $r = r_0$; $r_1$ may depend on $\mathbf{u}_0$.) 

Finally, let\/
\begin{math}
  \mathbf{u}\colon
  \overline{\mathfrak{X}}^{(r_1)} \times [t_0, t_0+T_1]\to \CC^M
\end{math}
be the continuous mapping obtained in\/
{\rm Proposition~\ref{prop-smooth_ext}} above, with
the number $T_1\in (0,T-t_0]$
depending on $r_1$ and $U$, but not on $t_0$.
Replace $T_0\in (0,T]$ by $\,\min\{ T_0, T_1\}\,$ if necessary, so that\/
$0 < T_0\leq T_1\leq T$ holds.
Then there exist constants\/
$\vartheta'\in (0, \vartheta_0]$ and\/
$T'\in (0,T_0]$, small enough, and a continuous mapping\/
\begin{math}
  \tilde{\mathbf{u}}\colon
  \overline{\mathfrak{X}}^{(r_1)} \times
  \left( t_0 + \Delta_{\vartheta'}^{T',T_1} \right)
  \to \CC^M
\end{math}
with the following properties:
\begin{itemize}
\item[{\rm (i)}]
For each $z_0\in \overline{\mathfrak{X}}^{(r_1)}$,
the (unique) weak solution
\begin{math}
  \mathbf{u}^{(z_0)}\in
  C\left( [t_0, t_0+T_1]\to \mathbf{B}^{s;p,p}(\RR^N)\right)
\end{math}
to the Cauchy problem \eqref{e:Cauchy+z_0} on the time interval\/
$[t_0, t_0+T_1]\subset [0,T]$ with the initial condition
$\mathbf{u}(\,\cdot\,,t_0) = \mathbf{u}_0^{(z_0)}$
at time $t = t_0$ satisfying\/
$z_0\in \overline{Q}^{(r_1)}_{\CC}$
possesses a unique holomorphic extension from
$(t_0, t_0+T_1)$ to
$t_0 + \Delta_{\vartheta'}^{T',T_1}$, such that
\begin{math}
    \mathbf{u}^{(z_0)}(\,\cdot\,, t_0 + s)
  = \tilde{\mathbf{u}} (\,\cdot\, + z_0, t_0 + s) \in U
\end{math}
holds for every\/
$s\in \Delta_{\vartheta'}^{T', T_1}$.
\item[{\rm (ii)}]
The complex function
\begin{math}
  \tilde{\mathbf{u}}\colon
  \mathfrak{X}^{(r_1)} \times
  \left( t_0 + \Delta_{\vartheta'}^{T',T_1} \right)
  \to \CC^M
\end{math}
is holomorphic (jointly) in all its variables,
$z = (z_1,z_2,\dots,z_N)\in \mathfrak{X}^{(r_1)}\subset \CC^N$
and\/
$t\in t_0 + \Delta_{\vartheta'}^{T',T_1} \subset \CC$.
\end{itemize}
\end{theorem}
\par\vskip 10pt

Let us recall that, by our notation in eq.~\eqref{glob:T^(r)},
for $\vartheta\in (0, \pi / 2)$,
$0 < t_0 < T < \infty$, and $0 < T'\leq T - t_0$, we have
\begin{equation}
\label{glob:t_0+T^(r)}
\begin{aligned}
  t_0 + \Delta_{\vartheta}^{T', T - t_0}
& = \left( t_0 + \Delta_{\vartheta}^{(T-t_0)} \right)
    \cap \{ t\in \CC\colon |\IM t| < T'\cdot \tan\vartheta \}
\\
& = \bigcup_{ t_0\leq \xi\leq T-T'}
    \{ \xi + t'\in \CC\colon t'\in \Delta_{\vartheta}^{(T')} \}
  = \bigcup_{ t_0\leq \xi\leq T-T'}
    \left( \xi + \Delta_{\vartheta}^{(T')} \right)
\end{aligned}
\end{equation}
with the closure
$t_0 + \overline{\Delta}_{\vartheta}^{T', T - t_0}$
in $\CC$.

\par\vskip 10pt
\begin{remark}\label{rem-smooth_ext}\nopagebreak
\begingroup\rm

{\rm (a)}$\;$
We remark that the main difference between our main result,
Theorem~\ref{thm-Main} (Section~\ref{s:Main}),
and Theorem~\ref{thm-smooth_ext} above is
the {\em temporally local\/} character of the latter stated for
the time interval 
$[t_0, t_0+T_1]\subset [0,T]$ with the additional
{\em analyticity\/} hypothesis on the initial data $\mathbf{u}_0$
(as in Part~{\bf (iii)} of Theorem~\ref{thm-Main}).

{\rm (b)}$\;$
Recalling our choice of the number $r\in (0,\infty)$
in {\rm Hypothesis} \eqref{hyp:analyt}
(before Remark~\ref{rem-hyp:analyt})
on the complex analy\-ti\-city of the initial data
\begin{math}
  \mathbf{u}_0 = \tilde{\mathbf{u}}_0
  \colon \mathfrak{X}^{(r)}\to \CC^M
\end{math}
extended to the complex strip
$\mathfrak{X}^{(r)}\subset \CC^N$,
we observe that the number $r_1\in (0,r)$, originally introduced
in the spatially ``shifted'' Cauchy problem \eqref{e:Cauchy+z_0},
is needed for sufficiently small perturbations (``shifts'')
$z_0\in \CC^N$ of the space variable
$z\in \mathfrak{X}^{(r_1)}$ in order to keep
$z + z_0\in \mathfrak{X}^{(r)}$.
As we have already mentioned after Remark~\ref{rem-hyp:analyt},
{\rm Hypotheses}
\eqref{hy:UnifEll}, \eqref{hy:Ellipt}, and\/ \eqref{hy:f_Hardy^2}
in Section~\ref{s:Main}
show that only the case $0 < r\leq r_0$ is useful.
We now recall from
Proposition~\ref{prop-smooth_ext} and Theorem~\ref{thm-smooth_ext}
that, in order to avoid excessive notation,
$r_1\in (0,r)$ must be chosen small enough, such that
\begin{math}
  \mathbf{u}_0^{(z_0)} = \mathbf{u}_0(\,\cdot\, + z_0)
  \in U_0\subset U\subset E_{ 1 - \frac{1}{p} ,\, p }
\end{math}
for every
\begin{math}
  z_0\in \overline{Q}^{(r_1)}_{\CC}
  \subset\overline{\mathfrak{X}}^{(r_1)} .
\end{math}
We recall that the sets $U_0$ and $U$ are defined in {\rm eqs.}\
\eqref{def:U_0} and \eqref{def:U}, respectively.
We stress that both, $U_0$ and $U$, are open in
$E_{ 1 - \frac{1}{p} ,\, p }$,
while being determined solely by the restriction to the real line $\RR$,
\begin{math}
  \mathbf{u}_0 = \tilde{\mathbf{u}}_0 \vert_{ \mathbb{R} }
  \in U_0\subset E_{ 1 - \frac{1}{p} ,\, p } \,,
\end{math}
of the initial data
\begin{math}
  \tilde{\mathbf{u}}_0\colon \mathfrak{X}^{(r)}\to \CC^M ,
\end{math}
i.e., by
$\mathbf{u}_0\colon \RR^N\to \CC^M$
as an element of the Besov space
$E_{ 1 - \frac{1}{p} ,\, p } = \mathbf{B}^{s;p,p}(\RR^N)$.
Consequently, the number
$r_1\equiv r_1(\mathbf{u}_0)\in (0,\infty)$
is determined by these initial data $\mathbf{u}_0\in U_0$; we have
$0 < r_1 < r\leq r_0$ where we may choose $r = r_0$, by
Remark~\ref{rem-hyp:analyt}.
Such a choice of $r_1\in (0,r)$ is possible thanks to the closed cube
$\overline{Q}^{(r_1)}$ being compact in $\RR^N$.
\hfill\Square
\endgroup
\end{remark}

\par\vskip 10pt
{\it Proof of\/} Theorem~\ref{thm-smooth_ext}.

Proof of~{\rm (i)}.$\;$
Let
$z_0 = x_0 + \ii y_0\in \overline{\mathfrak{X}}^{(r_1)}$
be arbitrary, but fixed, with
\begin{math}
  \max\{ |x_0|_{\infty} ,\, |y_0|_{\infty} \}\leq r_1 ,
\end{math}
i.e., $z_0\in \overline{Q}^{(r_1)}_{\CC}$.
We apply our time analyticity result in
Theorem~\ref{thm-Analytic_glob}
to the Cauchy problem \eqref{e:Cauchy+z_0} on the time interval\/
$[t_0, t_0+T_1]\subset [0,T]$ with the initial condition
$\mathbf{u}(\,\cdot\,,t_0) = \mathbf{u}_0^{(z_0)} \in U_0$
at time $t = t_0$ in order to derive the conclusion of
Part~{\rm (i)}.

Proof of~{\rm (ii)}.$\;$
The second part (Part~{\rm (ii)})
is obtained by combining Part~{\rm (i)} with
Proposition~\ref{prop-smooth_ext}, particularly Part~{\rm (d)}.
Finally, the joint time and space analyticity of the complex function
\begin{math}
  \tilde{\mathbf{u}}\colon
  \mathfrak{X}^{(r_1)} \times
  \left( t_0 + \Delta_{\vartheta'}^{T',T_1} \right)
  \to \CC^M
\end{math}
is obtained by applying the classical characterization of
holomorphic functions by the Cauchy\--Riemann equations
(Remark~\ref{rem-Cauchy-Riem}); see e.g.\
{\sc F.\ John} \cite[Theorem, p.~70]{John}
or
{\sc S.~G.\ Krantz} \cite[Definition~II, p.~3]{Krantz}.
\qed
\par\vskip 10pt

\section{Proofs of the main results, Theorem~\ref{thm-Main}
         and Proposition~\ref{prop-Main}}
\label{s:proof-Main}

Now we are ready to prove our main result, Theorem~\ref{thm-Main}.
Proposition~\ref{prop-Main} is then a consequence of
Theorem~\ref{thm-Main} and
Lemma~\ref{lem-perturb_MR:A<B}, ineq.~\eqref{est:B+A_sol};
we give its proof right after that of Theorem~\ref{thm-Main}.

\par\vskip 10pt
{\it Proof of\/} Theorem~\ref{thm-Main}.
Proof of~{\bf (i)}.$\;$
The existence and uniqueness of a weak solution
\begin{math}
  \mathbf{u}\in
  C\left( [0,T_1]\to \mathbf{B}^{s;p,p}(\RR^N) \right)
\end{math}
to the Cauchy problem \eqref{e:Cauchy},
local in time for $t\in [0,T_1]$ with some $T_1\in (0,T]$,
is obtained directly from the abstract result in
Theorem~\ref{thm-Clem-Li} where
\begin{math}
    E_{ 1 - \frac{1}{p} ,\, p }
  = \mathbf{B}^{s;p,p}(\RR^N) = [ B^{s;p,p}(\RR^N) ]^M .
\end{math}
The technical details in applying Theorem~\ref{thm-Clem-Li}
(an abstract result)
to problem~\eqref{e:Cauchy} have been given in
Section~\ref{s:Analyt_space}, right after problem~\eqref{e:Cauchy+z_0}.
The linear parabolic operator on the left\--hand side
in eq.~\eqref{e:Cauchy+z_0} is treated by the maximal $L^p$-regularity
described in Remark~\ref{rem-perturb}{\rm (a)}.
The special case of $p$ in this remark, $p_0 = 2$, is taken care of
by standard parabolic regularity making use of
{\rm G\r{a}rding's inequality} in Corollary~\ref{cor-hypos};
see, e.g., {\sc A.\ Friedman} \cite[Chapt.~10]{Friedman-64}.
If a weak solution
\begin{math}
  \mathbf{u}\in C\left( [0,T]\to \mathbf{B}^{s;p,p}(\RR^N) \right)
\end{math}
exists globally in time $t\in [0,T]$, then it is unique, by
Theorem~\ref{thm-Clem-Li_glob}, Part~{\rm (i)}.

Proof of~{\bf (ii)}.$\;$
The (unique) temporal extension of the function
$\mathbf{u}\colon \RR^N\times (0,T)\to \CC^M$
to a holomorphic function
\begin{math}
  \mathbf{u}^{\sharp}\colon
  \Delta_{\vartheta'}^{T',T}\to \mathbf{B}^{s;p,p}(\RR^N)
\end{math}
that possesses another extension to a continuous function on the closure
\begin{math}
  \overline{\Delta}_{\vartheta'}^{T',T} ,
\end{math}
denoted again by $\mathbf{u}^{\sharp}$, is derived from
Theorem~\ref{thm-smooth_ext}, Part~{\rm (i)}.
More precisely, Part~{\rm (i)} of Theorem~\ref{thm-smooth_ext}
is applied to the (global) weak solution
\begin{math}
  \mathbf{u}\in C\left( [0,T]\to \mathbf{B}^{s;p,p}(\RR^N)\right)
\end{math}
of the Cauchy problem \eqref{e:Cauchy}, which is assumed to exist,
in the temporal complex domain $t_0 + \Delta_{\vartheta'}^{(T')}$
with the initial value
$\mathbf{u}(\,\cdot\,,t_0)\in \mathbf{B}^{s;p,p}(\RR^N)$
at every initial time $t_0\in [0,T-T']$.
Here, we have used that
\begin{equation}
\label{glob:0+T^(r)}
  \Delta_{\vartheta'}^{T',T}
  = \bigcup_{ 0\leq t_0\leq T-T'}
    \left( t_0 + \Delta_{\vartheta'}^{(T')} \right) \,;
  \qquad\mbox{ cf.\ eq.~\eqref{glob:t_0+T^(r)}. }
\end{equation}
%

Proof of~{\bf (iii)}.$\;$
We remark that 
{\rm Hypothesis} \eqref{hyp:global} is satisfied with the function
\begin{math}
  \widehat{\mathbf{u}} = \mathbf{u}\in
  C\left( [0,T]\to \mathbf{B}^{s;p,p}(\RR^N)\right) ,
\end{math}
a globally defined weak solution
to the original Cauchy problem \eqref{e:Cauchy},
which exists by our hypothesis.
Let us recall the definitions of the bounded, open, and convex sets
$U_0$ and $U$ in $\mathbf{B}^{s;p,p}(\RR^N)$,
\begin{math}
  U_0\subset U = \tilde{U}\subset
  E_{ 1 - \frac{1}{p} ,\, p } = \mathbf{B}^{s;p,p}(\RR^N) ,
\end{math}
in eqs.\ \eqref{def:U_0} and~\eqref{def:U}, respectively,
where the radius $R_0\in (0,\infty)$ is an arbitrary positive number.
We recall also our hypothesis that the initial condition
\begin{math}
  \mathbf{u}(\,\cdot\,,0) = \mathbf{u}_0\in \mathbf{B}^{s;p,p}(\RR^N)
\end{math}
possesses a (unique) holomorphic extension
\begin{math}
  \tilde{\mathbf{u}}_0\colon \mathfrak{X}^{(\kappa_0)}\to \CC^M
\end{math}
from $\RR^N$ to the complex domain
\begin{math}
  \mathfrak{X}^{(\kappa_0)} = \RR^N + \ii Q^{(\kappa_0)}\subset \CC^N
\end{math}
(a tube), for some $\kappa_0\in (0,r_0]$,
that satisfies ineq.~\eqref{sup_y:u_0}.

We begin with a construction of the (unique) spatial extension of
the continuous function
$\mathbf{u}\colon \RR^N\times [0,T]\to \CC^M$
to a continuous function
\begin{math}
  \mathbf{u}^{\flat}\colon
  \overline{\mathfrak{X}}^{(\rho)}\times [0,T]\to \CC^M
\end{math}
that is holomorphic in the space variable
\begin{math}
  z = x + \ii y\in \mathfrak{X}^{(\rho)} = \RR^N + \ii Q^{(\rho)}
\end{math}
with some $\rho\in (0,\kappa_0]$ small enough.
Let us recall our notation with the ``shifted'' function
\begin{math}
  x\mapsto \mathbf{u}_0^{(z_0)}(x)\eqdef \mathbf{u}_0(x + z_0)
  \colon \RR^N\to \CC^M
\end{math}
introduced in the Cauchy problem \eqref{e:Cauchy+z_0}
spatially ``shifted'' by
$z_0 = x_0 + \ii y_0\in \mathfrak{X}^{(r)}\subset \CC^N$.
The constant $r\in (0,\infty)$ has been introduced in
{\rm Hypothesis} \eqref{hyp:analyt}; only the case
$0 < r < \kappa_0\leq r_0$ (${}< \infty$) is useful.
We wish to apply Proposition~\ref{prop-smooth_ext}
with the constant $r_1\in (0,r)$ specified there.
We choose $\rho\in (0,r_1)$ small enough, such that also
\begin{equation*}
  \left\Vert \mathbf{u}_0^{(\ii y)} - \mathbf{u}_0
  \right\Vert_{ B^{s;p,p}(\RR^N) }
  = \left\Vert \mathbf{u}_0(\,\cdot\, + \ii y) - \mathbf{u}_0
  \right\Vert_{ B^{s;p,p}(\RR^N) }
  < R_0 \quad\mbox{ holds for all }\, y\in \overline{Q}^{(\rho)} \,.
\end{equation*}
Here, we have used the (Lipschitz) continuity of the mapping
\begin{math}
  y\mapsto \mathbf{u}_0^{(\ii y)}\eqdef \mathbf{u}_0(\,\cdot\, + \ii y)
  \colon \overline{Q}^{(r_1)}\to \mathbf{B}^{s;p,p}(\RR^N) ,
\end{math}
thanks to $0 < r_1 < \kappa_0$ supplemented by
the Cauchy formula in a polydisc centered in
$\overline{\mathfrak{X}}^{(r_1)}$
(with radius ${}< \kappa_0 - r_1$)
and contained in the complex strip
$\mathfrak{X}^{(\kappa_0)}\subset \CC^N$.
From eq.~\eqref{def:U_0} we deduce that
$\mathbf{u}_0^{(\ii y)}\in U_0$ for all $y\in \overline{Q}^{(\rho)}$.
In analogy with our proof of Part~{\bf (i)} above, we apply
Theorem~\ref{thm-Clem-Li} to conclude that
the spatially shifted Cauchy problem \eqref{e:Cauchy+z_0},
with the shift $z_0 = \ii y$
($y\in \overline{Q}^{(\rho)}$),
possesses a unique weak solution
\begin{math}
  \mathbf{u}^{(\ii y)}\in
  C\left( [0,T_1]\to \mathbf{B}^{s;p,p}(\RR^N) \right) ,
\end{math}
local in time for $t\in [0,T_1]$ with some $T_1\in (0,T]$,
that satisfies
$\mathbf{u}^{(\ii y)}(\,\cdot\,,t)\in U$ for every $t\in [0,T_1]$.
We apply Part~{\rm (c)} or Part~{\rm (d)}
of Proposition~\ref{prop-smooth_ext} with $t_0 = 0$
to conclude that there is a number $R_1\in (0,R_0)$ small enough,
such that even
$\mathbf{u}^{(\ii y)}(\,\cdot\,,t)\in U_0\subset U$
holds for every $t\in [0,T_1]$, provided
$\rho\in (0,r_1)$ is chosen so small that also
\begin{equation*}
  \left\Vert \mathbf{u}_0^{(\ii y)} - \mathbf{u}_0
  \right\Vert_{ B^{s;p,p}(\RR^N) } < R_1 \quad (< R_0)
    \quad\mbox{ holds for all }\, y\in \overline{Q}^{(\rho)} \,.
\end{equation*}
Here, besides the (Lipschitz) continuity of the mapping
\begin{math}
  y\mapsto \mathbf{u}_0^{(\ii y)}\eqdef \mathbf{u}_0(\,\cdot\, + \ii y)
  \colon \overline{Q}^{(r_1)}\to \mathbf{B}^{s;p,p}(\RR^N)
\end{math}
mentioned above, we have used also
the continuous dependence of the solution
$\mathbf{u}^{(\ii y)}$ upon the initial data
$\mathbf{u}_0^{(\ii y)}\in \mathbf{B}^{s;p,p}(\RR^N)$
obtained in Theorem~\ref{thm-Clem-Li_glob}, Part~{\rm (ii)};
see also Remark~\ref{rem-Clem-Li}.
According to eq.~\eqref{def:u-ext},
we define the function
\begin{math}
  \mathbf{u}^{\flat}\colon
  \overline{\mathfrak{X}}^{(\rho)} \times [0,T_1]\to \CC^M
\end{math}
by the formula
\begin{equation}
\label{def:u^flat-ext}
  \mathbf{u}^{\flat}(x + \ii y,t)\eqdef \mathbf{u}^{(\ii y)}(x,t)
    \quad\mbox{ for all }\;
  (x,y,t)\in \RR^N\times \overline{Q}^{(\rho)}\times [0,T_1] \,.
\end{equation}
Clearly, by Proposition~\ref{prop-smooth_ext}, Part~{\rm (c)},
the function
\begin{math}
  \mathbf{u}^{\flat}\colon
  \overline{\mathfrak{X}}^{(\rho)} \times [0,T_1]\to \CC^M
\end{math}
is continuous and, by Proposition~\ref{prop-smooth_ext}, Part~{\rm (d)},
also holomorphic with respect to the complex variable
\begin{math}
  z = x + \ii y\in \mathfrak{X}^{(\rho)} = \RR^N + \ii Q^{(\rho)}
\end{math}
at every time $t\in [0,T_1]$.

Next, we set
\begin{math}
  \mathbf{u}_1^{(\ii y)}\eqdef \mathbf{u}^{(\ii y)}(\,\cdot\,,T_1)\in U_0
\end{math}
and repeat the procedure from above on the interval
$[T_1, 2 T_1]$ with the initial data
$\mathbf{u}_1^{(\ii y)}\in U_0$ at $t_0 = T_1$ in place of
$\mathbf{u}_0^{(\ii y)}\in U_0$ at $t_0 = 0$.
We stress that the interval length
$T_1\in (0,T-t_0]$ in Proposition~\ref{prop-smooth_ext}
is independent from the choice of the initial time $t_0\in (0,T)$ whenever
$[t_0, t_0+T_1]\subset [0,T]$.
Again, we apply Part~{\rm (c)} or Part~{\rm (d)}
of Proposition~\ref{prop-smooth_ext} with $t_0 = T_1$
(in place of $t_0 = 0$)
to conclude that there is a number $R_2\in (0,R_1)$ small enough,
such that even
$\mathbf{u}^{(\ii y)}(\,\cdot\,,t)\in U_0\subset U$
holds for every $t\in [0, 2 T_1]$, provided
$\rho\in (0,r_1)$ is chosen so small that also
\begin{equation*}
  \left\Vert \mathbf{u}_0^{(\ii y)} - \mathbf{u}_0
  \right\Vert_{ B^{s;p,p}(\RR^N) } < R_2 \quad (< R_1 < R_0)
    \quad\mbox{ holds for all }\, y\in \overline{Q}^{(\rho)} \,.
\end{equation*}
The desired function $\mathbf{u}^{\flat}$ is naturally extended from
the domain
\begin{math}
  \overline{\mathfrak{X}}^{(\rho)} \times [0,T_1]
\end{math}
to
\begin{math}
  \overline{\mathfrak{X}}^{(\rho)} \times [0, 2 T_1]
\end{math}
by setting
(cf.\ eq.~\eqref{def:u^flat-ext})
\begin{equation}
\label{def:u^flat-ext,2T_1}
  \mathbf{u}^{\flat}(x + \ii y,t)\eqdef \mathbf{u}^{(\ii y)}(x,t)
    \quad\mbox{ for all }\;
  (x,y,t)\in \RR^N\times \overline{Q}^{(\rho)}\times [T_1, 2 T_1] \,.
\end{equation}

We keep repeating this procedure (by ``induction'' on $k$)
with the initial data
\begin{math}
  \mathbf{u}_k^{(\ii y)}\eqdef \hfil\break
  \mathbf{u}^{(\ii y)}(\,\cdot\,, k T_1)\in U_0
\end{math}
successively for every $k = 0,1,2,\dots,m$
until reaching the inequalites
\begin{equation*}
  (m-1) T_1\leq T < m T_1 \quad\mbox{ at }\, k = m-1 \,.
\end{equation*}
In fact, setting
\begin{math}
  \mathbf{u}_{m-1}^{(\ii y)}\eqdef
  \mathbf{u}^{(\ii y)}(\,\cdot\,, (m-1) T_1)\in U_0
\end{math}
and repeating the procedure from above on the time interval
$[ (m-1) T_1, m T_1]$ with the initial data
$\mathbf{u}_{m-1}^{(\ii y)}\in U_0$ at $t_0 = (m-1) T_1$ in place of
$\mathbf{u}_0^{(\ii y)}\in U_0$ at $t_0 = 0$,
we can apply
Theorem~\ref{thm-Clem-Li} to conclude that
the spatially shifted Cauchy problem \eqref{e:Cauchy+z_0},
with the shift $z_0 = \ii y$
($y\in \overline{Q}^{(\rho)}$),
possesses a unique weak solution
\begin{math}
  \mathbf{u}^{(\ii y)}\in
  C\left( [ (m-1) T_1, m T_1]\to \mathbf{B}^{s;p,p}(\RR^N) \right) ,
\end{math}
local in time for $t\in [ (m-1) T_1, m T_1]$, that satisfies
\hfil\break
$\mathbf{u}^{(\ii y)}(\,\cdot\,,t)\in U$ for every
$t\in [ (m-1) T_1, m T_1]$.
Consequently, we may assume $T = m T_1$ instead of
$(m-1) T_1\leq T < m T_1$.
In this way we have constructed a finite set of numbers
$R_1, R_2, \dots, R_{m-1}, R_m$ such that
$0 < R_m < R_{m-1} < \dots < R_1 < R_0$ and, provided
$\rho\in (0,r_1)$ is chosen small enough, also
\begin{equation*}
  \left\Vert \mathbf{u}_k^{(\ii y)} - \mathbf{u}_0
  \right\Vert_{ B^{s;p,p}(\RR^N) }
  < R_k \quad\mbox{ holds for all }\, y\in \overline{Q}^{(\rho)} \,;\
    k = 0,1,2,\dots, m-1 \,,
\end{equation*}
together with
$\mathbf{u}^{(\ii y)}(\,\cdot\,,t)\in U_0\subset U$
for every $t\in [0, (m-1) T_1]$ and
$\mathbf{u}^{(\ii y)}(\,\cdot\,,t)\in U$
for every $t\in [ (m-1) T_1, m T_1]$.
Finally, the desired function $\mathbf{u}^{\flat}$ is defined successively
on the domains
\begin{math}
  \overline{\mathfrak{X}}^{(\rho)} \times [ (k-1) T_1, k T_1]
\end{math}
for each $k = 1,2,3,\dots,m$ by the formula
\begin{equation}
\label{def:u^flat-ext,T}
  \mathbf{u}^{\flat}(x + \ii y,t)\eqdef \mathbf{u}^{(\ii y)}(x,t)
    \quad\mbox{ for all }\;
  (x,y,t)\in \RR^N\times \overline{Q}^{(\rho)}\times [0,T] \,.
\end{equation}

To summarize the result of the procedure described above,
we have determined a constant $\rho\in (0,r_1)$, small enough,
such that for each shift $y\in \overline{Q}^{(\rho)}$
there is a unique weak solution
\begin{math}
  \mathbf{u}^{(\ii y)}\in
  C\left( [0,T]\to \mathbf{B}^{s;p,p}(\RR^N) \right)
\end{math}
to the spatially shifted Cauchy problem \eqref{e:Cauchy+z_0}
with the shift $z_0 = \ii y$, that satisfies
$\mathbf{u}^{(\ii y)}(\,\cdot\,,t)\in U_0\subset U$
for every $t\in [0, (m-1) T_1]$ and
$\mathbf{u}^{(\ii y)}(\,\cdot\,,t)\in U$
for every $t\in [0,T]$, where $T = m T_1$.
We apply Proposition~\ref{prop-smooth_ext}, Parts {\rm (c)} and~{\rm (d)},
once again to conclude that the function
\begin{math}
  \mathbf{u}^{\flat}\colon
  \overline{\mathfrak{X}}^{(\rho)}\times [0,T]\to \CC^M
\end{math}
constructed above in eq.~\eqref{def:u^flat-ext,T}
has the desired properties:
it is continuous and holomorphic in the space variable
\begin{math}
  z = x + \ii y\in \mathfrak{X}^{(\rho)} = \RR^N + \ii Q^{(\rho)}
\end{math}
with some $\rho\in (0,r_1)$ small enough, where
$0 < r_1 < \kappa_0\leq r_0$.

Now we are ready to finish our proof of Part~{\bf (iii)}
by further extending the (global) weak solution
\begin{math}
  \mathbf{u}\in C\left( [0,T]\to \mathbf{B}^{s;p,p}(\RR^N)\right)
\end{math}
of the Cauchy problem \eqref{e:Cauchy} from the domain 
$\overline{\mathfrak{X}}^{(\rho)}\times [0,T]$
of the (unique) spatial extension
\begin{math}
  \mathbf{u}^{\flat}\colon
  \overline{\mathfrak{X}}^{(\rho)}\times [0,T]\to \CC^M
\end{math}
to another continuous function
\begin{math}
  \tilde{\mathbf{u}}\colon \hfil\break
  \overline{\mathfrak{X}}^{(\rho)}\times
  \overline{\Delta}_{\vartheta'}^{T',T}\to \CC^M
\end{math}
which is holomorphic in
$\mathfrak{X}^{(\rho)}\times \Delta_{\vartheta'}^{T',T}$.
We recall that the solution
\begin{math}
  \mathbf{u}\in C\left( [0,T]\to \mathbf{B}^{s;p,p}(\RR^N)\right)
\end{math}
is assumed to exist by hypothesis (in Part~{\bf (iii)})
with the initial data
$\mathbf{u}_0\in \mathbf{B}^{s;p,p}(\RR^N)$
having  a (unique) holomorphic extension
\begin{math}
  \tilde{\mathbf{u}}_0\colon \mathfrak{X}^{(\kappa_0)}\to \CC^M
\end{math}
from $\RR^N$ to the complex domain
$\mathfrak{X}^{(\kappa_0)} \subset \CC^N$,
for some $\kappa_0\in (0,r_0]$.

We apply Theorem~\ref{thm-smooth_ext} to the function
\begin{math}
  \mathbf{u}^{\flat}\colon
  \overline{\mathfrak{X}}^{(\rho)}\times [0,T]\to \CC^M
\end{math}
on every time interval\hfil\break
$[t_0, t_0+T_1]\subset [0,T]$.
We remark that the number $T_1\in (0,T-t_0]$
depends on $r_1$ and $U$, but not on $t_0\in [0,T)$, provided
$[t_0, t_0+T_1]\subset [0,T]$.
In fact, making use of the same argument as above,
where we have extended the function $\mathbf{u}^{\flat}$ from the domain
$\overline{\mathfrak{X}}^{(\rho)}\times [0,T]$ to
$\overline{\mathfrak{X}}^{(\rho)}\times [0, m T_1]$
in case $(m-1) T_1\leq T < m T_1$,
we can extend $\mathbf{u}^{\flat}$ from the domain
$\overline{\mathfrak{X}}^{(\rho)}\times [t_0,T]$ to
$\overline{\mathfrak{X}}^{(\rho)}\times [t_0, t_0+T_1]$
in case $0\leq t_0 < T < t_0+T_1$.
Thus, if $T < t_0+T_1$ then we may replace $T$ by $T = t_0+T_1$
and, hence, assume that
$[t_0, t_0+T_1]\subset [0,T]$.
Consequently, by Theorem~\ref{thm-smooth_ext},
the (unique) weak solution
\begin{math}
  \mathbf{u}\in
  C\left( [0,T]\to \mathbf{B}^{s;p,p}(\RR^N)\right)
\end{math}
to the Cauchy problem \eqref{e:Cauchy}
possesses a unique holomorphic extension from the time interval
$(t_0, t_0+T_1)$ to the complex temporal domain
$t_0 + \Delta_{\vartheta'}^{T',T_1}$, such that
the continuous mapping
\begin{math}
  \tilde{\mathbf{u}}\colon
  \overline{\mathfrak{X}}^{(r_1)} \times
  \left( t_0 + \Delta_{\vartheta'}^{T',T_1} \right)
  \to \CC^M
\end{math}
constructed in Theorem~\ref{thm-smooth_ext} is holomorphic
in the space\--time domain
\begin{math}
  \mathfrak{X}^{(r_1)} \times
  \left( t_0 + \Delta_{\vartheta'}^{T',T_1} \right) .
\end{math}
Since the interval
$[t_0, t_0+T_1]\subset [0,T]$ is arbitrary, both statements
{\bf ($\mathbf{iii}_1$)} and {\bf ($\mathbf{iii}_2$)}
and the complex analyticity statement
{\bf ($\mathbf{iii}_3$)} in Part~{\bf (iii)} of Theorem~\ref{thm-Main}
follow from eq.~\eqref{glob:0+T^(r)}.
More precisely, the desired holomorphic extension
\begin{math}
  \tilde{\mathbf{u}}\colon
  \overline{\mathfrak{X}}^{(r_1)} \times \Delta_{\vartheta'}^{T',T}
  \to \CC^M
\end{math}
is obtained by shifting the temporal domain
$t_0 + \Delta_{\vartheta'}^{T',T_1}$ with the vertex $t_0$
ranging from the left to the right over the time interval $[0, T-T_1]$.
In this process, the uniqueness result
in Theorem~\ref{thm-smooth_ext}, Part~{\rm (i)},
guarantees that the function
\begin{math}
    \tilde{\mathbf{u}} (\,\cdot\, + \ii y,\, t)
  = \mathbf{u}^{(\ii y)}(t)
  = \tilde{\mathbf{u}} (\,\cdot\, + \ii y,\, t_0 + s)
  = \mathbf{u}^{(\ii y)}(t)\in \mathbf{B}^{s;p,p}(\RR^N) ,
\end{math}
with $0\leq s = t-t_0\leq T_1$, is well\--defined for all
$(y,t)\in Q^{(r_1)}\times \Delta_{\vartheta'}^{T',T}$
independently from the particular choice of the vertex $t_0\in [0, T-T_1]$
of the the complex temporal domain
$t_0 + \Delta_{\vartheta'}^{T',T_1}$.
Furthermore, in Part~{\bf (iii)}, {\bf ($\mathbf{iii}_1$)},
ineq.~\eqref{sup_(y,t):tilde_u}
holds with $\rho$ in place of $\kappa_0$, whereas
$r'\in (0,\kappa_0]$ has to be replaced by $\rho$, as well.

This concludes our proof of Theorem~\ref{thm-Main}.
\qed
\par\vskip 10pt

We conclude this section with the proof of the estimate in
ineq.~\eqref{eq:u_Hardy^2}.

\par\vskip 10pt
{\it Proof of\/} Proposition~\ref{prop-Main}.

We recall from Section~\ref{s:Analyt_space} that
the shifted continuous function
\begin{math}
    \mathbf{u}^{(\ii y)}\colon t\mapsto \mathbf{u}^{(\ii y)}(t)
  = \tilde{\mathbf{u}} (\,\cdot\, + \ii y,\, t)
    \colon [0,T]\to \mathbf{B}^{s;p,p}(\RR^N)
\end{math}
is a unique weak solution of
the spatially shifted Cauchy problem \eqref{e:Cauchy+z_0}
with the shift $z_0 = \ii y$ ($y\in Q^{(r_1)}$).
Consequently,
\begin{math}
  \mathbf{u}^{(\ii y)}\in
  C\left( [0,T]\to \mathbf{B}^{s;p,p}(\RR^N) \right)
\end{math}
is also a strict solution
(cf.\ Definition~\ref{def-strict_sol})
to the following abstract initial value problem,
for every $y\in Q^{(r_1)}$:
\begin{equation}
\label{eq:Cauchy_u^(iy)}
\left\{
\begin{alignedat}{2}
  \frac{\mathrm{d}}{\mathrm{d}t}\, \mathbf{u}^{(\ii y)}
  - A^{(\ii y)}(t) \mathbf{u}^{(\ii y)}
& {}
  = \mathbf{F}^{(\ii y)} \left( t,\, \mathbf{u}^{(\ii y)}(t) \right)
  \quad\mbox{ for a.e. } t\in (0,T) \,;
\\
  \mathbf{u}^{(\ii y)}(0)
& {}
  = \tilde{\mathbf{u}}_0(\,\cdot\, + \ii y)
  \in E_{ 1 - \frac{1}{p} ,\, p } = \mathbf{B}^{s;p,p}(\RR^N) \,,
\end{alignedat}
\right.
\end{equation}
cf.\ eqs.\ \eqref{e:Cauchy+z_0} and~\eqref{eq:Cauchy_dz}.
Here,
$A^{(z_0)}(t)\in \mathcal{L}(E_1\to E_0)$
is the bounded linear (partial differential) operator introduced in
eq.~\eqref{def:Cauchy_dz}, satisfying
$A^{(z_0)}(t)\in \mathrm{MR}_p(E_1\to E_0)$ for every $t\in [0,T]$, and
\begin{math}
  \mathbf{F}^{(z_0)}\colon [0,T]\times U\to E_0 = \mathbf{L}^p(\RR^N)
\end{math}
stands for the ``shifted'' Nemytskii operator defined in
Remark~\ref{rem-f_Lipschitz}, eq.~\eqref{def:F+z_0}.
We recall that both constants, $C_1\equiv C_1(C(U))$ and $L\equiv L(U)$
in inequalities \eqref{e:|df/dX|<C} and \eqref{est:F:v_1-v_2},
respectively, are independent from the shift by
$z_0\in \mathfrak{X}^{(r_1)}$ in case $x\in \RR^N$ is replaced by $x + z_0$;
with $z_0 = \ii y$ ($y\in Q^{(r_1)}$) in our case.

We now derive an estimate analogous to ineq.~\eqref{est:strict_sol}
for our shifted Cauchy problem \eqref{eq:Cauchy_u^(iy)}
in place of the (original) abstract problem \eqref{lin:abstr_Cauchy}.
Inspecting the proof of Theorem 2.1 in 
{\sc Ph.\ Cl\'ement} and {\sc S.\ Li} \cite[pp.\ 20--23]{Clem-Li}
and combining it with our linear perturbation result in
Lemma~\ref{lem-perturb_MR:A<B} and the estimate in
ineq.~\eqref{est:B+A_sol},
we arrive at the following estimate
for our shifted Cauchy problem \eqref{eq:Cauchy_u^(iy)}
in place of the abstract problem \eqref{lin:abstr_Cauchy},
\begin{equation}
\label{est:strict_sol^(iy)}
\begin{aligned}
&   \int_0^T
  \genfrac{\|}{\|}{}0{ \mathrm{d}\mathbf{u}^{(\ii y)} }{\mathrm{d}t}_{E_0}^p
    \,\mathrm{d}t
  + \int_0^T
    \left\| A^{(\ii y)}(0) \mathbf{u}^{(\ii y)}(t)
    \right\|_{E_0}^p \,\mathrm{d}t
\\
& {}
  \leq M_{p,T}
    \left(
      \left\| \mathbf{u}^{(\ii y)}(0)
      \right\|_{ E_{ 1 - \frac{1}{p} ,\, p } }^p
    + \int_0^T
    \left\| \mathbf{F}^{(\ii y)} \left( t,\, \mathbf{u}^{(\ii y)}(t)\right)
    \right\|_{E_0}^p \,\mathrm{d}t
    \right) \,,
\end{aligned}
\end{equation}
where $M_{p,T}\in (0,\infty)$
is a constant independent from the initial data
\begin{math}
  \mathbf{u}^{(\ii y)}(0) = \tilde{\mathbf{u}}_0(\,\cdot\, + \ii y)
  \in E_{ 1 - \frac{1}{p} ,\, p } = \mathbf{B}^{s;p,p}(\RR^N)
\end{math}
and the right\--hand side
\begin{math}
  \mathbf{F}^{(\ii y)} \left( t,\, \mathbf{u}^{(\ii y)}(t) \right)
\end{math}
of eq.~\eqref{eq:Cauchy_u^(iy)}, as well.
Since
\begin{math}
  A^{(\ii y)}(0)\in \mathrm{MR}_p(E_1\to E_0)
    \subset \mathrm{Hol}(E_1\to E_0)
\end{math}
holds by the proof of Proposition~\ref{prop-smooth_ext}, Part~{\rm (d)},
there is a number
$\lambda_0\in \RR_+ = [0,\infty)$, sufficiently large, such that
the bounded linear operator
$\lambda_0 I - A^{(\ii y)}(0)\colon E_1\to E_0$
is an isomorphism of $E_1$ onto $E_0$.
Hence, its inverse satisfies
\begin{math}
  ( \lambda_0 I - A^{(\ii y)}(0) )^{-1}\in \mathcal{L}(E_0\to E_1) .
\end{math}
We conclude that there are constants $c_1, C_1\in (0,\infty)$ and
$c_2, C_2\in \RR_+$
(both sufficiently large, depending on $\lambda_0\geq 0$)
such that both inequalities
\begin{equation*}
    c_1\, \| u\|_{E_1} - c_2\, \| u\|_{E_0}
  \leq \| A^{(\ii y)}(0) u\|_{E_0}\leq
    C_1\, \| u\|_{E_1} + C_2\, \| u\|_{E_0}
  \quad\mbox{ hold for all }\, u\in E_1 \,.
\end{equation*}
Consequently, we have also (respectively)
\begin{align}
\label{ineq:E_1<A^(iy)}
&   2^{-(p-1)}\, c_1\, \| u\|_{E_1}^p
  \leq
    \| A^{(\ii y)}(0) u\|_{E_0}^p + c_2^p\, \| u\|_{E_0}^p
  \quad\mbox{ and }
\\
\label{ineq:E_1>A^(iy)}
&   \| A^{(\ii y)}(0) u\|_{E_0}^p
  \leq 2^{p-1}
  \left( C_1^p\, \| u\|_{E_1}^p + C_2^p\, \| u\|_{E_0}^p \right)
  \quad\mbox{ for all }\, u\in E_1 \,.
\end{align}
Furthermore, for every $t\in [0,T]$, we split the expression
\begin{align*}
    \mathbf{F}^{(\ii y)} \left( t,\, \mathbf{u}^{(\ii y)}(t)\right)
  = \mathbf{F}^{(\ii y)} \left( t,\, \mathbf{0}\right)
  + \left[
    \mathbf{F}^{(\ii y)} \left( t,\, \mathbf{u}^{(\ii y)}(0)\right)
  - \mathbf{F}^{(\ii y)} \left( t,\, \mathbf{0}\right)
    \right]
\end{align*}
and apply Lemma~\ref{lem-f_Lipschitz} and Remark~\ref{rem-f_Lipschitz}
to derive the following analogue of ineq.~\eqref{est:F:v-0}
(where we insert $\mathbf{v} = \mathbf{u}^{(\ii y)}(t)$):
\begin{equation}
\label{est:F:u^(iy)-0}
    \left|
      \mathbf{F}^{(\ii y)} \left( t,\, \mathbf{u}^{(\ii y)}(t)\right) (x)
    \right|
  \leq
    \left|
      \mathbf{F}^{(\ii y)} (t,\mathbf{0})(x)
    \right|
  + L\sum_{|\beta|\leq m} \sum_{k=1}^M
    \left|
    \frac{ \partial^{|\beta|} }{ \partial x^{\beta} }\, u^{(\ii y)}_k(x,t)
    \right| \,,\quad x\in \RR^N \,,
\end{equation}
for all $t\in [0,T]$.
Here,
\begin{equation*}
    \mathbf{F}^{(\ii y)} (t,\mathbf{0})(x)
  = \mathbf{f}(x + \ii y,t; \vec{\mathbf{0}}) \,,\quad x\in \RR^N \,,\quad
    \vec{\mathbf{0}} = \left( 0\right)_{|\beta|\leq m}
  \equiv (0,\dots,0)\in \CC^{M\tilde{N}} \,,
\end{equation*}
satisfies
\begin{math}
  \mathbf{F}^{(\ii y)} (t,\mathbf{0})\in E_0 = \mathbf{L}^p(\RR^N) ,
\end{math}
by the $L^p$\--integrability condition in ineq.~\eqref{e:f_Hardy^2}, i.e.,
\begin{equation*}
    \| \mathbf{F}^{(\ii y)} (t,\mathbf{0}) \|_{E_0}
  = \left(
  \int_{\RR^N} |\mathbf{f} (x + \ii y,t; \vec{\mathbf{0}})|^p \,\mathrm{d}x
    \right)^{1/p}
  \leq K \quad\mbox{ for all $t\in [0,T]$, }
\end{equation*}
where $K\in (0,\infty)$ is a constant.
Each term
\begin{math}
  \frac{ \partial^{|\beta|} }{ \partial x^{\beta} }\,
    u^{(\ii y)}_k(\,\cdot\,,t) \in L^p(\RR^N)
\end{math}
on the right\--hand side of ineq.~\eqref{est:F:u^(iy)-0} above belongs to
the Besov space
$B^{s - |\beta|;p,p}(\RR^N)$
which, thanks to
\begin{math}
  |\beta|\leq m < s = 2m\left( 1 - \genfrac{}{}{}1{1}{p} \right) < 2m ,
\end{math}
is continuously imbedded into another Besov space,
\begin{math}
  B^{s - |\beta|;p,p}(\RR^N) \hookrightarrow
  B^{s-m;p,p}(\RR^N) \hookrightarrow L^p(\RR^N) .
\end{math}
Applying these estimates to the right\--hand side of
ineq.~\eqref{est:F:u^(iy)-0}, we thus obtain
\begin{equation}
\label{est^p:F:u^(iy)-0}
\begin{aligned}
    \left\| \mathbf{F}^{(\ii y)} \left( t,\, \mathbf{u}^{(\ii y)}(t)\right)
    \right\|_{E_0}^p
  \leq 2^{p-1}
    \left(
      K^p + \gamma_{N,m,M}\, L^p\cdot
      \| \mathbf{u}^{(\ii y)}(t)\|_{ B^{s;p,p}(\RR^N) }^p
    \right)
\end{aligned}
\end{equation}
for all $t\in [0,T]$ and every $y\in Q^{(r_1)}$, where
$\gamma_{N,m,M}\in (0,\infty)$
is a numerical constant depending only on $N$,$m$, and $M$.
Recalling
$\mathbf{u}^{(\ii y)}(t)\in U$ for all
$(y,t)\in Q^{(r_1)}\times [0,T]$ and the definition of the set
\begin{math}
  U\subset E_{ 1 - \frac{1}{p} ,\, p } = \mathbf{B}^{s;p,p}(\RR^N)
\end{math}
in eq.~\eqref{def:U},
we conclude that the right\--hand side of ineq.~\eqref{est^p:F:u^(iy)-0}
can be estimated from above by a constant
$C\equiv C(K,U)\in (0,\infty)$ independent from
$(y,t)\in Q^{(r_1)}\times [0,T]$:
\begin{equation}
\label{est^p:|F|<C}
\begin{aligned}
    \left\| \mathbf{F}^{(\ii y)} \left( t,\, \mathbf{u}^{(\ii y)}(t)\right)
    \right\|_{E_0}^p
  \leq C(K,U) \quad\mbox{ for all }\, (y,t)\in Q^{(r_1)}\times [0,T] \,.
\end{aligned}
\end{equation}

Finally, we combine this estimate with Theorem~\ref{thm-smooth_ext}
to arrive at
\begin{equation}
\label{est^p:|F|<C^tilde}
\begin{aligned}
    \left\| \mathbf{F}^{(\ii y)}
            \left( t,\, \tilde{\mathbf{u}} (\,\cdot\, + \ii y, t) \right)
    \right\|_{E_0}^p
  \leq \tilde{C}(K,U) \quad\mbox{ for all }\,
       (y,t)\in Q^{(r_1)}\times \Delta_{\vartheta'}^{T',T} \,,
\end{aligned}
\end{equation}
where
$\tilde{C}\equiv \tilde{C}(K,U)\in (0,\infty)$
is a constant independent from
$(y,t)\in Q^{(r_1)}\times \Delta_{\vartheta'}^{T',T}$.
Recall from our proof of Theorem~\ref{thm-Main}, Part~{\bf (iii)},
that the function
\begin{math}
  \tilde{\mathbf{u}}\colon
  \overline{\mathfrak{X}}^{(r_1)} \times \Delta_{\vartheta'}^{T',T}
  \to \CC^M
\end{math}
stands for the unique holomorphic temporal extension of the function
\begin{math}
  \mathbf{u}^{\flat}\colon
  \overline{\mathfrak{X}}^{(r_1)} \times (0,T)\to \CC^M
\end{math}
defined in formula \eqref{def:u^flat-ext,T}.
This extension, $\tilde{\mathbf{u}}$, satisfies
\begin{math}
  \tilde{\mathbf{u}}(x,y,t) = \mathbf{u}^{(\ii y)}(x,t)
\end{math}
for all
$(x,y,t)\in \RR^N\times \overline{Q}^{(r_1)}\times [0,T]$
and $\tilde{\mathbf{u}}(\,\cdot\, + \ii y,t)$ $\in U$ for all
\begin{math}
  (y,t)\in
  \overline{Q}^{(r_1)} \times \Delta_{\vartheta'}^{T',T} .
\end{math}
We employ the norm defined in eq.~\eqref{sup_y:u_0}
and ineq.~\eqref{est^p:|F|<C^tilde}
to estimate the right\--hand side of \eqref{est:strict_sol^(iy)}:
\begin{equation}
\label{est:sol^(iy)}
\begin{aligned}
&   \int_0^T
  \genfrac{\|}{\|}{}0{ \mathrm{d}\mathbf{u}^{(\ii y)} }{\mathrm{d}t}_{E_0}^p
    \,\mathrm{d}t
  + \int_0^T
    \left\| A^{(\ii y)}(0) \mathbf{u}^{(\ii y)}(t)
    \right\|_{E_0}^p \,\mathrm{d}t
\\
& {}
  \leq M_{p,T}
    \left(
      \sup_{ y\in Q^{(r_1)} }
    \| \tilde{\mathbf{u}}_0 (\,\cdot\, + \ii y)
    \|_{ B^{s;p,p}(\RR^N) }
    + \tilde{C}(K,U)\, T
    \right)
\end{aligned}
\end{equation}
for all $(y,t)\in Q^{(r_1)}\times [0,T]$.
However, in the norm
\begin{equation*}
  \mathfrak{N}^{(r_1)} (\tilde{\mathbf{u}}_0) \eqdef
    \sup_{ y\in Q^{(r_1)} }
    \| \tilde{\mathbf{u}}_0 (\,\cdot\, + \ii y)
    \|_{ B^{s;p,p}(\RR^N) } < \infty
\end{equation*}
we have
\begin{math}
  \mathbf{u}^{(\ii y)}(0) = \tilde{\mathbf{u}}_0(\,\cdot\, + \ii y)
  \in U
  \subset E_{ 1 - \frac{1}{p} ,\, p } = \mathbf{B}^{s;p,p}(\RR^N)
\end{math}
for every $y\in Q^{(r_1)}$ owing to our choice of the number
$r_1\in (0,r)$ being sufficiently small in (and before)
Proposition~\ref{prop-smooth_ext}.
Consequently, we can estimate the right\--hand side of
ineq.~\eqref{est:sol^(iy)} above by another constant
$\tilde{C}'\equiv \tilde{C}'(p,T,K,U)\in (0,\infty)$ independent from
$(y,t)\in Q^{(r_1)}\times [0,T]$:
\begin{equation*}
    \int_0^T
  \genfrac{\|}{\|}{}0{ \mathrm{d}\mathbf{u}^{(\ii y)} }{\mathrm{d}t}_{E_0}^p
    \,\mathrm{d}t
  + \int_0^T
    \left\| A^{(\ii y)}(0) \mathbf{u}^{(\ii y)}(t)
    \right\|_{E_0}^p \,\mathrm{d}t
  \leq \tilde{C}'(p,T,K,U) \,.
\end{equation*}
We estimate the left\--hand side of this inequality by
\eqref{ineq:E_1<A^(iy)} combined with
$\mathbf{u}^{(\ii y)}(t)\in U$ for all
$(y,t)\in Q^{(r_1)}\times [0,T]$, thus arriving at
\begin{equation}
\label{est:|u^(iy)|}
\begin{aligned}
&   \int_0^T
  \genfrac{\|}{\|}{}0{ \mathrm{d}\mathbf{u}^{(\ii y)} }{\mathrm{d}t}_{E_0}^p
    \,\mathrm{d}t
  + 2^{-(p-1)}\, c_1\int_0^T
    \left\| \mathbf{u}^{(\ii y)}(t) \right\|_{E_1}^p \,\mathrm{d}t
\\
& {}
  \leq \tilde{C}'(p,T,K,U)
  + c_2^p\int_0^T
    \left\| \mathbf{u}^{(\ii y)}(t) \right\|_{E_0}^p \,\mathrm{d}t
  \leq \hat{C}(p,T,K,U) \,,
\end{aligned}
\end{equation}
where
$\hat{C}\equiv \hat{C}(p,T,K,U)\in (0,\infty)$
is a constant independent from
$(y,t)\in Q^{(r_1)}\times [0,T]$.

Within the \underline{restriction} to the real time $t\in [0,T]$,
the desired estimate in ineq.~\eqref{eq:du/dt_Hardy^2}
is derived directly from ineq.~\eqref{est:|u^(iy)|} above
for all pairs
\begin{math}
  (y,t)\in Q^{(r')}\times [0,T]\subset
           Q^{(r')}\times \overline{\Delta}_{\vartheta'}^{T',T} .
\end{math}

In order to \underline{extend} ineq.~\eqref{eq:du/dt_Hardy^2}
to the complex time $t\in \overline{\Delta}_{\vartheta'}^{T',T}$
with $t = \sigma + \ii\tau$ ($\sigma, \tau\in \mathbb{R}$),
we take advantage of Theorem~\ref{thm-smooth_ext} once again.
We will consider the function
\begin{math}
  \tilde{\mathbf{u}}\colon
  \overline{\mathfrak{X}}^{(r_1)} \times \Delta_{\vartheta'}^{T',T}
  \to \CC^M
\end{math}
constructed in our proof of Theorem~\ref{thm-Main}, Part~{\bf (iii)},
along the complex temporal path
\begin{math}
  \tilde{\theta}\colon [0,T]\to \Delta_{\vartheta'}^{T',T}
  \colon
  s\mapsto \tilde{\theta}(s)\eqdef s + \ii\varsigma_1(s/T') \tau
\end{math}
which consists of two straight line segments,
\begin{alignat*}{2}
& \tilde{\theta}_1\colon [0,T']\to \Delta_{\vartheta'}^{T',T}
  \colon s\mapsto \tilde{\theta}_1(s)
      \eqdef \left( 1 + \ii\frac{\tau}{T'}\right) s
&&\quad\mbox{ with }\, 0\leq s\leq T' \,,
  \quad\mbox{ and }\quad
\\
& \tilde{\theta}_2\colon [T',T]\to \Delta_{\vartheta'}^{T',T}
  \colon s\mapsto \tilde{\theta}_2(s)\eqdef s + \ii\tau
&&\quad\mbox{ with }\, T'\leq s\leq T \,.
\end{alignat*}
Notice that
$\tilde{\theta}_1(0) = 0$,
$\tilde{\theta}_1(T') = \tilde{\theta}_2(T') = T'+ \ii\tau$, and
$\tilde{\theta}_2(T) = T + \ii\tau$.
We replace the (complex) time variable $t$
in the original Cauchy problem \eqref{e:Cauchy}
by the new (real) time variable $s\in [0,T]$,
thus obtaining two new abstract differential equations
with the time derivatives
\begin{alignat}{3}
\label{e:du/ds:s<T'}
    \frac{\partial\mathbf{u}}{\partial s}
& {}
  = \left( 1 + \ii\frac{\tau}{T'}\right)
    \frac{\partial\mathbf{u}}{\partial t} (x,t)
    \Big\vert_{ t = \tilde{\theta}_1(s) }
&&  \quad\mbox{ for }\, 0\leq s\leq T'
&&  \quad\mbox{ and }\quad
\\
\label{e:du/ds:s>T'}
    \frac{\partial\mathbf{u}}{\partial s}
& {}
  = \frac{\partial\mathbf{u}}{\partial t} (x,t)
    \Big\vert_{ t = \tilde{\theta}_2(s) }
&&  \quad\mbox{ for }\, T'\leq s\leq T \,,
&&  \quad\mbox{ respectively.}
\end{alignat}
Finally, we apply Lemma~\ref{lem-perturb_MR:A<B} and the estimate in
ineq.~\eqref{est:B+A_sol} to the new problem for $0\leq s\leq T'$,
thus arriving at the desired estimate in \eqref{eq:du/dt_Hardy^2}
for $0\leq \sigma\leq T'$.
For $T'\leq s\leq T$ we can use the definition of a strict solution
(Definition~\ref{def-strict_sol})
directly and combine it with the estimate in \eqref{est:strict_sol}
to obtain the estimate in \eqref{eq:du/dt_Hardy^2}
for $T'\leq \sigma\leq T$.
We conclude that \eqref{eq:du/dt_Hardy^2} is valid also for all pairs
\begin{math}
  (y,t)\in Q^{(r')}\times \overline{\Delta}_{\vartheta'}^{T',T}
\end{math}
with $t = \sigma + \ii\tau$ ($\sigma, \tau\in \mathbb{R}$).

The desired estimate in ineq.~\eqref{eq:u_Hardy^2}
now follows from ineq.~\eqref{eq:du/dt_Hardy^2}
by applying eqs.\ \eqref{norm:Y_p^p^} and~\eqref{norm:E_th^p}.
Also Proposition~\ref{prop-Main} is proved.
\qed
\par\vskip 10pt

\color{green}
$\maltese\maltese\maltese\maltese\maltese\maltese\maltese\maltese$
\color{black}

\section{An application to a Risk Model in Mathematical Finance}
\label{s:Appl}

Standard models in derivative pricing, including
the {\em Black\--Scholes model\/}
(see {\sc F.\ Black} and {\sc M.~S.\ Scholes} \cite{Black-Scholes}
 and {\sc R.~C.\ Merton} \cite{Merton}) and
the {\em Heston model\/}
(see {\sc S.~L.\ Heston} \cite{Heston})
take advantage of risk neutral valuation methods for
the arbitrage\--free (``fair'') price of the derivative.
The methods are economically justified by riskless hedging arguments
introduced in \cite{Black-Scholes, Merton}; see also
{\sc J.-P.\ Fouque}, {\sc G. Papanicolaou}, and {\sc K.~R.\ Sircar}
\cite{FouqPapaSir} and
{\sc J.~C.\ Hull} \cite{Hull-book}
for detailed explanations of these arguments.
An important assumption of these models,
which is used in most of the hedging arguments,
is the possibility to borrow and lend any amount of money at
a risk\--free interest rate.
This crucial conjecture has been questioned
as a consequence of the financial crisis
that started in 2007 and resulted in
the bankruptcy of major financial entities like Lehman Brothers.
{\sc Enron}'s bankruptcy in 2001 is briefly described in
\cite[p.~537]{Hull-book}, Business Snapshot 23.1.
Namely, traders have to take into consideration
the increased chance of a {\em\bfseries default\/}.
For this reason many trades contain
a {\em\bfseries collateral\/} against default and also
the pricing of non\--collateralized derivatives has to be adjusted.
A standard book on risk management has been written by
{\sc J.~C.\ Hull} \cite{Hull-risk}.
{\sc V.\ Piterbarg} \cite{Piterbarg} discusses
the differences or convexity adjustments between
the price processes of collateralized and non\--collaterlized contracts
which could result in
{\em\bfseries funding valuation adjustments\/} of the price processes.
It is only natural that traders have different
{\em funding costs\/} for transactions and try to include them
in the price of the contract.
{\sc J.\ Hull} and {\sc A.\ White} reasoned in \cite{Hull-White-short}
that there exists no theoretical basis for
such a {\em\bfseries funding value adjustment\/} (FVA).
Also {\sc C.\ Burgard} and {\sc M.\ Kjaer}
\cite{Burgard-Kjaer, Burgard-Kjaer-PDE-add}
came to a similar conclusion, by using different arguments.
However, since these theoretical arguments are not convincing
from a practitioner's point of view, but
traders make the adjustments anyway,
{\sc J.\ Hull} and {\sc A.\ White} studied the consequences of
funding value adjustment in a more practice\--oriented way in
\cite{Hull-White-3, Hull-White-XVA}.
Further common adjustments of the no\--default value of a derivative are
{\em\bfseries credit value adjustments\/} (CVA) and the related
{\em\bfseries  debit value adjustments\/} (DVA).
In a particular trade both parties have to take
the possibility of default of the counterparty into account
which is the {\em\bfseries bilateral counterparty risk\/}.
Price\--reducing credit value adjustments are made by
the trader to have a collateral against default of the counterparty
(e.g., a bank), whereas
debit value adjustments are
the corresponding adjustments made by the counterparty.
The {\em\bfseries sum of all\/} adjustments to the value of the derivative
evaluated in the absence of default is often refered to as {\bf XVA} with
\begin{equation*}
  \mathrm{XVA}\; =\; \mathrm{FVA}\; -\; \mathrm{CVA}\; +\;\mathrm{DVA} \,.
\end{equation*}
A general partial differential equation for
the {\em\bfseries adjusted value\/} under the bilateral counterparty risk
and funding value adjustments has been derived in
{\sc C.\ Burgard} and {\sc M.\ Kjaer} \cite[Section~3]{Burgard-Kjaer-PDE}
using hedging arguments.
This partial differential equation is {\em nonlinear\/}
if the mark\--to\--market value at default is considered to be
the total value of the derivative including all value adjustments
(see \cite[Section~4]{Burgard-Kjaer-PDE})
and {\em linear\/} if the mark\--to\--market value is given by
the no\--default value of the derivative
(see \cite[Section~5]{Burgard-Kjaer-PDE}).
In both cases, the partial differential equations are well\--suited
for numerical calculations of the adjusted value of the derivative;
see, e.g.,
{\sc I.\ Arregui}, {\sc B.\ Salvador}, and {\sc C.\ V\'{a}zquez}
\cite{ArrSalVaz} for recent results.
Following notation introduced in \cite{Burgard-Kjaer-PDE},
let us consider a derivative contract with payoff $H$ between
a {\em trader\/}, $B$, and a {\em counterparty\/}, $C$, on an asset~$S$,
e.g., a stock, that is not affected in case of a default by
one of the two counterparties and follows the stochastic dynamics,
\begin{equation}
\label{eq:stoch_dym}
  \td S_t = \mu(t) S_t\td t + \sigma(t) S_t\td W_t \,,
\end{equation}
where the drift $\mu\colon [0,\infty)\to \RR$ and the volatility
$\sigma\colon [0,\infty)\to \RR$
are positive deterministic (Borel measurable) functions and
$(W_t)_{t\geqslant 0}$ is a one\--dimensional Brownian motion.
Let $V$ denote the {\em fair price\/} (the ``risk\--less'' value) 
of the derivative in the setting without default and
let $\hat{V}$ denote the {\em adjusted price\/} (the ``risky'' value)
including {\em funding value adjustments\/} ($= \mathrm{FVA}$)
and {\em bilateral counterparty risk\/}
($= -\, \mathrm{CVA}\, +\, \mathrm{DVA}$),
\begin{equation*}
  \hat{V}\; =\; V\; +\; \mathrm{XVA}\;
  =\; V\; +\; \mathrm{FVA}\; -\; \mathrm{CVA}\; +\; \mathrm{DVA} \,.
\end{equation*}
By It\^{o}'s formula, the generator
$\mathcal{A}_t$ of the Markov process \eqref{eq:stoch_dym}
is the partial differential operator
\begin{equation}
\label{def:genA}
  \mathcal{A}_t
  = \frac{1}{2}\, \sigma^2 S^2\, \frac{\partial^2}{\partial S^2}
  + (q_S - \gamma_S) S\, \frac{\partial}{\partial S} \,,
\end{equation}
where $\gamma_S$ is the {\em dividend income rate\/} of $S$ and
$q_S$ represents the {\em financing costs\/} that depend on
the {\em risk\--free rate\/} $r$ and repo\--rate of the asset
(e.g., the Fed repo\--rate).
The decisive variable in the bilateral counterparty risk models studied in
\cite{ArrSalVaz, Burgard-Kjaer, Burgard-Kjaer-PDE, Burgard-Kjaer-PDE-add}
is the {\em\bfseries mark\--to\--market value\/}
(cf.\ ``close\--out''), $M$, introduced in
\cite[Sect.~3, Eq.~(24)]{Burgard-Kjaer-PDE}.
Only two different values of $M$ seem to be of significant interest,
namely, $M = \hat{V}$ and $M = V$ as described below:

If we set the mark\--to\--market value at default
$M = \hat{V}$, then the total value $\hat{V}$ satisfies
the {\em nonlinear\/} partial differential equation
\begin{equation}
\label{eq:Risk_nonlin}
    \frac{\partial}{\partial t}\, \hat{V}
  + \mathcal{A}_t\hat{V} - r\hat{V}
  = {}- (1-R_B)\lambda_B \hat{V}^{-}
      + (1-R_C)\lambda_C \hat{V}^{+} + s_F\hat{V}^{+} \,,
\end{equation}
with the final value $\hat{V}(S,T) = H(S)$ at maturity time $t = T$, by
\cite[Sect.~4, Eq.~(26)]{Burgard-Kjaer-PDE}.
Here, we have abbreviated
$x^{+}\eqdef \max\{ x, 0\}$ and
$x^{-}\eqdef \max\{-x, 0\}$ for $x\in \RR$; hence,
$x = x^{+} - x^{-}$.
We remark that the definition of the negative part $x^{-}$ of $x\in \RR$
often differs in the literature
(\cite{ArrSalVaz, Burgard-Kjaer, Burgard-Kjaer-PDE});
it may be used with the negative sign, i.e.,
$x^{-} = \min\{ x, 0\}$ ($\leq 0$), whence
$x = x^{+} + x^{-}$.
We will respect this convention below only when approximating the function
$x\mapsto x^{-}$ within the sum $x = x^{+} + x^{-}$.
Otherwise we use $x^{-} = \max\{-x, 0\}$ ($\geq 0$).
The parameters $\lambda_B$ and $\lambda_C$ are given by
$\lambda_B = r_B-r$ and $\lambda_C = r_C-r$, where
$r_B$ and $r_C$ are the {\em\bfseries yields on recovery\--less bonds\/}
for $B$ and $C$, respectively.
$R_B$ and $R_C$, respectively,
are {\em recovery rates\/} on the derivatives' mark\--to\--market value
at default and $s_F = r_F-r$ is the {\em funding spread\/} between
the {\em sellers funding rate\/} $r_F$ for borrowed cash and
the risk\--free rate~$r$.
We refer the interested reader to
{\sc C.\ Burgard} and {\sc M.\ Kjaer}
\cite[Sect.~2, pp.\ 2--4]{Burgard-Kjaer-PDE}
for further details concerning {\em\bfseries recovery\--less bonds\/}.

In contrast, if we assume the mark\--to\--market value
$M = V$, then the resulting partial differential equation for
$\hat{V}$ is {\em linear\/}, albeit inhomogeneous with
source terms on the right\--hand side,
\begin{equation}
\label{eq:Risk_lin}
    \frac{\partial}{\partial t}\, \hat{V}
  + \mathcal{A}_t\hat{V} - (r + \lambda_B + \lambda_C) \hat{V}
  = (R_B\lambda_B + \lambda_C) V^{-}
  - (\lambda_B + R_C\lambda_C) V^{+} + s_F V^{+} \,,
\end{equation}
with the final value $\hat{V}(S,T) = H(S)$, by
\cite[Sect.~5, Eq.~(46)]{Burgard-Kjaer-PDE}.
Of course, it is assumed that the fair price of the derivative, $V$,
is known.
It is claimed in \cite[Sect.~3]{Burgard-Kjaer-PDE}
that the vast majority of papers on valuation of conterparty risk uses
this choice ($M = V$).
From the mathematical point of view, also any {\it convex combination\/}
\begin{math}
  M = (1-\theta)\cdot V + \theta\cdot \hat{V}
    = V + \theta\cdot (\mathrm{XVA})
\end{math}
of $V$ and $\hat{V}$, with a constant $\theta\in [0,1]$,
might be of economic interest, as well.

We would like to investigate the question of {\em market completeness\/}
for the nonlinear model \eqref{eq:Risk_nonlin}
raised for related financial market models in
{\sc M.~H.~A.\ Davis} and {\sc J.\ Ob{\l}\'{o}j} \cite{Davis-Obloj}.
There, the authors have shown that the problem of
{\em\bfseries market completeness\/} in Mathematical Finance
is closely connected to (in fact, follows from)
the analyticity of the derivative price.
We refer to
{\sc P.\ Tak\'a\v{c}} \cite[Section 8, pp. 74--83]{Takac}
for a survey of results regarding the correlation between
market completeness and the analyticity of the solution and
an application of analyticity results
to the stochastic volatility model in
{\sc J.-P.\ Fouque}, {\sc G. Papanicolaou}, and {\sc K.~R.\ Sircar}
\cite[p, 47]{FouqPapaSir}.
The {\em Heston stochastic volatility model\/}
({\sc S.~L.\ Heston} \cite{Heston}, which is more popular)
is treated in 
{\sc B.\ Alziary} and {\sc P.\ Tak\'a\v{c}} \cite{Alziary-Tak}.
Market completeness for other stochastic volatility models
is discussed in \cite[Remark 8.7, pp. 82--83]{Takac}.
In our present work,
we are primarily interested in analyticity of the solution for
the nonlinear partial differential equation \eqref{eq:Risk_nonlin}
since the linear case \eqref{eq:Risk_lin} can be studied by applying
the results from \cite{Takac}.

The nonlinearities in {\rm eq.}~\eqref{eq:Risk_nonlin}
are uniformly Lipschitz\--continuous which enables us to apply
standard existence and uniqueness results for regular,
strongly parabolic semi\-linear Cauchy problems from, e.g.,
{\sc S.~D.\ Eidel$'$man} \cite{Eidel-69},
{\sc A.\ Friedman} \cite{Friedman, Friedman-64, Friedman-69},
{\sc A.\ Pazy} \cite[Chapt.~6, {\S}6.1, pp.\ 183--191]{Pazy},
or {\sc H.\ Tanabe} \cite{Tana-79}.
Due to the fact that the nonlinearities
$\hat{V}\mapsto \hat{V}^{\pm}\colon \RR\to \RR$
are not real\--analytic, we cannot expect any analyticity of the solution
\begin{math}
  (S,t)\mapsto \hat{V}(S,t)\colon (0,\infty)\times (0,\infty)\to \RR ,
\end{math}
neither in space nor in time.
In our approach we therefore modify the functions
$\hat{V}\mapsto \hat{V}^{\pm}\colon \RR\to \RR$ as follows:
We approximate them by real\--analytic functions with
complex\--analytic extensions to a domain ($\supset \mathbb{R}$)
in the complex plane $\CC$.
We attempt to justify this rather ``non\-rigorous'' step by arguing
that we deal with a model in Social Sciences (Economics)
where a precise nonlinear response
(i.e., the reaction function of type
 $\hat{V}\mapsto \hat{V}^{\pm}\colon \RR\to \RR$)
is hard to determine, while facing the dominant influence of
stochastic (and possibly also random) phenomena.
In our example with a single equation in one space dimension ($M=N=1$),
see eq.~\eqref{e:Cauchy:M=N=1},
we thus replace the nonlinearity $f(u) = f^{+}(u) + f^{-}(u)$
by a suitable linear combination of real\--analytic approximations of
the functions
$u\mapsto u^{+}\colon \RR\to \RR$ and $u\mapsto {}- u^{-}$
having the same asymptotic behavior at $\pm\infty$
(as $u\to \pm\infty$) and denoted by
$f^{(+)}(u)$ and $f^{(-)}(u)$, respectively,
with a complex variable $u\in \CC$.
We postpone explaining the details of this modification until
Example~\ref{exam-nontech} below.

It is a common approach to replace
$\hat{V}$ as a function of $(S,t)$ by the function
$\hat{v}(X,\tau) = \hat{V}(S,t)$ that depends on the logarithm of
the {\em asset price\/} $X = \ln S$ and
the {\em time to maturity\/} $\tau = T-t$.
Hence, by \eqref{eq:Risk_nonlin}, this function satisfies
the initial value problem
\begin{equation}
\label{eq:Risk_nonlinf2}
    \frac{\partial}{\partial \tau}\hat{v}
  - \tilde{\mathcal{A}}_t\hat{v} + r\hat{v}
  = {}- (1 - R_{B})\lambda_B f^{(-)}(\hat{v})
      - (1 - R_C)\lambda_C f^{(+)}(\hat{v}) - s_F f^{(+)}(\hat{v})
\end{equation}
with the initial value
$\hat{v}(X,0) = \hat{H}(X)\eqdef H(\exp(X))$
and the partial differential operator
\begin{equation}
\label{def:genA2}
    \tilde{\mathcal{A}}_t
  = \genfrac{}{}{}1{1}{2}\sigma^2\, \frac{\partial^2}{\partial X^2}
  + \left( q_S - \gamma_S + \genfrac{}{}{}1{1}{2}\sigma^2 \right)
    \frac{\partial}{\partial X} \,.
\end{equation}
As an alternative to the previous variable substitution it is also possible
to directly alter the stochastic processes by
\begin{equation*}
  \tilde{S}_t = \mathrm{e}^{r(T-t)} S_t \quad\mbox{ and }\quad
  \tilde{X}_t  = \ln\tilde{S}_t = X_t + r(T-t) \,,
\end{equation*}
which yields the same partial differential equation
\eqref{eq:Risk_nonlinf2} and allows for
a financial interpretation of the new variables.

Since the coefficients of the operator $\tilde{\mathcal{A}}_t$
defined in \eqref{def:genA2} are independent of the variables $X$ and~$\tau$,
the analyticity of the solution can be studied by means of the Green function;
see {\sc P.\ Tak\'{a}\v{c}} et al. \cite{TakacBoller}.
But if we replace the stochastic process \eqref{eq:stoch_dym}
that drives the value process of the asset by
a stochastic volatility process, e.g. the mean\--reverting process
from the classical paper of {\sc S.~L.\ Heston} \cite{Heston},
\begin{equation}
\label{eq:stoch_dym2}
\hspace*{-3.0cm}
\left\{\hspace*{0.5cm}
\begin{alignedat}{2}
& \mbox{ (\ref{eq:stoch_dym2}$_1$) \hspace*{2.00cm} }
& \td S_t & = \mu S_t\td t + \sqrt{V_t}\, S_t\td W_t^S \,,
\\
& \mbox{ (\ref{eq:stoch_dym2}$_2$) \hspace*{2.00cm} }
& \td V_t & = \kappa(\theta-V_t)\td t + \sigma^V\sqrt{V_t}\,\td W_t^V \,,
\\
& \mbox{ (\ref{eq:stoch_dym2}$_3$) \hspace*{2.00cm} }
& \rho & = \td W_t^S\td W_t^V \,,
\end{alignedat}
\right.
\end{equation}
the coefficients of the generator depend on the variables and
we can no longer calculate the Green function.
In the volatility process (\ref{eq:stoch_dym2}$_2$) above,
the parameters $\kappa$, $\theta$, and the volatility of volatility
$\sigma^V$ are positive constants and
$(W_t^S)_{t\geqslant 0}$ and $(W_t^V)_{t\geqslant 0}$
are one\--dimensional Brownian motions correlated by
a correlation factor $\rho\in[-1,1]$ through eq.\
(\ref{eq:stoch_dym2}$_3$).

\begin{remark}\label{rem-app:Op}\nopagebreak
\begingroup\rm
Our Hypotheses \eqref{hy:UnifEll} and \eqref{hy:Ellipt}
on the partial differential operator are consistent with
the hypotheses {\rm (H1)} and {\rm (H2)} in
{\sc P.\ Tak\'a\v{c}} \cite[p.~56]{Takac}.
As mentioned above, the operator connected to the stochastic volatility model
of {\sc J.-P.\ Fouque}, {\sc G. Papanicolaou}, and {\sc K.~R.\ Sircar}
\cite[p.~47]{FouqPapaSir}, which is parallel to
(but not more general than)
the Heston model \eqref{eq:stoch_dym2}, has been studied in detail in
\cite[Sect.~8, pp.\ 74--83]{Takac}.
Various other stochastic volatility models have been discussed in
\cite[Remark 8.7, pp.\ 82--83]{Takac}, as well.
Hence, \eqref{hy:UnifEll} and \eqref{hy:Ellipt} are satisfied for
these models; we refer the reader to \cite{Takac} for further details.
\hfill\Square
\endgroup
\end{remark}
\par\vskip 10pt

According to this remark (Remark \ref{rem-app:Op} above),
the hypotheses \eqref{hy:UnifEll} and \eqref{hy:Ellipt} are fullfilled for
{\rm eq.}~\eqref{eq:Risk_nonlinf2} even if we consider
a stochastic volatility model, e.g., like \eqref{eq:stoch_dym2},
instead of \eqref{eq:stoch_dym}.
We would like to give an example for suitable nonlinearities
$f^{(+)}$ and $f^{(-)}$ that satisfy the remaining hypothesis
\eqref{hy:f_Hardy^2} and approximate the functions
$u\mapsto u^{+}\colon \RR\to \RR$ and $u\mapsto {}- u^{-}$, respectively.

Our example is motivated by
{\sc P.\ Tak\'{a}\v{c}} \cite[Example 8.2, pp.\ 79--80]{Takac}.
We define the complex planar domains
\begin{align}
\label{theta>0:Z^(r)}
  \nabla_{\vartheta}^{(r)}
& \eqdef
  \{ \zeta = \xi\ee^{\ii\theta} + \ii\eta \in \CC\colon
     \xi\in \RR ,\ \eta\in (-r,r) ,\mbox{ and }
     |\theta| < \vartheta \} \,,
\\
\label{theta=0:Z^(r)}
  \nabla_0^{(r)}
& \eqdef
  \{ \zeta = \xi + \ii\eta \in \CC\colon \xi\in \RR ,\ \eta\in (-r,r) \}
  = \bigcap_{ 0 < \vartheta < \pi / 2 } \nabla_{\vartheta}^{(r)}
  = \RR + \ii (-r,r)
\end{align}
for $r\in (0,\infty)$ and $0 < \vartheta < \pi / 2$
with their respective closures in $\CC$ denoted by
$\overline{\nabla}_{\vartheta}^{(r)}$ and $\overline{\nabla}_0^{(r)}$;
both contain the origin $0\in \CC$.
For any given numbers
$r\in (0,\infty)$ and $0 < \vartheta < \pi / 2$,
the domain $\nabla_{\vartheta}^{(r)}$ is of the form
\begin{align*}
    \nabla_{\vartheta}^{(r)}
& = \nabla_{\vartheta}^{(0)} + \ii (-r,r)
  = \bigcup_{ \eta\in (-r,r) }
    \left( \ii\eta + \nabla_{\vartheta}^{(0)} \right)
  \subset \CC \,,\quad\mbox{ where }
\\
  \nabla_{\vartheta}^{(0)}
& \eqdef
  \{ \zeta = \xi\ee^{\ii\theta} \in \CC\colon
     \xi\in \RR\mbox{ and } |\theta| < \vartheta \}
  = (\Delta_{\vartheta}) \cup (- \Delta_{\vartheta}) \cup \{ 0\}
\end{align*}
is a {\em symmetric sector\/} in $\CC$ and the open sector
$\Delta_{\vartheta}$ is defined as in \eqref{loc:T^(infty)}.
We notice that $\nabla_0^{(r)}$ is a strip in $\CC$ and
$\mathfrak{X}^{(r)} = ( \nabla_0^{(r)} )^N \subset \CC^N$
for every $r\in(0,\infty)$.
At last, we define the domain
\begin{equation*}
  \mathcal{O}_1\eqdef
    \CC\setminus \{ \ii y\colon y\in (-\infty,-1]\cup [1,\infty) \}
\end{equation*}
that contains the closure
$\overline{\nabla}_{\vartheta}^{(r)}$ whenever
$0 < r < 1$ and $0 < \vartheta < \pi / 2$.
The definitions of these domains follow \cite[pp.\ 78--79]{Takac}.
\bigskip

We now give an example for functions $f^{(+)}$ and $f^{(-)}$
(approximating $v\mapsto v^{+}$ and $v\mapsto {}- v^{-}$, respectively)
that are analytic in $\mathcal{O}_1$ and whose first derivatives
are bounded in $\nabla_{\vartheta_0}^{(r_0)}$, whenever
$r_0\in (0,\infty)$ and $0 < \vartheta_0 < \pi / 2$.

\begin{example}\label{exam-nontech}\nopagebreak
\begingroup\rm
We consider the functions
\begin{align}
\label{ex:f}
  f^{(+)}(v) = v\left[ \frac{1}{2} + \frac{1}{\pi} \arctan(v)\right]
    \quad\mbox{and}\quad
  f^{(-)}(v) = v\left[ \frac{1}{2} - \frac{1}{\pi} \arctan(v)\right]
\end{align}
defined for every $v\in\RR$.
We have chosen $f^{(-)}$ such that
$f^{(-)}(v) = {}- f^{(+)}(-v)$ and
$f^{(+)}(v) + f^{(-)}(v) = v$ hold for all $v\in\RR$
since the nonlinearities $v^{+}$ and ${}- v^{-}$ in the original equation
\eqref{eq:Risk_nonlin} satisfy the same relations,
$v^{-} = (-v)^{+}$ and $v^{+} - v^{-} = v$, respectively.
In addition, by {\rm eq.}~\eqref{ex:f}, we have
\begin{align}
\label{exam:f}
  f^{(+)}(v) = \frac{1}{\pi}\,
               v\int_{-\infty}^v \frac{ \,\mathrm{d}t }{ 1 + t^2 }
    \quad\mbox{and}\quad
  f^{(-)}(v) = \frac{1}{\pi}\,
               v\int_v^{+\infty} \frac{ \,\mathrm{d}t }{ 1 + t^2 }
\end{align}
defined for every $v\in\RR$, which yields
\begin{equation*}
  f^{(+)}(v) > {}- \frac{1}{\pi} \;\mbox{ and }\;
  f^{(-)}(v) < \frac{1}{\pi}  \quad\mbox{ with the limits }\quad
  \lim_{v\to \mp\infty} f^{(\pm)}(v) = \mp\frac{1}{\pi} \,,
  \;\mbox{ respectively. }
\end{equation*}
We could immediately extend these two (real analytic) functions
to holomorphic (i.e., complex analytic) functions
\begin{math}
  f^{(+)}, f^{(-)}\colon \CC\setminus \{ -\ii,\, \ii\}
\end{math}
by replacing the Lebesgue integrals
$\int_{-\infty}^v \dots \,\mathrm{d}t$ and
$\int_v^{+\infty} \dots \,\mathrm{d}t$ in \eqref{exam:f}
over the real domains
$(-\infty,v]$ and $[v,+\infty)$ in \eqref{exam:f}, respectively,
by the complex path integrals along
some suitable (continuously differentiable) paths
\begin{equation*}
  \gamma_{+}\colon (-\infty,0]\to \CC\setminus \{ -\ii,\, \ii\}
  \quad\mbox{ and }\quad
  \gamma_{-}\colon [0,+\infty)\to \CC\setminus \{ -\ii,\, \ii\}
\end{equation*}
connecting the points $-\infty$ with $v$ and $v$ with $+\infty$,
respectively, whenever
$v\in \CC\setminus \{ -\ii,\, \ii\}$,
where the paths $\gamma_{+}$ and $\gamma_{-}$
do neither pass through nor wind around the points
$\pm\ii\in \CC$, i.e., they have the winding numbers
\begin{math}
    \mathrm{Ind}_{\gamma_{+}} (\pm\ii)
  = \mathrm{Ind}_{\gamma_{-}} (\pm\ii) = 0 \,.
\end{math}
As a consequence, this extension procedure could produce
{\em multi\--valued\/} analytic functions which is not desireable.
Therefore, we prefer to perform this holomorphic extension of the functions
$f^{(+)}$ and $f^{(-)}$ below, by formulas in \eqref{ex:f^tilde},
as they meet our current goals better.

We calculate the derivatives
\begin{align*}
  (f^{(+)}(v))' &= \frac{1}{2} + \frac{1}{\pi} \arctan(v)
                 + \frac{1}{\pi}\, \frac{v}{1+v^2} > 0 \,,
\\
  (f^{(-)}(v))' &= \frac{1}{2} - \frac{1}{\pi} \arctan(v)
                 - \frac{1}{\pi}\, \frac{v}{1+v^2} < 0 \,,
\end{align*}
with
\begin{math}
  (f^{(+)}(v))' + (f^{(-)}(v))' = 1
\end{math}
and
\begin{equation*}
    \lim_{v\to \infty}  (f^{(+)}(v))'
  = \lim_{v\to -\infty} (f^{(-)}(v))' = 1 \quad\mbox{and}\quad
    \lim_{v\to -\infty} (f^{(+)}(v))'
  = \lim_{v\to \infty}  (f^{(-)}(v))' = 0 \,.
\end{equation*}
By
\begin{align*}
  (f^{(+)}(v))''= \frac{1}{\pi}\, \frac{2}{ (1+v^2)^2 } > 0
    \quad\mbox{and}\quad
  (f^{(-)}(v))''= {}- \frac{1}{\pi}\, \frac{2}{ (1+v^2)^2 } < 0 \,,
\end{align*}
respectively,
$f^{(+)}$ is a strictly monotone increasing and strictly convex function,
whereas
$f^{(-)}$ is strictly monotone decreasing and strictly concave.
We can make use of
{\sc P.\ Tak\'{a}\v{c}} \cite[Example 8.2, pp.\ 79--80]{Takac}
and extend $f^{(+)}$ and $f^{(-)}$ to homolomorphic functions
$\tilde{f}^{(+)}$ and $\tilde{f}^{(-)}$ on the domain $\mathcal{O}_1$
via the following formulas,
\begin{align}
\label{ex:f^tilde}
\left\{\quad
\begin{aligned}
  \tilde{f}^{(+)}(z)
& = z\left[ \frac{1}{2} + \frac{\ii}{2\pi}\,
            \log\left( \frac{1 - \ii z}{1 + \ii z}\right)
     \right]
  = z\left[ \frac{1}{2} + \frac{1}{\pi} \arctan(z)\right]
    \quad\mbox{ and }\quad
\\
  \tilde{f}^{(-)}(z)
& = z\left[ \frac{1}{2} - \frac{\ii}{2\pi}\,
            \log\left( \frac{1 - \ii z}{1 + \ii z}\right)
     \right]
  = z\left[ \frac{1}{2} - \frac{1}{\pi} \arctan(z)\right]
\end{aligned}
\right.
\end{align}
for every
\begin{math}
    z\in \mathcal{O}_1
  = \CC\setminus \{ \pm\ii y\colon y\in [1,\infty) \} ,
\end{math}
thanks to the argument and logarithm formulas
\begin{equation*}
  \arg (1 + \ii y) = \arctan(y) \;\mbox{ for }\, y\in \RR
  \quad\mbox{ and }\quad
    \log\left( \frac{1 - \ii z}{1 + \ii z}\right)
  = {}- 2\ii\cdot \arctan(z) \;\mbox{ for }\, z\in \mathcal{O}_1 \,.
\end{equation*}
The extensions $\tilde{f}^{(+)}$ and $\tilde{f}^{(-)}$ have the restrictions
$\tilde{f}^{(\pm)}\vert_{\RR} = f^{(\pm)}$ to the real axis $\RR$,
respectively, and they are holomorphic on the domain $\mathcal{O}_1$
since the argument restriction
\begin{math}
  \arg \genfrac{(}{)}{}1{1 - \ii z}{1 + \ii z} \in (-\pi,\pi)
\end{math}
holds for $z\in \mathcal{O}_1$.
We refer to \cite[Eqs.\ (76)--(78), p.~79]{Takac}
for further discussion of the behavior of
$\log \genfrac{(}{)}{}1{1 - \ii z}{1 + \ii z}$.
As a consequence of the arguments in \cite[Example 8.2, pp.\ 79--80]{Takac},
we obtain
\begin{math}
  (\tilde{f}^{(+)}(z))' + (\tilde{f}^{(-)}(z))' = 1
\end{math}
for $z\in \mathcal{O}_1$ together with the limits
\begin{equation}
\label{z=infty:f(z)}
\left\{\qquad
\begin{alignedat}{2}
& (\tilde{f}^{(+)}(z))'\to 1
  \quad\mbox{ as }\; |z|\to \infty \;\mbox{ with $z\in \CC$, $\RE z > 0$, }
\\
& (\tilde{f}^{(-)}(z))'\to 0
  \quad\mbox{ as }\; |z|\to \infty \;\mbox{ with $z\in \CC$, $\RE z > 0$, }
\\
& (\tilde{f}^{(+)}(z))'\to 0
  \quad\mbox{ as }\; |z|\to \infty \;\mbox{ with $z\in \CC$, $\RE z < 0$, }
\\
& (\tilde{f}^{(-)}(z))'\to 1
  \quad\mbox{ as }\; |z|\to \infty \;\mbox{ with $z\in \CC$, $\RE z < 0$. }
\end{alignedat}
\right.
\end{equation}
The domain
$\mathcal{O}_1 = \CC\setminus \pm\ii [1,\infty)$
contains the strip
$\RR\times \ii (-r_0,r_0)$ for every $0 < r_0 < 1$ and
the imaginary parts of
$(\tilde{f}^{(+)}(z))'$ and $(\tilde{f}^{(-)}(z))'$
are uniformly bounded for $|\IM z| < r_0$.
Consequently, it suffices to verify eqs.\ \eqref{z=infty:f(z)}
only for $z = \RE z = x > 0$ and $x < 0$, respectively,
by Cauchy's integral theorem applied to the integral formula for
the function $\arctan(z)$ in $\mathcal{O}_1$.

Finally, in order to {\em approximate\/} the functions
$u\mapsto u^{+}\colon \RR\to \RR$ and $u\mapsto {}- u^{-}$, respectively,
we take $\eps\in (0,1)$ small enough and use the functions
\begin{align}
\label{ex:f^tilde^+-}
&
\left\{\qquad
\begin{aligned}
  \tilde{f}_{\eps}^{(+)}(z)\eqdef \eps\cdot
  \tilde{f}^{(+)}\genfrac{(}{)}{}1{z}{\eps}
& = z\left[ \frac{1}{2} + \frac{\ii}{2\pi}\,
            \log\left( \frac{ 1 - \ii (z/\eps) }{ 1 + \ii (z/\eps) }\right)
     \right] \,,
\\
  \tilde{f}_{\eps}^{(-)}(z)\eqdef \eps\cdot
  \tilde{f}^{(-)}\genfrac{(}{)}{}1{z}{\eps}
& = z\left[ \frac{1}{2} - \frac{\ii}{2\pi}\,
            \log\left( \frac{ 1 - \ii (z/\eps) }{ 1 + \ii (z/\eps) }\right)
     \right] \,,
\end{aligned}
\right.
\\
\nonumber
& \quad\mbox{ for every }\, z\in \mathcal{O}_{\eps} \eqdef
    \CC\setminus \{ \pm\ii y\colon y\in [\eps, \infty) \}
  = \eps\cdot \mathcal{O}_1 \,.
\end{align}
Notice that
\begin{math}
  \tilde{f}_{\eps}^{(+)}(z) + \tilde{f}_{\eps}^{(-)}(z) = z .
\end{math}
In particular, for every $u\in \RR$ we get
the {\em approximation\/}
\begin{align*}
  \tilde{f}_{\eps}^{(+)}(u) \,\longrightarrow\, u^{+}
    \quad\mbox{and}\quad
  \tilde{f}_{\eps}^{(-)}(u) \,\longrightarrow\, {}- u^{-}
    \quad\mbox{ as }\, \eps\to 0+ \,.
\end{align*}
This convergence is uniform on any compact interval
$[-R,R]\subset \RR$ with $0 < R < +\infty$.
\hfill\Square
\endgroup
\end{example}
\par\vskip 10pt

\begin{example}\label{exam-old}\nopagebreak
\begingroup\rm
Another example for real analytic functions $f^{(+)}$ and $f^{(-)}$
could be obtained by means of
{\sc P.\ Tak\'{a}\v{c}} \cite[Example 8.3, p.~80]{Takac}.
For this purpose, we could consider the real analytic functions
\begin{equation*}
  f^{(\pm)}(v) = \genfrac{}{}{}1{1}{2}\, v
    \pm \genfrac{}{}{}1{1}{2}\, \log (\cosh(v))
  \quad\mbox{ for every } v\in \RR \,,
\end{equation*}
which have similar properties as the functions defined in eq.~\eqref{ex:f}.
\hfill\Square
\endgroup
\end{example}
\par\vskip 10pt

\bigskip

We wish to apply our main result, Theorem \ref{thm-Main},
to the semilinear inital value problem in eq.~\eqref{eq:Risk_nonlinf2}
where we choose $f^{(+)}$ and $f^{(-)}$ as in eq.~\eqref{ex:f},
$\tilde{f}^{(+)}$ and $\tilde{f}^{(-)}$
as in eq.~\eqref{ex:f^tilde}, and
$\tilde{f}_{\eps}^{(+)}$ and $\tilde{f}_{\eps}^{(-)}$
as in eq.~\eqref{ex:f^tilde^+-}, respectively.
On the right\--hand side of eq.~\eqref{eq:Risk_nonlinf2}
we replace the (nonlinear) functions
$\hat{v}\mapsto \hat{v}^{+}\colon \RR\to \RR$ and
$\hat{v}\mapsto {}- \hat{v}^{-}$, respectively, by the pair of functions
$\hat{v}\mapsto \tilde{f}_{\eps}^{(+)}(\hat{v}) \colon \RR\to \RR$
and
$\hat{v}\mapsto \tilde{f}_{\eps}^{(-)}(\hat{v})$
defined in eq.~\eqref{ex:f^tilde^+-}.
Here, we take $\eps\in (0,1)$ small enough, but fixed.
The initial data in eq.~\eqref{eq:Risk_nonlinf2} are given by
a payoff function $\hat{H}\in B^{s;p,p}(\RR^N)$ with
$p > 2 + \frac{N}{m} = 3$ and
\begin{math}
  s = 2m\left( 1 - \frac{1}{p}\right) = 2\left( 1 - \frac{1}{p}\right)
\end{math}
($M = N = m = 1$).
We assume that these initial data possess a holomorphic extension
\begin{math}
  \tilde{H}\colon \mathfrak{X}^{(\kappa_0)}\to \CC^1
\end{math}
from $\RR^1$ to the complex domain
$\mathfrak{X}^{(\kappa_0)} \subset \CC^1$,
for some $\kappa_0\in (0,r_0]$, such that the function
\begin{equation*}
  \tilde{H} (\,\cdot\, + \ii y)
  \colon x\longmapsto \tilde{H} (x + \ii y)
  \colon \RR^1\longrightarrow \CC^1
\end{equation*}
belongs to $\mathbf{B}^{s;p,p}(\RR^1)$ 
for each $y\in Q^{(\kappa_0)}$ and has finite norm
$\mathfrak{N}^{(\kappa_0)} (\tilde{H}) < \infty$
which has been defined in eq.~\eqref{sup_y:u_0}.
In the case of a simple European call or put option, i.e.,
$\hat{H}(x) = (\ee^x - K)^{+}$ or
$\hat{H}(x) = (K - \ee^x)^{+}$, $x\in \RR^1$,
respectively, one may use the functions
$\tilde{f}_{\eps}^{(+)}$ and $\tilde{f}_{\eps}^{(-)}$
in order to find the desired holomorphic extension
$\tilde{H}$ of $\hat{H}$ that satisfies the hypotheses required
in Theorem \ref{thm-Main}, Part~{\bf (iii)}.

According to our Example~\ref{exam-nontech} above, the nonlinearity
\begin{equation*}
  f(\hat{v})
  = {}- (1 - R_{B})\lambda_B f^{(-)}(\hat{v})
      - (1 - R_C)\lambda_C f^{(+)}(\hat{v}) - s_F f^{(+)}(\hat{v})
\end{equation*}
on the right\--hand side of eq.~\eqref{eq:Risk_nonlinf2}
possesses a holomorphic extension
\begin{align}
\label{eq:Risk^tilde}
  \tilde{f}_{\eps}(\hat{v})
  = {}- (1 - R_{B})\lambda_B \tilde{f}_{\eps}^{(-)}(\hat{v})
      - (1 - R_C)\lambda_C \tilde{f}_{\eps}^{(+)}(\hat{v})
      - s_F \tilde{f}_{\eps}^{(+)}(\hat{v})
  \quad\mbox{ for all }\, \hat{v}\in \mathcal{O}_{\eps} \,,
\\
\nonumber
  \quad\mbox{ where }\,
  \mathcal{O}_{\eps} = \eps\cdot \mathcal{O}_1
  = \CC\setminus
    \{ \ii y\colon y\in (-\infty, -\eps]\cup [\eps,\infty) \} \,.
\end{align}
Now let us recall the definition of the complex planar domain
$\nabla_{\vartheta}^{(r)}\subset \CC$ in eq.~\eqref{theta>0:Z^(r)}
for $r\in (0,\infty)$ and $0 < \vartheta < \pi / 2$
with the closure $\overline{\nabla}_{\vartheta}^{(r)}$ in $\CC$.
We have
\begin{math}
  \overline{\nabla}_{\vartheta}^{(r)} \subset \mathcal{O}_{\eps}
\end{math}
whenever $0 < r < \eps$ and $0 < \vartheta < \pi / 2$.
Moreover, both $\tilde{f}_{\eps}$ and its complex derivative
$\tilde{f}_{\eps}'$ are uniformly bounded in
$\overline{\nabla}_{\vartheta}^{(r)}$.
Thus, fixing any number $\eps\in (0,1)$ and taking
$r\in (0,\eps)$, we observe that both
$\tilde{f}_{\eps}$ and $\tilde{f}_{\eps}'$ are uniformly bounded in
$\overline{\nabla}_{\vartheta}^{(r)}$ and, consequently,
also in the complex strip $\mathfrak{X}^{(r)}$,
\begin{math}
  \mathfrak{X}^{(r)} \subset
  \overline{\nabla}_{\vartheta}^{(r)} \subset \mathcal{O}_{\eps} .
\end{math}
We stress that the number $\eps\in (0,1)$ may be chosen arbitrarily small
in order to achieve a sufficiently precise approximation of
the reaction function $f = f(\hat{v})$ by the holomorphic function
$\tilde{f}_{\eps} = \tilde{f}_{\eps}(\hat{v})$
as desribed in Example~\ref{exam-nontech} above.
Naturally, the choice of a smaller number $\eps\in (0,1)$
diminishes the width of the strip $\mathfrak{X}^{(r)}$ according to
$0 < r < \eps$.

We have
$\tilde{f}_{\eps}^{(\pm)}\in A(\mathcal{O}_{\eps})$
and $f'$ is bounded in $\nabla_{\vartheta_0}^{(r_0)}$
for every $r_0\in (0,\, \eps / 2)$ and $0 < \vartheta_0 < \pi / 2$.
The technical estimate \eqref{e:f_Hardy^2} is trivially satisfied, owing to
$\tilde{f}_{\eps}^{(\pm)}(0) = 0$.

Following the discussion in Section \ref{s:Analyt_space},
we recall that the (unique) strict solution is restricted to
the bounded open set $U\subset B^{s;p,p}(\RR^N)$ defined in \eqref{def:U},
which is, due to the continuous Sobolev imbedding
$B^{s;p,p}(\RR^N)\hookrightarrow L^{\infty}(\RR^N)$,
bounded in $L^{\infty}(\RR^N)$, as well.
Hence, it is convenient to loosen Hypothesis \eqref{hy:f_Hardy^2}
in the sense that we replace the complex plane $\CC$
in the assumptions by smaller domains.
In particular, the function
$\tilde{f}_{\eps}$ in \eqref{eq:Risk^tilde} fulfills
Hypothesis \eqref{hy:f_Hardy^2} with such weakened assumptions.
Since all requirements are satisfied, we can apply
Theorem \ref{thm-Main} to the initial value problem \eqref{eq:Risk_nonlinf2}
and obtain the real analyticity of the solution.

Indeed, we apply our main result, Theorem \ref{thm-Main},
to the inital value problem \eqref{eq:Risk_nonlinf2}, where
we choose $f^{(+)}$ and $f^{(-)}$ as follows:
We replace the functions $f^{(+)}$ and $f^{(-)}$, respectively,
by their complexifications
$\tilde{f}_{\eps}^{(+)}(z)$ and
$\tilde{f}_{\eps}^{(-)}(z)$, respectively,
defined in formulas \eqref{ex:f^tilde^+-} for
$z\in \mathcal{O}_{\eps} = \eps\cdot \mathcal{O}_1$.
Here, $\eps > 0$ is as small as needed.
The initial data is given by the payoff function
$\hat{H}\in B^{s;p,p}(\RR^N)$ for $p>3$ and $s=2(1-1/p)$.
The partial differential operator $\tilde{\mathcal{A}}_t$
defined in \eqref{def:genA2} satisfies
Hypotheses \eqref{hy:UnifEll} and \eqref{hy:Ellipt} (with $N=1$).
If we replace $\tilde{\mathcal{A}}_t$ by the generator of
a stochastic volatility process like \eqref{eq:stoch_dym2},
then Hypotheses \eqref{hy:UnifEll} and \eqref{hy:Ellipt} (with $N=2$)
are still fulfilled, according to Remark \ref{rem-app:Op}.
For the nonlinearity $f(\hat{v})$ in \eqref{eq:Risk_nonlinf2},
extended in eq.~\eqref{eq:Risk^tilde} as
$\tilde{f}_{\eps}(\hat{v})$ for $\hat{v}\in \mathcal{O}_{\eps}$,
we have
$\tilde{f}_{\eps}\in A(\mathcal{O}_{\eps})$
and $\tilde{f}_{\eps}'$ is bounded in $\nabla_{\vartheta_0}^{(r_0)}$
for every $r_0\in (0,\, \eps / 2)$ and $0 < \vartheta_0 < \pi / 2$.
The technical estimate \eqref{e:f_Hardy^2} is trivially satisfied, owing to
$\tilde{f}_{\eps}^{(\pm)}(0) = 0$.

\section{Historical remarks and comments}
\label{s:Hist}

The questions we have studied in this paper are clearly related to
the classical Cauchy\--Kowalewski theorem
({\sc F.\ John} \cite{John}, Chapt.~3, Sect.\ 3(d), pp.\ 73--77).
It has been known since the work by
{\sc E.\ Holmgren} \cite{Holm-08}
that even the heat equation (in one space dimension(!))
has solutions that are {\em not\/} real analytic in the time variable
(cf.\ {\sc G.~G.\ Bilodeau} \cite[pp.\ 124--125]{Bilodeau}).
This phenomenon is due to a possibly very rapid growth of the solutions
as the spatial variable $x\in \RR$ escapes to $\pm\infty$;
to eliminate it one needs to restrict the function space,
where the solutions are considered at each time moment $t\in \RR_+$,
in order to prevent a too rapid growth of the solutions as
$x\to \pm\infty$.
This is precisely what has been done also in our present article.

Here, the emphasis is on the analytic dependence in time~$t$ and
the Cauchy problem~\eqref{e:Cauchy} is viewed as
an evolutionary equation in some suitable function space, e.g.,
$L^2(\RR)$ or $L^2(\RR^N)$.
Consequently, the solution is viewed as a vector\--valued function
$u\colon(0,T)\to L^2(\RR^N)$ and, thus, regularity results
(including analyticity results) have been obtained in this setting.
The interested reader is referred to
{\sc P.\ Tak\'{a}\v{c}} \cite[Sect.~9, pp.\ 83--85]{Takac}
for a number of pertinent references and their description;
for example,
{\sc T.\ Kato} and {\sc H.\ Tanabe} \cite{KatoTana-67},
{\sc H.\ Komatsu} \cite{Komatsu},
{\sc F.~J.\ Massey III} \cite{Massey},
{\sc K.\ Masuda} \cite{Masuda-72}, and in particular
{\sc H.\ Tanabe} \cite{Tana-79} and the references therein.

Investigation of
the {\it smoothing\/} (or {\it regularizing\/}) effect in
evolutionary equations of parabolic type has a long history; see e.g.\
{\sc S.~D.\ Eidel$'$man} \cite{Eidel-69},
{\sc A.\ Friedman} \cite{Friedman, Friedman-64, Friedman-69},
{\sc A.\ Pazy} \cite{Pazy}, and {\sc H.\ Tanabe} \cite{Tana-79}
and numerous references therein.
{\em\bfseries Analytic\/} smoothing (or regularizing) effects,
similar to those treated in our present article,
in the space ($x$) and/or time ($t$) variable(s),
have been obtained somewhat later, beginning with the theory of
analytic semigroups (in an abstract Banach space), see e.g.\
the monographs by
{\sc T.\ Kato} \cite{Kato}, {\sc J.-L.\ Lions} \cite{Lions-61},
{\sc A.\ Pazy} \cite{Pazy}, and {\sc H.\ Tanabe} \cite{Tana-79},
and applying (extending) it to
non\-autonomous analytic evolutionary equations, see e.g.\
{\sc T.\ Kato} and {\sc H.\ Tanabe} \cite{KatoTana-67},
{\sc H.\ Komatsu} \cite{Komatsu},
{\sc K.\ Masuda} \cite{Masuda-72}, and
{\sc H.\ Tanabe} \cite{Tana-79}.
Evolutionary equations exhibiting analytic smoothing effects
may be split into the following two classes:
{\em\bfseries dissipative\/} and {\em\bfseries dispersive\/}.
Again, we refer to
\cite[Sect.~9, pp.\ 83--85]{Takac}
for greater details about these two classes.
The results for dissipative evolutionary equations establish only
analyticity with respect to the time variable
$t\in (0,T)\subset \RR$.
{\sc N.\ Hayashi} and {\sc K.\ Kato} \cite{Hayashi-Kato}
establish an analogous time\--analyticity result for 
the nonlinear Schr\"odinger equation (NLS).
The early (general) treatments on the analytic smoothing effect 
with respect to the space variable $x\in \RR^N$ are given in
{\sc C.~S.\ Kahane} \cite{Kahane} and
{\sc C.\ Foias} and {\sc R.\ Temam} \cite{Foias-Tem-1, Foias-Tem-2}.

Finally, we mention the analyticity results by
{\sc G.\ Komatsu} \cite{Gen_Komatsu}
obtained for solutions to
elliptic and parabolic problems in a {\it bounded\/} spatial domain
$\Omega\subset \RR^N$
(with analytic boundary $\partial\Omega$).
Analyticity in the space variable $x$ and
$2$-nd Gevrey class regularity (weaker than analyticity)
in the time variable $t$ are established in
{\sc A.\ Cavallucci} \cite[Teorema 6.1, p.~166]{Cavallucci}
for linear parabolic equations.
Some results about the analyticity of solutions of
nonlinear parabolic systems, which are related to ours, are stated in
{\sc A.\ Friedman} \cite[Theorems 3 and 4]{Friedman}
without proofs, and for linear elliptic systems in
{\sc C.~B.\ Morrey, Jr.,} and {\sc L.\ Nirenberg} \cite{Morrey-Niren}.
For the Navier\--Stokes equations,
such analyticity results have been established in
{\sc K.\ Masuda} \cite{Masuda-80} and,
with respect to the space variable $x\in \RR^N$ only, earlier in
{\sc C.~S.\ Kahane} \cite{Kahane} and
{\sc K.\ Masuda} \cite{Masuda-67}.
These results state local analyticity of
infinitely differentiable solutions without any description of
their domain of holomorphy (i.e., domain of complex analyticity).
Our present article provides such description in
Theorem~\ref{thm-Main} and so do Refs.\
\cite{BonaGrujKal-2005, BonaGrujKal-2006}.
More results of global nature on the space analyticity can be found in
{\sc C.\ Bardos} and {\sc S.\ Benachour} \cite{BardosBen}
and
{\sc Z.\ Gruji\'c} and {\sc I.\ Kukavica} \cite{Grujic-Kuka}.

\section{Discussion and possible generalizations}
\label{s:Discuss}

In contrast to the analytic smoothing results established in
{\sc P.\ Tak\'{a}\v{c}} \cite[Theorem 3.3, p.~59]{Takac}
(for a linear parabolic problem)
and
{\sc P.\ Tak\'a\v{c}} et al.\ \cite[Theorem 2.1, p.~429]{TakacBoller}
(for a semi\-linear parabolic problem),
in the present work we have focused on {\it preserving\/}
the spatial analyticity of the initial data, $\mathbf{u}_0$,
for all times $t\in [0,T]$ as long as a (global) weak solution
$\mathbf{u}\in C\left( [0,T]\to \mathbf{B}^{s;p,p}(\RR)\right)$
to the Cauchy problem \eqref{e:Cauchy} exists, that is, loosely written,
$\mathbf{u}(\,\cdot\,,0) = \mathbf{u}_0$ is spatially analytic
(at $t=0$)
$\;\Longrightarrow\;$
$\mathbf{u}(\,\cdot\,,t)$ is spatially analytic at all times $t\in (0,T]$,
even for all $t\in (0, T+T_1]$ with some $T_1 > 0$ small enough.

However, also a spatial analytic smoothing result analogous to those in
\cite[Theorem 3.3, p.~59]{Takac}
and
\cite[Theorem 2.1, p.~429]{TakacBoller}
should hold in our present setting in the Besov space
$\mathbf{B}^{s;p,p}(\RR)$, by arguments similar to those used in
\cite[pp.\ 434--435]{TakacBoller}, proof of Lemma 3.4.
The Banach contraction principle can then be used in analogy with
\cite[pp.\ 437--438]{TakacBoller},
{\em Step~$4$\/} in the proof of Theorem 3.1.
This approach requires separation of the linear part of
the Cauchy problem \eqref{e:Cauchy}
(cf.\ \cite{Takac})
followed by an application of the Banach contraction principle
to the full semi\-linear parabolic problem in \eqref{e:Cauchy}
(cf.\ \cite[Theorem 2.1, p.~429]{TakacBoller}).

{\bf Acknowledgments.}$\;$
\begin{small}
This work was partially supported by a grant from
Deutsche Forschungs\-gemeinschaft (DFG, Germany),
Grant no.\ TA~213/16-1, and by
the Federal Ministry of Education and Research of Germany
under grant No.\ 57063847 (D.A.A.D.\ Program ``PPP''),
the latter within a common exchange program between
the Czech Republic and Germany.
The authors would like to express their sincere thanks to
Professors
Herbert Amann (University of Z\"urich, Switzerland),
Philippe Cl\'ement (Delft University of Technology, The Netherlands),
and
Alessandra Lunardi (University of Parma, Italy)
for their kind help with the literature on
the maximal regularity and analyticity questions
for linear parabolic problems.
\end{small}

\bigskip
\bigskip
\bigskip


%
%
\makeatletter \renewcommand{\@biblabel}[1]{\hfill#1.} \makeatother
%
%

\end{document}